\documentclass[reqno,11pt]{amsart}
\usepackage[english]{babel}
\renewcommand{\baselinestretch}{1.1}
\usepackage{graphicx,mathrsfs,tikz,latexsym,ifthen,amsmath,amsfonts,amssymb,amsthm,stmaryrd,fancyhdr}
\usepackage{enumerate,enumitem,empheq,mathtools,times}
\mathtoolsset{showonlyrefs,showmanualtags}
\usepackage[utf8x]{inputenc}
\usepackage[toc,page]{appendix}
\usepackage[pdftex]{hyperref}

\allowdisplaybreaks

\input xy
\xyoption{all}

\usepackage{color}
\definecolor{MyDarkBlue}{rgb}{0.15,0.25,0.45}
\usepackage{hyperref}
\hypersetup{
colorlinks=true,
citecolor=black,
linkcolor=black,
urlcolor=MyDarkBlue,
breaklinks=true
}

\newcommand{\hilb}[2]{\mathrm{Hilb}^{#1}(#2)}

\newcommand{\Acal}{{\mathcal A}}

\newcommand{\Ecal}{{\mathcal E}}

\newcommand{\Kcal}{{\mathcal K}}
\newcommand{\Ncal}{{\mathcal N}}
\newcommand{\Mcal}{{\mathcal{M}}}

\newcommand{\Zcal}{{\mathcal Z}}
\newcommand{\Ocal}{{\mathcal O}}
\newcommand{\Gcal}{{\mathcal G}}
\newcommand{\Fcal}{{\mathcal F}}
\newcommand{\Ical}{{\mathcal I}}

\newcommand{\Rcal}{{\mathcal R}}
\newcommand{\Lcal}{{\mathcal L}}

\newcommand{\Tcal}{{\mathcal T}}
\newcommand{\Pcal}{{\mathcal P}}

\newcommand{\Qfrak}{\mathfrak{Q}}
\newcommand{\hfrak}{\mathfrak{h}}
\newcommand{\slfrak}{\mathfrak{sl}}

\newcommand{\ufrak}{\mathfrak{u}}

\newcommand{\Uscr}{\mathscr{U}}
\newcommand{\Tscr}{\mathscr{T}}

\newcommand{\Xscr}{{\mathscr{X}}}
\newcommand{\Yscr}{\mathscr{Y}}
\newcommand{\Dscr}{{\mathscr{D}}}
\newcommand{\Bscr}{{\mathscr{B}}}

\newcommand{\IXscr}{{\mathcal{I}(\mathscr{X})}}
\newcommand{\IXkscr}{{\mathcal{I}(\mathscr{X}_k)}}

\newcommand{\qsf}{\mathsf{q}}
\newcommand{\Osf}{\mathsf{O}}

\newcommand{\Ebf}{\boldsymbol{E}}
\newcommand{\Ubf}{\boldsymbol{U}}
\newcommand{\ubf}{\boldsymbol{u}}
\newcommand{\Vbf}{\boldsymbol{V}}

\newcommand{\R}{{\mathbb{R}}}
\newcommand{\N}{{\mathbb{N}}}
\newcommand{\C}{{\mathbb{C}}}
\newcommand{\Z}{{\mathbb{Z}}}
\newcommand{\Q}{{\mathbb{Q}}}
\newcommand{\T}{{\mathbb{T}}}
\newcommand{\PP}{{\mathbb{P}}}
\newcommand{\A}{{\mathbb{A}}}
\newcommand{\G}{{\mathbb{G}}}
\newcommand{\F}{{\mathbb{F}}}

\def\DG{\operatorname{DG}}
\def\Pic{\operatorname{Pic}}
\def\Hom{\operatorname{Hom}}
\def\Ext{\operatorname{Ext}}
\def\Spec{\operatorname{Spec}}
\def\rk{\operatorname{rk}}
\def\ch{\operatorname{ch}}
\def\eu{\operatorname{Euler}}
\def\cone{\operatorname{Cone}}
\def\link{\operatorname{link}}
\def\boxx{\operatorname{Box}}
\def\crm{\operatorname{c}}

\newcommand{\splus}{S(\beta\alpha)^+}
\newcommand{\dplus}{d(\beta\alpha)^+}
\newcommand{\sminus}{S(\beta\alpha)^-}
\newcommand{\dminus}{d(\beta\alpha)^-}

\def\eulerclass{\prod_{\alpha,\beta=1}^r\limits\ \prod_{i=1}^k
  \limits\, m_{Y_{\alpha}^{i},
    {Y_{\beta}^{i}}}\big(\varepsilon_1^{(i)},\varepsilon_2^{(i)},
  a_{\beta\alpha}^{(i)} \big) \ \prod_{n=1}^{k-1} \limits \,
  \ell^{(n)}_{\vec{v}_{\beta\alpha}}\big(\varepsilon_1^{(n)},\varepsilon_2^{(n)},
  a_{\beta\alpha} \big)}

\def\adjointclass{\prod_{\alpha,\beta=1}^r\limits \ \prod_{i=1}^k
  \limits \, m_{Y_{\alpha}^{i},
    {Y_{\beta}^{i}}}\big(\varepsilon_1^{(i)},\varepsilon_2^{(i)},
  a_{\beta\alpha}^{(i)}+\mu\big) \ \prod_{n=1}^{k-1} \limits\,
  \ell^{(n)}_{\vec{v}_{\beta\alpha}}\big(\varepsilon_1^{(n)},\varepsilon_2^{(n)},
  a_{\beta\alpha}+\mu \big)}

\def\fundclass-s{\prod_{\alpha=1}^r\limits\ \prod_{i=1}^k\limits \,
  m_{Y_{\alpha}^{i}}\big(\varepsilon_1^{(i)},\varepsilon_2^{(i)},
  a_\alpha^{(i)}+\mu_s \big) \ \prod_{n=1}^{k-1}
  \limits\, \ell^{(n)}_{\vec{v}_{\alpha}} \big(\varepsilon_1^{(n)},\varepsilon_2^{(n)},
  a_\alpha+\mu_s \big)}

\def\e{{\,\rm e}\,}
\def\ii{{\,{\rm i}\,}}

\newcommand{\rootdivvec}{\sqrt[\vec{k}]{\vec{\Dscr}/\Xscr}}
\newcommand{\rootlines}{\sqrt[k]{(\mathcal{L},s)/\mathscr{X}}}
\newcommand{\rootline}{\sqrt[k]{\mathcal{L}/\mathscr{X}}}
\newcommand{\schemerootlines}{\sqrt[k]{(L,s)/X}}
\newcommand{\schemerootline}{\sqrt[k]{L/X}}
\newcommand{\rootdiv}{\sqrt[k]{\mathscr{D}/\mathscr{X}}}

\newcommand{\stackyfan}{\mathbf{\Sigma}}
\newcommand{\stackyfank}{\mathbf{\bar{\Sigma}}_k}

\newcommand{\triend}{\parbox{2mm}{\hfill} \hfill\mbox{\hspace{0.2mm}}\hfill$\triangle$}
\newcommand{\ocend}{\parbox{2mm}{\hfill} \hfill\mbox{\hspace{0.2mm}}\hfill$\oslash$}

\newtheorem{theorem}[equation]{Theorem}
\newtheorem{proposition}[equation]{Proposition}
\newtheorem{lemma}[equation]{Lemma}
\newtheorem{corollary}[equation]{Corollary}
\newtheorem{conjecture}[equation]{Conjecture}
\newtheorem*{theorem*}{Theorem}
\newtheorem*{proposition*}{Proposition}
\newtheorem*{corollary*}{Corollary}

\numberwithin{equation}{section}

\theoremstyle{remark}
\newtheorem{example}[equation]{Example}

\theoremstyle{remark}
\newtheorem{rem}[equation]{Remark}
\newenvironment{remark}{\begin{rem}}{\triend\end{rem}}

\theoremstyle{definition}
\newtheorem{defin}[equation]{Definition}
\newenvironment{definition}{\begin{defin}}{\ocend\end{defin}}

\begin{document}

\begin{flushright}
EMPG--13--21
\end{flushright}

\vskip 1cm

\title[Framed sheaves on root stacks and gauge theories 
    on ALE spaces]{\large{Framed sheaves on root stacks \\[5pt]
     and supersymmetric gauge theories on ALE spaces}}

\vskip 1cm

\maketitle \thispagestyle{empty}

\begin{center}
{\large \sc Ugo Bruzzo$^{\S\ddag}$, Mattia Pedrini$^{\S\ddag}$,\\[3pt]
  Francesco Sala$^{\P\star\bullet}$} \ and \ {\large \sc Richard J. Szabo$^{\P\star\circ}$} \\[8pt] 
$^\S$ Scuola Internazionale Superiore di Studi Avanzati {\sc (SISSA),}\\ Via Bonomea 265, 34136
Trieste, Italy; \\[3pt] $^\ddag$ Istituto Nazionale di Fisica Nucleare, Sezione di Trieste, Italy; \\[3pt]
$^\bullet$ Department of Mathematics, The University of Western Ontario,\\
Middlesex College, London N6A 5B7, Ontario, Canada;\\[3pt]
$^\P$ Department of Mathematics, Heriot-Watt University,\\ Colin
Maclaurin Building, Riccarton, Edinburgh EH14 4AS, United Kingdom; \\[3pt]
$^\star$ Maxwell Institute for Mathematical Sciences, Edinburgh,
United Kingdom; \\[3pt] $^\circ$ The Tait Institute, Edinburgh, United Kingdom
 \end{center}

\par\vfill

\noindent \begin{quote} 
{\sc Abstract.} \small
We develop a new approach to the study of supersymmetric gauge
theories on ALE spaces using the theory of framed sheaves on root toric
stacks, which illuminates relations with gauge theories on
$\R^4$ and with two-dimensional conformal field theory. We
construct a stacky compactification of the minimal resolution $X_k$ of
the $A_{k-1}$ toric singularity $\C^2/\Z_k$, which is a projective toric orbifold
$\Xscr_k$ such that $\Xscr_k\setminus X_k$ is a
$\Z_k$-gerbe. We construct   moduli spaces of torsion free
sheaves on $\Xscr_k$ which are framed along the compactification gerbe. We prove
that this moduli space is a smooth quasi-projective variety, compute
its dimension, and classify its fixed points under the natural induced toric action. We
use this construction to compute the partition functions
and correlators of chiral BPS operators for 
$\Ncal=2$ quiver gauge theories on $X_k$ with nontrivial
holonomies at infinity. The partition functions are computed with and
without couplings to bifundamental matter hypermultiplets and expressed in terms of toric
blowup formulas, which relate them to the
corresponding Nekrasov partition functions on the affine toric subsets
of $X_k$. We compare our new partition functions with previous
computations, explore their connections to the
representation theory of affine Lie algebras, and find new constraints on fractional instanton charges in the coupling to
fundamental matter. We show that the
partition functions in the low energy limit are characterised by the Seiberg-Witten
curves, and in some cases also by suitable blowup equations involving Riemann theta-functions
on the Seiberg-Witten curve with characteristics related to the
nontrivial holonomies.
\end{quote}

\par
\vfill
\parbox{.95\textwidth}{\hrulefill}\par
\noindent \begin{minipage}[c]{\textwidth}\parindent=0pt \renewcommand{\baselinestretch}{1.2}
\small
\emph{Date:} Revised  December 2015  \par 
 \par
\emph{E-Mail:} \texttt{bruzzo@sissa.it, mattia.pedro@gmail.com,
  salafra83@gmail.com,\\ R.J.Szabo@hw.ac.uk} \par
\end{minipage}

\newpage

\setlength{\parskip}{0.2ex}

\tableofcontents

\setlength{\parskip}{0.6ex plus 0.3ex minus 0.2ex}

\section{Introduction and summary}

\subsection{Background}

Gauge theories on ALE spaces of type $A_{k-1}$ have recently received
renewed interest as they provide a natural extension of the AGT
correspondence
\cite{art:aldaygaiottotachikawa2010,art:wyllard2009,art:aldaytachikawa2010}
to nontrivial four-manifolds. In this paper we aim to provide a new
sheaf-theoretic approach to the study of supersymmetric gauge theories on
$A_{k-1}$ ALE spaces which enables one to establish new relationships
between moduli spaces of sheaves and infinite-dimensional Lie
algebras. In this subsection we provide some physical and
mathematical background to the problems we shall address in this work.

An ALE space of type $A_{k-1}$ is a four-manifold $X$ which is
diffeomorphic to the minimal resolution $\varphi_k\colon
X_k\to\C^2/\Z_k$ of the simple Kleinian singularity $\C^2/\Z_k$, and
is equipped with a K\"ahler metric that is asymptotically locally
  Euclidean (ALE), i.e., there is a compact subset $K\subset X$ and a
diffeomorphism $X\setminus K\rightarrow \left(\mathbb{C}^2\setminus
  \mathrm{B}_r(0)\right)/\mathbb{Z}_k$ under which the metric is
approximated by the standard Euclidean metric on
$\mathbb{C}^2/\mathbb{Z}_k$. We call $X\setminus K$ the
\emph{infinity} of $X$. The ALE space $X$ can be realized as a
hyper-K\"ahler quotient \cite{art:kronheimer1989}, hence it depends on
a stability parameter $\xi:=(\xi_\R,\xi_\C)$ which belongs to an open subset $U$
of $H^2(X_k; \R)\times H^2(X_k; \C)\simeq \R^k\times\C^k$
\cite[Section~2]{art:joyce1999}; roughly speaking, a representative of
$\xi_\R$ is a K\"ahler form and a representative of $\xi_\C$ is a
holomorphic volume form on $X$. We denote this $\xi$-dependence of the
ALE space by $X_k(\xi)$. The natural toric structure on $\C^2/\Z_k$ lifts to its
resolution $X_k$. The McKay correspondence gives a bijection between the irreducible components of
$\varphi_k^{-1}(0)$, which are torus-invariant smooth projective
curves $D_i$ of genus zero for $i=1,\dots,k-1$, and the vertices of
the Dynkin diagram of type $A_{k-1}$; the intersection matrix of these
curves is exactly $-C$, where $C$ is the corresponding Cartan
matrix. Any ALE space inherits these properties, so the 
homology group $H_2(X_k(\xi); \Z)$ of an ALE space $X_k(\xi)$ can be
identified with the root lattice $\Qfrak$ of type $A_{k-1}$, while the
cohomology group $H^2(X_k(\xi); \R)$ is the real Cartan subalgebra of $\slfrak(k)$ \cite[Section~4]{art:kronheimer1989}.

The connection between gauge theories on ALE spaces and
two-dimensional conformal field theory in physics, or between moduli spaces and
infinite-dimensional Lie algebras from a mathematical perspective, goes back to pioneering works of Nakajima \cite{art:nakajima:1994-1,art:nakajima:1994-2,art:nakajima:1994-3}, who showed that one can construct highest weight representations of
affine Lie algebras using quiver varieties.   Nakajima's quiver
varieties arise from the ADHM construction of $U(r)$ gauge theory
instantons on ALE spaces. A $U(r)$ instanton on $X_k(\xi)$ is a pair
$(E,\nabla)$ consisting of a Hermitian vector bundle $E\to X_k(\xi)$
of rank $r$ and a unitary connection $\nabla$ on $E$ whose curvature is
anti-selfdual and square-integrable. The connection $\nabla$ is
additionally characterised by its behaviour at infinity: $\nabla$ is
flat but not necessarily trivial at infinity, so it is also
classified by its holonomy which takes values in the fundamental group
$\pi_1(X\setminus K)\simeq\pi_1(S^3/\Z_k)\simeq\Z_k$, and hence
corresponds to a homomorphism $\rho\colon \Z_k\to U(r)$. As described
in \cite{art:gochonakajima1992}, with each irreducible representation of
$\Z_k$ one can associate a Hermitian line bundle on $X_k(\xi)$, called a
tautological line bundle, which is endowed with an anti-selfdual
square-integrable connection whose holonomy at infinity corresponds
to the representation. The tautological line bundles form an integral
basis of the Picard group $\Pic(X_k(\xi))\simeq H^2(X_k(\xi);\Z)$;
their intersection matrix is minus the inverse of the Cartan matrix
$C$ of the $A_{k-1}$ Dynkin diagram. In
\cite{art:kronheimernakajima1990}, Kronheimer and Nakajima provided a
description of $U(r)$ instantons on $X_k(\xi)$ in terms of ADHM type
data depending on two vectors $\vec v=(v_0,v_1,\dots,v_{k-1}),\vec
w=(w_0,w_1,\dots,w_{k-1})\in\N^k$, which  
parameterize the Chern character $\ch(E)$, with $\crm_1(E)$ given
in terms of the first Chern classes of the tautological line bundles,
and the multiplicities of the $k$ one-dimensional irreducible
representations of $\Z_k$ in the holonomy at infinity $\rho$. By using
the ADHM description, moduli spaces parameterizing $U(r)$ instantons
on $X_k(\xi)$ can be realized as hyper-K\"ahler quotients
$\Mcal_\xi(\vec v,\vec w\,)$ depending on a \emph{chamber} containing $\xi$ in its closure. By perturbing the
real part of the moment map of the hyper-K\"ahler quotient
construction one obtains
\emph{larger} moduli spaces, which are Nakajima quiver varieties
associated with the affine extended Dynkin diagram of type $\hat
A_{k-1}$; they are smooth quasi-projective varieties. Nakajima
proved that the cohomology of these moduli spaces is a representation
of $\widehat{\slfrak}(k)$ acting at level $r$; the generators of the
affine Kac-Moody algebra are defined using geometric Hecke
correspondences. The McKay correspondence now gives a bijective
equivalence between the tautological line bundles on $X_k(\xi)$ and
integrable highest weight representations of the affine Lie algebra
$\widehat{\slfrak}(k)$ acting at level $r$; in the setting of
two-dimensional conformal field theory, the Chern characters
$\crm_1(E)$ and $\ch_2(E)$ are respectively identified
with the momentum $\vec p$ and the energy $L_0$ of the highest weight module over $\widehat{\slfrak}(k)_r$ corresponding to the unitary representation $\rho\colon \Z_k\to U(r)$.

A gauge theory realization of Nakajima's results first appeared in
the seminal work of Vafa and Witten~\cite{art:vafawitten1994} who
showed that the partition function of a certain topologically twisted version of $\Ncal=4$
supersymmetric Yang-Mills theory on $X_k(\xi)$ computes the
characters of highest weight representations of
$\widehat{\slfrak}(k)_r$. To extend this realization beyond the level of
characters, one should relate the instanton partition functions of
corresponding $\Ncal=2$ gauge theories on ALE spaces to correlators of vertex
operators of $\widehat{\slfrak}(k)_r$. To compute the instanton
partition functions, one needs some kind of smooth
completion of the moduli space of instantons. In the case of
instantons on $\R^4$ this completion can be obtained in three different but equivalent ways: 
by passing to a noncommutative deformation of $\R^4$
\cite{art:nekrasovschwarz1998}, by an algebro-geometric blowup, or by regarding the moduli space as a
hyper-K\"ahler quotient and perturbing the set of zeroes of the moment
map. In all three cases one obtains
a space that parameterizes framed torsion free sheaves on the
projective plane $\PP^2$. In the ALE case only the third technique is currently workable and yields Nakajima's quiver variety. To compute the instanton partition functions one can use the local description of quiver varieties around singular points (cf.\ \cite[Section 3]{art:nakajima2001} and \cite[Section 2.7]{art:nakajima2009})\footnote{We thank H. Nakajima for pointing this out.}. Another approach consists of realizing the quiver varieties as moduli spaces parameterizing framed torsion free sheaves. For example, there exists a chamber $C_0$ such that the associated Nakajima quiver varieties are isomorphic to moduli spaces of framed $\mathbb{Z}_k$-equivariant torsion free sheaves on $\mathbb{P}^2$ (see e.g. \cite[Section~2.3]{art:varagnolovasserot1999}). More generally, one can relate quiver varieties to framed sheaves on V-manifolds, as Nakajima did in \cite{art:nakajima2007}. In the present paper, we pursue a different approach relating quiver varieties to framed sheaves on Deligne-Mumford stacks (which form a category wider than the one of V-manifolds). The description of quiver varieties associated with $C_0$ in tems of $\mathbb{Z}_k$-equivariant framed sheaves on $\mathbb{P}^2$ provides a mathematical interpretation of the work \cite{art:fucitomoralespoghossian2004}, where the instanton partition functions of $\mathcal{N}=2$ supersymmetric gauge theories on  ALE spaces with parameters in the chamber $C_0$ were computed by considering a torus action (given by the torus $(\mathbb{C}^\ast)^2$ of $\mathbb{P}^2$ and the maximal torus of $GL(r, \mathbb{C})$) on the framed sheaves on $\mathbb{P}^2$ and taking into account only the fixed points that are invariant under the action of  $\mathbb{Z}_k$.

Recent developments in string theory have sparked new interest in the correspondence between four-dimensional
gauge theories and two-dimensional conformal field theories in view of the observation that the $\Ncal=2$ gauge theory on a four-manifold $X$ can also be
studied geometrically by embedding it in M-theory as the $(2,0)$
superconformal theory compactified on the corresponding Seiberg-Witten
curve $\Sigma$; this theory is the low-energy limit of the worldvolume
theory of a single M5-brane filling the six-dimensional manifold
$M_6=\Sigma\times X$. More generally, one can compactify on any
punctured Riemann surface $\Sigma'$ such that the Seiberg-Witten curve
$\Sigma$ is a branched cover of $\Sigma'$
\cite{art:losevmarshakovnekrasov2003,art:nekrasovokounkov2006}. The
chiral fields of the two-dimensional conformal field theory on
$\Sigma'$ which conjecturally describes the BPS sector of the
four-dimensional $\Omega$-deformed (i.e. equivariant) gauge theory arise as the zero
modes of the six-dimensional $(2,0)$ tensor multiplet. By applying
this machinery to the case $X={\mathbb{R}}^4$ one obtains the
correspondence conjectured
by Alday, Gaiotto and Tachikawa (AGT)
\cite{art:aldaygaiottotachikawa2010} (see also
\cite{art:wyllard2009}), which relates the Nekrasov
partition functions of $U(r)$ gauge theories
on $\mathbb{R}^4$ \cite{art:nekrasov2003} (see also \cite{art:flumepoghossian2003,art:bruzzofucitomoralestanzini2003}) to the conformal blocks of Toda conformal field
theories. From a mathematical perspective, the AGT correspondence
asserts a higher rank generalization of the celebrated result of
Nakajima \cite{art:nakajima1997} and Vasserot \cite{art:vasserot2001}
which relates the $(\C^\ast)^2$-equivariant cohomology of the Hilbert
schemes $\hilb{n}{\C^2}$ of $n$ points on $\C^2$ to representations of
the Heisenberg algebra over $H_{(\C^\ast)^2}^\ast(\C^2;\C)$. The higher
rank generalizations of $\hilb{n}{\C^2}$ are the moduli spaces
$\Mcal(r,n)$ of framed torsion free sheaves on $\PP^2$ of
rank $r$ and second Chern class $n$. The conjecture then supports the existence of a natural geometric action of the
$\mathcal{W}(\mathfrak{gl}_r)$-algebra on the direct sum over $n$ of
the equivariant cohomology groups of $\Mcal(r,n)$. This
correspondence was proven by Schiffmann and Vasserot
\cite{art:schiffmannvasserot2012}, and independently by Maulik and
Okounkov \cite{art:maulikokounkov2012}. The AGT correspondence allows
one to write the Nekrasov partition functions, defined as integrals over $\Mcal(r,n)$, as correlators
of vertex operators of the $\mathcal{W}(\mathfrak{gl}_r)$-algebra.

A generalization of the AGT correspondence to $\Ncal=2$ gauge theories
on ALE spaces could provide a better  gauge-theoretic explanation of 
Nakajima's work  than the Vafa-Witten theory. In
\cite{art:belavinfeigin2011,art:nishiokatachikawa2011,art:belavinbelavinbershtein2011}
the relevant   algebra is conjectured to be
\begin{equation}
\Acal(r,k)= \hfrak\oplus \widehat{\slfrak}(k)_r \oplus
\frac{\widehat{\slfrak}(r)_k \oplus
  \widehat{\slfrak}(r)_{\nu-k}}{\widehat{\slfrak}(r)_\nu} \ ,
\end{equation}
where $\hfrak$ is the infinite-dimensional Heisenberg algebra and the
parameter $\nu$ is related to the equivariant parameters of the torus
action. For $r=1$ this algebra is simply
$\Acal(1,k)=\hfrak\oplus\widehat{\slfrak}(k)_1$. For $k=2$ the algebra
$\Acal(r,2)$ is isomorphic to the sum of
$\widehat{\mathfrak{gl}}(2)_r$ acting at level $r$ and the
Neveu-Schwarz-Ramond algebra of supersymmetric Liouville theory. For
generic values of $(r,k)$, the algebra $\Acal(r,k)$ is the sum of the
affine Lie algebra $\widehat{\mathfrak{gl}}(k)_r$ and the
$\Z_k$-parafermionic $\mathcal{W}(\mathfrak{gl}_r)$-algebra. From this
perspective, the geometric representations of $\widehat{\slfrak}(k)_r$
constructed by Nakajima are only a part of the geometric
representations of $\Acal(r,k)$ that one expects to define on the
equivariant cohomology of the Nakajima quiver varieties. The
conjecture that has thus far emerged from field theory calculations
implies a new type of AGT correspondence which relates the equivariant
cohomology groups of moduli spaces of framed $\Z_k$-equivariant
torsion free sheaves on the projective plane $\PP^2$ and
representations of the algebra $\Acal(r,k)$; see e.g.\
\cite{art:bonellimaruyoshitanzini2011,art:bonellimaruyoshitanzini2012,art:wyllard2011,art:ito2011,art:belavinbershteinfeiginlitvinovtarnopolsky2011,art:bonellimaruyoshitanziniyagi2012,art:alfimovbelavintarnopolsky2013}.

There are three main outstanding issues that persist. Firstly, since the
Nakajima quiver varieties are all equivariantly diffeomorphic, but
isomorphic only when their stability parameters lie in the same
chamber, the computations of \cite{art:fucitomoralespoghossian2004}
provide instanton partition functions for $\Ncal=2$ gauge theories
with only adjoint matter hypermultiplets on \emph{any} ALE space
$X_k(\xi)$, but for $\Ncal=2$ gauge theories with fundamental matter
fields only on the ALE spaces with parameters in the chamber
$C_0$. Secondly, since the stability parameter of the surface $X_k$ does
\emph{not} belong to $C_0$ (see e.g. \cite{art:nagao2009}), for gauge
theories with only adjoint matter one expects a nontrivial
equivalence between the corresponding partition functions
\cite{art:belavinbershteinfeiginlitvinovtarnopolsky2011,art:bonellimaruyoshitanziniyagi2012,art:alfimovbelavintarnopolsky2013};
on the other hand, the relationship between the two partition
functions should be governed by a wall-crossing formula in the case of
gauge theories with fundamental matter
\cite{art:itomaruyoshiokuda2013}. Thirdly, one expects
\cite{art:bonellimaruyoshitanziniyagi2012} to find a blowup formula
for the instanton partition functions on ALE spaces in terms of
instanton partition functions on $\R^4$ depending on the equivariant
parameters of the torus action on the affine patches of $X_k$. This
generalizes the blowup formulas of \cite{art:nakajimayoshioka2005-I}
by taking into account that $X_k$ is obtained by a \emph{weighted blowup} at the singular point of $\C^2/\Z_k$.
In this paper, we address some of these problems. Instead of looking for possible smooth completions of the moduli spaces of instantons, we directly construct moduli spaces of framed sheaves which, as we will show, provide a suitable setting to compute instanton partition functions for gauge theories on $X_k$.

\subsection{Overview}

In this subsection we summarise the main problems we shall tackle in
some detail. The goal of this paper is to investigate a new geometrical approach to
the study of $\Ncal=2$ gauge theories on ALE spaces of type $A_{k-1}$
by finding a suitable compactification of $X_k$ on which to develop a
theory of framed sheaves. Our approach is the first attempt to use a
theory of framed sheaves in the study of gauge theories on ALE spaces
with stability parameters in the same chamber as the parameter of the
minimal resolution $X_k$. In particular, we aim to provide the
``correct" moduli space on which to construct natural geometric representations of the algebra $\mathcal{A}(r, k)$; this can be interpreted as a first step in the direction of  geometric realizations of representations of more complicated infinite-dimensional Lie algebras.

Note that any such compactification of a K\"ahler
surface $M$ is strongly constrained by a result of Bando
\cite{art:bando1993}: If $\bar M$ is a compactification of $M$ by a
smooth divisor with positive normal line bundle, which is a compact
K\"ahler surface, then holomorphic vector bundles on $\bar M$ that are
trivial on the compactification divisor $D=\bar M\setminus M$
correspond to holomorphic vector bundles on $M$ with anti-selfdual
square-integrable connections of trivial holonomy at infinity. Bando
also shows that the holonomy at infinity of the instantons should
correspond to a flat connection on the associated locally free sheaves
restricted to the compactification divisor; hence if $\bar M$ is obtained from $M$ by adding a projective line $D$, then only instantons on $M$ with trivial holonomy at infinity can be
described in terms of framed locally free sheaves on $\bar M$. When $M=X_k$, one possible way to avoid this restriction is to look for more general compactifications, that should properly allow
for the contributions of instantons with nontrivial holonomies at
infinity, which in gauge theory are sometimes called ``fractional'' instantons.

The first attempt at such a compactification of $X_k$ is due to
Nakajima \cite{art:nakajima2007}, who suggested a V-manifold
compactification; the compactification divisor $D$ carries a
$\Z_k$-action such that a framed torsion free sheaf restricted to $D$
is isomorphic to a $\Z_k$-equivariant locally free sheaf, which
should encode a fixed holonomy at infinity. Another approach was
pursued in \cite{art:bruzzopoghossiantanzini2011}, where torsion free
sheaves on Hirzebruch surfaces $\F_p$  framed along a divisor
$D_\infty$, with $D_\infty^2=p$, were used to compute the Vafa-Witten
partition functions of $\Ncal=4$ gauge theory on the total spaces
$\operatorname{Tot}(\Ocal_{\PP^1}(-p))$ of the line bundles
$\Ocal_{\PP^1}(-p)$, which can be realized as $\F_p\setminus
D_\infty$. The computations of \cite{art:bruzzopoghossiantanzini2011}
also make sense for \emph{fractional} Chern classes
$\crm_1\in\frac{1}{p}\,\Z$ which, although heuristic,
correctly incorporate the anticipated contributions to the Vafa-Witten
partition function from fractional instantons (see
e.g. \cite{art:fucitomoralespoghossian2006,art:griguoloseminaraszabotanzini2007});
hence the paper \cite{art:bruzzopoghossiantanzini2011} suggests that
the appropriate geometric arena for these computations should involve
a ``stacky'' compactification of
$\operatorname{Tot}(\Ocal_{\PP^1}(-p))$ which, as shown in  
\cite[Appendix~D]{art:bruzzosala2013}, is a root toric stack over $\F_p$. This reasoning was extended to compute the partition functions of $\Ncal=2$ gauge theories in \cite{art:bonellimaruyoshitanzini2011,art:bonellimaruyoshitanzini2012,art:bonellimaruyoshitanziniyagi2012}; other calculations of the fractional instanton contributions to supersymmetric gauge theory partition functions can be found in \cite{art:fucitomoralespoghossian2006,art:griguoloseminaraszabotanzini2007,art:ciraficikashanipoorszabo2011,art:szabo2009,art:ciraficiszabo2012}.
 
In the first part of this paper we address the problem of constructing
a suitable compactification of $X_k$. This is a projective toric
  orbifold $\Xscr_k$ over a  projective normal toric compactification
of $X_k$. The stacky structure of $\Xscr_k$ is concentrated at a
smooth effective Cartier divisor $\Dscr_\infty$, such that
$\Xscr_k\setminus \Dscr_\infty \simeq X_k$ and $\Dscr_\infty$ is a
$\Z_k$-gerbe over a \emph{football}, i.e., an orbifold curve over
$\PP^1$ with two orbifold points. One has $\pi_1(\Dscr_\infty^{\mathrm{top}})\simeq \Z_k$, where $\Dscr_\infty^{\mathrm{top}}$ is the underlying topological stack of $\Dscr_\infty$, so there exist $k$ complex line bundles endowed with flat connections associated with the $k$ irreducible unitary representations of $\Z_k$. Thus the locally free sheaf on $\Dscr_\infty$ which should encode the holonomy at infinity is the direct sum of these line bundles; we shall call it a \emph{framing sheaf}. It is associated with a homomorphism $\rho\colon \Z_k\to U(r)$.

By using the machinery developed in \cite{art:bruzzosala2013} we
construct a moduli space of torsion free sheaves on $\Xscr_k$ which
are isomorphic along $\Dscr_\infty$ to a given framing sheaf. By
defining a suitable toric action and characterizing the torus-fixed
point locus, we compute partition functions of gauge theories with
gauge groups $U(r)\times U(r'\,)$ on $X_k$, with and without matter in
the bifundamental representation; the extensions of our results to
general $A_n$-type quiver gauge theories  is carried out in \cite{art:bruzzosalaszabo2015}. By using this new geometrical approach, the
factorizations (blowup formulas)  are evident. In particular, we
will provide rigorous derivations of the partition functions for
$\Ncal=2$ gauge theories on $X_k$ which are conjecturally formulated
in \cite{art:bonellimaruyoshitanziniyagi2012} using heuristic
arguments. On the other hand, although our formulas can be shown to
match with those of \cite{art:bonellimaruyoshitanziniyagi2012} in
several nontrivial checks, our expressions are \emph{a priori}
different and we expect a nontrivial equivalence between the
corresponding partition functions. In this sense our partition functions may in fact yield a more transparent connection to conformal field theory in two dimensions under the AGT duality.

In addition to the partition functions, we also address the problem of
computing correlators of certain gauge invariant chiral BPS
operators. The topologically twisted $\Ncal=2$ gauge theory on $X_k$
has a natural set of holomorphic observables residing in the
nontrivial cohomology groups $H^p(X_k;\C)$
\cite{art:losevnekrasovshatashvili1997,art:marinomoore1998}. They are
labelled by an invariant polynomial $\Pcal$ on the Lie algebra
$\ufrak(r)$. The ``0-observables'' are then $\Osf^{(0)}:=\Pcal(\phi)$,
where $\phi$ is the complex scalar field of the $\Ncal=2$ vector
multiplet, while the ``2-observables'' $\Osf^{(2)}$ are obtained by
canonical descent equations from the invariant polynomial $\Pcal$;
there are similarly ``4-observables'', but their contributions can be
absorbed into the partition functions \cite{art:marinomoore1998}. The
generating functions for the correlators of nontrivial
$p$-observables are then of the schematic form
\begin{equation}
\Big\langle \exp\big(\Osf^{(0)}(P)+(\Osf^{(2)},S)\big)\Big\rangle_{X_k} \ ,
\end{equation}
where $P\in H_0(X_k;\C)$ and $S\in H_2(X_k;\C)$, while the braces
denote an expectation value in the quantum gauge theory on $X_k$. In
Donaldson-Witten theory, one takes $\Pcal$ to be linear in the
Casimir invariants of the gauge group; then the $p$-observables
correspond to single-trace chiral operators
$\Pcal_s(\phi)=\frac1{(2\pi\, {\rm i})^s\, s!}\,{\rm Tr}\, \phi^s$,
where $\operatorname{Tr}$ is the trace in the $r$-dimensional
representation, and in the $\Omega$-deformed gauge theory they are in
a bijective correspondence with characteristic classes
$\ch_s(\boldsymbol\Ecal)$ of the coherent rank $r$ universal framed
sheaf $\boldsymbol\Ecal$ associated to ``families'' of
$U(r)$ instantons parametrized by points of the moduli space of framed torsion free sheaves \cite{art:nekrasov2003,art:losevmarshakovnekrasov2003}. In particular, at the BRST fixed points of the topologically twisted gauge theory the scalar field $\phi$ is a certain $\ufrak(r)$-valued two-form on the moduli space, and the perturbation by arbitrary powers of holomorphic operators is then represented as
\begin{equation}
\Pcal(\phi) = \sum_{s=0}^\infty \, \frac{\tau_s}{s+1}\,
\ch_{s+1}(\boldsymbol\Ecal) \ ,
\end{equation}
where $\vec\tau=(\tau_0,\tau_1,\tau_2,\dots)$ is a collection of
complex parameters. In this case the generating functions for
$p$-observables are computed rigorously in~\cite{book:nakajimayoshioka2004} for gauge theories on $\R^4$ and on the blowup of $\R^4$ at the origin. 

Of particular interest is the case with $\Pcal(\phi)=\operatorname{Tr}\, \phi^2$, corresponding to observables obtained from the quadratic Casimir invariant; this corresponds to the choice
$\vec\tau=(0, -\tau_1,0,0,\dots)$ with $\Osf^{(0)}$ related to the
second Chern character of the universal sheaf $\boldsymbol\Ecal$. Then the
deformed partition functions can also be regarded as generating
functions for the equivariant cohomology version of the
Donaldson-Witten invariants of $X_k$. It is then natural to
investigate the structure of these correlators in the low energy limit where the
equivariant deformation ($\Omega$-background) is removed. In this
limit, we expect a relation between correlation functions on $X_k$ and
correlation functions on $X_k(\xi_0)$ analogous to the wall-crossing
formula that should relate the partition functions of the
corresponding $\Omega$-deformed gauge theories
\cite{art:itomaruyoshiokuda2013}; in the low energy limit, the exceptional set $\varphi_k^{-1}(0)$ of the minimal resolution $\varphi_k\colon X_k\to \C^2/\Z_k$ can be represented by an infinite number of local operators on $X_k(\xi_0)$. We therefore expect a blowup equation of the schematic form
\begin{equation}
\Big\langle \exp\big(\Osf^{(0)}(P)+t\,
(\Osf^{(2)},S)\big)\Big\rangle_{X_k} = 
\Big\langle\exp\big(\Osf^{(0)}\big(\varphi_k(P)\big) +{\sf
  F}[t|\Osf_2,\dots,\Osf_{r+1} ] \big)
\Big\rangle_{X_k(\xi_0)} \ , 
\end{equation}
where ${\sf F}$ is a holomorphic function and $\Osf_{\alpha}$,
$\alpha=2,\dots,r+1$ are the generators of the ring of local
BRST-invariant observables corresponding to the Casimir invariants of
$U(r)$; this formula generalizes the blowup equations of
\cite{art:losevnekrasovshatashvili1997,art:marinomoore1998} which
considered blowups of smooth points. Rigorous derivations of such
blowup equations for the Nekrasov partition functions were
obtained by Nakajima and Yoshioka in \cite{art:nakajimayoshioka2005-I,book:nakajimayoshioka2004} which
resemble the Fintushel-Stern formulas for Donaldson invariants
\cite{art:fintushelstern1996}; these equations were one of the
ingredients in the computations of Donaldson invariants from
Nekrasov partition functions in
\cite{art:gottschenakajimayoshioka2008,art:gottschenakajimayoshioka2011}. 
The low energy limit of $\Ncal=2$ supersymmetric Yang-Mills theory on
$X_k(\xi_0)$ is encoded in the Seiberg-Witten curve $\Sigma$ of genus
$r$; the blowup
equations are then expected
\cite{art:losevnekrasovshatashvili1997,art:marinomoore1998} to be
determined in terms of suitable $Sp(2r,\Z)$ modular forms. They thus unveil the modularity properties of the gauge
theory partition function and observables on $X_k$, in conjuction with
expectations from S-duality in gauge theory \cite{art:vafawitten1994}
and string theory.

\subsection{Summary of results}

In this article we use framed sheaves on the stacky compactification $\Xscr_k$ of $X_k$ to study supersymmetric gauge theories on $X_k$. In this subsection we summarise the main ingredients involved in this construction and our main results. 

\subsubsection{Stacky compactifications of ALE spaces}

The first problem we address in this paper consists of constructing a
stacky compactification of $X_k$. In Section
\ref{sec:stackycompactification} we describe in details a
compactification of $X_k$ given as a root toric stack; the theory of
root and toric stacks is recalled in Section
\ref{sec:preliminaries}. For this, we first compactify the ALE space
$X_k$ to a normal projective toric surface $\bar X_k$, with two
singular points of the same type, by adding a torus-invariant divisor
$D_\infty\simeq \PP^1$ such that for $k=2$ the surface $\bar X_2$
coincides with the second Hirzebruch surface $\F_2$. For $k\geq3$ the
surface $\bar X_k$ is singular, but one can associate with $\bar X_k$
its canonical toric stack $\Xscr_k^{\mathrm{can}}$ which is a
two-dimensional projective toric orbifold with Deligne-Mumford torus
$\C^\ast\times \C^\ast$ and coarse moduli space
$\pi_k^{\mathrm{can}}\colon \Xscr_k^{\mathrm{can}}\to \bar{X}_k$. By
``canonical'' we mean that the locus over which $\pi_k^{\mathrm{can}}$
is not an isomorphism has a nonpositive dimension; for $k=2$ one has
$\Xscr_2^{\mathrm{can}}\simeq \F_2$. Let us consider the
one-dimensional, torus-invariant, integral closed substack
$\tilde{\Dscr}_\infty:=(\pi_k^{\mathrm{can}})^{-1}(D_\infty)_{\mathrm{red}}\subset\Xscr_k^{\mathrm{can}}$. We
perform the $k$-th root construction on $\Xscr_k^{\mathrm{can}}$ along $\tilde{\Dscr}_\infty$ to extend the automorphism group of a generic point of $\tilde{\Dscr}_\infty$ by $\Z_k$; in this way we obtain a two-dimensional projective toric orbifold $\Xscr_k$ with Deligne-Mumford torus $\C^\ast\times \C^\ast$ and coarse moduli space $\pi_k\colon\Xscr_k\to\bar X_k$. The surface $X_k$ is isomorphic to the open subset $\Xscr_k \setminus\Dscr_\infty$ of $\Xscr_k$, where $\Dscr_\infty :=\pi_k^{-1}(D_\infty)_{\mathrm{red}}$. For $i=1,\ldots,k-1$ let $\Dscr_i:=\pi_k^{-1}(D_i)_{\mathrm{red}}$ be the divisors in $\Xscr_k$ corresponding to the exceptional divisors $D_i$ of the resolution of singularities $\varphi_k\colon X_k \to \C^2/\Z_k$; their intersection product is given by $-C$, where $C$ is the Cartan matrix of type $A_{k-1}$. Define the \emph{dual} classes
\begin{equation}
\omega_i := - \sum_{j=1}^{k-1}\, \big(C^{-1}\big)^{ij} \, \Dscr_j\ .
\end{equation}
We prove that these classes are integral. Let us denote by
$\Rcal_i := \Ocal_{\Xscr_k}(\omega_i)$ their associated line bundles
on $\Xscr_k$. Their restrictions to $X_k$ are precisely the
tautological line bundles of Kronheimer and Nakajima.
\begin{proposition*}[{Proposition \ref{prop:picardstack}}]
The Picard group $\Pic(\Xscr_k)$ of $\Xscr_k$ is freely generated over $\Z$ by $\Ocal_{\Xscr_k}(\Dscr_\infty)$ and $\Rcal_i$ with $i=1, \ldots, k-1$.
\end{proposition*}

We further provide a characterization of the divisor $\Dscr_\infty$ as a toric Deligne-Mumford stack with Deligne-Mumford torus $\C^\ast\times \Bscr\Z_k$ and coarse moduli space $r_k\colon \Dscr_\infty\to D_\infty$.
\begin{proposition*}[{Proposition \ref{prop:gerbestructure}}]
$\Dscr_\infty$ is isomorphic as a toric Deligne-Mumford stack to the toric global quotient stack
\begin{equation}
\left[\frac{\C^2\setminus\{0\}}{\C^\ast\times\Z_k}\right]\ ,
\end{equation}
where the group action is given in Equation \eqref{eq:Tactiongerbe}.
\end{proposition*}
Since line bundles on a global quotient stack $[X/G]$, with trivial $\operatorname{Pic}(X)$, are associated
with characters of the group $G$,  we find that the Picard
group $\Pic(\Dscr_\infty)$ is isomorphic to $\Z\oplus \Z_k$; it is
generated by the line bundles $\Lcal_1, \Lcal_2$ corresponding to the
characters $\chi_1,\chi_2\colon\C^\ast\times\Z_k \to \C^\ast$ given
respectively by the projections $(t,\omega)\mapsto t$ and
$(t,\omega)\mapsto\omega$, where $t\in\C^\ast$ and $\omega$ is a
primitive $k$-th root of unity. In particular, $\Lcal_2^{\otimes k}$
is trivial. As pointed out by \cite{art:eyssidieuxsala2013}, the
fundamental group of the underlying topological stack
$\Dscr_\infty^{\mathrm{top}}$ of $\Dscr_\infty$ is isomorphic to
$\Z_k$; in addition, for any $i=0, 1, \ldots, k-1$ the line bundle $\Lcal_2^{\otimes i}$ inherits a Hermitian metric and a unitary flat connection associated with the $i$-th irreducible representation of $\Z_k$.

\subsubsection{Moduli spaces of framed sheaves}

In order to construct moduli spaces of framed sheaves on $\Xscr_k$
which are needed for the formulation of supersymmetric gauge theories
on $X_k$, we first have to choose a suitable framing sheaf which
should encode the fixed holonomy at infinity of instantons. Since
the holonomy at infinity corresponds to a  representation of
$\Z_k$, the framing sheaf should have a Hermitian metric and a unitary flat connection associated with the   representation of $\Z_k$. Because of this, we choose as framing sheaf the locally free sheaf
\begin{equation}
\Fcal_{\infty}^{0, \vec w}:=\bigoplus_{i=0}^{k-1} \,
\big(\Ocal_{\Dscr_\infty}(0,i) \big)^{\oplus w_i}\ ,
\end{equation}
where $\Ocal_{\Dscr_\infty}(0,i)$ is $\Lcal_2^{\otimes i}$ for even
$k$ and $\Lcal_2^{\otimes i\, \frac{k+1}{2}}$ for odd $k$, and $\vec{w}:=(w_0,w_1, \ldots, w_{k-1})\in \mathbb{N}^{k}$ is a fixed vector. If we tensor $\Fcal_\infty^{0, \vec{w}}$ by a power $\Lcal_1^{\otimes s}$ of $\Lcal_1$, we obtain a more general framing sheaf $\Fcal_\infty^{s, \vec{w}}$; the degree of $\Fcal_\infty^{s, \vec{w}}$ is a rational multiple of $s$.

In Section \ref{sec:framedsheaves}, we apply the general theory of
framed sheaves on projective stacks developed in
\cite{art:bruzzosala2013} to construct a fine moduli space
$\Mcal_{r,\vec{u},\Delta}(\Xscr_k,\Dscr_\infty,\Fcal_\infty^{0,
  \vec{w}}\, )$ parameterizing $(\Dscr_\infty,\Fcal_\infty^{0,
  \vec{w}}\, )$-framed sheaves
$(\Ecal,\phi_\Ecal\colon \Ecal_{\vert \Dscr_\infty}\xrightarrow{\sim}
\Fcal_\infty^{0, \vec{w}}\, )$ on  $\Xscr_k$ with fixed rank
$r:=\sum_{i=0}^{k-1}\, w_i$, first Chern class
$\crm_1(\Ecal)=\sum_{i=1}^{k-1} \, u_i \, \omega_i$ and discriminant $\Delta(\Ecal)=\Delta$. The vector $\vec{u}=(u_1,\dots,u_{k-1})$ satisfies the constraint
\begin{equation}\label{eq:firstchernclass-intro}
k \, v_{k-1} = \sum_{i=0}^{k-1}\, i \, w_i \ \bmod{k}\ ,
\end{equation}
where $\vec{v}:=C^{-1}\vec{u}$
.
\begin{theorem*}[{Theorem \ref{thm:moduli}}]
The moduli space $\Mcal_{r,\vec{u},\Delta}(\Xscr_k,\Dscr_\infty,\Fcal_\infty^{0, \vec{w}})$ is a smooth quasi-projective variety of dimension 
\begin{equation}
\dim_{\C}\big(\Mcal_{r,\vec{u},\Delta}(\Xscr_k,\Dscr_\infty,\Fcal_\infty^{0,
  \vec{w}}\, )\big)=2r\, \Delta - \frac12\, \sum_{j=1}^{k-1}\, \big(C^{-1}\big)^{jj}\, \vec{w}\cdot\vec{w}(j)\ ,
\end{equation}
where for $j=1, \ldots, k-1$ the vector $\vec{w}(j)$ is $(w_j, \ldots, w_{k-1}, w_0, w_1, \ldots, w_{j-1})$. Moreover, the Zariski tangent
space of
$\Mcal_{r,\vec{u},\Delta}(\Xscr_k,\Dscr_\infty,\Fcal_\infty^{0,\vec{w}}\,
)$ at a point $[(\Ecal,\phi_{\Ecal})]$ is $\mathrm{Ext}^1(\Ecal,\Ecal\otimes \Ocal_{\Xscr_k}(-\Dscr_\infty))$.
\end{theorem*}
In the rank one case $r=1$, we show that there is an isomorphism of fine moduli spaces
\begin{equation}
\imath_{1,\vec u,n} \, \colon \, \hilb{n}{X_k} \ \xrightarrow{ \ \sim\ } \
\Mcal_{1,\vec{u},n}\big(\Xscr_k,\Dscr_\infty,\Fcal_\infty^{0, \vec{w}} \,
\big)\ ,
\end{equation}
where $\hilb{n}{X_k}$ is the Hilbert scheme of $n$ points on $X_k$.

In gauge theory, the tangent bundle $T
\Mcal_{r,\vec{u},\Delta}(\Xscr_k,\Dscr_\infty,\Fcal_\infty^{0,
  \vec{w}}\, )$ describes gauge fields and matter fields in the adjoint representation of the gauge group $U(r)$. Matter in the fundamental representation of the gauge group $U(r)$ is described by the coherent sheaf 
 \begin{equation}
\Vbf:=R^1 p_\ast\big({\boldsymbol\Ecal} \otimes
p_{\Xscr_k}^\ast(\Ocal_{\Xscr_k}(-\Dscr_\infty) )\big) \ ,
\end{equation}
where ${\boldsymbol\Ecal}$ is the universal sheaf of
$\Mcal_{r,\vec{u},\Delta}(\Xscr_k,\Dscr_\infty,\Fcal_\infty^{0,
  \vec{w}}\, )$, while $p$ and $p_{\Xscr_k}$ respectively denote the projections of ${\Mcal_{r,\vec{u},\Delta}(\Xscr_k,\Dscr_\infty,\Fcal_\infty^{0, \vec
    w}\, )} \times \mathscr X_k $ onto ${\Mcal_{r,\vec{u},\Delta}(\Xscr_k,\Dscr_\infty,\Fcal_\infty^{0, \vec
    w}\, )}$ and $\Xscr_k $. We
call $\Vbf$ the natural bundle of
$\Mcal_{r,\vec{u},\Delta}(\Xscr_k,\Dscr_\infty,\Fcal_\infty^{0, \vec
  w}\, )$.
\begin{proposition*}[{Proposition \ref{prop:naturalbundle}}]
$\Vbf$ is a locally free sheaf on
$\Mcal_{r,\vec{u},\Delta}(\Xscr_k,\Dscr_\infty,\Fcal_\infty^{0, \vec
  w}\, )$ of rank 
\begin{equation}
\rk(\Vbf) = \Delta +\frac{1}{2r}\, \vec v\cdot C\vec v - \frac{1}{2}\, \sum_{j=1}^{k-1}\, \big(C^{-1} \big)^{jj}\, w_j
\end{equation}
where $\vec{v}:=C^{-1} \vec{u}$.
\end{proposition*}
The computation of the rank of $\Vbf$ and the dimension of the moduli
spaces of framed sheaves is addressed in Appendix
\ref{sec:rankdimension}, where we use the T\"oen-Riemann-Roch theorem
and some   summation identities for complex roots of unity (derived in Appendix \ref{ap:rootsunity}) to obtain the explicit formulas.

There is a natural generalization of the vector bundles $T
\Mcal_{r,\vec{u},\Delta}(\Xscr_k,\Dscr_\infty,\Fcal_\infty^{0,
  \vec{w}}\, )$ and $\Vbf$ to a virtual vector bundle $\Ebf$ on the
product
$\Mcal_{r,\vec{u},\Delta}(\Xscr_k,\Dscr_\infty,\Fcal_\infty^{0, \vec
  w}\, )\times\Mcal_{r',\vec{u}\,',\Delta'}
(\Xscr_k,\Dscr_\infty,\Fcal_\infty^{0, \vec{w}'}\, )$ with fibre over a point $\big([(\Ecal,\phi_\Ecal)]\, ,\,
[(\Ecal',\phi_{\Ecal'})]\big)$ given by
\begin{equation}
\Ebf_{\left([(\Ecal,\phi_\Ecal)]\,,\, [(\Ecal',\phi_{\Ecal'})]\right)}=\Ext^1\big(\Ecal\,,\,\Ecal'\otimes \Ocal_{\Xscr_k}(-\Dscr_\infty) \big)\ ;
\end{equation}
it can be regarded as a higher rank generalization of the Carlsson-Okounkov bundle \cite{art:carlssonokounkov2012}. In gauge theory, the bundle $\Ebf$ describes matter in the bifundamental representation of $U(r)\times U(r'\,)$.

Let us now explain how these moduli spaces and their natural vector
bundles are used to compute partition functions of $\Ncal=2$ gauge
theories on $X_k$. The $\Ncal=2$ gauge theory on a four-dimensional
toric manifold $X$ in the $\Omega$-background is obtained as the
reduction of a six-dimensional $\Ncal=1$ gauge theory on a flat
$X$-bundle $M$ over $\T^2$ in the limit where the torus $\T^2$
collapses to a point~\cite[Section~3.1]{art:nekrasov2006}. The bundle $M$ can be realized as the quotient of $\C\times X$ by the $\Z^2$-action
\begin{equation}
(n_1,n_2)\triangleright (w,x)=\big(w+(n_1+\sigma \, n_2)\,,\, g_1^{n_1} \,
g_2^{n_2} (x) \big)\ ,
\end{equation}
where $x\in X$, $w\in\C$, $(n_1,n_2)\in\Z^2$, $g_1,g_2$ are two commuting isometries of $X$ and $\sigma$ is the complex structure modulus of $\T^2$. In the collapsing limit, fields of the gauge theory which are charged under the
R-symmetry group are sections of the pullback to $M$ of a flat
$T_t$-bundle over $\T^2$, where $T_t\subset GL(2,\C)$ is the torus of
$X$. As pointed out in \cite[Section~2.2.2]{art:nekrasovokounkov2006},
the chiral observables of the $\Ncal=2$ gauge theory in the
$\Omega$-background become closed forms on the moduli spaces of framed
instantons which are equivariant with respect to the action of the
torus $T:=T_t\times T_\rho$, where $T_\rho$ is the maximal torus of
the group $GL(r, \C)$ of constant gauge transformations which consists
of diagonal matrices. Thus correlation functions of chiral
BPS operators become integrals of equivariant classes over
the moduli spaces.

In our setting we therefore define gauge theory partition functions as
generating functions for $T$-equivariant integrals of suitable
characteristic classes over the moduli spaces
$\Mcal_{r,\vec{u},\Delta}(\Xscr_k,\Dscr_\infty,\Fcal_\infty^{0, \vec
  w}\, )$; the precise choice of cohomology classes depends on the
matter content of the gauge theory and on the chiral observables in
question. There is a natural $T$-action on
$\Mcal_{r,\vec{u},\Delta}(\Xscr_k,\Dscr_\infty,\Fcal_\infty^{0, \vec
  w}\, )$ given, on a point $[(\Ecal, \phi_{\Ecal})]$, by pullback
of $\Ecal$ via the natural automorphism of $\Xscr_k$ induced by an element
of $T_t$, and by ``rotation" of the framing $\phi_{\Ecal}$ by a diagonal
matrix $\rho$ of $T_\rho$. The $T$-fixed point locus of
$\Mcal_{r,\vec{u},\Delta}(\Xscr_k,\Dscr_\infty,\Fcal_\infty^{0,\vec
  w}\, )$ consists of a finite number of isolated points. 
\begin{proposition*}[{Proposition \ref{prop:fixedpoint}}]
A $T$-fixed point $[(\Ecal,\phi_\Ecal)]\in
\Mcal_{r,\vec{u},\Delta}(\Xscr_k,\Dscr_\infty,\Fcal_\infty^{0,\vec
  w}\, )^T$ decomposes as a direct sum of rank one framed sheaves
\begin{equation}
(\Ecal, \phi_{\Ecal})=\bigoplus_{\alpha=1}^r \, (\Ecal_\alpha, \phi_{\alpha})\ ,
\end{equation}
where for $i=0, 1, \ldots, k-1$ and $\sum_{j=0}^{i-1}\, w_j<\alpha\leq\sum_{j=0}^i\, w_j$:
\begin{itemize}
\item $\Ecal_\alpha$ is a tensor product
  $\imath_\ast(I_\alpha)\otimes\big(\bigotimes_{j=1}^{k-1}\,
  \Rcal_j^{\otimes (\vec{u}_\alpha)_j}\big)$, where $I_\alpha$ is an
  ideal sheaf of a $0$-dimensional subscheme $Z_\alpha$ of $X_k$ with
  length $n_\alpha$
  supported at the $T_t$-fixed points $p_1,\ldots,p_k$, while
  $\vec{u}_\alpha\in\Z^{k-1}$ obeys $\sum_{\alpha=1}^r\, \vec
  u_\alpha=\vec u$ and $\vec{v}_\alpha:=C^{-1}\vec{u}_\alpha$ satisfies
\begin{equation}
k\, (\vec{v}_\alpha)_{k-1}=  i\, \bmod{k} \ ;
\end{equation}
\item ${\phi_\alpha}\colon \Ecal_\alpha
  \xrightarrow{\sim}\Ocal_{\Dscr_\infty}(0,i)$ is induced by the
  canonical isomorphism $\bigotimes_{j=1}^{k-1}\, {\Rcal_j}^{\otimes
    (\vec{u}_\alpha)_j}_{\vert \Dscr_\infty}\simeq
  \Ocal_{\Dscr_\infty}(0,i)$;
\item $\displaystyle{ \Delta =\sum_{\alpha=1}^r \, n_\alpha+\frac{1}{2}\,
\sum_{\alpha=1}^r\, \vec{v}_\alpha\cdot C\vec{v}_\alpha-\frac{1}{2r}\,
\sum_{\alpha,\beta=1}^r \, \vec{v}_\alpha\cdot C\vec{v}_\beta \ \in \
\mbox{$\frac{1}{2\, r\, k}$} \, \mathbb{Z}} \ . $
\end{itemize}
\end{proposition*}
In the standard way, we can associate to the $T_t$-invariant ideal sheaf
$\imath_\ast(I_\alpha)$ a set of Young tableaux
$\vec{Y}_\alpha:=\{Y_\alpha^i\}_{i=1, \ldots, k}$. Hence to each fixed
point $[(\Ecal,\phi_\Ecal)]$ there corresponds a combinatorial datum $(\vec{\boldsymbol{Y}}, \vec{\boldsymbol{u}})$, where $\vec{\boldsymbol{Y}}:=(\vec{Y}_1, \ldots, \vec{Y}_r)$ and $\vec{\boldsymbol{u}}:=(\vec{u}_1, \ldots, \vec{u}_r)$.

\subsubsection{Gauge theory partition functions}

In Section \ref{sec:instantonpartitionfunctions} we compute partition functions and correlators of chiral observables for $\Ncal=2$ gauge theories on $X_k$, with and without bifundamental matter fields. Here we summarise the pertinent results only for the gauge theory with a single adjoint hypermultiplet
of mass $\mu$, i.e., the $\Ncal=2^\ast$ gauge theory; the pure $\Ncal=2$ theory is recovered formally in the limit $\mu\to\infty$, while the limit $\mu=0$ has enhanced maximal supersymmetry and reduces to the $\Ncal=4$ Vafa-Witten theory. This is the case that is conjectured to compute the correlation functions of chiral operators in a dual two-dimensional conformal field theory on a torus, obtained geometrically by a twisted compactification of the $(2,0)$ theory with defects on an elliptic curve. The more general cases are treated in the main text (see Sections \ref{sec:pure} and \ref{sec:fundamentalmatter}). 

Let $\varepsilon_1, \varepsilon_2, a_1, \ldots, a_r$ be the generators
of $H_T^\ast({\rm pt};\Q)$; in gauge theory on $X=\R^4$,
$a_1,\ldots,a_r$ are the expectation values of the complex scalar field $\phi$
of the $\Ncal=2$ vector multiplet and $\varepsilon_1, \varepsilon_2$
parameterize the holonomy of a flat connection on the $T_t$-bundle over
$\T^2$ used to define the $\Omega$-background.
For $i=1, \ldots, k$ define
\begin{equation}
\varepsilon_1^{(i)}:=(k-i+1)\, \varepsilon_1-(i-1)\, \varepsilon_2
\qquad \mbox{and} \qquad
\varepsilon_2^{(i)}:= -(k-i)\, \varepsilon_1+i\, \varepsilon_2\ ,
\end{equation}
and
\begin{equation}
\vec{a}^{(i)}:=\vec{a}+\varepsilon_1^{(i)}\, (\vec{v}\,)_i+\varepsilon_2^{(i)}\, (\vec{v}\,)_{i-1} \ ,
\end{equation}
where $(\vec{v}\,)_l:=((\vec{v}_1)_l, \ldots, (\vec{v}_{r})_l)$ for
$l=1, \ldots, k-1$ and $(\vec{v}\, )_0=(\vec{v}\, )_k=\vec{0}$. Roughly speaking, for fixed $i=1, \ldots, k$ the parameters
$\varepsilon_1^{(i)}, \varepsilon_2^{(i)}$ are related to the
coordinates of the affine toric neighbourhood of the fixed point $p_i$
of $X_k$ and the shift $\vec{a}^{(i)}$ of $\vec{a}$ is due to the
equivariant Chern character of the line bundles ${\mathcal
R}_j$, with $j=1, \ldots, k-1$, restricted to the open affine toric neighbourhood $U_i$ of $p_i$.

For any $T$-equivariant vector bundle $\boldsymbol G$ of rank $d$ on
$\Mcal_{r,\vec{u},\Delta}(\Xscr_k,\Dscr_\infty,\Fcal_\infty^{0,\vec
  w}\, )$, define the
$T$-equivariant characteristic class
\begin{equation}
{\rm E}_\mu(\boldsymbol G) := \sum_{j=0}^d \, (\crm_j)_T (\boldsymbol G) \, \mu^{d-j}\ .
\end{equation}
Let $\vec{v}\in \frac{1}{k}\, \mathbb{Z}^{k-1}$ such that $k\, v_{k-1}= \sum_{i=0}^{k-1}\, i\, w_i\ \mathrm{mod}\ k$. Then the generating function for $p$-observables (with $p=0,2$) of the gauge theory with massive adjoint matter on $X_k$, in the topological sector labelled by $\vec v=C^{-1}\vec u$, is defined by
\begin{multline}
\Zcal^\ast_{\vec{v}}\big(\varepsilon_1, \varepsilon_2, \vec{a}, \mu ; \qsf, \vec{\tau}, \vec{t}^{\;(1)}, \ldots,
\vec{t}^{\;(k-1)}\big) \\
\shoveleft{:=\sum_{\Delta\in \frac{1}{2\,r\, k}\mathbb{Z}}\,
  \qsf^{\Delta+\frac{1}{2r}\, \vec{v}\cdot C\vec{v}} \
  \int_{\Mcal_{r,\vec{u},\Delta}(\Xscr_k,\Dscr_\infty,\Fcal_\infty^{0,
      \vec w}\, )} \, {\rm
    E}_\mu\big(T\Mcal_{r,\vec{u},\Delta}(\Xscr_k,\Dscr_\infty,\Fcal_\infty^{0,
    \vec w}\, ) \big)} \\
\times \ \exp \bigg(\, \sum_{s=0}^\infty \ \Big(\, \sum_{i=1}^{k-1}\, t_s^{(i)} \, \big[\mathrm{ch}_{T}(\boldsymbol{\mathcal{E}})/[\mathscr{D}_i]\big]_s+\tau_s\, \big[\mathrm{ch}_T(\boldsymbol{\mathcal{E}})/[X_k]\big]_{s-1}\, \Big)\,\bigg) \ ,
\end{multline}
where $\tau_0:=\frac1{2\pi\,{\rm i}}\, \log\qsf$ is the \emph{bare
  complex gauge coupling}. Here $\boldsymbol{\mathcal{E}}$ is the
universal sheaf of $\Mcal_{r,\vec{u},\Delta}(\Xscr_k,$
$\Dscr_\infty,\Fcal_\infty^{0,\vec w}\, )$, while
$\ch_{T}(\boldsymbol{\mathcal{E}})/[\mathscr{D}_i]$ denotes the slant
product between $\ch_{T}(\boldsymbol{\mathcal{E}})$ and
$[\mathscr{D}_i]$, and the notation $[-]_p$ means to take the degree
$p$ part of the $T$-equivariant class; since $X_k$ is noncompact, the
class $\ch_T(\boldsymbol{\mathcal{E}})/[X_k]$ is defined   by
localization. In particular, setting $\vec\tau=\vec 0$ and
$\vec{t}^{\;(1)}= \cdots = \vec{t}^{\;(k-1)}= \vec 0$ we obtain the instanton partition function $\Zcal_{\vec{v}}^{\ast,\mathrm{inst}}$ in the topological sector labelled by $\vec v=C^{-1}\vec u$. By taking a weighted sum over all $\vec v$ we obtain the full generating function for $p$-observables for $\Ncal=2^\ast$ gauge theory on $X_k$ as 
\begin{multline}
\Zcal^{\ast}_{X_k}\big(\varepsilon_1, \varepsilon_2, \vec{a}, \mu;
\qsf, \vec{\xi},\vec{\tau}, \vec{t}^{\;(1)}, \ldots,
\vec{t}^{\;(k-1)}\big) \\
:= \sum_{\stackrel{\scriptstyle \vec{v}\in\frac{1}{k}\,
    \mathbb{Z}^{k-1} }{\scriptstyle k\, v_{k-1}= \sum_{i=0}^{k-1}\,
    i\, w_i \bmod{k}}} \, \vec{\xi}^{\ \vec{v}}\ \Zcal_{\vec{v}}^{\ast, \mathrm{inst}}(\varepsilon_1, \varepsilon_2, \vec{a}, \mu; \qsf, \vec{\tau}, \vec{t}^{\;(1)}, \ldots, \vec{t}^{\;(k-1)}\big)
\end{multline}
and the full instanton partition function as
\begin{equation}
\Zcal^{\ast, \mathrm{inst}}_{X_k}\big(\varepsilon_1, \varepsilon_2,
\vec{a}, \mu; \qsf, \vec{\xi}\ \big):=
\sum_{\stackrel{\scriptstyle\vec{v}\in\frac{1}{k}\, \Z^{k-1}
  }{\scriptstyle k\, v_{k-1}= \sum_{i=0}^{k-1}\, i\, w_i \bmod{k}}} \, \vec{\xi}^{\ \vec{v}}\ \Zcal_{\vec{v}}^{\ast,\mathrm{inst}}(\varepsilon_1, \varepsilon_2, \vec{a}, \mu; \qsf)\ ,
\end{equation}
where $\log\xi_i$ for $i=1,\dots,k-1$ are \emph{chemical potentials
  for the fractional instantons} and $\vec\xi^{\ \vec
  v}:=\prod_{i=1}^{k-1}\, \xi_i^{v_i}$. 

By using the localization theorem, we provide an explicit formula for
the generating function $\Zcal^{\ast}_{X_k}$ in Proposition
\ref{prop:deformed-adjoint}. It factorizes with respect to the
generating functions for $p$-observables of $\Ncal=2^\ast$ gauge
theory on the affine toric open subsets of $X_k$ only for $k=2$; for $k\geq
3$ the apparent lack of a factorization property for $\Zcal^{\ast}_{X_k}$ is
due to the fact that it involves terms which depend on pairs of
exceptional divisors intersecting at the fixed points of $X_k$, and
such terms do not split into terms each depending on a single affine
toric subset of $X_k$. This is a new and unexpected result which
extends previous generating functions for $p$-observables that depend only on one exceptional divisor \cite[Section~4]{book:nakajimayoshioka2004}.

For the instanton partition function we obtain the factorization formula
\begin{multline}
\Zcal^{\ast, \mathrm{inst}}_{X_k}\big(\varepsilon_1, \varepsilon_2, \vec{a}, \mu; \qsf, \vec{\xi}\ \big)\\
\shoveleft{= \sum_{\stackrel{\scriptstyle \vec{v}\in\frac{1}{k}\,
      \Z^{k-1}}{\scriptstyle k\, v_{k-1}= \sum_{i=0}^{k-1}\, i\, w_i
      \bmod{k}}} \, \vec{\xi}^{\ \vec{v}}\
  \sum_{\vec{\boldsymbol{v}}} \,\qsf^{\frac{1}{2}\,
    \sum\limits_{\alpha=1}^r \, \vec{v}_\alpha\cdot C\vec{v}_\alpha} \
  \prod_{\alpha,\beta=1}^r\ \prod_{n=1}^{k-1} \,
  \frac{\ell^{(n)}_{\vec{v}_{\beta\alpha}}\big(\varepsilon_1^{(n)},
    \varepsilon_2^{(n)},  a_{\beta\alpha}+ \mu \big)}{\ell^{(n)}_{\vec{v}_{\beta\alpha}}\big(\varepsilon_1^{(n)},  \varepsilon_2^{(n)}, a_{\beta\alpha}\big)}} \\ 
\times \ \prod_{i=1}^k\, \Zcal_{\C^2}^{\ast, \mathrm{inst}}\big(\varepsilon_1^{(i)}, \varepsilon_2^{(i)}, \vec{a}^{(i)}, \mu; \mathsf{q} \big)\ ,
\end{multline}
where $\vec{\boldsymbol{v}}=(\vec v_1,\dots,\vec v_r)$, $\vec
v_{\beta\alpha}=\vec v_\beta-\vec v_\alpha$ and
$a_{\beta\alpha}=a_\beta-a_\alpha$, while
$\Zcal_{\C^2}^{\ast,\mathrm{inst}}$ is the Nekrasov partition function
for the $\Ncal=2^\ast$ gauge theory on $\R^4$ with gauge group
$U(r)$. Here $\ell_{\vec v_{\beta\alpha}}^{(n)}$ are the \emph{edge contributions}; whilst their explicit expressions are somewhat complicated, we point out that they depend on the Cartan matrix (see Section \ref{sec:Eulerclasses} for the complete formulas and Appendix \ref{sec:edgecontribution} for their derivations). The analogous results in the case that fundamental hypermultiplets are included can be found in Section \ref{sec:fundamentalmatter}, while the pertinent formulas in the generic case with $U(r)\times U(r'\,)$ bifundamental matter is presented in Section \ref{sec:Eulerclasses}. 

Let us compare these results with those obtained in \cite{art:bonellimaruyoshitanziniyagi2012} where, based on a conjectural splitting of the full partition function on $X_k$ as a product of partition functions on the affine toric open subsets of $X_k$, the authors obtain an expression for the edge contributions which depends only on the combinatorial data of the fan of the toric variety $X_k$. The two sets of expressions appear to admit drastically different structures, and a general proof of their equivalence is not immediately evident. In the $k=2$ case, our instanton partition function assumes the form
\begin{multline}
\Zcal^{\ast, \mathrm{inst}}_{X_k}(\varepsilon_1, \varepsilon_2,
\vec{a}, \mu ; \qsf, \xi) \\
\shoveleft{=\sum_{\stackrel{\scriptstyle v\in\frac{1}{2}\, \Z }{\scriptstyle 2\,
    v= \, w_1 \bmod{2}}} \, \xi^{v} \ \sum_{\vec{\boldsymbol{v}}} \,\qsf^{\sum\limits_{\alpha=1}^r \, v_\alpha^2}\ \prod_{\alpha,\beta=1}^r\,\frac{\ell_{v_{\beta\alpha}}\big(2\varepsilon_1, \varepsilon_2 - \varepsilon_1, a_{\beta\alpha}+\mu \big)}{\ell_{v_{\beta\alpha}}\big(2\varepsilon_1, \varepsilon_2 - \varepsilon_1, a_{\beta\alpha} \big)}}\\  
  \times \ \Zcal_{\C^2}^{\ast, \mathrm{inst}}\big(2\varepsilon_1,
  \varepsilon_2 -\varepsilon_1, \vec{a}-2\varepsilon_1\,
  \vec{\boldsymbol{v}} , \mu ; \qsf\big) \
  \Zcal_{\C^2}^{\ast,\mathrm{inst}}\big(\varepsilon_1 - \varepsilon_2,
  2\varepsilon_2,\vec{a}-2\varepsilon_2\, \vec{\boldsymbol{v}} , \mu ; \qsf \big)\ ,
\end{multline}
where
\small
\begin{equation*}
\ell_{v_{\beta\alpha}}(2\varepsilon_1,\varepsilon_2 - \varepsilon_1, a_{\beta\alpha})=\left\{
\begin{array}{cl}
\prod_{i=0}^{\lfloor v_{\beta\alpha}\rfloor-1}\limits\hspace{0.2cm}
\prod_{j=0}^{2i+2\{v_{\beta\alpha}\}}\limits\,
\big(a_{\beta\alpha}+2i\, \varepsilon_1+j\, (\varepsilon_2-\varepsilon_1)\big) & \mbox{ for } \lfloor v_{\beta\alpha}\rfloor> 0\ ,\\[8pt]
1 & \mbox{ for } \lfloor v_{\beta\alpha}\rfloor= 0\ ,\\[8pt]
\prod_{i=1}^{-\lfloor v_{\beta\alpha}\rfloor}\limits\hspace{0.2cm}
\prod_{j=1}^{2i-2\{v_{\beta\alpha}\}-1}\limits \,
\big(a_{\beta\alpha}+2(2\{v_{\beta\alpha}\}-i)\, \varepsilon_1-j\,
(\varepsilon_2-\varepsilon_1) \big)  & \mbox{ for } \lfloor v_{\beta\alpha}\rfloor< 0 \ .
\end{array}
\right.
\end{equation*}
\normalsize
Since the computation of the edge contributions in this case is
equivalent to that of
\cite[Section~4.2]{art:bruzzopoghossiantanzini2011}, our formula
agrees with \cite[Equation~(6)]{art:bonellimaruyoshitanzini2011} and
\cite[Equation~(2.13)]{art:bonellimaruyoshitanzini2012} (which are
equivalent to \cite[Equation~(3.24)]{art:bonellimaruyoshitanziniyagi2012} for $k=2$). In \ref{app:k=3} we further check that, for $k=3$ and rank $r=2$, our
partition functions agree with those in
\cite[Appendix C]{art:bonellimaruyoshitanziniyagi2012} at leading
orders in the
$\mathsf{q}$-expansion for some choices of the holonomy at infinity
and the
first Chern class of the framed sheaves. The dictionary between our
notation and that of \cite{art:bonellimaruyoshitanziniyagi2012} is as
follows. Their $\vec{k}_\alpha$ are our $\vec{v}_\alpha$. On
the other hand, to fix the holonomy at infinity they use a vector
$\vec{I}=(I_1,\dots,I_r)$ with $I_\alpha\in\{0, 1, \ldots, k-1\}$ and
then their $\vec{k}_\alpha$ satisfy an equation depending on
$\vec{I}$: for $\vec{I}=(k-1, \ldots, k-1, k-2, \ldots, k-2, \ldots, 0,\ldots, 0)$,
where $k-i-1$ appears with multiplicity $w_i$ for $i=0, 1,\ldots, k-1$,
their constraint is equivalent to \eqref{eq:firstchernclass-intro}. The ``fugacities'' are related by
$\xi_l=\zeta_{l-1}\, \zeta_{l+1}\,/\,\zeta_{l}$ for $l=1, \ldots, k-1$
(where we set $\zeta_0=\zeta_k=1$).

In Section \ref{sec:fundamentalmatter} we make some new observations
concerning the consistency of gauge theories on ALE spaces with $N$
fundamental hypermultiplets. In this case, the generating function for
$p$-observables (and all other partition functions) is defined by
means of the integral of an equivariant characteristic
class depending on the natural bundle $\Vbf$ and the fundamental
masses; this integral makes sense only if the degree of the class is
nonnegative. For $U(r)$ gauge theories on $\R^4$ this constraint
implies $N\leq 2r$, i.e., that the gauge theory is asymptotically
free\footnote{For $\Ncal=2$ gauge theories on $\R^4$ with nonvanishing
  beta-function $\beta:=2r-N$, the effective expansion parameter is
  $\qsf\, \Lambda^\beta$ where $\Lambda$ is an \emph{energy scale parameter.} In this
  paper we assume for ease of notation that all quantities have been
  suitably rescaled and formally set $\Lambda=1$ in all nonconformal
  gauge theory partition functions.}. For $U(r)$ gauge theories on $X_k$ the constraint is
$\dim_\C\big(\Mcal_{r,\vec{u},\Delta}(\Xscr_k,\Dscr_\infty,\Fcal_\infty^{0,\vec
  w}\, )\big)\geq \rk(\Vbf)\, N$. Since we wish to factorize our
partition functions in terms of partition functions of gauge theories
on $\R^4$, by combining these two constraints we get an inequality on
the first Chern class depending on the holonomy at infinity; for $k=2$
this inequality reads $4v^2 \leq w_1^2$, where $v:=\frac{1}{2} \, {u}$
gives the first Chern class. In particular if the holonomy at infinity
is trivial, i.e., $\vec{w}=(r, 0)$, we get $v=0$, which is the case
considered in \cite{art:bonellimaruyoshitanzini2012}. Therefore to
define the partition functions for $\Ncal =2$ gauge theories on $X_k$
with $N\leq 2r$ fundamental matter fields, we restrict to moduli
spaces
$\Mcal_{r,\vec{u},\Delta}(\Xscr_k,\Dscr_\infty,\Fcal_\infty^{0,\vec
  w}\, )$ of framed sheaves whose first Chern class obeys an additional
constraint  (cf.\ Equation \eqref{eq:fundamental-set}); we call the corresponding gauge theory
\emph{asymptotically free}. In a similar way, one can define a
\emph{conformal} gauge theory on $X_k$ with $N=2r$ fundamental matter fields.

In Section \ref{sec:VW} we make another important check of our
results: we verify that the $\mu=0$ limit of the $\Ncal=2^\ast$ gauge theory partition function correctly reduces to the partition function of the Vafa-Witten topologically twisted $\Ncal=4$ gauge theory on $X_k$ \cite{art:vafawitten1994}. We show that \begin{equation}
\Zcal^{\mathrm{VW}}_{\mathrm{ALE}}\big(\qsf, \vec{\xi}\
\big):=\lim_{\mu \to 0}\,\Zcal^{\ast,\mathrm{inst}}_{\mathrm{ALE}}\big(\varepsilon_1, \varepsilon_2,
\vec{a}, \mu ; \qsf, \vec{\xi}\ \big)= \qsf^{\frac{r\,k}{24}}\,
\prod_{j=0}^{k-1}\, \left(\frac{\chi^{\widehat{\omega}_{j}}(\zeta, \tau_0 )}{\eta(\tau_0)}\right)^{w_j}\ ,
\end{equation}
where $\chi^{\widehat{\omega}_{j}}(\zeta, \tau_0)$ is the character of the
integrable highest weight representation of $\widehat{\slfrak}(k)$ at
level one with weight the $i$-th fundamental weight $\widehat{\omega}_i$ of type $\hat{A}_{k-1}$ for
$i=0,1,\dots,k-1$, $\eta(\tau_0)^{-1}$ is the character of the
Heisenberg algebra $\hfrak$, and $\xi_j=\exp(2\pi
\operatorname{i}(2\zeta_{j}-\zeta_{j-1}-\zeta_{j+1}))$ (we set
$\zeta_0=\zeta_k=0$). Thus in this case we correctly reproduce the
character of the representation of the affine Lie algebra
$\widehat{\mathfrak{gl}}(k)_r$, and hence confirm the modularity
(S-duality) of the partition function in the limit. One should be able
to reproduce the same result by computing the Poincar\'e polynomial of
the moduli spaces
$\Mcal_{r,\vec{u},\Delta}(\Xscr_k,\Dscr_\infty,\Fcal_\infty^{0, \vec
  w}\, )$.

Another useful limit of the generating function
$\Zcal^\ast_{X_k}$ is obtained by setting
$\vec{\tau}=(0,-\tau_1,0,\ldots)$ and $\vec{t}^{\;(1)}= \cdots =
\vec{t}^{\;(k-1)}= \vec 0$; we denote the resulting partition function by $\Zcal_{X_k}^{\ast,\circ}$. While this deformation simply has the effect of shifting the bare coupling $\tau_0\to\tau_0+\tau_1$ in the tree-level Lagrangian of the $\Ncal=2$ gauge theory, it enables one to combine the classical and instanton partition functions into a single correlation function. In this case we find
\begin{equation}
\Zcal_{X_k}^{\ast,\circ}\big(\varepsilon_1, \varepsilon_2, \vec{a},
\mu; \qsf, \vec{\xi},\tau_1
\big)=\Zcal_{\C^2}^{\mathrm{cl}}(\varepsilon_1,
\varepsilon_2,\vec{a};\tau_1)^{\frac{1}{k}} \ \Zcal_{X_k}^{\ast,
  \circ, \mathrm{inst}}\big(\varepsilon_1, \varepsilon_2, \vec{a}, \mu
; \qsf_{\mathrm{eff}}
,\vec{\xi},\tau_1 \big)\ ,
\end{equation}
where $\qsf_{\mathrm{eff}}:=\qsf\,\e^{\tau_1}$, $\Zcal_{\C^2}^{\mathrm{cl}}$ is the classical partition function on $\R^4$, and 
\begin{equation}
\Zcal_{X_k}^{\ast, \circ, \mathrm{inst}} \big(\varepsilon_1,
\varepsilon_2, \vec{a}, \mu ; \qsf_{\rm eff} , \vec{\xi},\tau_1 =0 \big)=
\Zcal_{X_k}^{\ast,\mathrm{inst}} \big(\varepsilon_1, \varepsilon_2,
\vec{a}, m; \qsf, \vec{\xi} \ \big) \ .
\end{equation}
The explicit expression for $\Zcal_{X_k}^{\ast, \circ,
  \mathrm{inst}}(\varepsilon_1, \varepsilon_2, \vec{a}, \mu; \qsf_{\rm
  eff} , \vec{\xi},\tau_1)$ is given in Equation
\eqref{eq:instantontau} for pure $\Ncal=2$ gauge theory; analogous
expressions can be written in the case of gauge theories with
matter. This factorization also holds in the pure case without matter (see Section \ref{sec:pure}) and in the case with fundamental matter (see Section \ref{sec:fundamentalmatter}).

To obtain a partition function that enjoys modularity properties
appropriate to its conjectural geometric description in the $(2,0)$
theory compactified on a two-torus, we should further multiply the
classical and instanton parts by the purely perturbative contribution,
which is independent of $\qsf$ and $\vec\xi$. In Section
\ref{sec:perturbativepartitionfunctions} we give a definition of
perturbative partition functions following \cite[Sections~3.1 and
3.2]{art:okounkov2006} and \cite[Section~6]{art:gasparimliu2010}. In
particular we define the perturbative part of the equivariant Chern
character of the bundle $\Ebf$: At a fixed point $([(\Ecal,
\phi_{\Ecal})]\,,\, [(\Ecal', \phi_{\Ecal'})])$ it is the sum of the
equivariant Euler characteristics of the fibre $\Ebf_{([(\Ecal,
  \phi_{\Ecal})]\,,\, [(\Ecal', \phi_{\Ecal'})])}$ respectively over
the two affine toric neighbourhoods of the two fixed points of the
compactification divisor $D_\infty$; in the main part of this paper
these neighbourhoods are denoted by $U_{{\infty,k}}$ and
$U_{{\infty,0}}$. We show that the computation of the perturbative
part reduces to the computation of the equivariant Euler
characteristics, over $U_{{\infty,k}}$ and $U_{{\infty,0}}$, of a Weil
divisor on $\bar{X}_k$ given by $D_\infty$, $D_0$ and $D_k$, which
depends only on the  holonomy at infinity, i.e., only on the
framing vector $\vec{w}$. By using the perturbative part of the
equivariant Chern character of the bundle $\Ebf$ we define the
pertubative partition functions. The general formulas are rather
complicated; in the case $k=2$ the perturbative partition function is
given by
\begin{multline}
\Zcal_{X_2}^{\ast, \mathrm{pert}}(\varepsilon_1,
\varepsilon_2,\vec{a},\mu ) \\
= \sum_{\vec{c}}\ \prod_{\alpha\neq \beta}\,
\frac{\exp\Big(\gamma_{\varepsilon_2-\varepsilon_1,
    2\varepsilon_2}\big(\mu +a_{\beta\alpha} +c_{\beta\alpha}\, (\varepsilon_2-\varepsilon_1)-2\varepsilon_2\big) +
\gamma_{\varepsilon_1-\varepsilon_2,
  2\varepsilon_1}\big(\mu+a_{\beta\alpha}
-2\varepsilon_1\big)\Big)}{\exp\Big(\gamma_{\varepsilon_2-\varepsilon_1
  , 2\varepsilon_2}\big(a_{\beta\alpha} +c_{\beta\alpha}\, (\varepsilon_2-\varepsilon_1)-2\varepsilon_2\big)+\gamma_{\varepsilon_1-\varepsilon_2, 2\varepsilon_1}\big(a_{\beta\alpha}-2\varepsilon_1\big)\Big) }
 \ ,
\end{multline}
where the sum runs over all $\vec{c}:=(c_1, \ldots, c_r)$ such that
$c_{\alpha}= 0$ if $\alpha\leq w_0$ and $c_\alpha=1$ if
$w_0<\alpha\leq w_0+w_1$, and for each $\alpha, \beta=1,\ldots, r$
with $\alpha\neq \beta$ the quantity $c_{\beta\alpha}\in\{0,1\}$ is the parity of $c_\beta-c_\alpha$. The function $\exp\big(\gamma_{\varepsilon_1,\varepsilon_2}(x) \big)$ is the double zeta-function regularization of the infinite product
\begin{equation}
\prod_{i,j=0}^\infty \, \big(x-i\, \varepsilon_1-j\, \varepsilon_2
\big) \ .
\end{equation}
This formula resembles the one-loop partition function computed in
\cite[Section~3.2]{art:bonellimaruyoshitanziniyagi2012}; however, the
definition of the perturbative partition function on a noncompact
space requires a choice and so a direct comparison is meaningless
without further input. In \cite{art:bonellimaruyoshitanzini2012} the
choice is made to match with the DOZZ (Dorn-Otto-Zamolodchikov-Zamolodchikov) three-point correlation functions of supersymmetric Liouville theory. The analogous results in the case that fundamental hypermultiplets are included can be found in Section \ref{sec:perturbative-fundamentalmatter}, while the pertinent formulas in the generic case with $U(r)\times U(r'\,)$ bifundamental matter is presented in Section \ref{sec:Eulerclasses-pert}. 

\subsubsection{Seiberg-Witten geometry}

In Section \ref{sec:SWgeometry} we study relations between the
supersymmetric gauge theory we have developed and Seiberg-Witten theory. The low energy limit
of $\Ncal=2$ gauge theories on $\R^4$
is completely characterised by the (punctured) Seiberg-Witten curve $\Sigma$ of
genus $r$ \cite{art:seibergwitten1994-I}. The
curve $\Sigma$ is equiped with a meromorphic differential
$\lambda_{\rm SW}$, called the \emph{Seiberg-Witten differential}, and
its periods determine the \emph{Seiberg-Witten prepotential} $\Fcal^\ast_{\C^2}(\vec
a,\mu;\qsf)$ which is a holomorphic function of all parameters. In a
symplectic basis $\{A_\alpha,B_\alpha\}_{\alpha=1,\ldots,r}\cup \{S\}$ for the
homology group $H_1(\Sigma;\Z)$, the periods of the Seiberg-Witten
differential determine the quantities
\begin{equation}
a_\alpha= \oint_{A_\alpha}\, \lambda_{\rm SW}\ ,\quad \frac{\partial\Fcal^\ast_{\C^2}}{\partial a_\alpha}(\vec a,\mu;\qsf)
= -2\pi\ii \,\oint_{B_\alpha}\, \lambda_{\rm SW}  \qquad \mbox{and}
\qquad  \mu=\oint_S \,\lambda_{\rm SW} \ .
\end{equation}
It follows that the period matrix $\tau=(\tau_{\alpha\beta})$ of the
Seiberg-Witten curve
$\Sigma$ is related to the prepotential by
\begin{equation}
\tau_{\alpha\beta} = -\frac1{2\pi\ii}\
\frac{\partial^2\Fcal^\ast_{\C^2}}{\partial a_\alpha\, \partial a_\beta}(\vec
a,\mu;\qsf)\ ,
\end{equation}
and it determines the infrared effective gauge couplings.

The Seiberg-Witten prepotential can be recovered from the partition
function for the $\Omega$-deformed $\Ncal=2^\ast$ gauge theory on $\R^4$ in the low energy limit
in which the equivariant parameters $\varepsilon_1,\varepsilon_2$
vanish; this result was originally conjectured by Nekrasov
\cite{art:nekrasov2003} and subsequently proven in
\cite{art:nakajimayoshioka2005-I,art:nekrasovokounkov2006}. In
Sections \ref{sec:instantonprepotential} and
\ref{sec:perturbativeprepotential} we prove analogous results for
gauge theory on $X_k$. Let $\tilde{k}=k/2$ for even $k$ and $\tilde{k}=k$ for odd $k$.

\begin{theorem*}[{Theorem \ref{thm:instantonprepotential-masses}}]
$F_{X_k}^{\ast, \mathrm{inst}}(\varepsilon_1,\varepsilon_2,\vec{a}, \mu;
\qsf, \vec{\xi} \ ) :=- \tilde{k}\, \varepsilon_1\, \varepsilon_2\, \log
\Zcal_{X_k}^{\ast, \mathrm{inst}}(\varepsilon_1,\varepsilon_2,\vec{a},
\mu; \qsf, \vec{\xi} \ )$ is analytic in $\varepsilon_1,\varepsilon_2$ near $\varepsilon_1=\varepsilon_2=0$ and
\begin{equation}
\lim_{\varepsilon_1,\varepsilon_2\to 0}\, F_{X_k}^{\ast,
  \mathrm{inst}}\big(\varepsilon_1,\varepsilon_2,\vec{a}, \mu; \qsf,
\vec{\xi} \ \big) = \frac{1}{k} \, \Fcal_{\C^2}^{\ast, \mathrm{inst}}(\vec{a}, \mu; \qsf)\ ,
\end{equation}
where $\Fcal_{\C^2}^{\ast, \mathrm{inst}}(\vec{a}, \mu; \qsf)$ is the instanton part of the Seiberg-Witten prepotential of $\Ncal=2^\ast$ gauge theory on $\R^4$.
\end{theorem*}
\begin{corollary*}
$F_{X_k}^{\ast, \circ}(\varepsilon_1,\varepsilon_2,\vec{a}, \mu; \qsf,
\vec{\xi}, \tau_1):=- \tilde{k}\, \varepsilon_1\, \varepsilon_2\, \log
\Zcal_{X_k}^{\ast, \circ}(\varepsilon_1,\varepsilon_2,\vec{a}, \mu ; \qsf, \vec{\xi}, \tau_1)$ is analytic in $\varepsilon_1,\varepsilon_2$ near $\varepsilon_1=\varepsilon_2=0$ and
\begin{equation}
\lim_{\varepsilon_1,\varepsilon_2\to 0}\, F_{X_k}^{\ast,
  \circ}\big(\varepsilon_1,\varepsilon_2,\vec{a}; \qsf,
\vec{\xi},\tau_1) = \frac1k\, \Big(\,
\frac{\tilde{k} \, \tau_1}{2}\, \sum_{\alpha=1}^r \, a_\alpha^2 +
\Fcal_{\C^2}^{\ast, \mathrm{inst}}(\vec{a}, \mu ; \qsf_{\rm eff})
\Big) \ .
\end{equation}
\end{corollary*}
\begin{theorem*}[{Theorem \ref{thm:perturbativeprepotential}}]
For any fixed holonomy vector at infinity $\vec{c}=(c_1, \ldots,  c_r)\in\{0,1,\ldots,k-1\}^r$ we have
\begin{equation}
\lim_{\varepsilon_1,\varepsilon_2\to 0} \, \varepsilon_1 \,
\varepsilon_2 \, F_{X_k}^{\ast, \mathrm{pert}}(\varepsilon_1,
\varepsilon_2, \vec{a}, \mu ; \vec{c}\, )=\frac{1}{k}\, \Fcal_{\C^2}^{\ast, \mathrm{pert}}(\vec{a}, \mu)\ ,
\end{equation}
where $\Zcal_{X_k}^{*,{\rm pert}}(\varepsilon_1,\varepsilon_2,\vec
a,\mu):=\sum_{\vec c}\, \exp\big(-F_{X_k}^{\ast, \mathrm{pert}}(\varepsilon_1,
\varepsilon_2, \vec{a}, \mu ; \vec{c}\, )\big)$ and $\Fcal_{\C^2}^{\ast, \mathrm{pert}}(\vec{a}, \mu)$ is the perturbative part of the Seiberg-Witten prepotential of $\Ncal=2^\ast$ gauge theory on $\R^4$.
\end{theorem*}
These results agree with the analogous derivations in \cite[Section~2.1]{art:bonellimaruyoshitanziniyagi2012}. They confirm that $\Ncal=2$ gauge theories on the ALE space $X_k$, like their counterparts on $\R^4$, are some sort of quantization of a Hitchin system with spectral curve $\Sigma$.

In Section \ref{sec:blowupequations} we derive blowup equations which relate
the generating functions for correlators of quadratic $2$-observables to the instanton
partition functions in the low energy limit. We focus
on the case $k=2$ in which the generating function
$\Zcal_{X_k}^{\ast}$ factorizes; in this
case, the computations are manageable because of the
factorization of the generating function for correlators of $2$-observables into two
copies of instanton partition functions. For $k\geq3$, due to this
apparent absence of factorization we are unable to derive analogous
blowup equations in the case of weighted blowups.

Let $\Zcal_{X_2}^{\ast,\bullet}(\varepsilon_1, \varepsilon_2, \vec{a},
\mu ; \qsf, \vec{\xi}, t)$ be the generating function
$\Zcal_{X_2}^{\ast}(\varepsilon_1,\varepsilon_2, \vec{a}, \mu ; \qsf,
\vec{\xi}, \vec{\tau}, \vec{t} \ )$ specialized at $\vec{\tau}=\vec 0$
and $\vec{t}:=(0, -t, 0, \ldots)$. Let
$\Theta\big[\genfrac{}{}{0pt}{}{\vec{\mu}}{\vec{\nu}}\big](\vec{\zeta}\,\vert\,
\tau)$ be the Riemann theta-function with characteristic $\big[\genfrac{}{}{0pt}{}{\vec{\mu}}{\vec{\nu}}\big]$ on the Seiberg-Witten curve
$\Sigma$ for $\Ncal=2^\ast$ gauge theory on $\R^4$. 
\begin{theorem*}[{Theorem \ref{thm:theta-adjoint}}]
$\Zcal_{X_2}^{\ast,\bullet}(\varepsilon_1, \varepsilon_2, \vec{a},
\mu; \qsf, \xi, t)/\Zcal^{\ast, \mathrm{inst}}_{X_2}(\varepsilon_1,
\varepsilon_2, \vec{a},\mu ; \qsf, \xi)$ is analytic in
$\varepsilon_1, \varepsilon_2$ near $\varepsilon_1=\varepsilon_2=0$, and
\begin{multline}
\lim_{\varepsilon_1,\varepsilon_2\to0}\ \frac{\Zcal_{X_2}^{\ast,\bullet}(\varepsilon_1, \varepsilon_2, \vec{a},
\mu; \qsf, \xi, t)}{\Zcal^{\ast, \mathrm{inst}}_{X_2}(\varepsilon_1,
\varepsilon_2, \vec{a},\mu ; \qsf, \xi)} \\
=\exp\bigg(\Big(\qsf \, \frac{\partial}{\partial \qsf}\Big)^2
\Fcal_{\C^2}^{\ast, \mathrm{inst}}(\vec{a}, \mu ; \qsf)\, t^2+2\pi
\operatorname{i}\, \sum_{\alpha=w_0+1}^r \, \zeta_\alpha \bigg)\,
\frac{\Theta\big[\genfrac{}{}{0pt}{}{0}{C\, \vec{\nu}}\big](C\,
  (\vec{\zeta}+\vec{\kappa})\,\vert\, C\,
  \tau)}{\Theta\big[\genfrac{}{}{0pt}{}{0}{C\, \vec{\nu}}\big](C\,
  \vec{\kappa}\,\vert\, C\, \tau)}\ ,
\end{multline}
where
\begin{equation}
\zeta_\alpha:=-\frac{t}{2\pi \operatorname{i}} \, \Big(a_\alpha+\qsf
\, \frac{\partial^2\Fcal_{\C^2}^{\ast,
    \mathrm{inst}}}{\partial\qsf\, \partial a_\alpha}(\vec{a}, \mu ;
\qsf) \Big) \ ,
\end{equation}
while $\kappa_\alpha:=\frac{1}{4\pi\operatorname{i}}\, \log(\xi)$ for
$\alpha=1,\ldots,r$ and
\begin{equation*}
\nu_\alpha:=\left\{
\begin{array}{ll}
\displaystyle
\sum_{\beta=w_0+1}^r\, \log
\left((a_{\beta}-a_{\alpha})^2-\mu^2\right)-\frac{2\pi\operatorname{i}w_1}{r}\,
\tau_0 & \\[5pt]
\displaystyle
\hskip20mm + \, \sum_{\beta=w_0+1}^r \, \frac{\partial^2 \Fcal_{\C^2}^{\ast,
    \mathrm{inst}}}{\partial a_\alpha \, \partial a_\beta}(\vec{a},
\mu ; \qsf) & \mbox{ for } \alpha=1, \ldots, w_0\ , \\[10pt]
\displaystyle
-\sum_{\beta=1}^{w_0}\, \log
\left((a_{\beta}-a_{\alpha})^2-\mu^2\right)+\frac{4\pi\operatorname{i}w_0}{r}\,
\tau_0 & \mbox{ for } \alpha=w_0+1, \ldots, r\ .
\end{array}\right.
\end{equation*}
\end{theorem*}
This blowup equation underlies the modularity properties of the
partition function and correlators of
quadratic $2$-observables on the Seiberg-Witten curve for $\Ncal=2^\ast$
gauge theory on $X_2$ with period matrix $\tau$ twisted by the $A_1$
Cartan matrix $C$, and it
generalizes the representation of the Vafa-Witten partition function
at $\mu=0$ in terms of modular forms.
If the fixed holonomy at infinity is trivial, i.e.,
$\vec{w}=(w_0,w_1)=(r,0)$, the characteristic vector $\vec{\nu}\in
\C^r$ vanishes and our result resembles
\cite[Theorem~8.1]{art:nakajimayoshioka2005-I} and
\cite[Equation~(2.25)]{art:bonellimaruyoshitanziniyagi2012}. In
general, the nontrivial holonomy at infinity is encoded in
$\vec{\nu}$. 

\subsection{Acknowledgements} 

We are grateful to G.\ Bonelli, P.\ Eyssidieux, A.\ Konechny, M.\ Kool, C.-C.\ M.\ Liu, K.\ Maruyoshi,
A.\ Tanzini and F.\ Yagi for helpful discussions and correspondence. This work was supported in part by PRIN ``Geometria delle variet\`a algebriche"  2010S47ARA, by GNSAGA-INDAM, by the Grant RPG-404 from the Leverhulme Trust, and by the Consolidated Grant ST/J000310/1 from the
UK Science and Technology Facilities Council.

\bigskip\section{Preliminaries on stacks}\label{sec:preliminaries}

\subsection{Deligne-Mumford stacks}

In this subsection we summarise the conventions about stacks that will
be used throughout this paper. Our main reference for the theory of stacks
is the book \cite{book:laumonmoretbailly2000}. In this paper all
schemes are defined over $\C$ and are Noetherian, unless otherwise
stated. A variety is an irreducible reduced separated scheme of finite
type over $\C$. The smooth locus of a variety $X$ is denoted by
$X_{\rm sm}$. 

By a \emph{Deligne-Mumford stack} we mean a separated Noetherian Deligne-Mumford stack $\Xscr$ of finite type over $\C$. An \emph{orbifold} is a smooth Deligne-Mumford stack with generically trivial stabilizer.

The \emph{inertia stack} $\IXscr$ of a Deligne-Mumford stack $\Xscr$
is the fibre product $\Xscr\times_{\Xscr\times\Xscr}\Xscr$ with
respect to the diagonal morphism $\Delta\colon\Xscr\to
\Xscr\times\Xscr$. For a scheme $T$, an object in $\IXscr(T)$ is a pair $(x,g)$ where $x$ is an object of $\Xscr(T)$ and $g\colon x
\xrightarrow{\sim} x$ is an automorphism. A morphism $(x,g)\to
(x',g'\, )$ is a morphism $f\colon x\to x'$ in $\Xscr(T)$ such that $f\circ g=g'\circ f$. Let $\varpi\colon \IXscr\to \Xscr$ be the \emph{forgetful morphism} which for any scheme $T$ sends a pair $(x,g)$ to $x$.

An \emph{\'etale presentation} of a Deligne-Mumford stack $\Xscr$ is a
pair $(U, u)$, where $U$ is a scheme and $u\colon U \to \Xscr$ is a
{representable \'etale surjective} morphism
\cite[Definition~4.1]{book:laumonmoretbailly2000}. A morphism between
two \'etale presentations $(U,u)$ and $(V,v)$ of $\Xscr$ is a pair
$(\varphi, \alpha)$, where $\varphi\colon U\to V$ is a morphism and
$\alpha\colon u\xrightarrow{\sim} v\circ\varphi$ is a
2-isomorphism. The \emph{\'etale groupoid} associated with the \'etale
presentation $u\colon U\to \Xscr$ is the groupoid
\begin{equation}
  \begin{tikzpicture}[xscale=2.8,yscale=-.7, ->, bend angle=25]
\node (A0_1) at (0,1) {$U\times_{\Xscr} U$};
\node (A1_1) at (1.15,1) {$U$}; 
\node (B1) at (0.35,0.9) {$ $};
\node (B2) at (1,0.9) {$ $};
\node (C1) at (0.35,1.1) {$ $};
\node (C2) at (1,1.1) {$ $};
\node (Comma) at (1.3,1.2) {.};
    \path (B1) edge [->]node [auto] {$\scriptstyle{}$} (B2);
        \path (C1) edge [->]node [below] {$\scriptstyle{}$} (C2);
  \end{tikzpicture}
\end{equation}
If $\sf P$ is a property of schemes which is local in the \'etale
topology (for example regular, normal, reduced, Cohen-Macaulay, etc.),
then $\Xscr$ has the property $\sf P$ if for one (and hence every)
\'etale presentation $u\colon U\to \Xscr$ the scheme $U$ has the
property $\sf P$.

A \emph{gerbe} over a Deligne-Mumford stack $\Xscr$ is a stack $\Yscr$ over $\Xscr$ which
\'etale locally admits a section and for which any two
local sections are locally 2-isomorphic. For any integer $k\geq2$, let
$\mu_k$ denote the group of complex $k$-th roots of unity. A gerbe
$\Yscr\rightarrow \Xscr$ is a \emph{$\mu_k$-banded gerbe}, or simply a
\emph{$\mu_k$-gerbe}, if for every \'etale presentation $U$ of $\Xscr$ and
every object $x\in\Yscr(U)$ there is an isomorphism $\alpha_x\colon
\mu_k \to \operatorname{Aut}_U(x)$ of sheaves of groups satisfying natural compatibility conditions~\cite[Section~6.1]{art:fantechimannnironi2010}.

A (quasi-)coherent sheaf $\Ecal$ on $\Xscr$ is a collection of pairs
$(\Ecal_{U,u}, \theta_{\varphi,\alpha})$, where for any \'etale
presentation $u\colon U\to \Xscr$, $\Ecal_{U,u}$ is a (quasi-)coherent
sheaf on $U$, and for any morphism $(\varphi,\alpha)\colon (U,u)\to
(V,v)$ between two \'etale presentations of $\Xscr$,
$\theta_{\varphi,\alpha}\colon \Ecal_{U,u} \xrightarrow{\sim}
  \varphi^\ast \Ecal_{V,v}$ is an isomorphism which satisfies a
  {cocycle condition} with respect to three \'etale presentations
  (see \cite[Lemma~12.2.1]{book:laumonmoretbailly2000} and \cite[Definition~7.18]{art:vistoli1989}). A torsion free (resp.\ locally free) sheaf on $\Xscr$ is a coherent sheaf $\Ecal$ such that all $\mathcal{E}_{U,u}$ are torsion free (resp.\ locally free).

If $\Xscr$ is a Deligne-Mumford stack, by
\cite[Corollary~1.3-(1)]{art:keelmori1997} there exists a \emph{coarse
  moduli space} $(X,\pi)$ (or simply $X$) where $X$ is a separated
algebraic space; amongst other properties, the morphism $\pi\colon \Xscr\to X$ is proper and quasi-finite. We shall always assume that $X$ is a scheme. We recall some properties of Deligne-Mumford stacks that we shall use in this paper:
\begin{itemize}\setlength{\itemsep}{2pt}
\item The functor $\pi_\ast\colon \mathrm{QCoh}(\Xscr)\to \mathrm{QCoh}(X)$ is exact and maps coherent sheaves to coherent sheaves \cite[Lemma~2.3.4]{art:abramovichvistoli2002}. In particular, $\Xscr$ is {tame} \cite[Definition~3.1]{art:abramovicholssonvistoli2008}.
\item $H^\bullet(\Xscr, \Ecal)\simeq H^\bullet (X, \pi_\ast\Ecal)$ for any quasi-coherent sheaf $\Ecal$ on $\Xscr$ \cite[Lemma~1.10]{art:nironi2008}.
\end{itemize} 

\subsubsection*{Notation}

We use the symbols $\Ecal$, $\Gcal$, $\Fcal$, $\ldots$ for sheaves on a Deligne-Mumford stack, and the symbols $E$, $F$, $G$, $\ldots$ for sheaves on a scheme. For any coherent sheaf $\Fcal$ on a Deligne-Mumford stack $\Xscr$ we denote by $\Fcal^\vee$ its \emph{dual} $\mathcal{H}om(\Fcal,\Ocal_{\Xscr})$. We denote in the same way the dual of a coherent sheaf on a scheme. By $\A^n$ we denote the $n$-dimensional affine space
over $\C$, and by $\G_m$ the multiplicative group $\C^\ast$. The projection morphism  $T\times Y \to Y$  is written as  $p_Y$ or $p_{T\times Y,Y}$.

\subsection{Root stacks}\label{sec:rootstack}

In this subsection, we give a brief survey of the theory of root
stacks as presented in \cite{art:cadman2007} (see also \cite[Appendix
B]{art:abramovichgrabervistoli2008}).

Let $\Xscr$ be an algebraic stack. We use the standard fact that there
is an equivalence between the category of line bundles on $\Xscr$ and
the category of morphisms $\Xscr\to \Bscr\G_m$, where the morphisms in
the former category are taken to be isomorphisms of line
bundles. There is also an equivalence between the category of pairs
$(\Lcal,s)$, with $\Lcal$ a line bundle on $\Xscr$ and $s\in \Gamma(\Xscr,\Lcal)$, and the category of morphisms $\Xscr\to [\A^1/\G_m]$, where $\G_m$ acts on $\A^1$ by multiplication \cite[Example~5.13]{art:olsson2003}.

Throughout this subsection, $\Xscr$ will be an algebraic stack, $\Lcal$ a line bundle on $\Xscr$, $s\in \Gamma(\Xscr,\Lcal)$ a global section, and $k$ a positive integer.
The pair $(\Lcal,s)$ defines a morphism $\Xscr\to [\A^1/\G_m]$ as above. Let $\theta_k\colon [\A^1/\G_m]\to [\A^1/\G_m]$ be the morphism induced by the morphisms
\begin{equation}
x\in\A^1 \ \longmapsto \ x^k\in\A^1 \qquad \mbox{and} \qquad
t\in\G_m \ \longmapsto \ t^k\in\G_m\ ,
\end{equation}
which sends a pair $(\Lcal,s)$ to its $k$-th tensor power $(\Lcal^{\otimes k}, s^{\otimes k})$.
\begin{definition}
Let $\rootlines$ be the algebraic stack obtained as the fibre product
\begin{equation}
  \begin{tikzpicture}[xscale=1.5,yscale=-1.2]
    \node (A0_0) at (0, 0) {$\rootlines$};
    \node (A0_2) at (2, 0) {$[\A^1/\G_m]$};
    \node (A1_1) at (1, 1) {$\square$};
    \node (A2_0) at (0, 2) {$\Xscr$};
    \node (A2_2) at (2, 2) {$[\A^1/\G_m]$};
    \node (Comma) at (2.7, 1) {$.$};
    \path (A0_0) edge [->]node [auto] {$\scriptstyle{}$} (A2_0);
    \path (A0_0) edge [->]node [auto] {$\scriptstyle{}$} (A0_2);
    \path (A0_2) edge [->]node [auto] {$\scriptstyle{\theta_k}$} (A2_2);
    \path (A2_0) edge [->]node [auto] {$\scriptstyle{}$} (A2_2);
  \end{tikzpicture} 
\end{equation}
We say that $\rootlines$ is the \emph{root stack} obtained from $\Xscr$ by the $k$-th root construction.
\end{definition}
 
\begin{remark}
By \cite[Example~2.4.2]{art:cadman2007}, if $s$ is a nowhere vanishing
section then $\rootlines\simeq \Xscr$. This shows that all of the ``new" stacky structure in $\rootlines$ is concentrated at the zero locus of $s$.
\end{remark}
\begin{definition}
Let $\rootline$ be the algebraic stack obtained as the fibre product
\begin{equation}
  \begin{tikzpicture}[xscale=1.5,yscale=-1.2]
    \node (A0_0) at (0, 0) {$\rootline$};
    \node (A0_2) at (2, 0) {$\Bscr\G_m$};
    \node (A1_1) at (1, 1) {$\square$};
    \node (A2_0) at (0, 2) {$\Xscr$};
    \node (A2_2) at (2, 2) {$\Bscr\G_m$};
    \node (Comma) at (2.6, 1) {$,$};
    \path (A0_0) edge [->]node [auto] {$\scriptstyle{}$} (A2_0);
    \path (A0_0) edge [->]node [auto] {$\scriptstyle{}$} (A0_2);
    \path (A0_2) edge [->]node [auto] {$\scriptstyle{}$} (A2_2);
    \path (A2_0) edge [->]node [auto] {$\scriptstyle{}$} (A2_2);
  \end{tikzpicture} 
\end{equation}
where $\Xscr\to \Bscr\G_m$ is determined by $\Lcal$ and $\Bscr\G_m\to \Bscr\G_m$ is given by the map $t\in\G_m\mapsto t^k\in\G_m$.
\end{definition}

As described in \cite[Example~2.4.3]{art:cadman2007}, $\rootline$ is a
closed substack of $\sqrt[k]{(\Lcal,0)/\Xscr}.$ In general, let
$\Dscr$ be the vanishing locus of $s\in \Gamma(\Xscr,\Lcal)$.  One has a chain of inclusions of closed substacks
\begin{equation}
\sqrt[k]{\Lcal_{\vert \Dscr}/\Dscr}\ \subset \
\sqrt[k]{(\Lcal_{\vert \Dscr},0)/\Dscr}\ \subset \ \rootlines\ .
\end{equation}
In addition, $\sqrt[k]{\Lcal_{\vert \Dscr}/\Dscr}$ is isomorphic to the
reduced stack $\big(\, \sqrt[k]{(\Lcal_{\vert \Dscr},0)/\Dscr}\,
\big)_{\mathrm{red}}$. Locally, $\rootline$ is a quotient of $\Xscr$
by a trivial action of $\mu_k$, though this is not true globally. In
general, $\rootline$ is a $\mu_k$-gerbe over $\Xscr$. Its cohomology
class in the \'etale cohomology group $H^2(\Xscr;\mu_k)$ is obtained
from the class $[\Lcal] \in H^1(\Xscr; \G_m)$ via the boundary
homomorphism $\delta\colon H^1(\Xscr; \G_m)\to H^2(\Xscr; \mu_k)$ given by the Kummer exact sequence
\begin{equation}
1\ \longrightarrow \ \mu_k\ \longrightarrow \ \G_m\ \xrightarrow{(-)^k}
\ \G_m\ \longrightarrow \ 1\ . 
\end{equation}
Since the class $\delta([\Lcal])$ has trivial image in $H^2(\Xscr; \G_m)$, the gerbe $\rootline$ is called \emph{essentially trivial} \cite[Definition~2.3.4.1 and Lemma~2.3.4.2]{art:lieblich2007}.

By \cite[Corollaries~2.3.6 and 2.3.7]{art:cadman2007} we get the following result.
\begin{proposition}
The projection $\rootlines\to\Xscr$ is faithfully flat and
quasi-compact. If $X$ is a scheme and $L$ is a line bundle on $X$ with
global section $s$, then $X$ is a coarse moduli space for both $\schemerootlines$ and $\schemerootline$ under the projections to $X$.
\end{proposition}
\begin{theorem}[{\cite[Theorem~2.3.3]{art:cadman2007}}]
If $\Xscr$ is a Deligne-Mumford stack, then $\rootlines$ is also a
Deligne-Mumford stack. 
\end{theorem}

\subsubsection{Roots of an effective Cartier divisor on a smooth algebraic stack}\label{sec:rootdivisor}

\begin{definition}
Let $\Xscr$ be a smooth algebraic stack, $\Dscr\subset \Xscr$ an effective Cartier divisor and $k$ a positive integer. We denote by $\rootdiv$ the root stack $\sqrt[k]{(\Ocal_\Xscr(\Dscr),s_\Dscr)/\Xscr}$, where $s_\Dscr$ is the tautological section of $\Ocal_\Xscr(\Dscr)$ which vanishes along $\Dscr$. Let now $\Dscr_i\subset \Xscr$ be effective Cartier divisors and $k_i$ positive integers for $i=1,\ldots,n$. We denote by $\rootdivvec$ the fibre product
\begin{equation}
\sqrt[k_1]{\Dscr_1/\Xscr} \times_{\Xscr} \cdots \times_{\Xscr} \sqrt[k_n]{\Dscr_n/\Xscr}\ .
\end{equation}
\end{definition}
There is an equivalent definition of $\rootdivvec$ in \cite[Definition~2.2.4]{art:cadman2007}, which gives rise to a morphism $\rootdivvec \to [\A^1/\G_m]^n$. Each of the components $\rootdivvec \to [\A^1/\G_m]$ corresponds to an effective divisor $\tilde{\Dscr}_i$, i.e., the reduced closed substack $\pi^{-1}(\Dscr_i)_{\mathrm{red}}$, where $\pi\colon\rootdivvec \to \Xscr$ is the natural projection morphism. Moreover
\begin{equation}
\pi^\ast\Ocal_{\Xscr}(\Dscr_i)\simeq \Ocal_{\rootdivvec}(k_i\, \tilde{\Dscr}_i) \ .
\end{equation}
As explained in \cite[Section~2.1]{art:bayercadman2010}, since $\Xscr$
is smooth, each $\Dscr_i$ is smooth and has simple normal crossing.
Moreover, $\rootdivvec$ is smooth and each $\tilde{\Dscr}_i$ has simple
normal crossing; each $\tilde{\Dscr}_i$ is the root stack
$\sqrt[k_i]{\Ocal_{\Xscr}(\Dscr_i)_{\vert \Dscr_i}/\Dscr_i}$ and hence
is a $\mu_{k_i}$-gerbe over $\Dscr_i$. 

In closing this subsection we provide a useful characterization of the
Picard group of $\rootdivvec$. There is an exact sequence of groups \cite[Corollary~3.1.2]{art:cadman2007}
\begin{equation}\label{eq:rootstack}
0\ \longrightarrow \ \Pic(\Xscr) \ \xrightarrow{\pi^\ast} \
\Pic\Big(\rootdivvec\, \Big)\ \xrightarrow{ \ q \ } \ \prod_{i=1}^n\,
\mu_{k_i} \ \longrightarrow \ 0\ .
\end{equation}
Every line bundle $\Lcal\in \Pic\Big(\rootdivvec\, \Big)$ can
therefore be written as $\Lcal\simeq \pi^\ast\Mcal\otimes
\bigotimes_{i=1}^n\, \Ocal_{\rootdivvec}(m_i\, \tilde{\Dscr}_i)$, where
$\mathcal{M}\in \Pic(\Xscr)$ and $0\leq m_i<k_i$ for $i=1, \ldots,
n$; each integer $m_i$ is unique and $\Mcal$ is unique up to
isomorphism. The morphism $q$ maps $\Lcal$ to $(m_1,\ldots,
m_n)$.

In the following we shall denote by $\lfloor x\rfloor\in\Z$
  the integer part (floor function) and by $\{x\}:=x-\lfloor
  x\rfloor\in[0,1)$ the fractional part  of a rational 
  number $x$.

\begin{lemma}[{\cite[Theorem 3.1.1]{art:cadman2007}}]
Let $\Xscr$ be a smooth algebraic stack, $\Mcal$ a coherent sheaf on $\Xscr$ and $m_1, \ldots, m_n\in \Z$. Then
\begin{equation}
\pi_\ast\Big(\pi^{\ast}\Mcal\otimes\bigotimes_{i=1}^n\,
\Ocal_{\rootdivvec}(m_i\, \tilde{\Dscr}_i)\Big)\simeq \Mcal \otimes
\bigotimes_{i=1}^n \, \mathcal{O}_\Xscr\Big(\left\lfloor
  \frac{m_i}{k_i}\right\rfloor\, \Dscr_i\Big) \ .
\end{equation}
\label{lem:cadman3.1.1}\end{lemma}


\subsection{Toric stacks}\label{sec:toricstacks}

In this subsection we shall collect some results about toric stacks. We
describe two equivalent approaches to the theory. The first   is due to
Fantechi, Mann and Nironi \cite{art:fantechimannnironi2010}, and is
based on a generalization of the notions of tori and toric
varieties. The second   is of a combinatorial nature and is based on the
notion of \emph{stacky fan} due to Borisov, Chen and Smith
\cite{art:borisovchensmith2004}. 

Our reference for the geometry of
toric varieties is the book \cite{book:coxlittleschenck2011}. 
Let $T$ be a torus. Denote by $M:=\Hom(T,\C^\ast)$ the lattice of
characters of $T$ and by $N:=\Hom_{\Z}(M, \Z)$ the dual lattice of
one-parameter subgroups; for an element $m\in M$ we denote the
corresponding character by $e^m$. A \emph{toric variety} is a variety
$X$ containing a torus $T$ as an open subset, such that the action of
$T$ on itself extends to an algebraic action on the whole variety
$X$. It is characterized by combinatorial data encoded by a \emph{fan}
$\Sigma$ in the vector space $N_\Q:=N\otimes_\Z \Q$. We say that $\Sigma$ is a
\emph{rational simplicial fan} in $N_{\Q}$ if every cone $\sigma\in
\Sigma$ is generated by linearly independent vectors over $\Q$. We
denote by $\Sigma(j)$ the set of $j$-dimensional cones of $\Sigma$ for
$j=0, 1, \ldots, \dim_\C(X)$.  A \emph{ray} is a one-dimensional cone
of $\Sigma$. By the orbit-cone correspondence, rays of $\Sigma$
correspond to $T$-invariant divisors of $X$. We denote by $D_\rho$ the
$T$-invariant divisor associated with the ray $\rho$. In order to
avoid overly cumbersome notation, when the ray has an index, say
$\rho_\alpha$, we denote the associated divisor by $D_\alpha$. We
denote by $p_\sigma$ the torus-fixed point associated with a cone
$\sigma\in \Sigma$ of maximal dimension, and by $U_\sigma$ the affine
toric variety $\Spec(\C[\sigma^\vee\cap M])$ which is an open subset
of $X$. When the cone has an index, say $\sigma_i$, we denote the
associated fixed point by $p_i$ and by $U_i$ the associated affine variety.

\begin{definition}
A \emph{Deligne-Mumford torus} $\Tscr$ is a product $T\times \Bscr G$
where $T$ is an ordinary torus and $G$ is a finite abelian group.
\end{definition}

Fantechi, Mann and Nironi give a definition of Deligne-Mumford tori in
terms of \emph{Picard stacks}, depending on suitable morphisms of
finitely generated abelian groups
\cite[Definition~2.4]{art:fantechimannnironi2010}. They also show \cite[Proposition~2.6]{art:fantechimannnironi2010} that their definition is equivalent to the present one. So a \emph{morphism} between Deligne-Mumford tori is simply a morphism as Picard stacks.

\begin{definition}
A \emph{toric Deligne-Mumford stack} is a smooth Deligne-Mumford stack $\Xscr$ with coarse moduli space $\pi\colon \Xscr\to X$, together with an open immersion of a Deligne-Mumford torus $\imath\colon \Tscr\hookrightarrow \Xscr$ with dense image such that the action of $\Tscr$ on itself extends to an action $a\colon \Tscr\times \Xscr\to \Xscr$. A \emph{morphism} of toric Deligne-Mumford stacks $f\colon \Xscr \to \Xscr'$ is a morphism of stacks between $\Xscr$ and $\Xscr'$ which extends a morphism of the corresponding Deligne-Mumford tori $\Tscr\to \Tscr'$.
\end{definition}

A toric Deligne-Mumford stack is an orbifold if and only if its Deligne-Mumford torus is an ordinary torus. By \cite[Proposition~3.6]{art:fantechimannnironi2010} the toric structure on a toric Deligne-Mumford stack $\Xscr$ with Deligne-Mumford torus $\Tscr$ induces a toric structure on the coarse moduli space $X$ with torus $T$, which is a simplicial toric variety. In particular, $\Xscr$ is irreducible.

Let $\Xscr$ be a toric Deligne-Mumford stack of dimension $d$ with coarse moduli space $\pi\colon \Xscr\to X$. By \cite[Proposition~5.1, Theorem 5.2, Theorem 6.25 and Corollary 6.26]{art:fantechimannnironi2010}, $\pi$ factorizes as 
\begin{equation}\label{eq:factorization}
\begin{aligned}
  \begin{tikzpicture}[xscale=2.5,yscale=-1.3,->, bend angle=45]
  \node (A0_0) at (0, 0) {$\Xscr$};
\node (A1_0) at (1, 0) {$\Xscr^{\mathrm{rig}}$};
\node (A2_0) at (2, 0) {$\Xscr^{\mathrm{can}}$};
\node (A3_0) at (3, 0) {$X$};
\node (A0_3) at (3.2, 0) {$,$};
    \path (A0_0) edge [bend right]node [auto] {$\scriptstyle{\pi}$} (A3_0);
        \path (A0_0) edge [->]node [above] {$\scriptstyle{r}$} (A1_0);
        \path (A1_0) edge [->]node [above] {$\scriptstyle{f^{\mathrm{rig}}}$} (A2_0);
        \path (A2_0) edge [->]node [above] {$\scriptstyle{\pi^{\mathrm{can}}}$} (A3_0);
        \path (A1_0) edge [bend left] node [below] {$\scriptstyle{\pi^{\mathrm{rig}}}$} (A3_0);
  \end{tikzpicture}
  \end{aligned}
\end{equation}
with the following properties:
\begin{itemize}
\item $\Xscr^{\mathrm{can}}$ is the \emph{canonical toric orbifold} of
  $X$, i.e., the unique (up to isomorphism) smooth $d$-dimensional
  toric Deligne-Mumford stack such that the locus where
  $\pi^{\mathrm{can}}$ is not an isomorphism has dimension $\leq d-2$.
\item $\Xscr^{\mathrm{rig}}$ is the \emph{rigidification} of $\Xscr$ with
  respect to the generic stabilizer\footnote{The \emph{generic
      stabilizer} $G$ of $\Xscr$ is defined as the union, inside the
    inertia stack $\IXscr$, of all the components of maximal dimension
    and it is a subsheaf of groups of $\IXscr$. Intuitively, the
    {rigidification} of $\Xscr$ by its generic stabilizer $G$ is the
    stack with the same objects as $\Xscr$ and the automorphism group
    of an object $x$ of $\Xscr$ is the quotient
    $\mathrm{Aut}_\Xscr(x)/G$. Rigidifications can also be defined for any central subgroup of the generic stabilizer. For the general construction we refer to \cite[Appendix~A]{art:abramovicholssonvistoli2008} (see also \cite[Section~5.1]{art:abramovichcortivistoli2003}).}.
\item Let $D_1, \ldots, D_m$ be torus-invariant divisors in
  $X$, and define the smooth integral closed substacks $\tilde{\Dscr}_i:=(\pi^{\mathrm{can}})^{-1}(D_i)_{\mathrm{red}}$ of codimension one in $\Xscr^{\mathrm{can}}$. Then there exist integers $k_1,\ldots,k_m\in\N$ such that
\begin{equation}
\Xscr^{\mathrm{rig}} \simeq
\sqrt[\vec{k}]{\vec{\tilde{\Dscr}}/\Xscr^{\rm can}}\ .
\end{equation}
\item There exist $\ell$ line bundles $\Lcal_1,\ldots,\Lcal_\ell\in\Pic(\Xscr^{\mathrm{rig}})$ and integers $b_i\in\N$, $i=1,\ldots, \ell$ such that
\begin{equation}
\Xscr \simeq \sqrt[b_1]{\Lcal_1/\Xscr^{\mathrm{rig}}} \times_{\Xscr^{\mathrm{rig}}} \cdots \times_{\Xscr^{\mathrm{rig}}} \sqrt[b_\ell]{\Lcal_\ell/\Xscr^{\mathrm{rig}}}\ ;
\end{equation}
  thus $r\colon \Xscr \to \Xscr^{\mathrm{rig}}$ is an essentially
trivial $\prod_{i=1}^\ell\, \mu_{b_i}$-gerbe.
\end{itemize}

\subsubsection{Gale duality with torsion}\label{sec:galedual}

Here we follow the presentation of \emph{generalized Gale duality} in
\cite[Section~2]{art:borisovchensmith2004}; this paper extends the
classical Gale duality construction (see
e.g. \cite[Section~14.3]{book:coxlittleschenck2011}) to a larger class
of maps. Let $N$ be a finitely generated abelian group and
$\beta\colon \Z^n \to N$ a group homomorphism. 
Define the \emph{Gale dual} map $\beta^\vee\colon
(\Z^n)^\ast \to \DG(\beta)$ as follows (here we denote
$(-)^\ast:=\Hom(-,\Z)$). Take projective resolutions $E^\bullet$ and
$F^\bullet$ for $\Z^n$ and $N$ respectively. By
\cite[Theorem~2.2.6]{book:weibel1994}, $\beta$ lifts to a morphism
$E^\bullet\to F^\bullet$, and by
\cite[Subsection~1.5.8]{book:weibel1994} there is a short exact
sequence of cochain complexes $0\to F^\bullet \to \cone(\beta) \to
E^\bullet[1] \to 0$, where $\cone(\beta)$ is the mapping cone of $\beta$. Since
$E^\bullet$ is projective, it gives an exact sequence of cochain complexes
\begin{equation}
0 \ \longrightarrow \ E^\bullet[1]^\ast \ \longrightarrow \
\cone(\beta)^\ast \ \longrightarrow \ (F^\bullet)^\ast \
\longrightarrow \ 0
\end{equation}
which induces a long exact sequence in cohomology   containing the segment
\begin{equation}\label{eq:longexactgale}
N^\ast \ \xrightarrow{\beta^\ast} \ (\Z^n)^\ast \ \longrightarrow \
H^1\big(\cone(\beta)^\ast \big) \ \longrightarrow \ \Ext^1_\Z(N,\Z) \
\longrightarrow \ 0\ .
\end{equation} 
Define $\DG(\beta):=H^1(\cone(\beta)^\ast)$ and $\beta^\vee\colon (\Z^n)^\ast \to \DG(\beta)$ to be the second map in \eqref{eq:longexactgale}. By this definition, it is evident that the construction is natural.

There is also an explicit description of the   map $\beta^\vee$. If
$d$ is the rank of $N$, one can choose a projective re\-so\-lution of
$N$ of the form $0 \to \Z^r \xrightarrow{Q} \Z^{d+r} \to N \to 0$,
where $Q$ is a matrix of integers. Then $\beta\colon\Z^n \to N$ lifts to
a map $\Z^n \xrightarrow{B} \Z^{d+r}$ where $B$ is a matrix of integers. Thus the mapping cone $\cone(\beta)$ is the complex $0 \to
\Z^{n+r} \xrightarrow{[B\, Q]} \Z^{d+r} \to 0$, hence $\DG(\beta)=
(\Z^{n+r})^\ast / \mathrm{Im}([B\, Q]^\ast)$ and $\beta^\vee$ is the
composition of the inclusion map $(\Z^n)^\ast \to (\Z^{n+r})^\ast$
with the quotient map $(\Z^{n+r})^\ast \to \DG(\beta)$. If $N$ is
free, then $Q=0$ and $\mathrm{DG}(\beta)=(\Z^n)^*/\mathrm{Im}(B^*)$,
and moreover the kernel of $\beta^\vee$ is $N^*$ \cite[Proposition~2.2]{art:borisovchensmith2004}.

We give a property of the generalized Gale dual that will be useful in the following.
\begin{lemma}[{\cite[Lemma~2.3]{art:borisovchensmith2004}}]\label{lem:galedualsequences}
A morphism of short exact sequences
\begin{equation}
\begin{aligned}
  \begin{tikzpicture}[xscale=2.8,yscale=-0.8]
    \node (A0_0) at (0.4, 0) {$0$};
    \node (A0_1) at (1, 0) {$\Z^{n_1}$};
    \node (A0_2) at (2, 0) {$\Z^{n_2}$};
    \node (A0_3) at (3, 0) {$\Z^{n_3}$};
    \node (A0_4) at (3.6, 0) {$0$};
    \node (A1_5) at (3.77, 1) {$,$};
    \node (A2_0) at (0.4, 2) {$0$};
    \node (A2_1) at (1, 2) {$N_1$};
    \node (A2_2) at (2, 2) {$N_2$};
    \node (A2_3) at (3, 2) {$N_3$};
    \node (A2_4) at (3.6, 2) {$0$};    
    \path (A0_0) edge [->]node [left] {$\scriptstyle{}$} (A0_1);
    \path (A0_1) edge [->]node [left] {$\scriptstyle{}$} (A0_2);
    \path (A0_2) edge [->]node [left] {$\scriptstyle{}$} (A0_3);
    \path (A0_3) edge [->]node [left] {$\scriptstyle{}$} (A0_4);

    \path (A2_0) edge [->]node [left] {$\scriptstyle{}$} (A2_1);
    \path (A2_1) edge [->]node [above] {$\scriptstyle{}$} (A2_2);
    \path (A2_2) edge [->]node [left] {$\scriptstyle{}$} (A2_3);
    \path (A2_3) edge [->]node [left] {$\scriptstyle{}$} (A2_4);

    \path (A0_1) edge [->]node [auto] {$\scriptstyle{\beta_1}$} (A2_1);
    \path (A0_2) edge [->]node [auto] {$\scriptstyle{\beta_2}$} (A2_2);    
    \path (A0_3) edge [->]node [auto] {$\scriptstyle{\beta_3}$} (A2_3);    
  \end{tikzpicture} 
\end{aligned}
\end{equation}
in which the columns have finite cokernels, induces a morphism of short exact sequences
\begin{equation}
\begin{aligned}
  \begin{tikzpicture}[xscale=2.8,yscale=-0.8]
    \node (A0_0) at (0.4, 0) {$0$};
    \node (A0_1) at (1, 0) {$(\Z^{n_3})^\ast$};
    \node (A0_2) at (2, 0) {$(\Z^{n_2})^\ast$};
    \node (A0_3) at (3, 0) {$(\Z^{n_1})^\ast$};
    \node (A0_4) at (3.6, 0) {$0$};
    \node (A1_5) at (3.75, 1) {$.$};
    \node (A2_0) at (0.4, 2) {$0$};
    \node (A2_1) at (1, 2) {$\DG(\beta_3)$};
    \node (A2_2) at (2, 2) {$\DG(\beta_2)$};
    \node (A2_3) at (3, 2) {$\DG(\beta_1)$};
    \node (A2_4) at (3.6, 2) {$0$};    
    \path (A0_0) edge [->]node [left] {$\scriptstyle{}$} (A0_1);
    \path (A0_1) edge [->]node [left] {$\scriptstyle{}$} (A0_2);
    \path (A0_2) edge [->]node [left] {$\scriptstyle{}$} (A0_3);
    \path (A0_3) edge [->]node [left] {$\scriptstyle{}$} (A0_4);

    \path (A2_0) edge [->]node [left] {$\scriptstyle{}$} (A2_1);
    \path (A2_1) edge [->]node [above] {$\scriptstyle{}$} (A2_2);
    \path (A2_2) edge [->]node [left] {$\scriptstyle{}$} (A2_3);
    \path (A2_3) edge [->]node [left] {$\scriptstyle{}$} (A2_4);

    \path (A0_1) edge [->]node [auto] {$\scriptstyle{\beta_3^\vee}$} (A2_1);
    \path (A0_2) edge [->]node [auto] {$\scriptstyle{\beta_2^\vee}$} (A2_2);    
    \path (A0_3) edge [->]node [auto] {$\scriptstyle{\beta_1^\vee}$} (A2_3);    
  \end{tikzpicture} 
\end{aligned}
\end{equation}
\end{lemma}

\subsubsection{Stacky fans}\label{sec:stackyfans}

Let $N$ be a finitely generated abelian group of rank $d$. We write
$\overline{N}$ for the lattice generated by $N$ in the $d$-dimensional
$\Q$-vector space $N_{\Q}$. The natural map $N\to \overline{N}$ is
denoted by $b\mapsto \overline{b}$. Let $\Sigma$ be a rational
simplicial fan in $N_{\Q}$. We assume that the rays $\rho_1, \ldots,
\rho_n$ of $\Sigma$
span $N_\Q$ and   fix an element $b_i\in N$ such that $\bar{b}_i$
generates the cone $\rho_i$ for $i=1,\ldots,n$. The set $\{b_1,
\ldots, b_n\}$ defines a homomorphism of groups $\beta\colon \Z^n \to
N$ with finite cokernel. We call the triple $\stackyfan:=
(N,\Sigma,\beta)$ a \emph{stacky fan}.

A stacky fan $\stackyfan$ encodes a group action on a quasi-affine
variety $Z_\Sigma$. To describe this action, let $\C[z_1, \ldots, z_n]$ be
the coordinate ring of $\A^n$. The quasi-affine variety $Z_\Sigma$ is
$\A^n\setminus {\mathbb V}(J_\Sigma)$, where $J_\Sigma$ is the ideal generated
by the monomials $\prod_{\rho_i\not\subseteq \sigma} \, z_i$ for
$\sigma\in \Sigma$. The $\C$-valued points of $Z_\Sigma$ are the points $z\in
\A^n$ such that the cone generated by the set $\{\rho_i\, \vert\,
z_i=0\}$ belongs to $\Sigma$. We equip $Z_\Sigma$ with an action of the group
$G_{\stackyfan}:=\Hom_{\Z}(\DG(\beta), \C^\ast)$ as follows. By
applying the functor $\Hom_{\Z}(-, \C^\ast)$ to the map $\beta^\vee\colon (\Z^n)^*
\to \DG(\beta)$ we obtain a homomorphism $\alpha\colon G_{\stackyfan}
\to (\C^\ast)^n$. Then the natural action of $(\C^\ast)^n$ on $\A^n$
induces an action of $G_{\stackyfan}$ on $\A^n$. Since ${\mathbb V}(J_\Sigma)$
is a union of coordinate subspaces, the variety $Z_\Sigma$ is $G_{\stackyfan}$-invariant.

Let us define the global quotient stack $\Xscr(\stackyfan):=[Z_\Sigma/G_{\stackyfan}]$. 
By \cite[Lemma~7.15 and Theorem~7.24]{art:fantechimannnironi2010} one
has the following result.

\begin{proposition}\label{prop:stackyfan}
$\Xscr(\stackyfan)$ is a toric Deligne-Mumford stack of dimension $d$
with coarse moduli space the simplicial toric variety $X(\Sigma)$
associated with $\Sigma$. Conversely, let $\Xscr$ be a toric
Deligne-Mumford stack with coarse moduli space $X$. Denote by $\Sigma$
a fan of $X$ in $N_\Q$, and assume that the rays of $\Sigma$ span
$N_\Q$. Then there exists a finitely generated abelian group $N$ and a
group homomorphism $\beta\colon \Z^n\to N$ such that the stack $\Xscr(\stackyfan)$ associated with $\stackyfan=(N,\Sigma,\beta)$ is isomorphic to $\Xscr$ as toric Deligne-Mumford stacks.
\end{proposition}

\begin{remark}\label{rem:pictoric}
Let $\stackyfan=(N,\Sigma,\beta)$ be a stacky fan. By
\cite[Remark~7.19-(1)]{art:fantechimannnironi2010}, the Picard group $\Pic(\Xscr(\stackyfan))$ is isomorphic to $\DG(\beta)$.
\end{remark}

We will now reformulate the factorization \eqref{eq:factorization} in
terms of stacky fans. Let $\stackyfan=(N,\Sigma,\beta)$ be a stacky
fan. Denote by $\beta^{\mathrm{rig}}$ the composition of $\beta$
with  the quotient map $N\to N/N_{\mathrm{tor}}$, where
$N_{\mathrm{tor}}$ is the torsion subgroup of $N$. For $i=1,\ldots,
n$, denote by $v_i$ the unique generator of $\rho_i\cap
(N/N_{\mathrm{tor}})$. Then there exist unique integers $a_i$ such
that $\beta^{\mathrm{rig}}(e_i)=a_i\, v_i$, where $e_i$ with
$i=1,\dots,n$ is the standard basis of $\Z^n$. Therefore there exists
a unique group homomorphism $\beta^{\mathrm{can}}\colon\Z^n \to N/N_{\mathrm{tor}}$ such that the following diagram commutes
\begin{equation}
  \begin{tikzpicture}[xscale=1.5,yscale=-1.2]
    \node (A0_0) at (0, 0) {$\Z^n$};
    \node (A0_2) at (2, 0) {$N$};
    \node (A1_1) at (2.5, 1) {$.$};
    \node (A2_0) at (0, 2) {$\Z^n$};
    \node (A2_2) at (2, 2) {$N/N_{\mathrm{tor}}$};
    \path (A0_0) edge [->]node [left] {$\scriptstyle{\mathrm{diag}(a_1,\ldots,a_n)}$} (A2_0);
 \path (A0_0) edge [->]node [right=0.2cm] {$\scriptstyle{\beta^{\mathrm{rig}}}$} (A2_2);   
    \path (A0_0) edge [->]node [above] {$\scriptstyle{\beta}$} (A0_2);
    \path (A0_2) edge [->]node [auto] {$\scriptstyle{}$} (A2_2);
    \path (A2_0) edge [->]node [below] {$\scriptstyle{\beta^{\mathrm{can}}}$} (A2_2);
  \end{tikzpicture} 
\end{equation}
Set
$\stackyfan^{\mathrm{can}}:=(N/N_{\mathrm{tor}},\Sigma,\beta^{\mathrm{can}})$
and
$\stackyfan^{\mathrm{rig}}:=(N/N_{\mathrm{tor}},\Sigma,\beta^{\mathrm{rig}})$. By
\cite[Lemma~7.15, Theorem~7.17 and Theorem
7.24]{art:fantechimannnironi2010} we have the following result.
\begin{proposition}
Let $\stackyfan=(N,\Sigma,\beta)$ be a stacky fan. Assume that the
rays of $\Sigma$ span $N_\Q$. Then:
\begin{itemize}
\item $\Xscr(\stackyfan)^{\mathrm{rig}}\simeq \Xscr(\stackyfan^{\mathrm{rig}})$, and $\Xscr(\stackyfan)$ is a toric orbifold if and only if $N$ is free, or equivalently,  if and only if $\stackyfan=\stackyfan^{\mathrm{rig}}$.
\item $\Xscr(\stackyfan)^{\mathrm{can}}\simeq \Xscr(\stackyfan^{\mathrm{can}})$, and $\Xscr(\stackyfan)$ is canonical if and only if $\stackyfan=\stackyfan^{\mathrm{can}}$.
\end{itemize}
\end{proposition}

\subsubsection{Closed and open substacks}\label{sec:toricsubstacks}

Fix a cone $\sigma$ in the fan $\Sigma$. Let $N_\sigma$ be the
subgroup of $N$ generated by the set $\{b_i\,\vert\, \rho_i\subseteq
\sigma\}$ and let $N(\sigma)$ be the quotient group $N/N_\sigma$. By
extending scalars, the quotient map $N\to N(\sigma)$ becomes a
surjection $N_{\mathbb{Q}}\to N(\sigma)_{\mathbb{Q}}$. The quotient
fan $\Sigma/\sigma$ in $N(\sigma)_{\mathbb{Q}}$ is the set
$\{\tilde{\tau}=\tau + N(\sigma)_{\mathbb{Q}}\,\vert\, \sigma\subseteq
\tau \mbox{ and } \tau\in \Sigma\}$, and the \emph{link} of $\sigma$
is the set $\link(\sigma):=\{\tau\,\vert\, \tau+\sigma\in\Sigma,
\tau\cap\sigma=0\}$. For each ray $\rho_i$ in $\link(\sigma)$, we
write $\tilde{\rho}_i$ for the ray in $\Sigma/\sigma$ and
$\tilde{b}_i$ for the image of $b_i$ in $N(\sigma)$. Let $\ell$ be the
number of rays in $\link(\sigma)$ and let $\beta(\sigma)\colon
\Z^\ell\to N(\sigma)$ be the map determined by the list
$\{\tilde{b}_i\,\vert\, \rho_i\in\link(\sigma)\}$. The \emph{quotient
  stacky fan} $\bf \Sigma/\sigma$ is the triple $(N(\sigma),
\Sigma/\sigma, \beta(\sigma))$. By using
\cite[Proposition~4.2]{art:borisovchensmith2004} and the
characterization of the invariant divisors of a toric variety as
universal geometric quotients \cite{art:cox1995} one can prove the
following result.
\begin{proposition}\label{prop:closedsubstacks}
Let $\sigma$ be a cone in the stacky fan $\bf \Sigma$. Then $\Xscr(\bf \Sigma/\sigma)$ defines a closed substack of $\Xscr(\stackyfan)$. In particular, if $\rho$ is a ray in $\Sigma$, the closed substack $\Dscr_\rho=\Xscr(\bf \Sigma/\rho)$ is the effective Cartier divisor $\pi^{-1}(D_\rho)_{\mathrm{red}}$ of $\Xscr(\stackyfan)$ with coarse moduli space $D_\rho$, where $\pi\colon \Xscr(\stackyfan)\to X(\Sigma)$ is the coarse moduli space morphism.
\end{proposition}
In the following we shall call $\Dscr_\rho$ the torus-invariant
divisor $\Xscr(\stackyfan /\rho)$ associated with the ray $\rho$.
\begin{remark}\label{rem:restrction}
One can explicitly construct a map $\DG(\beta)\to \DG(\beta(\sigma))$
that gives both the restriction morphism $\Pic(\Xscr(\stackyfan)) \to
\Pic(\Xscr(\stackyfan/\sigma))$ and the group homomorphism
$G_{\stackyfan/\sigma}\to G_{\stackyfan}$. The method was first
described in the proof of \cite[Proposition
4.2]{art:borisovchensmith2004} but it had a gap, so we refer here to
\cite[Section~5.1]{art:jiangtseng2010} where the gap was fixed. Let
$\ell$ be the number of rays in $\link(\sigma)$ and $d$ the dimension
of $\sigma$. Consider the morphisms of short exact sequences
\begin{equation}
\begin{aligned}
  \begin{tikzpicture}[xscale=2.6,yscale=-0.8]
    \node (A0_0) at (0.4, 0) {$0$};
    \node (A0_1) at (1, 0) {$\Z^{\ell+d}$};
    \node (A0_2) at (2, 0) {$\Z^{n}$};
    \node (A0_3) at (3, 0) {$\Z^{n-\ell-d}$};
    \node (A0_4) at (3.6, 0) {$0$};
    \node (A1_5) at (3.77, 1) {$$};
    \node (A2_0) at (0.4, 2) {$0$};
    \node (A2_1) at (1, 2) {$N$};
    \node (A2_2) at (2, 2) {$N$};
    \node (A2_3) at (3, 2) {$0$};
    \node (A2_4) at (3.6, 2) {$0$};    
    \path (A0_0) edge [->]node [left] {$\scriptstyle{}$} (A0_1);
    \path (A0_1) edge [->]node [left] {$\scriptstyle{}$} (A0_2);
    \path (A0_2) edge [->]node [left] {$\scriptstyle{}$} (A0_3);
    \path (A0_3) edge [->]node [left] {$\scriptstyle{}$} (A0_4);

    \path (A2_0) edge [->]node [left] {$\scriptstyle{}$} (A2_1);
    \path (A2_1) edge [->]node [above] {$\scriptstyle{\simeq}$} (A2_2);
    \path (A2_2) edge [->]node [left] {$\scriptstyle{}$} (A2_3);
    \path (A2_3) edge [->]node [left] {$\scriptstyle{}$} (A2_4);

    \path (A0_1) edge [->]node [auto] {$\scriptstyle{\tilde{\beta}}$} (A2_1);
    \path (A0_2) edge [->]node [auto] {$\scriptstyle{\beta}$} (A2_2);    
    \path (A0_3) edge [->]node [auto] {$\scriptstyle{}$} (A2_3);    
  \end{tikzpicture} 
\end{aligned}
\end{equation}
and
\begin{equation}
\begin{aligned}
  \begin{tikzpicture}[xscale=2.6,yscale=-0.8]
    \node (A0_0) at (0.4, 0) {$0$};
    \node (A0_1) at (1, 0) {$\Z^{d}$};
    \node (A0_2) at (2, 0) {$\Z^{\ell+d}$};
    \node (A0_3) at (3, 0) {$\Z^{\ell}$};
    \node (A0_4) at (3.6, 0) {$0$};
    \node (A1_5) at (3.77, 1) {$$};
    \node (A2_0) at (0.4, 2) {$0$};
    \node (A2_1) at (1, 2) {$N_\sigma$};
    \node (A2_2) at (2, 2) {$N$};
    \node (A2_3) at (3, 2) {$N(\sigma)$};
    \node (A2_4) at (3.6, 2) {$0$};    
    \path (A0_0) edge [->]node [left] {$\scriptstyle{}$} (A0_1);
    \path (A0_1) edge [->]node [left] {$\scriptstyle{}$} (A0_2);
    \path (A0_2) edge [->]node [left] {$\scriptstyle{}$} (A0_3);
    \path (A0_3) edge [->]node [left] {$\scriptstyle{}$} (A0_4);

    \path (A2_0) edge [->]node [left] {$\scriptstyle{}$} (A2_1);
    \path (A2_1) edge [->]node [above] {$\scriptstyle{}$} (A2_2);
    \path (A2_2) edge [->]node [left] {$\scriptstyle{}$} (A2_3);
    \path (A2_3) edge [->]node [left] {$\scriptstyle{}$} (A2_4);

    \path (A0_1) edge [->]node [auto] {$\scriptstyle{\beta_\sigma}$} (A2_1);
    \path (A0_2) edge [->]node [auto] {$\scriptstyle{\tilde{\beta}}$} (A2_2);    
    \path (A0_3) edge [->]node [auto] {$\scriptstyle{\beta(\sigma)}$} (A2_3);    
  \end{tikzpicture} 
\end{aligned}
\end{equation}
Here we see $\mathbb{Z}^n$ as the group freely generated over $\mathbb{Z}$ by the rays of $\Sigma$ and $\mathbb{Z}^{\ell+d}$ as the group freely generated over $\mathbb{Z}$ by the rays in the link of $\sigma$ and by the rays of $\sigma$. Moreover, $\tilde{\beta}$ is the restriction of $\beta$ to $\mathbb{Z}^{\ell+d}$ and $\beta_\sigma\colon{\mathbb{Z}}^d\to N_\sigma$ is the group
homomorphism determined by the set $\{b_i\,\vert\, \rho_i\subseteq
\sigma\}$. By applying Lemma \ref{lem:galedualsequences} to these diagrams, one
obtains from the first diagram a morphism $\DG(\beta)\to
\DG(\tilde{\beta})$ and from the second   an isomorphism $\DG(\beta(\sigma))\to \DG(\tilde{\beta})$. Composing the inverse of the latter with the former, one has a morphism $\DG(\beta)\to \DG(\beta(\sigma))$.
\end{remark}

One can give a characterization of the inertia stack $\Ical
(\Xscr(\stackyfan))$ in terms of these closed substacks of
$\Xscr(\stackyfan)$. Recall that if $\Xscr$ is a
global quotient stack of the form $[Z/G]$, then
$\Ical(\Xscr)=\bigsqcup_{g\in G}\, [Z^g/G]$ where $Z^g$ is the fixed
point locus of $Z$ with respect to the element $g\in G$ (see e.g. \cite[Section~4]{art:borisovchensmith2004}).
After fixing a stacky fan $\mathbf{\Sigma}=(N,\Sigma,\beta)$, for every maximal cone $\sigma\in\Sigma(d)$ define the set
\begin{equation}
\boxx(\sigma) := \Big\{ v \in N\ \Big\vert\
\bar{v}=\mbox{$\sum\limits_{\rho_i \subseteq \sigma}$}\, q_i \,
\bar{b}_i \ \mbox{ for } \ 0\leq q_i < 1 \Big\}\ .
\end{equation}
The set $\boxx(\sigma)$ is in a one-to-one correspondence with the finite
group $N(\sigma)$. Define $\boxx(\Sigma) = \bigcup_{\sigma \in
  \Sigma(d) } \, \boxx(\sigma)$, and for every $v \in N$ call
$\sigma(v)$ the unique maximal cone containing $\bar{v}$. By
\cite[Lemma 4.6 and Theorem 4.7]{art:borisovchensmith2004} we get the
following result.
\begin{theorem}\label{thm:inertiatoric}
If $\Sigma$ is a complete fan, the elements $v\in \boxx(\Sigma)$ are in one-to-one correspondence with the elements $g \in G_\mathbf{\Sigma}$ which fix a point in $Z_\Sigma$ and
\begin{equation}
\Xscr\big(\mathbf{\Sigma} / \sigma(v) \big) \simeq \big[ Z_\Sigma^g /
G_\mathbf{\Sigma} \big] \ .
\end{equation}
Moreover, the inertia stack can be characterized as
\begin{equation}
\Ical\big(\Xscr(\mathbf{\Sigma})\big) = \bigsqcup_{v \in
  \boxx(\Sigma)} \, \Xscr\big(\mathbf{\Sigma} / \sigma(v)\big) \ .
\end{equation}
\end{theorem}

Viewing a $d$-dimensional cone $\sigma\in \Sigma$ as the fan
consisting of the cone $\sigma$ and all of its faces, we can identify
$\sigma$ with an open substack of $\Xscr(\stackyfan)$. The induced
stacky fan ${\boldsymbol \sigma}$ is the triple $(N,\sigma,
\beta_\sigma)$, where $\beta_\sigma\colon \Z^d\to N$ is the group
homomorphism introduced in Remark~\ref{rem:restrction}.
\begin{proposition}[{\cite[Proposition 4.3]{art:borisovchensmith2004}}]\label{prop:opensubstack}
Let $\sigma$ be a $d$-dimensional cone in the fan $\Sigma$. Then $\Xscr({\boldsymbol \sigma})$ is an open substack of $\Xscr(\stackyfan)$ of the form  $[V(\sigma)/N(\sigma)]$, where $V(\sigma)\simeq \C^d$ and $N(\sigma)$ is a finite abelian group acting on it, whose coarse moduli space is the open affine toric subset $U_\sigma$ of $X(\Sigma)$.
\end{proposition}
\begin{remark}
By varying the $d$-dimensional cones $\sigma$ of $\Sigma$, the open substacks $\Xscr(\bf \sigma)$ form an open cover of $\Xscr(\stackyfan)$. 
\end{remark}

\bigskip\section{Stacky compactification of the ALE space $X_k$}\label{sec:stackycompactification}

\subsection{Minimal resolution of $\C^2/\Z_k$}\label{sec:minimalresolution}

Let $k\geq 2$ be an integer and denote by $\mu_k$ the group of $k$-th roots of unity in $\C$. A choice of a primitive $k$-th root of unity $\omega$ defines an isomorphism of groups $\mu_k\simeq \Z_k$. We define an action of $\mu_k\simeq \Z_k$ on $\C^2$ as
\begin{equation}
\omega\triangleright (z_1,z_2):= (\omega\, z_1, \omega^{-1}\, z_2)\ .
\end{equation}
The quotient $\C^2/\Z_k$ is a normal affine toric surface. To describe
its fan we need to introduce some notation. Let $N\simeq \Z^2$ be the
lattice of one-parameter subgroups of the torus $T_t:=\C^\ast\times
\C^\ast$. Fix a $\Z$-basis $\{ \vec{e}_1, \vec{e}_2\}$ of $N$ and define the vectors
$\vec{v}_i:=i\, \vec{e}_1-(i-1)\, \vec{e}_2\in N$ for any integer $i\geq 0$. Then the fan of $\C^2/\Z_k$ consists of the two-dimensional cone $\sigma:=\cone(\vec{v}_0,\vec{v}_k)\subset N_\Q$ and its subcones.
The origin is the unique singular point of $\C^2/\Z_k$, and is a particular case of a \emph{rational double point} or \emph{Du~Val singularity} \cite[Definition~10.4.10]{book:coxlittleschenck2011}.

By \cite[Example~10.1.9 and Corollary~10.4.9]{book:coxlittleschenck2011}, the minimal resolution of singularities of $\C^2/\Z_k$ is the smooth toric surface $\varphi_k\colon X_k\to\C^2/\Z_k$ defined by the fan $\Sigma_k\subset N_{\Q}$ where
\begin{align}
\Sigma_k(0)&:=\big\{\{0\}\big\}\ ,\\[4pt]
\Sigma_k(1)&:=\big\{\rho_i:=\cone(\vec{v}_i)\,\big\vert\, i=0,1,2, \ldots,
k \big\}\ ,\\[4pt]
\Sigma_k(2)&:=\big\{\sigma_i:=\cone(\vec{v}_{i-1}, \vec{v}_i)\,\big\vert\, i=1,2,
\ldots, k\big\} \ . 
\end{align}
The vectors $\vec v_i$ are the minimal generators of the rays $\rho_i$ for
$i=0, 1, \ldots, k$.

\begin{remark}\label{rem:mckayALE}
Recall the \emph{McKay correspondence}: There is a one-to-one correspondence
between the irreducible representations of $\mu_k$ and the irreducible
components of the exceptional divisor $\varphi_k^{-1}(0)$ of the
minimal resolution $\varphi_k\colon X_k\to \C^2/\Z_k$, which are 
$T_t$-invariant rational curves $D_i$ for $i=1, \ldots, k-1$
\cite[Corollary~10.3.11]{book:coxlittleschenck2011}. By
\cite[Equation~(10.4.3)]{book:coxlittleschenck2011}, the intersection
matrix $(D_i\cdot D_j)_{1\leq i,j\leq k-1}$ is given by minus the
symmetric Cartan matrix $C$ of the root system of type $A_{k-1}$, i.e., one has
\begin{equation}
\left( D_i \cdot D_j \right)_{1\leq i,j\leq k-1} = -C= 
\begin{pmatrix}
-2 & 1 & \cdots & 0 \\
1 & -2 & \cdots & 0 \\
\vdots & \vdots & \ddots & \vdots\\
0 & 0 & \cdots & -2
\end{pmatrix} \ .
\end{equation}
\end{remark}

Let $U_i$ be the toric affine open subset of $X_k$
corresponding to the two-dimensional cone $\sigma_i$ for $i=1, \ldots,
k$. Then ${\mathbb{C}}[U_i]:={\mathbb{C}}[\sigma_i^\vee\cap
M]={\mathbb{C}}[T_1^{2-i}\,
T_2^{1-i},T_1^{i-1}\, T_2^i]$ for $i=1, \ldots, k$. By the relations
\begin{equation}\label{eq:variablerelations}
T_1=t_1^k\qquad\mbox{and}\qquad T_2=t_2 \, t_1^{1-k}\ .
\end{equation}
we have ${\mathbb{C}}[U_i]={\mathbb
{C}}[t_1^{k-i+1}\,
t_2^{1-i},t_1^{i-k}\, t_2^{i}]$.

After identifying the characters of $T_t$ with the one-dimensional
$T_t$-modules, we denote by $\varepsilon_1$ and $\varepsilon_2$ the
equivariant first Chern class of $t_1$ and $t_2$, respectively. Define
\begin{equation*}
\chi^i_1(t_1,t_2)=t_1^{k-i+1}\, t_2^{1-i}\qquad
\mbox{and}\qquad\chi_2^i(t_1,t_2)=t_1^{i-k}\, t_2^{i}\ .
\end{equation*}
Let $\varepsilon_j^{(i)}$ be the equivariant first Chern class of
$\chi^i_j$ for $i=1,\dots,k$ and $j=1,2$. Then
\begin{equation*}
\varepsilon_1^{(i)}(\varepsilon_1,\varepsilon_2)=(k-i+1)\,
\varepsilon_1-(i-1)\,
\varepsilon_2\qquad\mbox{and}\qquad\varepsilon_2^{(i)}(\varepsilon
_1,\varepsilon_2)=
-(k-i)\, \varepsilon_1+i\, \varepsilon_2\ .
\end{equation*}
%
\begin{lemma}\label{lem:charactertangent}
Let $i\in\{1, 2, \ldots, k\}$. The character of the tangent space to
$X_k$ at the torus-invariant point $p_i$ is given by
\begin{equation*}
\operatorname{ch}_{T_t}(T_{p_i}X_k)=\chi^i_1+\chi_2^i\ .
\end{equation*}
\end{lemma}

By \cite[Section~2.3]{art:brionvergne1997} one gets the
following result.
%
\begin{lemma}\label{lem:characterline}
Let $i\in\{0, 1, 2,\ldots, k\}$ and $j\in\{1, 2, \ldots, k\}$. The equivariant Chern
character of the line bundle ${\mathcal O}_{X_k}(D_i)$ restricted to $U_j$ is
\begin{equation*}
\operatorname{ch}_{T_t}\big({\mathcal O}_{X_k}(D_i)_{\vert U_j} \big)=\left\{
\begin{array}{l@{\quad}l}
{\,\mathrm{e}\,}^{-\varepsilon_1^{(i)}} \ , & j=i \ , \\[4pt]
{\,\mathrm{e}\,}^{-\varepsilon_2^{(i+1)}} \ , & j=i+1 \ , \\[4pt]
0 & \mbox{otherwise} \ .
\end{array}
\right.
\end{equation*}
\end{lemma}

\subsection{Normal compactification}

Let us define the vector $\vec{b}_\infty:=-\vec{v}_0-\vec{v}_k=-k\, \vec{e}_1+(k-2)\, \vec{e}_2$ in
$N$. Denote by $\rho_\infty$ the ray $\cone(\vec{b}_\infty)\subset N_\Q$ and
by $\vec{v}_\infty$ its minimal generator. For even $k$ we have
$\vec{v}_\infty=\frac{1}{2}\, \vec{b}_\infty$, while for odd $k$ we have $\vec{v}_\infty=\vec{b}_\infty$. Let $\sigma_{\infty,k}$ and $\sigma_{\infty,0}$ be the two-dimensional cones $\cone(\vec{v}_{k},\vec{v}_\infty)\subset N_\Q$ and $\cone(\vec{v}_{0}, \vec{v}_\infty)\subset N_\Q$, respectively.

Let $\bar{X}_k$ be the normal projective toric surface\footnote{The
  completeness of the fan is equivalent to the completeness of the
  surface \cite[Theorem~3.4.6]{book:coxlittleschenck2011}. In two
  dimensions, the completeness of a surface is equivalent to its
  projectivity \cite[Proposition~6.3.25]{book:coxlittleschenck2011}.}
defined by the complete fan $\bar{\Sigma}_k\subset N_\Q$ with
\begin{align}
\bar{\Sigma}_k(0)&:=\big\{\{0\} \big\} \ = \ \Sigma_k(0) \ ,\\[4pt]
\bar{\Sigma}_k(1)&:=\{\rho_i\,\vert\,i=0,1,2, \ldots,
k\}\cup\{\rho_\infty\} \ = \ \Sigma_k(1)\cup\{\rho_\infty\}\ ,\\[4pt]
\bar{\Sigma}_k(2)&:=\{\sigma_i\,\vert\, i=1,2, \ldots,
k\}\cup\{\sigma_{\infty,k},\sigma_{\infty,0}\} \ = \ \Sigma_k(2)\cup\{\sigma_{\infty,k}, \sigma_{\infty,0}\}\ .
\end{align}
The surface $X_k$ is an open dense subset of $\bar{X}_k$. 

Henceforth we will denote by ${\tilde k}\in\Z$ the integer $k/2$ if
$k$ is even, $k$ if $k$ is odd. By
\cite[Example~1.3.20]{book:coxlittleschenck2011}, the affine toric
open subsets $U_{\sigma_{\infty,k}}$ and $U_{\sigma_{\infty,0}}$ are
isomorphic to $\C^2/\Z_{\tilde k}$. In particular, for $k=2$ one has
$\tilde k=1$, and hence the toric surface $\bar{X}_2$ is smooth;
indeed, $\bar X_2$ is the second Hirzebruch surface~$\mathbb{F}_2$. 

\begin{proposition}\label{prop:intersections}
The intersection matrix $\left(D_i \cdot D_j \right)_{i,j=0, 1,\ldots, k, \infty}$ is given by
\begin{equation}\label{eq:intersectionamtrix}
\left(D_i \cdot D_j \right)_{i,j=0, 1,\ldots, k, \infty} = 
\begin{pmatrix}
\displaystyle\frac{2-k}{k}\; & 1 & 0 & \cdots & 0 & 0 & \displaystyle\frac{1}{\tilde k} \\[8pt]
1 & \; -2 \; & 1 & \cdots & 0& 0 & 0 \\[8pt]
0 & 1 & \; -2 \; & \cdots & 0& 0 & 0 \\[8pt]
\vdots & \vdots& \vdots & \ddots &\vdots & \vdots \\[8pt]
0& 0& 0& \cdots & -2& 1& 0 \\[8pt]
0 & 0 & 0 & \cdots& 1& \displaystyle\frac{2-k}{k} & \displaystyle\frac{1}{\tilde k} \\[8pt]
\displaystyle\frac{1}{\tilde k} & 0 & 0 & \cdots & 0& \displaystyle\frac{1}{\tilde k} &\; \displaystyle\frac{k}{\tilde{k}^2}
\end{pmatrix} \ .
\end{equation}
\end{proposition}
\proof
By \cite[Proposition~6.4.4-(a)]{book:coxlittleschenck2011} we have directly $D_\infty\cdot D_i=0$ for $i=1, \ldots, k-1$. On the other hand by  \cite[Lemma~6.4.2]{book:coxlittleschenck2011} we get
\begin{equation}
D_\infty\cdot D_0=\frac{\mathrm{mult}(\rho_0)}{\mathrm{mult}(\sigma_{\infty,0})}\ , 
\end{equation}
where $\mathrm{mult}(\rho_0)$ is the index of the sublattice $\Z \vec{v}_0$
in $\Q \vec{v}_0\cap N$, so that $\mathrm{mult}(\rho_0)=1$;
$\mathrm{mult}(\sigma_{\infty,0})$ is the index of $\Z \vec{v}_0+\Z
\vec{v}_\infty$ in $\left(\Q \vec{v}_0+\Q \vec{v}_\infty\right)\cap N$. Since $\Z \vec{v}_0+\Z
\vec{v}_\infty=\tilde{k}\, \Z \vec{e}_1+\Z \vec{e}_2$, we have $\mathrm{mult}(\sigma_{\infty,0})=\tilde{k}$. In the same way one computes $
D_\infty\cdot D_k$.

Suppose that $k$ is odd. We have $\vec{v}_0+\vec{v}_\infty+\vec{v}_k=0$, hence by using \cite[Proposition~6.4.4]{book:coxlittleschenck2011} we get
\begin{equation}
D_\infty\cdot
D_\infty=\frac{\mathrm{mult}(\rho_\infty)}{\mathrm{mult}(\sigma_{\infty,k})}=\frac{1}{k}
= \frac k{\tilde k^2} \ .
\end{equation}
In the same way, using the relation $k\, \vec{v}_1-(k-2)\, \vec{v}_0+\vec{v}_\infty=0$
we get $D_0\cdot D_0=-\frac{k-2}{k}$, and analogously $k\,
\vec{v}_{k-1}-(k-2)\, \vec{v}_k+\vec{v}_\infty=0$ so we have $D_k\cdot
D_k=-\frac{k-2}{k}$. For even $k$ one uses analogous relations. By
\cite[Corollary~6.4.3]{book:coxlittleschenck2011}, one has $D_i\cdot
D_j=1$ for $i=0, j=1$ and $i=k-1, j=k$. The result now follows by Remark \ref{rem:mckayALE}. 
\endproof

\begin{rem}
By \cite[Theorem~4.2.8]{book:coxlittleschenck2011} the divisors
$\tilde k\,D_\infty$, $\tilde k\,D_0$ and $\tilde k\, D_k$ are
Cartier, while by \cite[Theorem~6.3.12]{book:coxlittleschenck2011} the
divisor $\tilde k\,D_\infty$ is nef. Since $\big(\tilde k\,
D_\infty\big)^2=k$, the divisor $\tilde k\, D_\infty$ is   big as well.\hfill$\triangle$
\end{rem}
 

\subsection{Canonical stack}

Let $\pi_k^{\mathrm{can}}\colon\Xscr_k^{\mathrm{can}}\to \bar{X}_k$ be
the two-dimensional \emph{canonical} projective toric
orbifold\footnote{A Deligne-Mumford stack is \emph{projective} if it
  is a global quotient stack and it has a projective scheme as coarse
  moduli space \cite[Definition~2.20]{art:nironi2008}.} with
Deligne-Mumford torus $T_t$ and with coarse moduli space
$\bar{X}_k$. Since $\Xscr_k^{\mathrm{can}}$ is canonical, the
locus where $\pi_k^{\mathrm{can}}$ is not an isomorphism has a nonpositive
dimension. In particular, $\pi_k^{\mathrm{can}}$ is an isomorphism
precisely over the smooth locus $(\bar{X}_k)_{\rm sm}$ of $\bar{X}_k$. 

The stacky fan of $\Xscr_k^{\mathrm{can}}$ is
$\stackyfank^{\mathrm{can}}:=(N,\bar{\Sigma}_k,\beta^{\mathrm{can}})$,
where $\beta^{\mathrm{can}}\colon\Z^{k+2}\to N\simeq \Z^2$ is the map
sending the $i$-th coordinate vector of $\Z^{k+2}$ to $\vec{v}_i$ for
$i=0,1, \ldots,k,\infty$. As explained in Section \ref{sec:stackyfans}, the stacky fan determines the structure of $\Xscr_k^{\mathrm{can}}$ as a global quotient stack, and $\Xscr_k^{\mathrm{can}}$ is the quotient stack $[Z_{\bar{\Sigma}_k}/G_{\stackyfank^{\mathrm{can}}}]$, where $Z_{\bar{\Sigma}_k}$ is the union over all cones $\sigma\in\bar{\Sigma}_k$ of the open subsets
\begin{equation}\label{eq:Zsigma}
Z_\sigma:=\big\{\vec{z}\in\C^{k+2}\, \big\vert \, z_i\neq0 \ \mbox{ if
} \ \rho_i\notin\sigma \big\} \ \subset \ \C^{k+2}\ .
\end{equation}
The group $G_{\stackyfank^{\mathrm{can}}}$ is given by
\begin{equation}
G_{\stackyfank^{\mathrm{can}}}=\Hom_{\Z}\big(\DG(\beta^{\mathrm{can}}),\C^\ast\big)\ ,
\end{equation}
where $\DG(\beta^{\mathrm{can}})$ is the cokernel of the map
$(\beta^{\mathrm{can}})^\ast\colon \Z^2 \to\Z^{k+2}$ dual to the map
$\beta^{\mathrm{can}}$. Thus $\DG(\beta^{\mathrm{can}})$ $\simeq\Z^{k}$
and $G_{\stackyfank^{\mathrm{can}}}\simeq (\C^\ast)^k$. The action of
$G_{\stackyfank^{\mathrm{can}}}$ on
$Z_{\bar{\Sigma}_k}\subset\C^{k+2}$ is given by applying the functor
$\Hom_\Z(- ,\C^\ast)$ to the quotient map
$(\beta^{\mathrm{can}})^\vee\colon \Z^{k+2}\to
\DG(\beta^{\mathrm{can}})\simeq\Z^k$ to obtain an injective group morphism
\begin{equation}
\imath_{G_{\stackyfank^{\mathrm{can}}}}\, \colon\,
G_{\stackyfank^{\mathrm{can}}}=\Hom_{\Z}\big(\DG(\beta^{\mathrm{can}}),\C^\ast
\big) \ \longrightarrow \ \Hom_\Z\big(\Z^{k+2},\C^\ast \big)\simeq(\C^\ast)^{k+2}\ .
\end{equation}
By restricting the standard action of $(\C^\ast)^{k+2}$ on $\C^{k+2}$
with respect to this morphism, we find that the action of
$G_{\stackyfank^{\mathrm{can}}}$ on $Z_{\bar{\Sigma}_k}$ is given by
\begin{equation*}
(t_1,\dots,t_k)\triangleright(z_1,\dots,z_{k+2})=\left\{
\begin{array}{ll}
\Big(\,\prod\limits_{i=1}^{k-1\, }t_i^i\, t_k^{2-k}\,
z_1\,,\,\prod\limits_{i=1}^{k-1}\, t_i^{-(i+1)}\, t_k^k\, z_2\,,\,t_1\,
z_3\,,\,\dots\,,\,t_k\, z_{k+2}\, \Big) \ , & \ \mbox{odd $k$}\\[8pt]
\Big(\,\prod\limits_{i=1}^{k-1}\, t_i^i\, t_k^{1-\tilde{k}}\,
z_1\,,\,\prod\limits_{i=1}^{k-1}\, t_i^{-(i+1)}\, t_k^{\tilde{k}}\,
z_2\,,\,t_1\, z_3\,,\,\dots\,,\,t_k\, z_{k+2}\, \Big) \ , & \ \mbox{ even $k$}
\end{array}\right.
\end{equation*}
for $(t_1,\ldots,t_k)\in G_{\stackyfank^{\mathrm{can}}}$ and
$(z_1,\ldots,z_{k+2})\in Z_{\bar{\Sigma}_k}$.

The boundary divisor $\Xscr_k^{\mathrm{can}}\setminus T_t$ is a simple
normal crossing divisor with $k+2$ irreducible components which we denote by
$\tilde{\Dscr}_0,\tilde{\Dscr}_1, \ldots,\tilde{\Dscr}_k,\tilde{\Dscr}_\infty$. For $i=0, 1,\ldots, k, \infty$, the effective divisor $\tilde{\Dscr}_i$ is Cartier and   is $T_t$-invariant, hence it corresponds to the line bundle $\Ocal_{\Xscr_k^{\mathrm{can}}}(\tilde{\Dscr}_i)$, and the latter has a  canonical section $s_i$. Since $\Xscr_k^{\mathrm{can}}$ is a global quotient stack with trivial $\operatorname{Pic}(Z_{\bar{\Sigma}_k})$, the line bundles on $\Xscr_k^{\mathrm{can}}$ are in one-to-one correspondence with characters of $G_{\stackyfank^{\mathrm{can}}}$ \cite[Remark~1.1]{art:fantechimannnironi2010}. In particular, $\Ocal_{\Xscr_k^{\mathrm{can}}}(\tilde{\Dscr}_i)$ is associated with the character 
\begin{equation}
G_{\stackyfank^{\mathrm{can}}} \
\xrightarrow{\imath_{G_{\stackyfank^{\mathrm{can}}}}} \
(\C^\ast)^{k+2} \ \xrightarrow{ \ p_i \ } \ \C^\ast\ ,
\end{equation}
where $p_i$ is the $i$-th projection. The canonical section $s_i$ is the $i$-th coordinate of $Z_{\bar{\Sigma}_k}$. 

\begin{remark}
For $i=0,1, \ldots, k, \infty$ the divisor $\tilde{\Dscr}_i$ is a
projective toric orbifold with Deligne-Mumford torus $\C^\ast$ and
with coarse moduli space $D_i$. It is the torus-invariant divisor
associated with the ray $\rho_i$. We shall describe in detail later on the stacky fan of $\tilde{\Dscr}_\infty$.
\end{remark}

We can now characterize the Picard group
$\Pic(\Xscr_k^{\mathrm{can}})$ of $\Xscr_k^{\mathrm{can}}$. Firstly,
by \cite[Remark~4.14-(b)]{art:fantechimannnironi2010}, the group $\Pic(\Xscr_k^{\mathrm{can}})$ fits into the short exact sequence
\begin{equation}\label{eq:picardcanonical}
0 \ \longrightarrow \ M \ \longrightarrow \
\mathrm{Div}_{T_t}(\bar{X}_k) \ \longrightarrow \
\Pic(\Xscr_k^{\mathrm{can}}) \ \longrightarrow \ 0\ ,
\end{equation}
where the map $\mathrm{Div}_{T_t}(\bar{X}_k)\to
\Pic(\Xscr_k^{\mathrm{can}})$ is given by the composition of the
natural map $\mathrm{Div}_{T_t}(\bar{X}_k)\to A^1(\bar{X}_k)$ and the
pullback morphism $(\pi_k^{\mathrm{can}})^\ast\colon A^1(\bar{X}_k)\to
\Pic(\Xscr_k^{\mathrm{can}})$ which sends $[D]\in A^1(\bar{X}_k)$ to
the line bundle
$\Ocal_{\Xscr_k^{\mathrm{can}}}((\pi_k^{\mathrm{can}})^{-1}(D\cap
(\bar{X}_k)_{\rm sm}))$. By Remark \ref{rem:pictoric}, one has
$\Pic(\Xscr^{\mathrm{can}}_k)\simeq \DG(\beta^{\mathrm{can}})$ and so
the morphism $\mathrm{Div}_{T_t}(\bar{X}_k)\to
\Pic(\Xscr_k^{\mathrm{can}})$ is the quotient map
$(\beta^{\mathrm{can}})^\vee$. 

\subsubsection*{Characterization of $\tilde{\Dscr}_\infty$}

The Cartier divisor $\tilde{\Dscr}_\infty$ is a one-dimensional
projective toric orbifold with Deligne-Mumford torus $\C^\ast$ and
with coarse moduli space $D_\infty$. Its stacky fan is
$\stackyfank^{\mathrm{can}}/\rho_\infty:=(N(\rho_\infty),$
$\bar{\Sigma}_k/\rho_\infty,$ $\beta^{\mathrm{can}}(\rho_\infty))$
where the quotient group $N(\rho_\infty):=N/\Z \vec{v}_\infty$ is isomorphic
to $\Z$ and the quotient fan $\bar{\Sigma}_k/\rho_\infty\subset
N(\rho_\infty)\otimes_{\Z}\Q\simeq\Q$ is given by
\begin{equation}\label{eq:stackyfaninfty}
\begin{aligned}
\begin{array}{lll}
\bar{\Sigma}_k/\rho_\infty(0)&:=&\big\{\{0\}\big\}\ ,\\[4pt]
\bar{\Sigma}_k/\rho_\infty(1)&:=&\big\{\rho'_0:=\cone(1) \ , \
\rho'_\infty:=\cone(-1) \big\}\ ,
\end{array}
\end{aligned}
\end{equation}
while the map $\beta^{\mathrm{can}}(\rho_\infty)\colon\Z^2\to
N(\rho_\infty)$ is defined as multiplication by $(\tilde k,-\tilde
k)$. As described in Section \ref{sec:toricstacks}, toric orbifolds
are obtained by performing {root stack constructions} over 
canonical toric orbifolds along the torus-invariant divisors associated
with the rays of the stacky fan of the orbifold. In particular,
$\tilde{\Dscr}_\infty$ is obtained from $D_\infty\simeq \PP^1$ by
performing a $(\tilde k,\tilde k)$-root stack construction on the
torus-fixed points of $D_\infty$, which we denote by $0,\infty$, so that
\begin{equation}
\tilde{\Dscr}_\infty\simeq\sqrt[(\tilde k,\tilde
k)]{(0,\infty)/\mathbb{P}^1} \ \xrightarrow{ \ \tilde{\pi}_k \ } \ D_\infty\simeq \mathbb{P}^1\ ,
\end{equation}
where $\tilde{\pi}_k:=(\pi_k^{\mathrm{can}})_{\vert \tilde{\Dscr}_\infty}$; indeed, 
by \cite[Remark~7.19-(2)]{art:fantechimannnironi2010} the orders of
the root stack construction over $D_\infty$, relative to the points 0
and $\infty$,  are given by the absolute values of the components of
the matrix of $\beta^{\mathrm{can}}(\rho_\infty)$, so they are both
equal to $\tilde{k}$. 
\begin{proposition}\label{prop:quotientcanonicalinfty}
The toric orbifold $\tilde{\Dscr}_\infty$ is isomorphic as a global quotient stack to
\begin{equation}
\left[\frac{\C^2\setminus\{0\}}{\C^\ast\times\mu_{\tilde{k}}}\right]\ ,
\end{equation}
where the action of $\C^\ast\times\mu_{\tilde{k}}$ on $\C^2\setminus\{0\}$ is given by 
\begin{equation}\label{eq:actiontildeinfinity}
(t,\omega)\triangleright (z_1,z_2)=(t\, \omega\, z_1,t \, z_2)
\end{equation}
for $(t,\omega)\in\C^\ast\times\mu_{\tilde{k}}$ and $(z_1,z_2)\in\C^2\setminus\{0\}$.
\end{proposition}
\proof
By Proposition \ref{prop:closedsubstacks} we have $\tilde{\Dscr}_\infty\simeq [Z_{\bar{\Sigma}_k/\rho_\infty}/G_{\stackyfank^{\mathrm{can}}/\rho_\infty}]$, where $Z_{\bar{\Sigma}_k/\rho_\infty}:=\C^2\setminus\{0\}$ and $G_{\stackyfank^{\mathrm{can}}/\rho_\infty}:=\Hom_{\Z}(\DG(\beta^{\mathrm{can}}(\rho_\infty)),\C^\ast)$. It remains to prove that $G_{\stackyfank^{\mathrm{can}}/\rho_\infty}$ is isomorphic to $\C^\ast\times\mu_{\tilde{k}}$ and the action of $G_{\stackyfank^{\mathrm{can}}/\rho_\infty}$ on $Z_{\bar{\Sigma}_k/\rho_\infty}$ is given by \eqref{eq:actiontildeinfinity}.

As described in Section \ref{sec:galedual}, the abelian group $\DG(\beta^{\mathrm{can}}(\rho_\infty))$ can be realized as the cokernel of the map 
\begin{equation}
\beta^{\mathrm{can}}(\rho_\infty)^\ast\, \colon \, \Z \
\longrightarrow \ \Z^2\ , \qquad
m \ \longmapsto \ m\, \tilde{k}\, \vec{e}_1-m\, \tilde{k}\, \vec{e}_2\ .
\end{equation}
Hence $\DG(\beta^{\mathrm{can}}(\rho_\infty))\simeq\Z\oplus \Z_{\tilde{k}}$ and therefore $G_{\stackyfank^{\mathrm{can}}/\rho_\infty}\simeq \C^\ast\times\mu_{\tilde{k}}$. The quotient map $\beta^{\mathrm{can}}(\rho_\infty)^{\vee}\colon\Z^2\to \DG(\beta^{\mathrm{can}}(\rho_\infty))$ is given by the matrix
\begin{equation}\label{eq:betacaninfty}
\begin{pmatrix}
1 & 1\\
1 & 0
\end{pmatrix} \ .
\end{equation}
The action of $G_{\stackyfank^{\mathrm{can}}/\rho_\infty}$ on $\C^2\setminus\{0\}$ is the restriction of the standard action of $(\C^\ast)^2$ on $\C^2\setminus\{0\}$ via the immersion
\begin{equation}
\imath_{G_{\stackyfank^{\mathrm{can}}/\rho_\infty}}\, \colon \,
\C^\ast\times\mu_{\tilde{k}}\simeq
G_{\stackyfank^{\mathrm{can}}/\rho_\infty}\simeq\Hom_\Z
\big(\DG(\beta^{\mathrm{can}}(\rho_\infty)),\C^\ast \big) \
\longrightarrow \ \Hom_\Z\big(\Z^2,\C^\ast \big)\simeq(\C^\ast)^2
\end{equation}
obtained by applying the functor $\Hom_\Z(- , \C^\ast)$ to the quotient map $\beta^{\mathrm{can}}(\rho_\infty)^{\vee}$. Therefore the action of $G_{\stackyfank^{\mathrm{can}}/\rho_\infty}$ on $Z_{\bar{\Sigma}_k/\rho_\infty}$ is exactly \eqref{eq:actiontildeinfinity}.
\endproof
This characterization of $\tilde{\Dscr}_\infty$ as a global quotient stack yields the following result.
\begin{corollary}\label{cor:linebundletilde}
The Picard group $\Pic(\tilde{\Dscr}_\infty)$ of $\tilde{\Dscr}_\infty$ is isomorphic to $\Z\oplus\Z_{\tilde{k}}$. It is generated by the line bundles $\tilde{\Lcal}_1$ and $\tilde{\Lcal}_2$ corresponding, respectively, to the characters
\begin{equation}
\tilde{\chi}_1\, \colon\, (t,\omega)\in\C^\ast\times\mu_{\tilde{k}} \
\longmapsto \ t\in\C^\ast\qquad\mbox{and}\qquad \tilde{\chi}_2\,
\colon\, (t,\omega)\in\C^\ast\times\mu_{\tilde{k}} \ \longmapsto \ \omega\in\C^\ast\ .
\end{equation}
\end{corollary}
We give a geometrical interpretation of the line bundles
$\tilde{\Lcal}_1$ and $\tilde{\Lcal}_2$. Recall that the boundary
divisor $\tilde{\Dscr}_\infty\setminus \C^\ast$ is a simple normal
crossing divisor with two irreducible components $\tilde{p}_0$ and
$\tilde{p}_\infty$. These are effective Cartier divisors associated
with the rays $\rho'_0$ and $\rho'_\infty$, respectively, and they coincide with the closed substacks $\tilde{\pi}_k^{-1}(0)_{\mathrm{red}}$ and $\tilde{\pi}_k^{-1}(\infty)_{\mathrm{red}}$. Define the Cartier divisor $\tilde{p}:=\tilde{p}_0-\tilde{p}_\infty$.
\begin{proposition}
The line bundle $\Ocal_{\tilde{\Dscr}_\infty}(\tilde{p}_\infty)$ is isomorphic to $\tilde{\Lcal}_1$ and the line bundle $\Ocal_{\tilde{\Dscr}_\infty}(\tilde{p})$ is isomorphic to $\tilde{\Lcal}_2$.
\end{proposition}
\proof
Recall that the quotient map $\beta^{\mathrm{can}}(\rho_\infty)^{\vee}\colon\Z^2\to \DG(\beta^{\mathrm{can}}(\rho_\infty))\simeq \Z\oplus \Z_{\tilde{k}}$ sends the coordinate vector $\vec{e}_1$ of $\Z^2$ to $\vec{f}_1+\vec{f}_2$ and the coordinate vector $\vec{e}_2$ of $\Z^2$ to $\vec{f}_1$, where $\{\vec{f}_1,\vec{f}_2\}$ is the standard basis of $\Z\oplus \Z_{\tilde{k}}$. By Remark \ref{rem:pictoric}, $\Pic(\tilde{\Dscr}_\infty)$ is isomorphic to $\DG(\beta^{\mathrm{can}}(\rho_\infty))\simeq \Z\oplus \Z_{\tilde{k}}$, hence the vector $\vec{f}_i$ corresponds to $\tilde{\Lcal}_i$ for $i=1,2$ by Corollary \ref{cor:linebundletilde}. So the quotient map $\beta^{\mathrm{can}}(\rho_\infty)^{\vee}\colon\Z^2\to \DG(\beta^{\mathrm{can}}(\rho_\infty))$ can be interpreted as a map from $\mathrm{Div}_{T_t}(\tilde{\Dscr}_\infty)\simeq \Z \rho_0'\oplus\Z\rho_\infty'$ to $\Pic(\tilde{\Dscr}_\infty)$, hence we have $\Ocal_{\tilde{\Dscr}_\infty}(\tilde{p}_0) \simeq \tilde{\Lcal}_1\otimes\tilde{\Lcal}_2$ and $\Ocal_{\tilde{\Dscr}_\infty}(\tilde{p}_\infty)\simeq\tilde{\Lcal}_1$. Then $\Ocal_{\tilde{\Dscr}_\infty}(\tilde{p})\simeq\Ocal_{\tilde{\Dscr}_\infty}(\tilde{p}_0)\otimes\Ocal_{\tilde{\Dscr}_\infty}(-\tilde{p}_\infty)\simeq\tilde{\Lcal}_2$.
\endproof

The following result will be useful later on.
\begin{lemma}\label{lem:restrictioncanonical}
The line bundle $\Ocal_{\Xscr_k^{\mathrm{can}}}(\tilde{\Dscr}_\infty)_{\vert\tilde{\Dscr}_\infty}$ assumes the following form with respect to the generators $\Ocal_{\tilde{\Dscr}_\infty}(\tilde{p}_\infty)$ and $\Ocal_{\tilde{\Dscr}_\infty}(\tilde{p})$ of $\Pic(\tilde{\Dscr}_\infty)$:
\begin{equation}\label{eq:decompositiontilde}
\Ocal_{\Xscr_k^{\mathrm{can}}}(\tilde{\Dscr}_\infty)_{\vert\tilde{\Dscr}_\infty}\simeq\left\{
\begin{array}{cc}
\Ocal_{\tilde{\Dscr}_\infty}(2\tilde{p}_\infty)\otimes\Ocal_{\tilde{\Dscr}_\infty}(\tilde{p}) & \mbox{ for even $k$}\ ,\\[5pt]
\Ocal_{\tilde{\Dscr}_\infty}(\tilde{p}_\infty)\otimes\Ocal_{\tilde{\Dscr}_\infty}\big(\frac{k+1}{2}\,
\tilde{p}\big) & \mbox{ for odd $k$}\ .
\end{array}\right.
\end{equation}
Moreover 
\begin{equation}
\Ocal_{\Xscr_k^{\mathrm{can}}}(\tilde{\Dscr}_0)_{\vert\tilde{\Dscr}_\infty}\simeq\Ocal_{\tilde{\Dscr}_\infty}(\tilde{p}_0)\qquad\mbox{and}\qquad \Ocal_{\Xscr_k^{\mathrm{can}}}(\tilde{\Dscr}_k)_{\vert\tilde{\Dscr}_\infty}\simeq\Ocal_{\tilde{\Dscr}_\infty}(\tilde{p}_\infty)\ .
\end{equation}
\end{lemma}
\proof
We need to compute the restriction map
$\Pic(\Xscr_k^{\mathrm{can}})\to \Pic(\tilde \Dscr_\infty)$. Since by
Remark \ref{rem:pictoric}, $\Pic(\Xscr_k^{\mathrm{can}})\simeq
\DG(\beta^{\mathrm{can}})$ and $\Pic(\tilde \Dscr_\infty)\simeq
\DG(\beta^{\mathrm{can}}(\rho_\infty))$, we need only to determine a
\emph{restriction} map from $\DG(\beta^{\mathrm{can}})$ to
$\DG(\beta^{\mathrm{can}}(\rho_\infty))$. By applying the procedure
described in Remark \ref{rem:restrction} we obtain morphisms of short
exact sequences
\begin{equation}\label{eq:diagchinaI}
\begin{aligned}
  \begin{tikzpicture}[xscale=3.3,yscale=-0.8]
    \node (A0_0) at (0.4, 0) {$0$};
    \node (A0_1) at (1, 0) {$\Z^{k-1}$};
    \node (A0_2) at (2, 0) {$\Z^{k+2}$};
    \node (A0_3) at (3, 0) {$\Z^{3}$};
    \node (A0_4) at (3.6, 0) {$0$};
    \node (A1_5) at (3.75, 1) {$$};
    \node (A2_0) at (0.4, 2) {$0$};
    \node (A2_1) at (1, 2) {$\Z^{k-1}$};
    \node (A2_2) at (2, 2) {$\DG(\beta^{\mathrm{can}})$};
    \node (A2_3) at (3, 2) {$\DG(\tilde{\beta}^{\mathrm{can}})$};
    \node (A2_4) at (3.6, 2) {$0$};    
    \path (A0_0) edge [->]node [left] {$\scriptstyle{}$} (A0_1);
    \path (A0_1) edge [->]node [left] {$\scriptstyle{}$} (A0_2);
    \path (A0_2) edge [->]node [left] {$\scriptstyle{}$} (A0_3);
    \path (A0_3) edge [->]node [left] {$\scriptstyle{}$} (A0_4);

    \path (A2_0) edge [->]node [left] {$\scriptstyle{}$} (A2_1);
    \path (A2_1) edge [->]node [above] {$\scriptstyle{}$} (A2_2);
    \path (A2_2) edge [->]node [left] {$\scriptstyle{}$} (A2_3);
    \path (A2_3) edge [->]node [left] {$\scriptstyle{}$} (A2_4);

    \path (A0_1) edge [->]node [auto] {$\scriptstyle{\mathrm{id}}$} (A2_1);
    \path (A0_2) edge [->]node [auto] {$\scriptstyle{(\beta^{\mathrm{can}})^\vee}$} (A2_2);    
    \path (A0_3) edge [->]node [auto] {$\scriptstyle{(\tilde{\beta}^{\mathrm{can}})^\vee}$} (A2_3);    
  \end{tikzpicture} 
\end{aligned}
\end{equation}
and
\begin{equation}\label{eq:diagchinaII}
\begin{aligned}
  \begin{tikzpicture}[xscale=3.3,yscale=-0.8]
    \node (A0_0) at (0.4, 0) {$0$};
    \node (A0_1) at (1, 0) {$\Z^2$};
    \node (A0_2) at (2, 0) {$\Z^{3}$};
    \node (A0_3) at (3, 0) {$\Z$};
    \node (A0_4) at (3.6, 0) {$0$};
    \node (A1_5) at (3.75, 1) {$,$};
    \node (A2_0) at (0.4, 2) {$0$};
    \node (A2_1) at (1, 2) {$\DG(\beta^{\mathrm{can}}(\rho_\infty))$};
    \node (A2_2) at (2, 2) {$\DG(\tilde{\beta}^{\mathrm{can}})$};
    \node (A2_3) at (3, 2) {$\DG(\beta^{\mathrm{can}}_{\rho_\infty})$};
    \node (A2_4) at (3.6, 2) {$0$};    
    \path (A0_0) edge [->]node [left] {$\scriptstyle{}$} (A0_1);
    \path (A0_1) edge [->]node [left] {$\scriptstyle{}$} (A0_2);
    \path (A0_2) edge [->]node [left] {$\scriptstyle{}$} (A0_3);
    \path (A0_3) edge [->]node [left] {$\scriptstyle{}$} (A0_4);

    \path (A2_0) edge [->]node [left] {$\scriptstyle{}$} (A2_1);
    \path (A2_1) edge [->]node [above] {$\scriptstyle{\simeq}$} (A2_2);
    \path (A2_2) edge [->]node [left] {$\scriptstyle{}$} (A2_3);
    \path (A2_3) edge [->]node [left] {$\scriptstyle{}$} (A2_4);

    \path (A0_1) edge [->]node [auto] {$\scriptstyle{(\beta^{\mathrm{can}}(\rho_\infty))^\vee}$} (A2_1);
    \path (A0_2) edge [->]node [auto] {$\scriptstyle{(\tilde{\beta}^{\mathrm{can}})^\vee}$} (A2_2);    
    \path (A0_3) edge [->]node [auto] {$\scriptstyle{(\beta^{\mathrm{can}}_{\rho_\infty})^\vee}$} (A2_3);    
  \end{tikzpicture} 
 \end{aligned}
\end{equation}
where $\tilde{\beta}^{\mathrm{can}}\colon \Z^3\to N$ is the
restriction of $\beta^{\mathrm{can}}\colon \Z^{k+2}\to N$ to the
subgroup $\Z^3\subset\Z^{k+2}$ generated by the rays $\rho_0, \rho_k,
\rho_\infty$ (here we identify $\Z^{k+2}$ with the free abelian group generated by the rays of $\bar{\Sigma}_k$). Since $N_{\rho_\infty}$ is the subgroup of $N$ generated by $\vec{v}_\infty$, the map $\beta_{\rho_\infty}^{\mathrm{can}}\colon \Z\to N_{\rho_\infty}$ sends $1$ to $\vec{v}_\infty$.

Let us denote by $\phi$ the isomorphism $\DG(\beta^{\mathrm{can}}_{\rho_\infty})\xrightarrow{\sim}\DG(\tilde{\beta}^{\mathrm{can}})$. One can explicitly compute the map $(\tilde{\beta}^{\mathrm{can}})^\vee\colon\Z^3\to\DG(\tilde{\beta}^{\mathrm{can}})\simeq \Z\oplus\Z_{\tilde{k}}$ in the commutative diagram \eqref{eq:diagchinaI}, obtaining the matrix
\begin{equation}
\begin{pmatrix}
1 & 1 & 2\\
-1 & 0 & -1
\end{pmatrix} \ \mbox{ for even $k$}\qquad \mbox{and} \qquad 
\begin{pmatrix}
1 & 1 & 1\\
-1 & 1 & 0
\end{pmatrix} \ \mbox{ for odd $k$}\ .
\end{equation}
Since $\beta^{\mathrm{can}}(\rho_\infty)^\vee$ is given by the matrix \eqref{eq:betacaninfty}, the map $\phi$ in the commutative diagram \eqref{eq:diagchinaII} is represented by the matrix 
\begin{equation}
\begin{pmatrix} 1 & 0 \\ 0 & -1 \end{pmatrix} \ \mbox{ for even
  $k$}\qquad \mbox{and} \qquad \begin{pmatrix} 1 & 0 \\ 1 &
  -2 \end{pmatrix} \ \mbox{ for odd $k$}\ .
\end{equation}
Then its inverse is represented by
\begin{equation}
\begin{pmatrix} 1 & 0 \\ 0 & -1 \end{pmatrix} \ \mbox{ for even
  $k$}\qquad \mbox{and} \qquad \begin{pmatrix} 1 & 0 \\ \frac{k+1}{2}
  & \frac{k-1}{2} \end{pmatrix} \ \mbox{ for odd $k$}\ .
\end{equation}
Thus the restriction map $\Pic(\Xscr_k)\to\Pic(\tilde{\Dscr}_\infty)$
is given by the composition of the map $\DG(\beta^{\mathrm{can}})\to
\DG(\tilde{\beta}^{\mathrm{can}})$ in the commutative diagram
\eqref{eq:diagchinaI} with $\phi^{-1}$. Since the line bundle
$\Ocal_{\Xscr_k^{\mathrm{can}}}(\tilde{\Dscr}_\infty)$ is the image of
$(0, \ldots, 0,1)\in \Z^{k+2}$ via the map
$(\beta^{\mathrm{can}})^\vee$, the first statement follows. 

The results for the restrictions of the line bundles $\Ocal_{\Xscr_k^{\mathrm{can}}}(\tilde{\Dscr}_0)$ and $\Ocal_{\Xscr_k^{\mathrm{can}}}(\tilde{\Dscr}_k)$ follow in the same way.
\endproof
\begin{remark}\label{rem:ampleness}
By this Lemma, it is easy to see that the line bundle $\Ocal_{\Xscr_k^{\mathrm{can}}}(\tilde{\Dscr}_\infty)$ is \emph{$\pi_k^{\mathrm{can}}$-ample}, i.e., the representation of its fibre at any geometric point of $\Xscr_k^{\mathrm{can}}$ is faithful. 
\end{remark}

\subsection{Root toric stack}

Let
$\Xscr_k:=\sqrt[k]{\tilde{\Dscr}_\infty/\Xscr_k^{\mathrm{can}}}\xrightarrow{\phi_k}
\Xscr_k^{\mathrm{can}}$ be the stack obtained from
$\Xscr_k^{\mathrm{can}}$ by performing a $k$-th {root construction} along the divisor $\tilde{\Dscr}_\infty$. It is a two-dimensional projective toric orbifold with Deligne-Mumford torus $T_t$ and with coarse moduli space $\pi_k:=\pi_k^{\mathrm{can}}\circ \phi_k\colon \Xscr_k\to\bar{X}_k$ (cf.\ Section \ref{sec:toricstacks}). Its stacky fan is
$\stackyfank:=(N,\bar{\Sigma}_k,\beta)$, where $\beta\colon\Z^{k+2}\to
N$ is given by $\{\vec{v}_0,\vec{v}_1, \dots,\vec{v}_k,k\, \vec{v}_\infty\}$. 

As a global quotient stack, $\Xscr_k$ is isomorphic to
$[Z_{\bar{\Sigma}_k}/G_{\stackyfank}]$ where $Z_{\bar{\Sigma}_k}$ is
the same as for $\Xscr_k^{\mathrm{can}}$, since $\bar{\Sigma}_k$ is
the fan in both stacky fans. The group $G_{\stackyfank}$ is
$\Hom_\Z(\DG(\beta),\C^\ast)$, where $\DG(\beta)$ is the cokernel of
the group homomorphism $\beta^\ast\colon\Z^{2}\to\Z^{k+2}$. Then we
obtain $\DG(\beta)\simeq\Z^{k}$ and therefore
$G_{\stackyfank}\simeq(\C^\ast)^k$. By applying the functor $\Hom_\Z(-
,\C^\ast)$ to the quotient map $\beta^\vee\colon\Z^{k+2}\to
\DG(\beta)$ we obtain an injective morphism
$G_{\stackyfank}\to(\C^\ast)^{k+2}$ which is given by
\begin{equation*}
(t_1, \ldots, t_k) \ \longmapsto \ \left\{\begin{array}{ll}
\Big(\, \prod\limits_{i=1}^{k-1}\, t_i^i\,
t_k^{2k-k^2}\,,\,\prod\limits_{i=1}^{k-1}\, t_i^{-(i+1)}\,
t_k^{k^2}\,,\,t_1\,,\,\dots\,,\,t_k\, \Big)&\mbox{ for odd $k$}\ ,\\[8pt]
\Big(\, \prod\limits_{i=1}^{k-1}\, t_i^i\, t_k^{k-k\,
  \tilde{k}}\,,\,\prod\limits_{i=1}^{k-1}\, t_i^{-(i+1)}\, t_k^{k\,
  \tilde{k}}\,,\,t_1\,,\,\dots\,,\,t_k\, \Big)&\mbox{ for even $k$}\ .
\end{array}\right.
\end{equation*}
By restricting the standard action of $(\C^\ast)^{k+2}$ on $Z_{\bar{\Sigma}_k}\subset\C^{k+2}$ we obtain an action of $G_{\stackyfank}$ on $Z_{\bar{\Sigma}_k}$.  

As before, the boundary divisor $\Xscr_k\setminus T_t$ is a simple
normal crossing divisor with $k+2$ irreducible components
$\Dscr_0,\Dscr_1, \ldots,\Dscr_k,\Dscr_\infty$, which are the
effective Cartier divisors corresponding to the rays $\rho_0,\rho_1, \ldots,\rho_k,\rho_\infty$. By construction (cf.\  Section \ref{sec:rootstack}), we have
\begin{equation*}
\phi_k^\ast\big(\Ocal_{\Xscr_k^{\mathrm{can}}}(\tilde{\Dscr}_i) \big)\simeq\left\{
\begin{array}{ll}
\Ocal_{\Xscr_k}(\Dscr_i) & \mbox{   for   }i=0,1, \ldots,k\ ,\\[8pt]
\Ocal_{\Xscr_k}(k\, \Dscr_\infty) & \mbox{   for   }i=\infty\ .
\end{array}\right.
\end{equation*}

\begin{remark}\label{rem:intersectiontheory}
Since $\Xscr_k$ is a quotient stack, there is a well defined integral
intersection theory \cite{art:edidingraham1998}. As $\Xscr_k$ is
smooth, its rational Chow groups are isomorphic to the rational Chow
groups of $\bar{X}_k$ via ${\pi_k}_{\ast}$ by
\cite[Proposition~6.1]{art:vistoli1989}. In particular,
${\pi_k}_\ast(\Dscr_i)=D_i$ for $i=0, 1, \ldots, k$ and ${\pi_k}_\ast(\Dscr_\infty)= \frac{1}{k}\,D_\infty$.
\end{remark}
Recall that by Remark \ref{rem:mckayALE} the intersection matrix
$(D_i\cdot D_j)_{1\leq i,j\leq k-1}$ is minus the Cartan matrix $C$ of
the Dynkin diagram of type $A_{k-1}$. The matrix $C$ is not unimodular
and its inverse has matrix elements
\begin{equation}
\big(C^{-1}\big)^{ij}=\min(i,j)-\frac{i\, j}{k}\ .
\end{equation}
In $\Pic(\Xscr_k)_\Q$ we can define the classes 
\begin{equation}\label{eq:tautologicalclasses}
\omega_i:=-\sum_{j=1}^{k-1} \, \big(C^{-1}\big)^{ij}\, \Dscr_j
\end{equation}
for $i=1, \ldots, k-1.$ 
\begin{lemma}\label{lem:tautological}
The classes $\omega_i$ are integral combinations of $\Dscr_j$ for $i=1, \ldots, k-1$, $j=0,1, \ldots, k$ and $\Dscr_\infty$ in $\Pic(\Xscr_k)$.
\end{lemma}
\proof
We argue as in \cite[Section~5.2]{art:ciraficikashanipoorszabo2011}.
Let $\vec{v}_\infty=-\tilde{k}\, \vec{e}_1+a \, \vec{e}_2$ be the minimal generator of
the ray
$\rho_\infty$, where $a:=\tilde{k}-1\in\Z$ if $k$ is even and  $a:=k-2$ if $k$ is odd. Let us consider the following relations in $\Pic(\Xscr_k)$:
\begin{align}\label{eq:condition1}
0&=\mathrm{div}\big(\chi^{(1,0)} \big)=\Dscr_1+2\,\Dscr_2+\cdots +k\,\Dscr_k-\tilde{k}\, k\,\Dscr_\infty\ ,\\[4pt] \label{eq:condition2}
0&=\mathrm{div}\big(\chi^{(0,1)} \big)=\Dscr_0-\Dscr_2-\cdots - (k-1)\,\Dscr_k+a\, k\,\Dscr_\infty\ ,
\end{align}
where $\chi^{(1,0)},\chi^{(0,1)}$ are the characters of $T_t$ associated with $(1,0), (0,1)\in M$, respectively.
Since by definition
\begin{equation}
\omega_1=-\sum_{j=1}^{k-1}\, \frac{(k-j)}{k}\,
\Dscr_j\qquad\mbox{and}\qquad\omega_{k-1}=-\sum_{j=1}^{k-1}\,
\frac{j}{k}\, \Dscr_j\ ,
\end{equation}
we get $\omega_1=\Dscr_0+(-\tilde{k}\, (k-1)+a \, k)\, \Dscr_\infty$
and $\omega_{k-1}=\Dscr_k-\tilde{k}\, \Dscr_\infty$. For $i=2, \ldots,
k-2$ we have $\omega_{i}=\omega_{i-1}-\omega_{k-1}-\sum_{j=i}^{k-1}\,
\Dscr_j$, and the assertion follows.
\endproof
\begin{definition}
For $i=1, \ldots, k-1$, the $i$\emph{-th tautological line bundle} $\Rcal_i$ on $\Xscr_k$ is the line bundle associated to the Cartier divisor $\omega_i$.
\end{definition}
By the proof of Lemma \ref{lem:tautological}, the tautological line bundles $\Rcal_i$ for $i=1, \ldots, k-1$ can be written as
\begin{equation}\label{eq:tautologicalintegralrelations}
\Rcal_i=
\left\{
\begin{array}{ll}
 \Ocal_{\Xscr_k}\Big(\Dscr_0-\sum\limits_{j=1}^{i-1}\, (j-1)\,
 \Dscr_j-(i-1)\, \sum\limits_{j=i}^k\, \Dscr_j + (i-2)\, \tilde{k}\, \Dscr_\infty\Big) & \mbox{for } i=1,\ldots, k-2\ , \\[8pt]
\Ocal_{\Xscr_k}\big(\Dscr_k-\tilde{k}\, \Dscr_\infty\big) &  \mbox{for } i=k-1\ .
\end{array}\right.
\end{equation}
\begin{remark}
We call these line bundles ``tautological'' as they play the same role as the tautological line bundles considered by Kronheimer and Nakajima \cite{art:kronheimernakajima1990}. Indeed, their restrictions to $X_k\subset \Xscr_k$ yield exactly their tautological line bundles. Moreover, note that
\begin{equation}\label{eq:tautological}
\int_{\Xscr_k}\, \mathrm{c_1}(\Rcal_i)\cdot
\mathrm{c_1}(\Rcal_j)=\int_{\Xscr_k}\, \omega_i\cdot \omega_j =
-\big(C^{-1} \big)^{ij}\qquad\mbox{for } \ i,j=1,\ldots, k-1\ ,
\end{equation}
which is the same result as in \cite[Theorem A.7]{art:kronheimernakajima1990}.
\end{remark}
\begin{proposition}\label{prop:picardstack}
The Picard group $\Pic(\Xscr_k)$ of $\Xscr_k$ is freely generated over $\Z$\, by $\Ocal_{\Xscr_k}(\Dscr_\infty)$ and $\Rcal_i$\, with $i=1, \ldots, k-1$.
\end{proposition}
\proof
By Section \ref{sec:rootdivisor}, any line bundle $\Lcal$ on $\Xscr_k$
is of the form $\phi_k^\ast(\Mcal)\otimes \Ocal_{\Xscr_k}(m\,
\Dscr_\infty)$ for $\Mcal$ a line bundle on $\Xscr^{\mathrm{can}}_k$
and $m$ an integer such that $0\leq m\leq k-1$. Moreover $m$ is unique
and $\Mcal$ is unique up to isomorphism. By the short exact sequence
\eqref{eq:picardcanonical}, the line bundle $\Mcal$ is given as an
integral combination of $\Ocal_{\Xscr_k}(\tilde{\Dscr}_i)$ for $i=1,
\ldots, k-1$ and $\Ocal_{\Xscr_k}(\tilde{\Dscr}_\infty).$\footnote{The
  divisors $\tilde{\Dscr}_0$ and $\tilde{\Dscr}_k$ can be expressed in
  terms of the divisors $\tilde{\Dscr}_i$ for $i=1,\ldots, k-1,
  \infty$ in $\Pic(\Xscr_k^{\mathrm{can}})$. For example, it is enough
  to use the explicit expressions for $\mathrm{div}(\chi^{(0,1)})$ and
  $\mathrm{div}(\chi^{(1,1)})$ in $\mathrm{Div}_{T_t}(\bar{X}_k)$ in
  terms of the divisors $D_i$ for $i=0, 1, \ldots, k,\infty$.} Therefore the line bundle $\Lcal$ is an integral combination of $\Ocal_{\Xscr_k}(\Dscr_i)$ for $i=1, \ldots, k-1$ and $\Ocal_{\Xscr_k}(\Dscr_\infty)$. Since the line bundles $\Ocal_{\Xscr_k}(\Dscr_j)$ for $j=0, \ldots, k$ can be given as integral combinations of $\Rcal_i$ for $i=1, \ldots, k-1$, the claim follows.
\endproof

\begin{remark}\label{rem:rootlattice-cohomology}  
There is a relation between line bundles on $\Xscr_k$ and elements of
the \emph{root lattice $\Qfrak$ of type $A_{k-1}$}. As explained in
\cite[Section~4]{art:kronheimer1989}, the cohomology group $H^2(X_k;
\R)\simeq \Pic(X_k)\otimes_{\Z}\R$ can be identified with the
\emph{real Cartan subalgebra} $\hfrak$ associated with the Dynkin
diagram of type $A_{k-1}$. In this picture, $H_2(X_k; \Z)$ can be
identified with the {root lattice $\Qfrak$ of type $A_{k-1}$}. Under
this correspondence, the classes $[D_1], \ldots, [D_{k-1}]$ are the
simple roots. Since $\Pic(\Xscr_k)$ has no torsion, the map
$\jmath\colon\Pic(\Xscr_k)\to \Pic(\Xscr_k)\otimes_{\Z}\R$ is
injective. Consider the restriction map
$\imath^\ast\colon\Pic(\Xscr_k)\otimes_{\Z}\R\to
\Pic(X_k)\otimes_{\Z}\R$ with respect to the inclusion morphism
$\imath \colon X_k \to \Xscr_k$. The map $\imath^\ast$ is surjective because of Equation \eqref{eq:tautologicalclasses}.

For $j=1, \ldots, k-1$, let $\gamma_j$ be a simple root in $\Qfrak$,
and $\gamma:=\sum_{i=1}^{k-1}\, y_i\, \gamma_i$ an element of the root
lattice. Under the correspondence described above, this gives a linear
combination of the divisors $\sum_{i=1}^{k-1}\, y_i\, [D_i]$ and
singles out a unique line bundle $\Ocal_{X_k}\big(\sum_{i=1}^{k-1}\,
y_i \, D_i)$. After fixing an integer $u_\infty\in\Z$ and setting
$\vec{u}=-C\vec{y}$, where $\vec{y}:=(y_1, \ldots, y_{k-1})$, we can
define a line bundle $\bigotimes_{i=1}^{k-1}\, \Rcal_i^{\otimes u_i}\otimes\Ocal_{\Xscr_k}(\Dscr_\infty)^{\otimes u_\infty}$ such that
\begin{equation}
\bigotimes_{i=1}^{k-1}\, \Rcal_i^{\otimes
  u_i}\otimes\Ocal_{\Xscr_k}(u_\infty \, \Dscr_\infty)_{\vert X_k}
\simeq \Ocal_{X_k}\Big(\, \mbox{$\sum\limits_{i=1}^{k-1}$}\, y_i\, D_i\, \Big)\ .
\end{equation}
\end{remark}

\subsubsection*{Characterization of $\Dscr_\infty$}

All results presented below hold for every $k\geq 2$, but the proofs
and computations are given only for the cases $k>2$. For $k=2$ there
are some differences; since $\bar X_2$ is the second Hirzebruch surface,
for that case one can use the calculations done in  \cite[Appendix~D]{art:bruzzosala2013}
for the  stacky Hirzebruch surfaces.

By construction (cf.\ Section \ref{sec:rootstack}), $\Dscr_\infty$ is
isomorphic to the root stack
$\sqrt[k]{\Ocal_{\Xscr_k^{\mathrm{can}}}(\tilde{\Dscr}_\infty)_{\vert\tilde{\Dscr}_\infty}/\tilde{\Dscr}_\infty}$. By
\cite[Theorem~6.25-(1)]{art:fantechimannnironi2010} the divisor
$\Dscr_\infty$ is a projective toric Deligne-Mumford stack with
Deligne-Mumford torus $\Tscr\simeq T_t\times\Bscr\mu_k$. Its stacky
fan is the quotient stacky fan
$\stackyfank/\rho_\infty:=(N(\rho_\infty),\bar{\Sigma}_k/\rho_\infty,\beta(\rho_\infty))$,
where $N(\rho_\infty):=N/k\, \Z \vec{v}_\infty\simeq\Z\oplus \Z_k$ and the
quotient fan $\bar{\Sigma}_k/\rho_\infty\subset
N(\rho_\infty)\otimes_{\Z}\mathbb{Q}\simeq\mathbb{Q}$ is the same as
the quotient fan \eqref{eq:stackyfaninfty} of $\tilde{\Dscr}_\infty$. The quotient map $N\to N(\rho_\infty)\simeq\Z\oplus\Z_k$ is given by
\begin{equation}
\begin{pmatrix}
  1-\tilde{k} & -\tilde{k}\\
  -1 & -1
  \end{pmatrix} \ 
\mbox{ for even $k$}\qquad \mbox{and} \qquad \begin{pmatrix}
  k-2 & k\\
  -\frac{k-1}{2} & -\frac{k+1}{2}
  \end{pmatrix} \ \mbox{  for odd $k$}\ .
\end{equation}
On the other hand, the map $\beta(\rho_\infty)\colon\Z^2\to N(\rho_\infty)\simeq\Z\oplus\Z_k$ is given by the matrix
\begin{equation}
M\big(\beta(\rho_\infty)\big)=\begin{pmatrix}
  \tilde{k} & -\tilde{k}\\
  -1 & -1
  \end{pmatrix} \ 
\mbox{ for even $k$}\qquad \mbox{and} \qquad M\big(\beta(\rho_\infty)\big)=\begin{pmatrix}
  k & -k\\
  \frac{k-1}{2} & \frac{k-1}{2}
  \end{pmatrix} \ \mbox{  for odd $k$}\ .
\end{equation}

By Section \ref{sec:toricstacks}, the toric stack $\Dscr_\infty$ is an essentially trivial gerbe with banding group $\Hom_\Z(N(\rho_\infty)_{\mathrm{tor}},$ $\C^*)\simeq\mu_k$ over its rigidification $\Dscr_\infty^{\mathrm{rig}}$, which is $\tilde{\Dscr}_\infty$. Let $\tilde{\phi}_k:=(\phi_k)_{\vert\Dscr_\infty}\colon\Dscr_\infty\to\tilde{\Dscr}_\infty$ be the $\mu_k$-gerbe structure morphism. Then $r_k:=\tilde{\pi}_k\circ\tilde{\phi}_k\colon\Dscr_\infty\to D_\infty\simeq\PP^1$ is the coarse moduli space of $\Dscr_\infty$.
\begin{proposition}\label{prop:gerbestructure}
The toric stack $\Dscr_\infty$ is isomorphic to the global quotient stack
\begin{equation}
\left[\frac{\C^2\setminus\{0\}}{\C^\ast\times\mu_k}\right]\ ,
\end{equation}
where the action is given by
\begin{equation}\label{eq:Tactiongerbe}
(t,\omega)\triangleright (z_1,z_2)=\left\{
\begin{array}{ll}
\big(t^{\tilde{k}}\, \omega \, z_1\,,\,t^{\tilde{k}}\, \omega^{-1}\,
z_2 \big) & \mbox{ for even $k$}\ ,\\
\big(t^k\, \omega^{\frac{k+1}{2}}\, z_1\,,\,t^k\,
\omega^{\frac{k-1}{2}}\, z_2\big) & \mbox{ for odd $k$}\ ,
\end{array}
\right.
\end{equation}
for $(t,\omega)\in\C^\ast\times\mu_k$ and $(z_1,z_2)\in\C^2\setminus\{0\}$.
\end{proposition}
\proof
By arguing along the lines of the analogous characterization for $\tilde{\Dscr}_\infty$ (cf.\ Proposition \ref{prop:quotientcanonicalinfty}), we get that $\Dscr_\infty$ is isomorphic to $[Z_{\bar{\Sigma}_k/\rho_\infty}/G_{\stackyfank/\rho_\infty}]$, where $Z_{\bar{\Sigma}_k/\rho_\infty}=\C^2\setminus\{0\}$ is the same as for $\tilde{\Dscr}_\infty$. The group $G_{\stackyfank/\rho_\infty}$ is $\Hom_\Z(\DG(\beta(\rho_\infty)),\C^\ast)$. Since $N(\rho_\infty)$ has torsion, $\DG(\beta(\rho_\infty))$ is obtained in the following way. Consider a free resolution of $N(\rho_\infty)$
\begin{equation}
0\ \longrightarrow\ \Z \ \xrightarrow{ \ Q \ } \ \Z^2 \
\longrightarrow \ N(\rho_\infty)\simeq\Z\oplus \Z_k \ \longrightarrow \ 0
\end{equation}
where $Q$ is the map which sends $1\in\Z$ to $k\,
\vec{e}_2\in\Z^2$. Consider also a lifting $B\colon\Z^2\to\Z^2$ of
$\beta(\rho_\infty)$, so that $B$ can be represented by the matrix
$M(\beta(\rho_\infty))$. Define the map $[B\, Q]:\Z^3\to\Z^2$ by
adding the column $Q$ to the matrix of $B$. Then
$\DG(\beta(\rho_\infty))=\mathrm{Coker}([B\, Q]^\ast)$ and $[B\, Q]^\ast$ is given by the matrix
\begin{equation}
H=
\begin{pmatrix}
\tilde{k} & -1\\
-\tilde{k} & -1\\
0 & k
\end{pmatrix} \ 
\mbox{ for even $k$}\qquad \mbox{and} \qquad
H=\begin{pmatrix}
k & \frac{k-1}{2}\\
-k & \frac{k-1}{2}\\
0 & k
\end{pmatrix} \ 
\mbox{ for odd $k$}\ .
\end{equation}
In both cases, $H$ is equivalent to the matrix
\begin{equation}
K=\begin{pmatrix}
1 & 0\\
0 & k\\
0 & 0
\end{pmatrix} \ ,
\end{equation}
in the sense that there exist two unimodular matrices $T\in GL(3,\Z)$
and $P\in GL(2,\Z)$ such that $H=T\, K\, P$. Hence we have
$\DG(\beta(\rho_\infty))\simeq\Z\oplus \Z_k$ and
$G_{\stackyfank/\rho_\infty}\simeq\C^\ast\times\mu_k$. The action of
$\C^\ast\times\mu_k$ on $\C^2\setminus\{0\}$ is given by composition
of the standard $(\C^\ast)^2$-action with the map
$\C^\ast\times\mu_k\to(\C^\ast)^2$ obtained by applying the functor
$\Hom_\Z(-,\C^\ast)$ to the composition $\Z^2\hookrightarrow\Z^3\to
\DG(\beta(\rho_\infty))\simeq\Z\oplus \Z_k$, where the second map is the quotient map. This gives the assertion.
\endproof
Recall that $\Pic(\Dscr_\infty)\simeq \DG(\beta(\rho_\infty))$ (cf.\
Remark \ref{rem:pictoric}). The explicit characterization of the
global quotient structure of $\Dscr_\infty$ then yields the following result.
\begin{corollary}
The Picard group $\Pic(\Dscr_\infty)$ of $\Dscr_\infty$ is isomorphic
to $\Z\oplus\Z_k$. It is generated by the line bundles $\Lcal_1$ and
$\Lcal_2$ corresponding respectively to the two characters of
$G_{\stackyfank/\rho_\infty}\simeq\C^*\times\mu_k$ given by
\begin{equation}
\chi_1\,\colon\, (t,\omega)\in\C^\ast\times\mu_k \ \longmapsto \
t\in\C^\ast\qquad \mbox{and}\qquad \chi_2\, \colon\,
(t,\omega)\in\C^\ast\times\mu_k\ \longmapsto\ \omega\in\C^\ast\ .
\end{equation}
In particular, $\Lcal_2^{\otimes k}$ is trivial.
\end{corollary}

By the commutative diagram of \cite[Equation (6.28)]{art:fantechimannnironi2010} we also know that $\Pic(\Dscr_\infty)$ fits into a commutative diagram
\begin{equation}\label{eq:shortexactpic}
\begin{aligned}
  \begin{tikzpicture}[xscale=2.9,yscale=-0.8]
    \node (A0_0) at (0.4, 0) {$0$};
    \node (A0_1) at (1, 0) {$\Z$};
    \node (A0_2) at (2, 0) {$\Z$};
    \node (A0_3) at (3, 0) {$\Z_k$};
    \node (A0_4) at (3.6, 0) {$0$};
    \node (A1_5) at (3.75, 1) {$,$};
    \node (A2_0) at (0.4, 2) {$0$};
    \node (A2_1) at (1, 2) {$\Pic(\tilde{\Dscr}_\infty)$};
    \node (A2_2) at (2, 2) {$\Pic(\Dscr_\infty)$};
    \node (A2_3) at (3, 2) {$\Z_k$};
    \node (A2_4) at (3.6, 2) {$0$};    
    \path (A0_0) edge [->]node [left] {$\scriptstyle{}$} (A0_1);
    \path (A0_1) edge [->]node [above] {$\scriptstyle{\cdot\ k}$} (A0_2);
    \path (A0_2) edge [->]node [left] {$\scriptstyle{}$} (A0_3);
    \path (A0_3) edge [->]node [left] {$\scriptstyle{}$} (A0_4);

    \path (A2_0) edge [->]node [left] {$\scriptstyle{}$} (A2_1);
    \path (A2_1) edge [->]node [above] {$\scriptstyle{\tilde{\phi}_k^\ast}$} (A2_2);
    \path (A2_2) edge [->]node [left] {$\scriptstyle{}$} (A2_3);
    \path (A2_3) edge [->]node [left] {$\scriptstyle{}$} (A2_4);

    \path (A0_1) edge [->]node [auto] {$\scriptstyle{\tilde{f}}$} (A2_1);
    \path (A0_2) edge [->]node [auto] {$\scriptstyle{f}$} (A2_2);    
    \path (A0_3) edge [->]node [auto] {$\scriptstyle{\mathrm{id}}$} (A2_3);    
  \end{tikzpicture} 
 \end{aligned}
\end{equation}
where the morphisms $\tilde{f}$ and $f$ send $1\mapsto
\Ocal_{\Xscr_k^{\mathrm{can}}}(\tilde{\Dscr}_\infty)_{\vert\tilde{\Dscr}_\infty}$
and $1\mapsto \Ocal_{\Xscr_k}(\Dscr_\infty)_{\vert\Dscr_\infty}$,
respectively. From this commutative diagram, one can argue that every
line bundle $\Lcal$ on $\Dscr_\infty$ can be written as a tensor
product $\tilde{\phi}_k^\ast(\Ncal)\otimes\Ocal_{\Xscr_k}(\ell\, \Dscr_\infty)_{\vert\Dscr_\infty}$ for a line bundle $\Ncal$ on $\tilde{\Dscr}_\infty$ and $0\leq \ell<k$ an integer.

Now we characterize the restrictions of line bundles from $\Xscr_k$ to $\Dscr_\infty$.
\begin{lemma}\label{lem:restrictionDinfty}
The restriction of $\Ocal_{\Xscr_k}(\Dscr_\infty)$ to $\Dscr_\infty$
is isomorphic to $\Lcal_1$. For even $k$ one has
\begin{equation}
\Ocal_{\Xscr_k}(\Dscr_0)_{\vert\Dscr_\infty}\simeq\Lcal_1^{\otimes\tilde{k}}\otimes\Lcal_2\qquad\mbox{and}\qquad\Ocal_{\Xscr_k}(\Dscr_k)_{\vert\Dscr_\infty}\simeq\Lcal_1^{\otimes\tilde{k}}\otimes\Lcal_2^{\otimes-1}\ ,
\end{equation}
while for odd $k$ one has
\begin{equation}
\Ocal_{\Xscr_k}(\Dscr_0)_{\vert\Dscr_\infty}\simeq\Lcal_1^{\otimes k}\otimes\Lcal_2^{\otimes\frac{k+1}{2}}\qquad\mbox{and}\qquad\Ocal_{\Xscr_k}(\Dscr_k)_{\vert\Dscr_\infty}\simeq\Lcal_1^{\otimes k}\otimes\Lcal_2^{\otimes\frac{k-1}{2}}\ .
\end{equation}
\end{lemma}
\proof
The proof is analogous to that of Lemma \ref{lem:restrictioncanonical}. We need only to point out that the analogue of the map $\phi$ in this case 
is not uniquely determined by just imposing the commutativity of the
diagram analogous to the diagram \eqref{eq:diagchinaII}. To compute that map,  one has to follow the proof of Lemma \ref{lem:galedualsequences}. 
In particular, in our case we find that it is the identity.
\endproof
If one works out some of the details of the proof of this
Lemma and in particular writes   explicitly some of the maps in the
required commutative diagrams, one can for example obtain that $\Ocal_{\Xscr_k}(\Dscr_i)_{\vert\Dscr_\infty}\simeq\Ocal_{\Dscr_\infty}$ for $i=1,\ldots,k-1$. Moreover, if we denote by $p_0$ and $p_\infty$ respectively the divisors $\tilde{\phi}_k^{-1}(\tilde{p}_0)_{\mathrm{red}}$ and $\tilde{\phi}_k^{-1}(\tilde{p}_\infty)_{\mathrm{red}}$ in $\Dscr_\infty$ corresponding to $\rho_0'$ and $\rho_\infty'$, from the explicit form of the Gale dual of $\beta(\rho_\infty)$ one has the following result.
\begin{corollary}\label{cor:decomposition}
For even $k$ we have
\begin{equation}
\Ocal_{\Dscr_\infty}(p_0)\simeq\Lcal_1^{\otimes\tilde{k}}\otimes\Lcal_2
\qquad\mbox{and}\qquad
\Ocal_{\Dscr_\infty}(p_\infty)\simeq\Lcal_1^{\otimes\tilde{k}}\otimes\Lcal_2^{\otimes-1}\ .
\end{equation}
For odd $k$ we have
\begin{equation}
\Ocal_{\Dscr_\infty}(p_0)\simeq\Lcal_1^{\otimes k}\otimes\Lcal_2^{\otimes\frac{k+1}{2}}
\qquad\mbox{and}\qquad
\Ocal_{\Dscr_\infty}(p_\infty)\simeq\Lcal_1^{\otimes k}\otimes\Lcal_2^{\otimes\frac{k-1}{2}}.
\end{equation}
In particular, for any $k\geq2$ we have $\Ocal_{\Dscr_\infty}(p_0)\simeq\Ocal_{\Xscr_k}(\Dscr_0)_{\vert\Dscr_\infty}$ and $\Ocal_{\Dscr_\infty}(p_\infty)\simeq\Ocal_{\Xscr_k}(\Dscr_k)_{\vert\Dscr_\infty}$.
\end{corollary}
\begin{remark}
This corollary makes it clear that the line bundles associated with
the torus-invariant divisors are not sufficient to generate the whole Picard group of the gerbe $\Dscr_\infty$. This is evident if we consider the exact sequence \eqref{eq:longexactgale}, which in our case becomes
\begin{equation}
\Z \ \xrightarrow{\beta(\rho_\infty)^\ast} \ \Z^2 \
\xrightarrow{\beta(\rho_\infty)^\vee} \ \Pic(\Dscr_\infty) \
\longrightarrow \ \Ext^1_\Z\big(N(\rho_\infty),\Z \big)\simeq\Z_k \
\longrightarrow \ 0\ ,
\end{equation}
and our statement is equivalent to the fact that  $\beta(\rho_\infty)^\vee$ is not surjective. 
\end{remark}
\begin{corollary}\label{cor:restrictiontautological}
The restrictions of the tautological line bundles $\Rcal_i$ on
$\Xscr_k$  to $\Dscr_\infty$ are given by
\begin{equation}
{\Rcal_i}_{\vert\Dscr_\infty}\simeq\Lcal_2^{\otimes i}\quad\mbox{for
  even $k$}\qquad \mbox{and} \qquad
{\Rcal_i}_{\vert\Dscr_\infty}\simeq\Lcal_2^{\otimes i\, \frac{k+1}{2}}\quad\mbox{for odd $k$}\ .
\end{equation}
\end{corollary}
\proof 
By using Equation \eqref{eq:tautologicalintegralrelations} we get
\begin{align}
{\Rcal_i}_{\vert\Dscr_\infty}&\simeq 
\Ocal_{\Xscr_k}\big(\Dscr_0-(i-1)\, \Dscr_k+(i-2)\, \tilde{k}\,
\Dscr_\infty\big)_{\vert\Dscr_\infty} \qquad \mbox{ for } \
i=1,\ldots,k-2 \ , \\[4pt]
{\Rcal_{k-1}}_{\vert\Dscr_\infty}&\simeq 
\Ocal_{\Xscr_k}\big(\Dscr_k-\tilde{k}\,
\Dscr_\infty\big)_{\vert\Dscr_\infty} \ .
\end{align}
The result now follows from Lemma \ref{lem:restrictionDinfty}. 
\endproof
Finally, we need to relate the Picard groups of $\tilde{\Dscr}_\infty$
and $\Dscr_\infty$, in particular making the map $\tilde{\phi}_k^*$ in
the diagram \eqref{eq:shortexactpic} explicit. At this stage, from the commutativity of the diagram we only know that
\begin{equation}
\tilde{\phi}_k^\ast\Ocal_{\Xscr_k^{\mathrm{can}}}(\tilde{\Dscr}_\infty)_{\vert\tilde{\Dscr}_\infty}\simeq\Ocal_{\Xscr_k}(k\,
\Dscr_\infty)_{\vert\Dscr_\infty}\simeq\Lcal_1^{\otimes k} \ .
\end{equation}
\begin{proposition}
One has isomorphisms
\begin{equation}
\tilde{\phi}_k^\ast\Ocal_{\tilde{\Dscr}_\infty}(\tilde{p}_0)\simeq\Ocal_{\Dscr_\infty}(p_0)
\qquad\mbox{and}\qquad
\tilde{\phi}_k^\ast\Ocal_{\tilde{\Dscr}_\infty}(\tilde{p}_\infty)\simeq\Ocal_{\Dscr_\infty}(p_\infty)\ .
\end{equation}
For odd $k$ we have
$\tilde{\phi}_k^\ast\Ocal_{\tilde{\Dscr}_\infty}(\tilde{p})\simeq\Lcal_2$
while for even $k$ we have $\tilde{\phi}_k^\ast\Ocal_{\tilde{\Dscr}_\infty}(\tilde{p})\simeq\Lcal_2^{\otimes 2}$.
\end{proposition}
\proof
To give an explicit form to the map $\tilde{\phi}_k^\ast$ in the
diagram \eqref{eq:shortexactpic} we use the commutative diagram in \cite[Equation (7.21)]{art:fantechimannnironi2010}:
\begin{equation}\label{eq:diagramfantechi}
\begin{aligned}
  \begin{tikzpicture}[xscale=3.7,yscale=-1.4]
    \node (A3_5) at (3.75, 3) {$,$};
    
    \node (A1_1) at (1, 1) {$0$};
    \node (A1_2) at (2, 1) {$0$};
    \node (A1_3) at (3, 1) {$0$};
    
    \node (A2_0) at (0.4, 2) {$0$};
    \node (A2_1) at (1, 2) {$\Z$};
    \node (A2_2) at (2, 2) {$\Z^2$};
    \node (A2_3) at (3, 2) {$\Z$};
    \node (A2_4) at (3.6, 2) {$0$};    

    \node (A3_0) at (0.4, 3) {$0$};
    \node (A3_1) at (1, 3) {$\Z^2$};
    \node (A3_2) at (2, 3) {$\Z^3$};
    \node (A3_3) at (3, 3) {$\Z$};
    \node (A3_4) at (3.6, 3) {$0$};

    \node (A4_0) at (0.4, 4) {$0$};
    \node (A4_1) at (1, 4) {$\DG\big(\beta^{\mathrm{can}}(\rho_\infty)
      \big)$};
    \node (A4_2) at (2, 4) {$\DG\big(\beta(\rho_\infty) \big)$};
    \node (A4_3) at (3, 4) {$\Z_k$};
    \node (A4_4) at (3.6, 4) {$0$};

    \node (A5_1) at (1, 5) {$0$};
    \node (A5_2) at (2, 5) {$0$};
    \node (A5_3) at (3, 5) {$0$};

    \path (A2_0) edge [->]node [left] {$\scriptstyle{}$} (A2_1);
    \path (A2_1) edge [->]node [above] {$\scriptstyle{}$} (A2_2);
    \path (A2_2) edge [->]node [left] {$\scriptstyle{}$} (A2_3);
    \path (A2_3) edge [->]node [left] {$\scriptstyle{}$} (A2_4);

    \path (A3_0) edge [->]node [left] {$\scriptstyle{}$} (A3_1);
    \path (A3_1) edge [->]node [above] {$\scriptstyle{}$} (A3_2);
    \path (A3_2) edge [->]node [left] {$\scriptstyle{}$} (A3_3);
    \path (A3_3) edge [->]node [left] {$\scriptstyle{}$} (A3_4);

    \path (A4_0) edge [->]node [left] {$\scriptstyle{}$} (A4_1);
    \path (A4_1) edge [->]node [above] {$\scriptstyle{\tilde{\phi}_k^\ast}$} (A4_2);
    \path (A4_2) edge [->]node [left] {$\scriptstyle{}$} (A4_3);
    \path (A4_3) edge [->]node [left] {$\scriptstyle{}$} (A4_4);

    \path (A1_1) edge [->]node [auto] {$\scriptstyle{}$} (A2_1);
    \path (A1_2) edge [->]node [auto] {$\scriptstyle{}$} (A2_2);    
    \path (A1_3) edge [->]node [auto] {$\scriptstyle{}$} (A2_3);    

    \path (A2_1) edge [->]node [auto] {$\scriptstyle{\beta^{\mathrm{can}}(\rho_\infty)^\ast}$} (A3_1);
    \path (A2_2) edge [->]node [auto] {$\scriptstyle{[B\, Q]^\ast}$} (A3_2);    
    \path (A2_3) edge [->]node [auto] {$\scriptstyle{}$} (A3_3);    

    \path (A3_1) edge [->]node [auto] {$\scriptstyle{\beta^{\mathrm{can}}(\rho_\infty)^\vee}$} (A4_1);
    \path (A3_2) edge [->]node [auto] {$\scriptstyle{\pi_\beta}$} (A4_2);    
    \path (A3_3) edge [->]node [auto] {$\scriptstyle{}$} (A4_3);    

    \path (A4_1) edge [->]node [auto] {$\scriptstyle{}$} (A5_1);
    \path (A4_2) edge [->]node [auto] {$\scriptstyle{}$} (A5_2);    
    \path (A4_3) edge [->]node [auto] {$\scriptstyle{}$} (A5_3);    
  \end{tikzpicture} 
 \end{aligned}
\end{equation}
where $\pi_\beta$ is the quotient projection
$\Z^3\to\mathrm{Coker}([B\, Q]^\ast)$. By commutativity of the diagram, the map $\tilde{\phi}_k^\ast\colon \DG(\beta^{\mathrm{can}}(\rho_\infty))\simeq\Z\oplus\Z_{\tilde{k}}\to \DG(\beta(\rho_\infty))\simeq\Z\oplus\Z_k$ is represented by the matrix
\begin{equation}
\begin{pmatrix}
\tilde{k} & 0 \\
-1 & 2
\end{pmatrix} \ \mbox{for even $k$}\qquad \mbox{and} \qquad
\begin{pmatrix}
k & 0 \\
\frac{k-1}{2} & 1
\end{pmatrix} \ \mbox{for odd $k$}\ .
\end{equation}
The result now follows by taking the images of the vectors
$(1,0),(1,1),(0,1)$ in
$\Pic(\tilde{\Dscr}_\infty)\simeq\Z\oplus\Z_{\tilde{k}}$ which
correspond respectively to the line bundles $\Ocal_{\tilde{\Dscr}_\infty}(\tilde{p}_0),\Ocal_{\tilde{\Dscr}_\infty}(\tilde{p}_\infty),\Ocal_{\tilde{\Dscr}_\infty}(\tilde{p})$.
\endproof
\begin{remark}$ $
\begin{itemize}
\item For odd $k$, since $\tilde{\phi}_k^\ast\Ocal_{\tilde{\Dscr}_\infty}(\tilde{p})\simeq\Lcal_2$, the line bundle $\tilde{\phi}_k^\ast\Ocal_{\tilde{\Dscr}_\infty}(\tilde{p})$ generates the torsion part of the Picard group $\Pic(\Dscr_\infty)$ of $\Dscr_\infty$. For even $k$ this is not true since $\tilde{\phi}_k^\ast\Ocal_{\tilde{\Dscr}_\infty}(\tilde{p})\simeq\Lcal_2^{\otimes 2}$ is not sufficient to generate   the torsion part of $\Pic(\Dscr_\infty)$.
\item Following the proof of
  \cite[Proposition~7.20]{art:fantechimannnironi2010}, the last short
  exact sequence in the diagram \eqref{eq:diagramfantechi} is an
  element of $\Ext_\Z^1(N_{\mathrm{tor}},\Pic(\tilde{\Dscr}_\infty))$,
  which by \cite[Proposition~6.9]{art:fantechimannnironi2010} induces
  an element $[\Ocal_{\Xscr_k^{\rm can}}(\tilde{\Dscr}_\infty)_{\vert \tilde{\Dscr}_\infty}]\in\Pic(\tilde{\Dscr}_\infty)/k\Pic(\tilde{\Dscr}_\infty)$. The last column of the diagram is a free (hence projective) resolution of $\Z_k$, so we can lift the identity map of $\Z_k$ to obtain a morphism of short exact sequences
\begin{equation}\label{eq:morshortexact}
\begin{aligned}
  \begin{tikzpicture}[xscale=2.9,yscale=-0.8]
    \node (A0_0) at (0.4, 0) {$0$};
    \node (A0_1) at (1, 0) {$\Z$};
    \node (A0_2) at (2, 0) {$\Z$};
    \node (A0_3) at (3, 0) {$\Z_k$};
    \node (A0_4) at (3.6, 0) {$0$};
    \node (A1_5) at (3.75, 1) {$.$};
    \node (A2_0) at (0.4, 2) {$0$};
    \node (A2_1) at (1, 2) {$\Pic(\tilde{\Dscr}_\infty)$};
    \node (A2_2) at (2, 2) {$\Pic(\Dscr_\infty)$};
    \node (A2_3) at (3, 2) {$\Z_k$};
    \node (A2_4) at (3.6, 2) {$0$};    
    \path (A0_0) edge [->]node [left] {$\scriptstyle{}$} (A0_1);
    \path (A0_1) edge [->]node [above] {$\scriptstyle{\cdot\ k}$} (A0_2);
    \path (A0_2) edge [->]node [left] {$\scriptstyle{}$} (A0_3);
    \path (A0_3) edge [->]node [left] {$\scriptstyle{}$} (A0_4);

    \path (A2_0) edge [->]node [left] {$\scriptstyle{}$} (A2_1);
    \path (A2_1) edge [->]node [above] {$\scriptstyle{\tilde{\phi}_k^\ast}$} (A2_2);
    \path (A2_2) edge [->]node [left] {$\scriptstyle{}$} (A2_3);
    \path (A2_3) edge [->]node [left] {$\scriptstyle{}$} (A2_4);

    \path (A0_1) edge [->]node [auto] {$\scriptstyle{\tilde{f}}$} (A2_1);
    \path (A0_2) edge [->]node [auto] {$\scriptstyle{f}$} (A2_2);    
    \path (A0_3) edge [->]node [auto] {$\scriptstyle{\mathrm{id}}$} (A2_3);    
  \end{tikzpicture} 
 \end{aligned}
\end{equation}
The choices of liftings $\tilde{f}$ and $f$ are not unique. In
particular the choice of $\tilde{f}$ corresponds to a choice of a line
bundle in the class
$[\Ocal_{\Xscr_k^{\rm can}} (\tilde{\Dscr}_\infty)_{\vert \tilde{\Dscr}_\infty}]\in\Pic(\tilde{\Dscr}_\infty)/k\Pic(\tilde{\Dscr}_\infty)$,
while the choice of $f$ is equivalent to a choice of a line bundle in
the class $[\Ocal_{\Xscr_k}
(\Dscr_\infty)_{\vert \Dscr_\infty}]\in\Pic(\Dscr_\infty)/k\Pic(\Dscr_\infty)$. Note that by choosing
$\Ocal_{\Xscr_k^{\rm
    can}}(\tilde{\Dscr}_\infty)_{\vert \tilde{\Dscr}_\infty}$ and
$\Ocal_{\Xscr_k} (\Dscr_\infty)_{\vert \Dscr_\infty}$, respectively,  the diagram \eqref{eq:morshortexact} coincides with 
 the diagram \eqref{eq:shortexactpic}.
\end{itemize}
\end{remark}

We conclude this section by computing the \emph{degree}, i.e., the
integral of the first Chern class, of line bundles on
$\Dscr_\infty$. This result will be used in the Section~\ref{sec:framedsheaves}.
\begin{lemma}\label{lem:degree}
The degree of a line bundle $\Lcal=\Lcal_1^{\otimes
  a}\otimes\Lcal_2^{\otimes b}$ on $\Dscr_\infty$ with
$a,b\in\mathbb{Z}$ is given by
\begin{equation}
\int_{\Dscr_\infty}\, \crm_1(\Lcal)=\frac{a}{k\,\tilde{k}^2}\ ,
\end{equation}
where by $\int_{\Dscr_\infty}\colon A^\ast(\Dscr_\infty)_\Q\to\Q$ we denote the pushforward $A^\ast(\Dscr_\infty)_\Q\to A^\ast(\Spec(\C))_\Q\simeq\Q$.
\end{lemma}
\proof
First observe that for any $a,b\in\Z$ we have
\begin{equation}
\Lcal^{\otimes k\, \tilde{k}}\simeq\big(\Lcal_1^{\otimes
  a}\otimes\Lcal_2^{\otimes b} \big)^{\otimes k\,
  \tilde{k}}\simeq\Lcal_1^{\otimes a\, k\,
  \tilde{k}}\simeq\Ocal_{\Dscr_\infty}(a\, k\, p_\infty)\ .
\end{equation}
Since $\Dscr_\infty$ is smooth, by  \cite[Proposition
6.1]{art:vistoli1989} the structure map $r_k\colon\Dscr_\infty\to
D_\infty$ induces an isomorphism ${r_k}_\ast\colon
A^\ast(\Dscr_\infty)_{\Q}\xrightarrow{\sim}A^\ast(D_\infty)_{\Q}\simeq
A^\ast(\PP^1)_{\Q}$, and therefore
\begin{equation}
\int_{\Dscr_\infty}\, \crm_1\big(\Lcal^{\otimes k\, \tilde{k} } \big)=
\int_{\Dscr_\infty}\, \crm_1\big(\Ocal_{\Dscr_\infty}(a\, k\,
p_\infty) \big)=\int_{D_\infty}\,
{r_k}_\ast\crm_1\big(\Ocal_{\Dscr_\infty}(a\, k \, p_\infty) \big) \ .
\end{equation}
By  \cite[Example 6.7]{art:vistoli1989}, we obtain
${r_k}_\ast[p_\infty]=\frac{1}{d_\infty}\, [\infty]$, where $d_\infty$
is the order of the stabilizer of the point $p_\infty$. By using the
quotient presentation of $\Dscr_\infty$ in Proposition
\ref{prop:gerbestructure}, we find $d_\infty=k\, \tilde{k}$ and hence
\begin{equation}\label{eq:degreeDinfty}
\int_{\Dscr_\infty}\, \crm_1(\Lcal)=\frac{1}{k\, \tilde{k}}\,
\int_{\PP^1}\, \frac{1}{k\, \tilde{k}}\, \crm_1\big(\Ocal_{\PP^1}(a\, k)
\big)=\frac{a}{k\, \tilde{k}^2}\ .
\end{equation}
\endproof

\bigskip\section{Moduli spaces of framed sheaves\label{sec:framedsheaves}}

\subsection{Preliminaries on framed sheaves}

We begin by introducing some notation. Given a vector
$\vec{u}= (u_1,\ldots,u_{k-1}) \in\Z^{k-1}$, we denote by $\Rcal^{\vec{u}}$ the line bundle
$\bigotimes_{i=1}^{k-1}\, \Rcal_i^{\otimes u_i}$, and by $\Rcal_0$ the trivial line bundle $\Ocal_{\Xscr_k}$.

Let us fix $s\in\Z$. For $i=0,1, \ldots, k-1$ define the line bundles
\begin{equation*}
\Ocal_{\Dscr_\infty}(s,i):=\left\{
\begin{array}{ll}
\Lcal_1^{\otimes s}\otimes \Lcal_2^{\otimes i} & \mbox{for even $k$}\ ,\\[8pt]
\Lcal_1^{\otimes s}\otimes \Lcal_2^{\otimes i\, \frac{k+1}{2}} & \mbox{for odd $k$}\ .
\end{array}
\right.
\end{equation*}
In addition, let us fix $\vec{w}:=(w_0, w_1\ldots, w_{k-1})\in \N^{k}$ and $r:=\sum_{i=0}^{k-1}\, w_i$. Define 
\begin{equation}
\Fcal_{\infty}^{s,\vec{w}}:=\bigoplus_{i=0}^{k-1} \, \Ocal_{\Dscr_\infty}(s,i)^{\oplus w_i}\ .
\end{equation}
We shall call $\Fcal_\infty^{s,\vec{w}}$ a \emph{framing sheaf}. It is
a locally free sheaf on $\Dscr_\infty$ of degree $\frac{s\, r}{k\,
  \tilde{k}^2}$ (cf.\ Lemma \ref{lem:degree}).
  
\begin{remark}
The line bundles $\Ocal_{\Dscr_\infty}(0,j)$ are endowed with a unitary flat connection associated with the $j$-th 
irreducible unitary representation $\rho_j$ of $\Z_k$ for $j=0, \ldots, k-1$ (cf.\ \cite[Remark 6.5]{art:eyssidieuxsala2013}). So the framing 
sheaf $\Fcal_\infty^{0,\vec{w}}$ has a unitary flat connection associated with the unitary representation $\bigoplus_{j=0}^
{k-1}\, \rho_j^{\oplus w_j}$.
\end{remark}

\begin{definition}
A \emph{$(\Dscr_\infty,\Fcal_\infty^{s,\vec{w}}\, )$-framed sheaf} on
$\Xscr_k$ is a pair $(\Ecal, \phi_{\Ecal})$, where $\Ecal$ is a
torsion-free sheaf on $\Xscr_k$, locally free in a neighbourhood of
$\Dscr_\infty$, and $\phi_{\Ecal}\colon
\Ecal_{\vert\Dscr_\infty}\xrightarrow{\sim}
\Fcal^{s,\vec{w}}_\infty$ is an isomorphism. We call $\phi_{\Ecal}$ a
\emph{framing} of $\Ecal.$ A \emph{morphism} between
$(\Dscr_\infty,\Fcal_\infty^{s,\vec{w}}\, )$-framed sheaves $(\Ecal,\phi_\Ecal)$ and $(\Gcal,\phi_\Gcal)$ is a morphism  
$f\colon \Ecal\to\Gcal$ such that $\phi_\Gcal\circ f_{\vert\Dscr_\infty} = \phi_\Ecal$.
\end{definition}
\begin{remark}
Let $\Xscr$ be a toric Deligne-Mumford stack with coarse moduli space
a projective toric variety $\pi\colon \Xscr \to X$. The rank of a
coherent sheaf $\Ecal$ on $\Xscr$ is the degree zero part of the Chern
character $\ch(\Ecal)$ of its class $[\Ecal]$ in the K-theory group $K(\Xscr)$, which is
isomorphic to the Grothendieck group generated by the locally free
sheaves on $\Xscr$ because $\Xscr$ is projective and smooth. Roughly
speaking\footnote{A more precise discussion involves $K(\IXscr)$ and the
  T\"oen-Riemann-Roch theorem
  \cite[Remark~2.28]{art:bruzzosala2013}.}, the pushforward $\pi_\ast$
preserves the rank of the invariant part of a sheaf $\Ecal$ with
respect to the action of the generic stabilizer of $\Xscr$ on it. If
$\Xscr=\Xscr_k$, then since $\Xscr_k$ is an orbifold, its generic
stabilizer is trivial and the rank of a coherent sheaf $\Ecal$
coincides with the rank of its pushforward ${\pi_k}_\ast(\Ecal)$. By
applying the same argument to $\Dscr_\infty$, one finds that the rank
of a coherent sheaf $\Fcal$ on $\Dscr_\infty$ does not coincide in general with
the rank of ${r_k}_\ast(\Fcal)$, and indeed the latter is only the
rank of ${r_k}_\ast(\Fcal^0)$ where $\Fcal^0$ is the $\mu_k$-invariant
part of $\Fcal$ with respect to the action of $\mu_k$
\cite[Section~3.2]{art:nironi2008}. Since the K-theory groups
$K(\Xscr_k)$ and $K(\Dscr_\infty)$ are both generated by line bundles
(cf.\ \cite[Sections~4 and 6]{art:borisovhorja2006}),
the rank is preserved under the restriction to $\Dscr_\infty$. In
particular, if $(\Ecal, \phi_{\Ecal})$ is a
$(\Dscr_\infty,\Fcal_\infty^{s,\vec{w}}\, )$-framed sheaf on
$\Xscr_k$, we get $\rk(\Ecal)=\rk(\Fcal_\infty^{s,\vec{w}}\, )=r$.

Moreover, the Picard group of $\Xscr_k$ is isomorphic to its second
singular cohomology group with integral coefficients via the
\emph{first Chern class map} \cite[Section
3.1.2]{art:iritani2009}.\footnote{This result is a generalization of
  the analogous result for toric varieties \cite[Theorem
  12.3.2]{book:coxlittleschenck2011}.} Thus fixing the determinant
line bundle of a coherent sheaf $\Ecal$ on $\Xscr_k$ is equivalent to fixing its first Chern class. 
\end{remark}
\begin{lemma}\label{lem:firstchernclass}
Let $(\Ecal, \phi_{\Ecal})$ be a
$(\Dscr_\infty,\Fcal_\infty^{s,\vec{w}}\, )$-framed sheaf on
$\Xscr_k$. Then the determinant $\det(\Ecal)$ of $\Ecal$ is of the
form $\Rcal^{\vec{u}}\otimes \Ocal_{\Xscr_k}(s\, r\, \Dscr_\infty)$, where the vector $\vec{u}\in \Z^{k-1}$ satisfies the condition
\begin{equation}\label{eq:condition-determinant}
\sum_{j=1}^{k-1}\, j\, u_j=   \sum_{i=0}^{k-1}\, i\, w_i \ \bmod{k}  \ .
\end{equation}
\end{lemma}
\proof
Let $\det(\Ecal)=\bigotimes_{j=1}^{k-1}\, {\Rcal_j^{\otimes
    u_j}}\otimes \Ocal_{\Xscr_k}(u_\infty \, \Dscr_\infty)$ be the
determinant line bundle of $\Ecal$ for some integer $u_\infty$ and
some vector $\vec{u}:=(u_1, \ldots, u_{k-1})\in\Z^{k-1}$. Since
$\det(\Fcal_\infty^{s,\vec{w}}\, )\simeq\det(\Ecal_{\vert \Dscr_\infty})$, we get
\begin{equation}
\bigotimes_{i=0}^{k-1}\, \Ocal_{\Dscr_\infty}(s,i)^{\otimes w_i}\simeq
\bigotimes_{j=1}^{k-1}\, {\Rcal_j^{\otimes u_j}}_{\vert
  \Dscr_\infty}\otimes \Ocal_{\Xscr_k}(u_\infty \, \Dscr_\infty)_{\vert \Dscr_\infty}\ .
\end{equation}
By Corollary \ref{cor:restrictiontautological} we have ${\Rcal_i}_{\vert\Dscr_\infty}\simeq\Ocal_{\Dscr_\infty}(0,i)$ for $i=1,\ldots,k-1$ and $\Ocal_{\Xscr_k}(\Dscr_\infty)_{\vert\Dscr_\infty}\simeq\Ocal_{\Dscr_\infty}(1,0)$, hence we get the assertion.
\endproof
\begin{remark}\label{rem:v-y}
Set $\vec{v}:=C^{-1}\vec{u}$. Then Equation \eqref{eq:condition-determinant} implies the relation
\begin{equation}\label{eq:v-condition}
k\, v_j=-j\, \sum_{i=0}^{k-1}\, i\, w_i \ \bmod{k}
\end{equation}
 for $j=1, \ldots, k-1$. Note also that for $i\in\{1, \ldots, k-1\}$ a
 component $v_i$ is integral if and only if every component is. Let
 $c\in\{0,1,\ldots, k-1\}$ be the equivalence class modulo $k$ of
 $\sum_{i=0}^{k-1}\, i\, w_i$ and define $\vec{y}:=C^{-1}\vec{e}_{c}-\vec{v}$ if $c>0$, otherwise $\vec{y}:=-\vec{v}$. Then $\vec{y}\in \Z^{k-1}$. We identify $\vec{y}$ with an element of the root lattice $\Qfrak$ (cf.\ Remark \ref{rem:rootlattice-cohomology}).
\end{remark}
Now we are ready to apply the machinery developed in
\cite{art:bruzzosala2013} to construct moduli spaces of framed
sheaves on root toric orbifolds. In particular, note that $\tilde{k}\,
D_\infty$ is a big and nef Cartier divisor, which contains the
singular points of $\bar{X}_k$; the line bundle
$\Ocal_{\Xscr_k^{\mathrm{can}}}(\tilde{\Dscr}_\infty)$ is a
$\pi_k^{\mathrm{can}}$-ample sheaf (cf.\ Remark
\ref{rem:ampleness}). Moreover, since $\Fcal_\infty^{s,\vec{w}}$ is
given as direct sums of line bundles of the same degree, it is a
\emph{good framing sheaf}
\cite[Definition~5.7]{art:bruzzosala2013}. Then by
\cite[Theorem~6.2]{art:bruzzosala2013} there exists a \emph{fine}
moduli space
$\Mcal_{r,\vec{u},\Delta}(\Xscr_k,\Dscr_\infty,\Fcal_\infty^{s,\vec{w}}\,
)$
parametrizing isomorphism classes of
$(\Dscr_\infty,\Fcal_\infty^{s,\vec{w}}\, )$-framed sheaves
$(\Ecal,\phi_\Ecal)$ on $\Xscr_k$, where $\Ecal$ is a torsion-free
sheaf of rank $r$, determinant $\Rcal^{\vec{u}}\otimes
\Ocal_{\Xscr_k}(s\, r \, \Dscr_\infty)$, and discriminant 
\begin{equation}
\Delta:=\Delta(\Ecal)=\int_{\Xscr_k}\,
\Big(\crm_2(\Ecal)-\frac{r-1}{2r}\, \crm_1^2(\Ecal)\Big)\ ,
\end{equation}
and $\vec{u}\in\mathbb{Z}^{k-1}$ satisfies Equation
\eqref{eq:condition-determinant}.

\begin{rem}\label{rem:universalsheaf}
By ``fine'' one means that there exists a \emph{universal framed
  sheaf} $(\boldsymbol{\Ecal}, \boldsymbol{\phi}_{\boldsymbol{\Ecal}})$, where
$\boldsymbol{\Ecal}$ is a coherent sheaf on $\Xscr_k\times
\Mcal_{r,\vec{u},\Delta}(\Xscr_k,\Dscr_\infty,\Fcal_\infty^{s,\vec{w}}\,
)$ which is flat over
$\Mcal_{r,\vec{u},\Delta}(\Xscr_k,\Dscr_\infty,\Fcal_\infty^{s,\vec{w}}\,
)$,
and $\boldsymbol{\phi}_{\boldsymbol{\Ecal}}\colon \boldsymbol{\Ecal}\to
p_{\Xscr_k}^\ast(\Fcal_\infty^{s,\vec{w}}\, )$ is a
morphism such that $(\boldsymbol{\phi}_{\boldsymbol{\Ecal}})_{\vert
  \Dscr_\infty\times
  \Mcal_{r,\vec{u},\Delta}(\Xscr_k,\Dscr_\infty,\Fcal_\infty^{s,\vec{w}}\,
  )}$ is an isomorphism; the fibre over $[(\Ecal,\phi_\Ecal)]\in\Mcal_{r,\vec{u},\Delta}(\Xscr_k,\Dscr_\infty,\Fcal_\infty^{s,\vec{w}}\,
)$ is itself the $(\Dscr_\infty,\Fcal_\infty^{s,\vec{w}}\, )$-framed sheaf $(\Ecal,\phi_\Ecal)$ on $\Xscr_k$. In the following we shall call $\boldsymbol{\Ecal}$ a \emph{universal sheaf}.\hfill$\triangle$
\end{rem}

\subsection{Smoothness of the moduli space}

\begin{proposition}
For any $s\in \mathbb{Z}$ and $i=0,1, \ldots, k-1$, the pushforward
${r_k}_\ast(\Ocal_{\Dscr_\infty}(s,i))$ of $\Ocal_{\Dscr_\infty}(s,i)$
is given by:
\begin{itemize}
\item ${r_k}_\ast(\Ocal_{\Dscr_\infty}(s,i)) = 0$ if $s$ and $i$ do not satisfy the following conditions:
\begin{align}
s+i\, \tilde{k}&=   0 \ \bmod{k} \qquad\mbox{for even $k$}\
, \label{eq:decompositioneven} \\[4pt]
s&=   0 \ \bmod{k}  \qquad\mbox{for odd $k$}\ . \label{eq:decompositionodd}
\end{align}
\item Otherwise
\begin{equation*}
{r_k}_\ast\big(\Ocal_{\Dscr_\infty}(s,i) \big)\simeq \left\{
\begin{array}{ll}
\Ocal_{\mathbb{P}^1}\Big(\left\lfloor \frac{s+i\, \tilde{k}}{k\,
    \tilde{k}}\right\rfloor+\left\lfloor \frac{s-i\, \tilde{k}}{k\, \tilde{k}}\right\rfloor\Big) & \mbox{for even $k$}\ ,\\[8pt]
 \Ocal_{\mathbb{P}^1}\Big(\left\lfloor \frac{\frac{1-k}{2}\,
     \frac{s-i\, k}{k}}{k}\right\rfloor+\left\lfloor
   \frac{\frac{1+k}{2}\, \frac{s-i\, k}{k}+i}{k}\right\rfloor\Big) & \mbox{for odd $k$}\ .
\end{array}
\right.
\end{equation*}
\end{itemize}
\end{proposition}
\proof
Let $s\in \mathbb{Z}$ and $i=0,1, \ldots, k-1$. First recall that the
banding group of the gerbe $\tilde{\phi}_k\colon \Dscr_\infty\to
\tilde{\Dscr}_\infty$ is $\mu_k$, which fits into the exact sequence
given by:
\begin{itemize}
\item for even $k$:
\begin{equation}
1 \ \longrightarrow \ \mu_k \ \xrightarrow{i_{\rm ev}} \
\mathbb{C}^\ast\times \mu_k \ \xrightarrow{q_{\rm ev}} \
\mathbb{C}^\ast\times \mu_{\tilde{k}} \ \longrightarrow \ 1\ ,
\end{equation}
where $i_{\rm ev}\colon \eta\mapsto (\eta,\eta^{\tilde{k}})$ and
$q_{\rm ev}\colon (t,\omega)\mapsto (t^{\tilde{k}}\, \omega^{-1}, \omega^2)$;
\item for odd $k$:
\begin{equation}
1 \ \longrightarrow \ \mu_k \ \xrightarrow{i_{\mathrm{odd}}} \
\mathbb{C}^\ast\times \mu_k \ \xrightarrow{q_{\mathrm{odd}}} \
\mathbb{C}^\ast\times \mu_{k} \ \longrightarrow 1\ ,
\end{equation}
where $i_{\mathrm{odd}}\colon \eta\mapsto (\eta,1)$ and
$q_{\mathrm{odd}}\colon (t,\omega)\mapsto (t^{k}\, \omega^{\frac{k-1}{2}}, \omega)$.
\end{itemize}
Any coherent sheaf on $\Dscr_\infty$ decomposes as a direct sum of
eigensheaves with respect to the characters of $\mu_k$. The
pushforward of $\tilde{\phi}_k$ preserves only the $\mu_k$-invariant
part of a coherent sheaf on $\Dscr_\infty$. In view of the previous exact sequences
we find that the pushforward
${\tilde{\phi}}_k{}_\ast(\Ocal_{\Dscr_\infty}(s,i))$ is nonzero if and
only if Equations
\eqref{eq:decompositioneven}--\eqref{eq:decompositionodd} are
satisfied. 
For even $k$, and for $s$ and $i$ satisfying Equation \eqref{eq:decompositioneven}, we get
\begin{equation}
\Ocal_{\Dscr_\infty}(s,i)\simeq
\tilde{\phi}_k^\ast\Big(\Ocal_{\tilde{\Dscr}_\infty}\big(\mbox{$\frac{s+i\,
    \tilde{k}}{k}$} \, \tilde{p}_0\big)\otimes
\Ocal_{\tilde{\Dscr}_\infty}\big(\mbox{$\frac{s+i\, \tilde{k}-i\,
    k}{k}$} \, \tilde{p}_\infty\big)\Big)\ .
\end{equation}
By the projection formula, which holds for the rigidification morphism
$\tilde{\phi}_k$ \cite{question:vistoli2013}, we thus have
\begin{equation}
{\tilde{\phi}_k}{}_\ast\big(\Ocal_{\Dscr_\infty}(s,i) \big)\simeq \Ocal_{\tilde{\Dscr}_\infty}\big(\mbox{$\frac{s+i\,
    \tilde{k}}{k}$} \, \tilde{p}_0\big)\otimes
\Ocal_{\tilde{\Dscr}_\infty}\big(\mbox{$\frac{s+i\, \tilde{k}-i\,
    k}{k}$} \, \tilde{p}_\infty\big) \ .
\end{equation}
Recall that $\tilde{\Dscr}_\infty$ is obtained from $D_\infty$ by performing a $(\tilde{k},\tilde{k})$-root construction at the points $0,\infty\in D_\infty\simeq \PP^1$.
The asserted result for ${r_k}_\ast(\Ocal_{\Dscr_\infty}(s,i))$ now
follows by Lemma \ref{lem:cadman3.1.1}. In the same
way, for odd $k$ and for $s$ satisfying Equation
\eqref{eq:decompositionodd}, we get the asserted result for
${r_k}_\ast(\Ocal_{\Dscr_\infty}(s,i))$.
\endproof
\begin{corollary}\label{cor:vanishing}
For $i\in\{0,1, \ldots, k-1\}$ and any negative integer $s$ we have 
\begin{equation}
H^0\big(\Dscr_\infty,\Ocal_{\Dscr_\infty}(s,i)\big)=H^0\big(\PP^1,{r_k}_\ast(\Ocal_{\Dscr_\infty}(s,i))
\big)=0\ .
\end{equation}
\end{corollary}
From now we shall assume that $s$ is a negative integer. Thanks to this corollary, we can argue exactly as in the proof of the analogous result for framed sheaves on smooth toric surfaces \cite[Proposition~2.1]{art:gasparimliu2010} and obtain easily the following result. Note that the proof involves Serre duality for stacks as stated in \cite[Theorem~B.7]{art:bruzzosala2013}.
\begin{proposition}\label{prop:vanishing}
The group $\Ext^i(\Ecal,\Ecal'\otimes\Ocal_{\Xscr_k}(-\Dscr_\infty))$
vanishes for $i=0,2$ and for any pairs consisting of a
$(\Dscr_\infty,\Fcal_\infty^{s,\vec{w}}\, )$-framed sheaf
$(\Ecal,\phi_{\Ecal})$ and a
$(\Dscr_\infty,\Fcal_\infty^{s,\vec{w}'}\, )$-framed sheaf $(\Ecal',\phi_{\Ecal'})$ on $\Xscr_k$. If in addition $s=0$ then $H^i(\Xscr_k, \Ecal\otimes \Ocal_{\Xscr_k}(-\Dscr_\infty))=0$ for $i=0,2$.
\end{proposition}

\begin{theorem}\label{thm:moduli}
The moduli space
$\Mcal_{r,\vec{u},\Delta}(\Xscr_k,\Dscr_\infty,\Fcal_\infty^{s,\vec{w}}\,
)$ is a smooth quasi-projective variety of dimension
\begin{equation}
\dim_{\mathbb{C}}\big(\Mcal_{r,\vec{u},\Delta}(\Xscr_k,\Dscr_\infty,\Fcal_\infty^{s,\vec{w}}\,
) \big)=2\, r\, \Delta- \frac12\, \sum_{j=1}^{k-1}\, \big(C^{-1} \big)^{jj}\, \vec{w}\cdot\vec{w}(j)\ ,
\end{equation}
where for $j=1, \ldots, k-1$ the vector $\vec{w}(j)$ is $(w_j, \ldots, w_{k-1}, w_0, w_1, \ldots, w_{j-1})$ and $C$ is the Cartan matrix of the Dynkin diagram of type $A_{k-1}$. Moreover, the Zariski tangent
space of
$\Mcal_{r,\vec{u},\Delta}(\Xscr_k,\Dscr_\infty,\Fcal_\infty^{s,\vec{w}}\,
)$ at a point $[(\Ecal,\phi_{\Ecal})]$ is $\mathrm{Ext}^1(\Ecal,\Ecal\otimes \Ocal_{\Xscr_k}(-\Dscr_\infty))$.
\end{theorem}
\proof
The moduli space
$\Mcal_{r,\vec{u},\Delta}(\Xscr_k,\Dscr_\infty,\Fcal_\infty^{s,\vec{w}}\,
)$ is a separated quasi-projective scheme of finite type over $\C$ by
\cite[Theorem~6.2]{art:bruzzosala2013}. By
\cite[Theorem~4.17-(ii)]{art:bruzzosala2013}, the group
$\Ext^2(\Ecal,\Ecal\otimes\Ocal_{\Xscr_k}(-\Dscr_\infty))$ contains the
obstruction to the smoothness of
$\Mcal_{r,\vec{u},\Delta}(\Xscr_k,\Dscr_\infty,\Fcal_\infty^{s,\vec{w}}\,
)$ at $[(\Ecal,\phi_{\Ecal})]$. By Proposition \ref{prop:vanishing},
the group $\Ext^2(\Ecal,\Ecal\otimes\Ocal_{\Xscr_k}(-\Dscr_\infty))$
vanishes for all points $[(\Ecal,\phi_{\Ecal})]$ of
$\Mcal_{r,\vec{u},\Delta}(\Xscr_k,\Dscr_\infty,\Fcal_\infty^{s,\vec{w}}\,
)$, and so
$\Mcal_{r,\vec{u},\Delta}(\Xscr_k,\Dscr_\infty,\Fcal_\infty^{s,\vec{w}}\,
)$ is everywhere smooth. Thus
$\Mcal_{r,\vec{u},\Delta}(\Xscr_k,\Dscr_\infty,\Fcal_\infty^{s,\vec{w}}\,
)$ is a smooth quasi-projective variety over $\mathbb{C}$.

By \cite[Theorem~4.17-(i)]{art:bruzzosala2013}, the Zariski tangent
space of
$\Mcal_{r,\vec{u},\Delta}(\Xscr_k,\Dscr_\infty,\Fcal_\infty^{s,\vec{w}}\,
)$ at a point $[(\Ecal,\phi_{\Ecal})]$ is $\mathrm{Ext}^1(\Ecal,\Ecal\otimes \Ocal_{\Xscr_k}(-\Dscr_\infty))$. Hence
\begin{equation}
\dim_{\mathbb{C}}\big(\Mcal_{r,\vec{u},\Delta}(\Xscr_k,\Dscr_\infty,\Fcal_\infty^{s,\vec{w}}\,
)\big)=\dim_{\mathbb{C}}\big(\mathrm{Ext}^1(\Ecal,\Ecal\otimes
\Ocal_{\Xscr_k}(-\Dscr_\infty)) \big)\ ,
\end{equation}
and the latter dimension is computed in Appendix
\ref{sec:rankdimension} (see Corollary \ref{cor:dimension}).
\endproof 

By \cite[Theorem 6.9]{art:eyssidieuxsala2013} the moduli space $\Mcal_{r, \vec{u}, \Delta}(\Xscr_k, \Dscr_\infty, \Fcal_\infty^{0, \vec{w}}\,)$ contains as an open subset the moduli space of $U(r)$-instantons on $X_k$ with first Chern class $\sum_{i=1}^{k-1}\,u_i\,\crm_1({\Rcal_i}_{\vert X_k})$, discriminant $\Delta$ and holonomy at infinity 
associated with the unitary representation $\bigoplus_{j=0}^{k-1}\, \rho_j^{\oplus w_j}$. Because of this fact, we state the 
following conjecture:
\begin{conjecture}
The moduli spaces $\Mcal_{r, \vec{u}, \Delta}(\Xscr_k, \Dscr_\infty, \Fcal_\infty^{0, \vec{w}}\,)$ are isomorphic to 
Nakajima quiver varieties $\Mcal_{\xi_k}(\vec{s}, \vec{w})$ for a suitable choice of $\vec{s}\in\N^k$, where $\xi_k$ is 
the stability parameter of $X_k$.
\end{conjecture}
As we shall see in the following section, in the rank one case our moduli spaces are isomorphic to Hilbert schemes of 
points on $X_k$ and the latter moduli spaces are Nakajima quiver varieties by \cite{art:kuznetsov2007}.

\subsection{Rank one}\label{sec:rankonecase}

Let $\hilb{n}{X_k}$ be the Hilbert scheme of $n$ points of $X_k$,
i.e., the fine moduli space parameterizing zero-dimensional subschemes of
$X_k$ of length $n$. It is a smooth quasi-projective variety of dimension $2n$. If $Z$ is a point of $\hilb{n}{X_k}$,  the
pushforward $\imath_\ast I_Z$ of the ideal sheaf $I_Z$ with respect to
the inclusion morphism $\imath\colon X_k \to \Xscr_k$ is a rank one
torsion-free sheaf on $\Xscr_k$ with $\det(\imath_\ast I_Z)\simeq
\Ocal_{\Xscr_k}$ and $\int_{\Xscr_k}\, \crm_2(\imath_\ast(I_Z))=n$. The morphism $\imath$ induces an isomorphism $\imath\colon X_k\xrightarrow{\sim}\Xscr_k\setminus\Dscr_\infty$, hence $Z\subset X_k$ is disjoint from $\Dscr_\infty$ so that $\imath_\ast I_Z$ is locally free in a neighbourhood of $\Dscr_\infty$. 

Let $\vec{u}\in \Z^{k-1}$ and fix $i\in\{0,1, \ldots, k-1\}$ such that
\begin{equation}
i = \sum_{j=1}^{k-1} \, j \, u_j \ \bmod{k}\ .
\end{equation}
Let $s\in \Z$. Then  $\Ecal:=\imath_\ast (I_Z)\otimes
\Rcal^{\vec{u}}\otimes \Ocal_{\Xscr_k}(s \, \Dscr_\infty)$ is a rank
one torsion-free sheaf on $\Xscr_k$, locally free in a neighbourhood
of $\Dscr_\infty$, with a framing $\phi_{\Ecal}\colon
\Ecal_{\vert\Dscr_\infty}\xrightarrow{\sim} \Ocal_{\Dscr_\infty}(s,i)$
induced canonically by the isomorphism ${\Rcal^{\vec{u}}\otimes
  \Ocal_{\Xscr_k}(s \, \Dscr_\infty)}_{\vert\Dscr_\infty}\simeq
\Ocal_{\Dscr_\infty}(s,i)$ (cf.\ Corollary
\ref{cor:restrictiontautological}). So we get a
$(\Dscr_\infty,\Ocal_{\Dscr_\infty}(s,i))$-framed sheaf $(\Ecal,
\phi_{\Ecal})$ on $\Xscr_k$ (as the line bundle
$\Ocal_{\Dscr_\infty}(s,i)$ coincides with $\Fcal_\infty^{s,\vec{w}}$
for the vector $\vec{w}$ such that $w_i=1$ and $w_j=0$ for $j\neq
i$). Moreover, $\det(\Ecal)\simeq \Rcal^{\vec{u}}\otimes
\Ocal_{\Xscr_k}(s \, \Dscr_\infty)$ and 
\begin{equation}
\int_{\Xscr_k}\,  \ch_2(\Ecal)=\frac{1}{2}\, \int_{\Xscr_k}\, \crm_1
\big(\Rcal^{\vec{u}}\otimes \Ocal_{\Xscr_k}(s \, \Dscr_\infty) \big)^2-n \ .
\end{equation}
This singles out a point $[(\Ecal,\phi_{\Ecal})]$ in $\Mcal_{1,\vec{u},n}(\Xscr_k,\Dscr_\infty,\Ocal_{\Dscr_\infty}(s,i))$, so that an
injective morphism of fine moduli spaces
\begin{equation}\label{eq:inclu}
\imath_{(1,\vec{u},n)}\, \colon \, \hilb{n}{X_k} \ \longrightarrow \
\Mcal_{1,\vec{u},n}\big(\Xscr_k,\Dscr_\infty,\Ocal_{\Dscr_\infty}(s,i)
\big)
\end{equation}
is defined. This argument extends straightforwardly to families of zero-dimensional subschemes of $X_k$ of length $n$. 
\begin{proposition}\label{prop:rankone}
The inclusion morphism \eqref{eq:inclu} is an isomorphism of fine moduli spaces.
\end{proposition}
\proof
We can define an inverse morphism $\jmath_{(1,\vec{u},n)}\colon \Mcal_{1,\vec{u},n}(\Xscr_k,\Dscr_\infty,\Ocal_{\Dscr_\infty}(s,i)) \to \hilb{n}{X_k}$ in the following way. Let $[(\Ecal,\phi_{\Ecal})]$ be a point in $\Mcal_{1,\vec{u},n}(\Xscr_k,\Dscr_\infty,\Ocal_{\Dscr_\infty}(s,i))$. The torsion-free sheaf $\Ecal$ fits into the exact sequence
\begin{equation}
0 \ \longrightarrow \ \Ecal \ \longrightarrow \ \Ecal^{\vee\vee} \
\longrightarrow \ \mathcal{Q} \ \longrightarrow \ 0\ ,
\end{equation}
where $\Ecal^{\vee\vee}$ is the line bundle $\Rcal^{\vec{u}}\otimes
\Ocal_{\Xscr_k}(s \, \Dscr_\infty)$ and $\mathcal{Q}$ is a zero-dimensional sheaf whose support has length $n.$ Since $\Ecal$ is locally free in a neighbourhood of $\Dscr_\infty$, the support of $\mathcal{Q}$ is disjoint from $\Dscr_\infty.$ Hence the quotient
\begin{equation}
\Ocal_{\Xscr_k}\simeq \Ecal^{\vee\vee}\otimes \Rcal^{-\vec
  u}\otimes\Ocal_{\Xscr_k}(-s\, \Dscr_\infty) \ \longrightarrow \
\mathcal{Q}\otimes \Rcal^{-\vec u}\otimes\Ocal_{\Xscr_k}(-s \,
\Dscr_\infty)\ \longrightarrow \ 0
\end{equation}
defines a zero-dimensional subscheme $Z\subset \Xscr_k$ of length $n$
which is disjoint from $\Dscr_\infty$, and the quotient
$\Ocal_{X_k}\to \imath^\ast(\Ocal_Z)\to 0$ defines a point
$Z\in\hilb{n}{X_k}$ with $\Ecal\simeq \imath_\ast (I_Z)\otimes
\Rcal^{\vec{u}}\otimes \Ocal_{\Xscr_k}(s\, \Dscr_\infty)$. It is easy
to see that this argument can be generalized to families of framed
sheaves. Moreover, one has $\imath_{(1,\vec{u},n)}\circ\jmath_{(1,\vec{u},n)} = \mathrm{id}$ and $\jmath_{(1,\vec{u},n)}\circ\imath_{(1,\vec{u},n)}=\mathrm{id}$.
\endproof
\begin{remark}\label{rem:rankone}
A consequence of this proposition is that after fixing $i\in \{0, 1,
\ldots,$ $ k-1\}$, a vector $\vec{u}\in \Z^{k-1}$ such that
$\sum_{j=1}^{k-1}\, j \, u_j= i \bmod{k} $, and an integer $s\in \Z$,
for any $(\Dscr_\infty,\Ocal_{\Dscr_\infty}(s,i))$-framed sheaf
$(\Ecal, \phi_{\Ecal})$ of rank one on $\Xscr_k$ the torsion-free
sheaf $\Ecal$ is isomorphic to $\imath_\ast(I)\otimes
\Rcal^{\vec{u}}\otimes \Ocal_{\Xscr_k}(s \, \Dscr_\infty)$, where $I$
is the ideal sheaf of some zero-dimensional subscheme of $X_k$, and
$\phi_{\Ecal}$ is canonically induced by the isomorphism
${\Rcal^{\vec{u}}\otimes \Ocal_{\Xscr_k}(s \, \Dscr_\infty)}_{\vert\Dscr_\infty}\simeq \Ocal_{\Dscr_\infty}(s,i)$.
\end{remark}
Since $\imath_{(1,\vec{u},n)}$ is an isomorphism between fine moduli
spaces, we also obtain an isomorphism between the corresponding
universal objects. Let us denote by $\boldsymbol{Z}\subset
\hilb{n}{X_k}\times X_k$ the universal subscheme of $\hilb{n}{X_k}$,
whose fibre over $Z\in\hilb{n}{X_k}$ is the zero-dimensional subscheme
$Z$ itself. Consider the commutative diagram
\begin{equation}
  \begin{tikzpicture}[xscale=6.5,yscale=-1.2]
    \node (A0_0) at (0, 0) {$\hilb{n}{X_k}\times X_k$};
    \node (A0_1) at (1, 0)
    {$\Mcal_{1,\vec{u},n}\big(\Xscr_k,\Dscr_\infty,\Ocal_{\Dscr_\infty}(s,i)
      \big)\times \Xscr_k$};
    \node (A) at (1.6, 0.5) {$.$};  
    \node (A1_0) at (0, 1) {$\hilb{n}{X_k}$};
    \node (A1_1) at (1, 1)
    {$\Mcal_{1,\vec{u},n}\big(\Xscr_k,\Dscr_\infty,\Ocal_{\Dscr_\infty}(s,i)
      \big)$};
    \path (A0_0) edge [->]node [auto] {$\scriptstyle{(\imath_{(1,\vec{u},n)},\imath)}$} (A0_1);
    \path (A0_1) edge [->]node [auto] {$\scriptstyle{}$} (A1_1);
    \path (A1_0) edge [->]node [auto] {$\scriptstyle{\imath_{(1,\vec{u},n)}}$} (A1_1);
    \path (A0_0) edge [->]node [auto] {$\scriptstyle{}$} (A1_0);
  \end{tikzpicture}
\end{equation}
Then $(\imath_{(1,\vec{u},n)},\imath )^\ast\boldsymbol{\Ecal}$ is the ideal sheaf of
$\boldsymbol{Z}$ and $(\imath_{(1,\vec{u},n)},\imath )^\ast\boldsymbol{\phi}_{\boldsymbol{\Ecal}}=0$, where $(\boldsymbol{\Ecal}, \boldsymbol{\phi}_{\boldsymbol{\Ecal}})$ is the universal framed sheaf on $\Mcal_{1,\vec{u},n}(\Xscr_k,\Dscr_\infty,\Ocal_{\Dscr_\infty}(s,i))\times \Xscr_k$ introduced in Remark \ref{rem:universalsheaf}.

\subsection{Natural bundle}

In this subsection we set $s=0$. The \emph{natural bundle} $\Vbf$ on
$\Mcal_{r,\vec{u},\Delta}(\Xscr_k,\Dscr_\infty,\Fcal_\infty^{0,\vec
  w}\, )$
is defined in terms of  the universal sheaf 
$\boldsymbol{\Ecal}$ as
 \begin{equation}
\Vbf:=R^1 p_\ast\Big(\boldsymbol{\Ecal}\otimes
p_{\Xscr_k}^\ast\big(\Ocal_{\Xscr_k}(-\Dscr_\infty) \big)\Big)\ ,
\end{equation}
where $p\colon
{\Mcal_{r,\vec{u},\Delta}(\Xscr_k,\Dscr_\infty,\Fcal_\infty^{0,\vec
    w}\, )} \times \mathscr X_k \to
{\Mcal_{r,\vec{u},\Delta}(\Xscr_k,\Dscr_\infty,\Fcal_\infty^{0,\vec
    w}\, )}$ is the projection.
\begin{proposition}\label{prop:naturalbundle}
$\Vbf$ is a locally free sheaf on
$\Mcal_{r,\vec{u},\Delta}(\Xscr_k,\Dscr_\infty,\Fcal_\infty^{0,\vec
  w}\, )$ of rank
\begin{equation}
\rk(\Vbf) = \Delta +\frac{1}{2r}\, \vec v\cdot C\vec v - \frac{1}{2}
\, \sum_{j=1}^{k-1}\, \big(C^{-1}\big)^{jj}\, w_j
\end{equation}
where $\vec{v}:=C^{-1} \vec{u}$.
\end{proposition}
\proof
First note that 
\begin{equation}
\mathrm{id}\times \pi_k\, \colon \,
\Mcal_{r,\vec{u},\Delta}(\Xscr_k,\Dscr_\infty,\Fcal_\infty^{0,\vec
  w}\, )\times \Xscr_k \ \longrightarrow \
\Mcal_{r,\vec{u},\Delta}(\Xscr_k,\Dscr_\infty,\Fcal_\infty^{0,\vec
  w}\, )\times \bar{X}_k
\end{equation}
is the coarse moduli space of
$\Mcal_{r,\vec{u},\Delta}(\Xscr_k,\Dscr_\infty,\Fcal_\infty^{0,\vec
  w}\, )\times \Xscr_k$. Moreover, the projection morphism 
\begin{equation}
\Mcal_{r,\vec{u},\Delta}(\Xscr_k,\Dscr_\infty,\Fcal_\infty^{0,\vec
  w}\, )\times \bar{X}_k \ \longrightarrow \
\Mcal_{r,\vec{u},\Delta}(\Xscr_k,\Dscr_\infty,\Fcal_\infty^{0,\vec w}
\, )
\end{equation}
is proper. Since  the stack
$\Mcal_{r,\vec{u},\Delta}(\Xscr_k,\Dscr_\infty,\Fcal_\infty^{0,\vec
  w}\, )\times \Xscr_k$ is tame over
$\Mcal_{r,\vec{u},\Delta}(\Xscr_k,\Dscr_\infty,\Fcal_\infty^{0,\vec
  w}\, )$, we can apply the cohomology and base change theorem \cite[Theorem~1.7]{art:nironi2008}. By using  the vanishing results of Proposition \ref{prop:vanishing} one can prove that 
$\Vbf$ is a locally free sheaf of rank $\dim_{\C} H^1(\Xscr_k,
\Ecal\otimes \Ocal_{\Xscr_k}(-\Dscr_\infty))$, where $\Ecal$ is the
underlying sheaf of a point
$[(\Ecal,\phi_{\Ecal})]\in\Mcal_{r,\vec{u},\Delta}(\Xscr_k,\Dscr_\infty,\Fcal_\infty^{0,\vec
  w}\, )$. The computation of this dimension is given in Appendix
\ref{sec:rankdimension} (see Theorem \ref{thm:rank}).
\endproof

\subsection{Carlsson-Okounkov bundle}\label{sec:carlssonokounkov}

In this subsection we introduce the Carlsson-Okounkov bundle for the
moduli spaces
$\Mcal_{r,\vec{u},\Delta}(\Xscr_k,\Dscr_\infty,\Fcal_\infty^{s,\vec
  w}\, )$, generalizing the definition given in
\cite{art:carlssonokounkov2012} for the rank one case. Let
\begin{equation}
\boldsymbol{\Ecal}_i:={p_{i 3}}^\ast(\boldsymbol{\Ecal}) \ \in \
K\big(\Mcal_{r,\vec{u},\Delta}(\Xscr_k,\Dscr_\infty,\Fcal_\infty^{s,\vec
  w}\, )
\times\Mcal_{r',\vec{u}\,',\Delta'}(\Xscr_k,\Dscr_\infty,\Fcal_\infty^{s,\vec
  w\,'}\, )\times \Xscr_k \big)
\end{equation}
for $i=1, 2$, where $p_{ij}$ is the projection of $\Mcal_{r,\vec{u},\Delta}(\Xscr_k,\Dscr_\infty,\Fcal_\infty^{s,\vec
  w}\, )
\times\Mcal_{r',\vec{u}\,',\Delta'}(\Xscr_k,\Dscr_\infty,\Fcal_\infty^{s,\vec
  w\,'}\, )\times \Xscr_k$ on the product of the $i$-th and $j$-th factors. Denote by $p_3$ the
projection of the same product onto $\Xscr_k$.
\begin{definition}
The \emph{Carlsson-Okounkov bundle} is the element
\begin{equation}
\Ebf:={p_{12}}_\ast\big(-\boldsymbol{\Ecal}_1^\vee\cdot \boldsymbol{\Ecal}_2\cdot p_3^\ast(\Ocal_{\Xscr_k}(-\Dscr_\infty))\big)
\end{equation}
in the K-theory
$K\big(\Mcal_{r,\vec{u},\Delta}(\Xscr_k,\Dscr_\infty,\Fcal_\infty^{s,\vec
  w}\, )
\times\Mcal_{r',\vec{u}\,',\Delta'}(\Xscr_k,\Dscr_\infty,\Fcal_\infty^{s,\vec
  w\,'}\, ) \big)$.
\end{definition}
The fibre of $\Ebf$ over a pair of points $\big([(\Ecal,\phi_\Ecal)]\,,\, [(\Ecal',\phi_{\Ecal'})]\big)$ is
\begin{equation}\label{eq:fibre}
\Ebf_{\left([(\Ecal,\phi_\Ecal)]\,,\,
    [(\Ecal',\phi_{\Ecal'})]\right)}=-\chi_{\Xscr_k}\big(\Ecal,\Ecal'\otimes
\Ocal_{\Xscr_k}(-\Dscr_\infty) \big)\ ,
\end{equation}
where for any pair of coherent sheaves $\Fcal$ and $\Gcal$ on
$\Xscr_k$ we set $\chi_{\Xscr_k}(\Fcal,\Gcal):=\sum_i\, (-1)^i\, \Ext^i(\Fcal,\mathcal{G})$. By Proposition \ref{prop:vanishing} and the fact that all the Ext groups vanish in degree greater than two,
the fibre \eqref{eq:fibre} reduces to
\begin{equation}
\Ebf_{\left([(\Ecal,\phi_\Ecal)]\,,\,
    [(\Ecal',\phi_{\Ecal'})]\right)}=\Ext^1\big(\Ecal,\Ecal'\otimes
\Ocal_{\Xscr_k}(-\Dscr_\infty) \big)\ .
\end{equation}
\begin{rem}\label{rem:carlssonokounkovbundle}
Let us take
$\Mcal_{r',\vec{u}\,',\Delta'}(\Xscr_k,\Dscr_\infty,\Fcal_\infty^{s,\vec
  w\,'}\,
)=\Mcal_{r,\vec{u},\Delta}(\Xscr_k,\Dscr_\infty,\Fcal_\infty^{s,\vec
  w}\, )$. Then the fibre of $\Ebf$ at a point
$\big([(\Ecal,\phi_\Ecal)]\,,\, [(\Ecal,\phi_{\Ecal})]\big)$ is
isomorphic to the tangent space of
$\Mcal_{r,\vec{u},\Delta}(\Xscr_k,\Dscr_\infty,\Fcal_\infty^{s,\vec
  w}\, )$ at $[(\Ecal,\phi_{\Ecal})]$.\footnote{We expect that the
  restriction of $\Ebf$ to the diagonal of
  $\Mcal_{r,\vec{u},\Delta}(\Xscr_k,\Dscr_\infty,\Fcal_\infty^{s,\vec
    w}\, )
  \times\Mcal_{r,\vec{u},\Delta}(\Xscr_k,\Dscr_\infty,\Fcal_\infty^{s,\vec
    w}\, )$ is isomorphic to the tangent bundle of
  $\Mcal_{r,\vec{u},\Delta}(\Xscr_k,\Dscr_\infty,\Fcal_\infty^{s,\vec
    w}\, )$ thanks to a framed version of the Kodaira-Spencer map
  which is established in \cite{art:sala2012} for moduli spaces of
  framed sheaves on a surface. We expect that this construction can be
  generalized to our moduli spaces
  $\Mcal_{r,\vec{u},\Delta}(\Xscr_k,\Dscr_\infty,\Fcal_\infty^{s,\vec
    w}\, )$.}  \hfill\mbox{\hspace{0.2mm}}\hfill$\triangle$
\end{rem}

\subsection{Torus action and fixed points}

Let us begin by recalling some definitions which will be used in the
combinatorial expressions below. Let
$Y\subset\N_{>0} \times\N_{>0} $ be a Young tableau, which we think of as sitting ``in the first quadrant".
Define the \emph{arm length} and \emph{leg length} of a box $s=(i,j)\in Y$  as $A(s)=A_Y(s):= \lambda_i-j$ and $L(s)=L_Y(s):= \lambda_j' - i$ respectively, where $\lambda_i$ is the length of the $i$-th column of $Y$ and $\lambda_j'$ is the length of the $j$-th row of $Y$. The \emph{arm colength} and \emph{leg colength} are  given by $A'(s)= A_Y'(s):=j-1$ and $L'(s)= L_Y'(s):=i-1$, respectively. We also define the \emph{weight} $\vert Y\vert$ of a Young tableau as the number of boxes $s\in Y$.

Let $T_\rho$ be the maximal torus of $GL(r,\C)$ consisting of diagonal
matrices and set $T:=T_t\times T_\rho$. We define an action of $T$ on
$\Mcal_{r,\vec{u},\Delta}(\Xscr_k,\Dscr_\infty,\Fcal_\infty^{s,\vec
  w}\, )$ as follows. For any element $(t_1, t_2)\in T_t$, let
$F_{(t_1, t_2)}$ be the automorphism of $\Xscr_k$ induced by the torus
action of $T_t$ on $\Xscr_k$. Define an action of the torus $T=
T_t\times T_\rho$ on
$\Mcal_{r,\vec{u},\Delta}(\Xscr_k,\Dscr_\infty,\Fcal_\infty^{s,\vec
  w}\, )$ by
\begin{equation}
(t_1,t_2,\vec\rho\, )\triangleright \big[(\Ecal,\phi_{\Ecal}) \big]
:=\big[\big( (F_{(t_1,t_2)}^{-1})^\ast(\Ecal)\,,\, \phi_\Ecal'\big) \big] \ ,
\end{equation}
where $\vec\rho=(\rho_1,\dots,\rho_r)\in T_\rho$ and $\phi_\Ecal'$ is the composition of isomorphisms
\begin{equation}
\phi_\Ecal' \, \colon\,  \big(F_{(t_1,t_2)}^{-1} \big)^\ast\Ecal_{|_{\Dscr_\infty}} \ \xrightarrow{(F_{(t_1,t_2)}^{-1})^\ast(\phi_{\Ecal})} \
\big(F_{(t_1,t_2)}^{-1} \big)^\ast \Fcal_\infty^{s,\vec w} \
\longrightarrow \ \Fcal_\infty^{s,\vec w} \ \xrightarrow{ \vec\rho \ \cdot} \ \Fcal_\infty^{s,\vec w} \ ;
\end{equation}
here the middle arrow is given by the $T_t$-equivariant structure
induced on $\Fcal_\infty^{s,\vec{w}}$ by restriction of the torus
action of $\Xscr_k$ to $\mathscr D_\infty$.

\begin{proposition}\label{prop:fixedpoint} 
A $T$-fixed point $[(\Ecal,\phi_\Ecal)]\in
\Mcal_{r,\vec{u},\Delta}\big(\Xscr_k,\Dscr_\infty,\Fcal_\infty^{s,\vec
  w}\, \big)^T$  decomposes as a direct sum of rank one framed sheaves
\begin{equation}
(\Ecal, \phi_{\Ecal})=\bigoplus_{\alpha=1}^r \, (\Ecal_\alpha, \phi_{\alpha})\ ,
\end{equation}
where for $\sum_{j=0}^{i-1}\, w_j<\alpha\leq\sum_{j=0}^i\,  w_j$ with
$i=0, 1, \ldots, k-1$ we have that:
\begin{itemize}
\item $\Ecal_\alpha$ is a tensor product
  $\imath_\ast(I_\alpha)\otimes\Rcal^{\vec{u}_\alpha}\otimes\Ocal_{\Xscr_k}(s
  \, \Dscr_\infty)$, where $I_\alpha$ is an ideal sheaf of a zero-dimensional subscheme $Z_\alpha$ of $X_k$ supported at the $T_t$-fixed points $p_1,\ldots,p_k$ and $\vec{u}_\alpha\in\Z^{k-1}$ is such that 
\begin{equation}\label{eq:decompositionalpha}
\sum_{j=1}^{k-1}\, j \, (\vec{u}_\alpha)_j=  i \ \bmod{k} \ ;
\end{equation}
\item The framing $\phi_\alpha\colon
  \Ecal_\alpha \xrightarrow{\sim}\Ocal_{\Dscr_\infty}(s,i)$ is induced
  by the isomorphism $\Rcal^{\vec{u}_\alpha}\otimes\Ocal_{\Xscr_k}(s
  \, \Dscr_\infty)_{\vert \Dscr_\infty}\simeq \Ocal_{\Dscr_\infty}(s,i)$.
\end{itemize}
\end{proposition}
\proof
By using the same arguments as in the proof of an analogous result for framed sheaves on smooth projective surfaces \cite[Proposition~3.2]{art:bruzzopoghossiantanzini2011}, we obtain a decomposition 
\begin{equation}
\Ecal=\bigoplus_{\alpha=1}^r \, \Ecal_\alpha
\end{equation}
where each $\Ecal_\alpha$ is a $T$-invariant rank one torsion-free
sheaf on $\Xscr_k$. The restriction ${\phi_\Ecal}_{\vert\Ecal_\alpha}$
gives a canonical framing to a direct summand of $\Fcal_\infty^{s,\vec
  w}$. By reordering the indices $\alpha$ if necessary, for $i=0,1,
\ldots,k-1$ and for each $\alpha$ such that $\sum_{j=0}^{i-1}\,
w_j<\alpha\leq\sum_{j=0}^i \, w_j$ we get an induced framing on $\Ecal_\alpha$
\begin{equation}
\phi_\alpha:={\phi_\Ecal}_{\vert\Ecal_\alpha}\, \colon\, \Ecal_\alpha
\ \xrightarrow{ \ \sim \ } \ \Ocal_{\Dscr_\infty}(s,i)\ .
\end{equation}
Thus $(\Ecal_\alpha,\phi_\alpha)$ is a $(\Dscr_\infty,
\Ocal_{\Dscr_\infty}(s,i))$-framed sheaf of rank one on $\Xscr_k$. As
explained in Remark \ref{rem:rankone}, the torsion-free sheaf
$\Ecal_\alpha$ is the tensor product of an ideal sheaf $I_\alpha$ of a
zero-dimensional subscheme $Z_\alpha$ of length $n_\alpha$ supported
on $X_k$ and the line bundle
$\Rcal^{\vec{u}_\alpha}\otimes\Ocal_{\Xscr_k}(s\, \Dscr_\infty)$ for a
vector $\vec{u}_\alpha\in \Z^{k-1}$ that satisfies Equation
\eqref{eq:decompositionalpha} in view  of Lemma
\ref{lem:firstchernclass}. Since the torsion-free sheaf $\Ecal$ is
fixed by the $T_t$-action, $Z_\alpha$ is fixed as well and so it is supported at the $T_t$-fixed points $p_1,\ldots,p_k$.
\endproof
Let $\big[(\Ecal, \phi_{\Ecal}) \big]=\big[\bigoplus_{\alpha=1}^r \,
(\Ecal_\alpha, \phi_{\alpha}) \big]$ be a $T$-fixed point in
$\Mcal_{r,\vec{u},\Delta}(\Xscr_k,\Dscr_\infty,\Fcal_\infty^{s,\vec
  w}\, )$. Then 
\begin{equation}
\Rcal^{\vec{u}}\otimes \Ocal_{\Xscr_k}(s\,r \, \Dscr_\infty)\simeq
\det(\Ecal)\simeq\bigotimes_{\alpha=1}^r \, \det(\Ecal_\alpha)\simeq
\bigotimes_{\alpha=1}^r\, \big(\Rcal^{\vec{u}_\alpha}\otimes
\Ocal_{\Xscr_k}(s\, \Dscr_\infty) \big)\ ,
\end{equation}
hence $\sum_{\alpha=1}^r \, \vec{u}_\alpha=\vec{u}$. On the other
hand, $I_\alpha$ is the ideal sheaf of a $T_t$-fixed zero-dimensional
subscheme $Z_\alpha$ of length $n_\alpha$ for $\alpha\in\{1, \ldots,
r\}$. Hence it is a disjoint union of zero-dimensional subschemes
$Z_\alpha^i$ supported at the $T_t$-fixed points $p_i$ for $i=1,
\ldots, k$; each $Z_\alpha^i$ corresponds to a Young tableau
$Y_\alpha^i$ \cite{art:ellingsrudstromme1987}, and $Z_\alpha$
corresponds to the set of Young tableaux
$\vec{Y}_\alpha=\{Y_\alpha^i\}_{i=1,\ldots, k}$ such that
$\sum_{i=1}^k\, \vert Y_\alpha^i\vert = n_\alpha$.

Thus we can parametrize the point $[(\Ecal, \phi_{\Ecal})]$ by the pair $(\vec{\boldsymbol{Y}},\vec{\ubf})$, where 
\begin{itemize}
\item $\vec{\boldsymbol{Y}}=(\vec{Y}_1, \ldots, \vec{Y}_r)$, where
  $\vec{Y}_\alpha=\{Y^i_\alpha\}_{i=1, \ldots, k}$ for any $\alpha=1,
  \ldots, r$ is a set of Young tableaux such that $\sum_{i=1}^k\, \vert Y_\alpha^i\vert = n_\alpha$;
\item $\vec{\ubf}=(\vec{u}_1, \ldots, \vec{u}_r)$, where
  $\vec{u}_\alpha=((\vec{u}_\alpha)_1, \ldots,
  (\vec{u}_\alpha)_{k-1})$ for any $\alpha=1, \ldots, r$ is an integer
  vector such that the relation \eqref{eq:decompositionalpha} holds
  and $\sum_{\alpha=1}^r\, \vec{u}_\alpha=\vec{u}$.
\end{itemize}
If we set $\vec{v}_\alpha:=C^{-1}\vec{u}_\alpha$ for $\alpha=1, \ldots, r$, we denote the same point by $(\vec{\boldsymbol{Y}},\vec{\boldsymbol{v}})$ where $\vec{\boldsymbol{v}}=(\vec{v}_1, \ldots, \vec{v}_r)$. In this case any $\vec{v}_\alpha$ satisfies the relation
\begin{equation}
k\, (\vec{v}_\alpha)_{k-1}= i \ \bmod{k}
\end{equation}
if $\sum_{j=0}^{i-1}\, w_j<\alpha\leq\sum_{j=0}^i\, w_j$ for $i=0,1, \ldots, k-1$.
We shall call these pairs the \emph{combinatorial data} of the torus-fixed point $[(\Ecal, \phi_{\Ecal})]$.

\begin{remark}
It is easy to see that
\begin{align}
\int_{\Xscr_k}\, \ch_2(\Ecal)&= \sum_{\alpha=1}^r \ \int_{\Xscr_k}\,
\ch_2 \big(\imath_\ast
I_\alpha\otimes\Rcal^{\vec{u}_\alpha}\otimes\Ocal_{\Xscr_k}(s\,
\Dscr_\infty) \big) \\[4pt]
&=\sum_{\alpha=1}^r \, \Big(\, \int_{\Xscr_k}\, \frac{1}{2}\,
\crm_1\big(\Rcal^{\vec{u}_\alpha}\otimes\Ocal_{\Xscr_k}(s\,
\Dscr_\infty)\big)^2-n_\alpha\, \Big)
=\sum_{\alpha=1}^r\, \Big(\, \frac{s^2}{2k\,
  \tilde{k}^2}-\frac{1}{2}\, \vec{v}_\alpha\cdot C\vec{v}_\alpha
-n_\alpha\, \Big) \\
&=\frac{r\, s^2}{2k\, \tilde{k}^2}-\frac{1}{2}\, \sum_{\alpha=1}^r\,
\vec{v}_\alpha\cdot C\vec{v}_\alpha -\sum_{\alpha=1}^r\,  n_\alpha \
\in \ \mbox{$\frac{1}{2k\, \tilde{k}^2}$}\, \Z\ .
\end{align}
Then
\begin{equation}
\int_{\Xscr_k}\, \crm_2(\Ecal)= \frac{r\, (r-1)\, s^2}{2k\,
  \tilde{k}^2}+\sum_{\alpha=1}^r\, n_\alpha -\frac{1}{2}\,
\sum_{\alpha\neq \beta} \, \vec{v}_\alpha\cdot C\vec{v}_\beta \ \in \
\mbox{$\frac{1}{2k\, \tilde{k}^2}$}\, \Z\ ,
\end{equation}
and therefore
\begin{equation}
\Delta=\sum_{\alpha=1}^r \, n_\alpha+\frac{1}{2r}\,
\sum_{\alpha=1}^r\, \vec{v}_\alpha\cdot C\vec{v}_\alpha-\frac{1}{2r}\,
\sum_{\alpha, \beta=1}^r\,  \vec{v}_\alpha\cdot C\vec{v}_\beta \ \in \
\mbox{$\frac{1}{2r\, k}$}\, \Z\ .
\end{equation}
As a by-product, this computation shows that the discriminant of any
$(\Dscr_\infty, \Fcal_\infty^{s, \vec{w}}\, )$-framed sheaf on
$\Xscr_k$ takes values in $ \frac{1}{2r\, k}\, \Z$.
\end{remark}

\subsection{Euler classes\label{sec:Eulerclasses}}

We introduce the equivariant parameters of the torus $T=T_t\times T_\rho$. For $\alpha=1, \ldots, r$, let $e_\alpha$ be the one-dimensional $T_\rho$-module corresponding to the projection $(\mathbb{C}^\ast)^r\to \mathbb{C}^\ast$ onto the $\alpha$-th factor and $a_\alpha$ its equivariant first Chern class.
The corresponding $T_t$-module parameters $t_j$ and $\varepsilon_j$
for $j=1,2$ were introduced in Section
\ref{sec:minimalresolution}. Then $H_{T}^\ast({\rm pt};\mathbb{Q})=\mathbb{Q}[\varepsilon_1, \varepsilon_2,a_1, \ldots, a_r]$.

On $\Mcal_{r,\vec{u},\Delta}(\Xscr_k,\Dscr_\infty,\Fcal_\infty^{s,\vec
  w}\, )\times
\Mcal_{r',\vec{u}\,',\Delta'}(\Xscr_k,\Dscr_\infty,\Fcal_\infty^{s,\vec
  w\,'}\, )$ there is a natural action of the extended torus $\tilde T=T_t\times T_{\rho}\times T_{\rho'}$, which acts as $T_t \times T_{\rho}$ on the first factor ($T_{\rho'}$ acting trivially) and as $T_t\times T_{\rho'}$ on the second factor ($T_{\rho}$ acting trivially). We want to compute the character 
$
\ch_{\tilde{T}}\Ebf_{\left([(\Ecal,\phi_\Ecal)]\,,\, [(\Ecal',\phi_{\Ecal'})]\right)}
$
of the Carlsson-Okounkov bundle at a fixed point $$\big([(\Ecal,\phi_\Ecal)]\,,\,
[(\Ecal',\phi_{\Ecal'})]\big)\in
\big(\Mcal_{r,\vec{u},\Delta}(\Xscr_k,\Dscr_\infty,\Fcal_\infty^{s,\vec
  w}\, )\times \Mcal_{r',\vec{u}\,',\Delta'}(\Xscr_k,\Dscr_\infty
\Fcal_\infty^{s,\vec w\,'}\, ) \big)^{\tilde T}\ .$$
 Let
$\big((\vec{\boldsymbol{Y}},\vec{\ubf})\,,\,(\vec{\boldsymbol{Y}}{}',\vec{\ubf}{}')\big)$
be the corresponding combinatorial data. Since the torsion-free
sheaves $\Ecal$ and $\Ecal'$ decompose as
\begin{equation}
\Ecal= \bigoplus_{\alpha=1}^{r}\,
\imath_\ast(I_\alpha)\otimes\Rcal^{\vec{u}_\alpha}\otimes\Ocal_{\Xscr_k}(s\,
\Dscr_\infty) \qquad\mbox{and}\qquad \Ecal'= \bigoplus_{\beta=1}^{r'}
\,
\imath_\ast(I_\beta')\otimes\Rcal^{\vec{u}_\beta{}'}\otimes\Ocal_{\Xscr_k}(s\,
\Dscr_\infty) \ ,
\end{equation}
we get
\begin{align}
\ch_{\tilde{T}}\Ebf_{\left([(\Ecal,\phi_\Ecal)]\,,\, [(\Ecal',\phi_{\Ecal'})]\right)} &= \ch_{\tilde{T}} \Ext^1\big(\Ecal,\Ecal'\otimes \Ocal_{\Xscr_k}(-\Dscr_\infty)\big)\\[4pt]
=-\sum_{\alpha=1}^r \ & \sum_{\beta=1}^{r'} \, \ch_{\tilde{T}}
\Ext^\bullet\big(\imath_\ast(I_\alpha)\otimes\Rcal^{\vec{u}_\alpha},\imath_\ast(I_\beta')\otimes\Rcal^{\vec{u}_\beta{}'}\otimes
\Ocal_{\Xscr_k}(-\Dscr_\infty) \big)\\[4pt]
 =-\sum_{\alpha=1}^r \ \sum_{\beta=1}^{r'}\, & e_\beta\,
 e_\alpha^{\prime\, -1}\, \ch_{T_t}
 \Ext^\bullet\big(\imath_\ast(I_\alpha) \otimes\Rcal^{\vec{u}_\alpha}\,,\,\imath_\ast(I_\beta')
 \otimes\Rcal^{\vec{u}_\beta{}'}\otimes\Ocal_{\Xscr_k}(-\Dscr_\infty)
 \big)\ .
\end{align}
Let
\begin{align}\label{eq:edge}
L_{\alpha\beta}(t_1,t_2):=& -\ch_{T_t}
\Ext^\bullet\big(\Rcal^{\vec{u}_\alpha},\Rcal^{\vec{u}_\beta{}'}\otimes
\Ocal_{\Xscr_k}(-\Dscr_\infty) \big)\\[4pt] 
=& -\chi_{T_t}\big(\Xscr_k,\Rcal^{\vec{u}_\beta{}'-\vec{u}_\alpha}\otimes
\Ocal_{\Xscr_k}(-\Dscr_\infty) \big)\ ,
\end{align}
and
\begin{align}\label{eq:vertex}
M_{\alpha\beta}(t_1,t_2):=&\ch_{T_t}
\Ext^\bullet\big(\Rcal^{\vec{u}_\alpha},\Rcal^{\vec{u}_\beta{}'}\otimes
\Ocal_{\Xscr_k}(-\Dscr_\infty) \big)\\ &
-\,\ch_{T_t}
\Ext^\bullet\big(\imath_\ast(I_\alpha)\otimes\Rcal^{\vec{u}_\alpha},\imath_\ast(I_\beta')\otimes\Rcal^{\vec{u}_\beta{}'}\otimes
\Ocal_{\Xscr_k}(-\Dscr_\infty) \big)\ .
\end{align}
Then
\begin{equation}\label{eq:character-carlssonokounkovbundle}
\ch_{\tilde{T}}\Ebf_{\left([(\Ecal,\phi_\Ecal)]\,,\,
    [(\Ecal',\phi_{\Ecal'})]\right)} =\sum_{\alpha=1}^r\
\sum_{\beta=1}^{r'}\, e_\beta\, e_\alpha^{\prime\, -1}\,
\big(M_{\alpha\beta}(t_1,t_2)+L_{\alpha\beta}(t_1,t_2) \big)\ .
\end{equation}

Here and in the following we use an approach similar to that of
\cite{art:gasparimliu2010}, which computes the character of the
tangent bundle at a fixed point of the moduli space of framed sheaves
on a smooth projective toric surface, and so we shall borrow some of
their terminology. In particular, we call $M_{\alpha\beta}(t_1,t_2)$ a
\emph{vertex contribution} to $\ch_{\tilde{T}}\Ebf_{
  \left([(\Ecal,\phi_\Ecal)]\,,\, [(\Ecal',\phi_{\Ecal'})]\right)}$
and $L_{\alpha\beta}(t_1,t_2)$ an \emph{edge contribution} to
$\ch_{\tilde{T}}\Ebf_{\left([(\Ecal,\phi_\Ecal)]\,,\,
    [(\Ecal',\phi_{\Ecal'})]\right)}$. The paper
\cite{art:gasparimliu2010} uses this terminology because
$M_{\alpha\beta}(t_1,t_2)$ will depend on the torus-fixed points $p_i$
of
$X_k$ (which one can match to vertices of the toric graph) and
$L_{\alpha\beta}(t_1,t_2)$ will depend on the torus-invariant divisors $D_i$
of $X_k$ (which one can match to edges of the toric graph such that an edge joins two vertices if and only if the corresponding fixed points lie in the divisor corresponding to the edge).

\subsubsection*{Edge contribution}

By definition we have
\begin{equation}
L_{\alpha\beta}(t_1,t_2)= -\chi_{T_t}\Big(\Xscr_k\,,\,
\mbox{$\bigotimes\limits_{j=1}^{k-1}$}\, \Rcal_j^{\otimes (u_\beta')_j-(u_\alpha)_j}\otimes \Ocal_{\Xscr_k}(-\Dscr_\infty)\Big)\ .
\end{equation}
The computation of this character is done in Appendix \ref{sec:edgecontribution-result}.

\subsubsection*{Vertex contribution}

\begin{proposition}\label{prop:vertexcontribution}
\begin{equation}\label{eq:vertexcontribution}
M_{\alpha\beta}(t_1,t_2)=\sum_{i=1}^k\,
(\chi_1^i)^{(\vec{v}_\beta{}'\,)_{i}-(\vec{v}_\alpha)_{i}}\,
(\chi_2^i)^{(\vec{v}_\beta{}'\,)_{i-1}-(\vec{v}_\alpha)_{i-1}} \,
M_{Y_{\alpha}^{i}, {Y_{\beta}^{i}}{}'}(\chi_1^i,\chi_2^i)
\end{equation}
where $\chi_1^i$ and $\chi_2^i$ are introduced in Section \ref{sec:minimalresolution}, and given two Young tableaux $Y, Y'$ ,we set
\begin{equation}
M_{Y,Y'}(x,y):=\sum_{s\in Y}\, x^{-L_{Y'}(s)}\,
y^{A_Y(s)+1}+\sum_{s'\in Y'}\, x^{L_Y(s'\, )+1}\, y^{-A_{Y'}(s'\, )}\ .
\end{equation}
\end{proposition}
Before proving this Proposition, we need some preliminary results. By
Proposition \ref{prop:opensubstack} the open substack $\Uscr_i$ of
$\Xscr_k$ corresponding to $\sigma_i$ is $[V_i/N(\sigma_i)]\simeq U_i
\simeq \C^2$ for $i=1, \ldots, k$, and the open substacks $\Uscr_{k+1},\Uscr_{k+2}$ of
$\Xscr_k$ corresponding to $\sigma_{\infty, j}$ for
$j=0,k$ are
$[V_{k+1}/N(\sigma_{\infty,0})], [V_{k+2}/N(\sigma_{\infty,k})]\simeq [\C^2/\mu_{k\, \tilde{k}}]$.
Set $\Ubf=\bigsqcup_{i=1}^{k+2}\, V_i$. Since the morphisms $\Ubf \to
\bigsqcup_{i=1}^{k+2} \, \Uscr_i$ and $\bigsqcup_{i=1}^{k+2}\, \Uscr_i
\to \Xscr_k$ are \'etale and surjective, the composition $\ubf\colon
\Ubf\to \Xscr_k$ is \'etale and surjective as well, hence the pair
$(\Ubf,\ubf)$ is an \'etale presentation of $\Xscr_k$. Denote by
$\Ubf_\bullet\to \Xscr_k$ the \emph{strictly} simplicial algebraic
  space associated to the simplicial algebraic space which is obtained by taking the
0-coskeleton of $(\Ubf,\ubf)$ \cite[Section 4.1]{art:olsson2007}; for
any $n\geq 0$ we have
\begin{equation}
\Ubf_n=\bigsqcup_{\stackrel{\scriptstyle i_0, i_1, \ldots, i_k
    \in\{1,\ldots, k+2\} }{\scriptstyle i_0< i_1<\cdots <i_n}}\,  V_{i_0}\times_{\Xscr_k}V_{i_1}\times_{\Xscr_k}\cdots\times_{\Xscr_k} V_{i_n}\ .
\end{equation}
By \cite[Proposition~6.12]{art:olsson2007}, the category of coherent
sheaves on $\Xscr_k$ is equivalent to the category of
\emph{simplicial} coherent sheaves on $\Ubf_\bullet$ (see
\cite{art:olsson2007} for the definition of simplicial coherent sheaf on a strictly simplicial algebraic space).
\proof[Proof of Proposition \ref{prop:vertexcontribution}]
As explained in \cite[Section~6]{art:olsson2007}, one has an isomorphism between the Ext groups of coherent sheaves on $\Xscr_k$ and the Ext groups of the corresponding simplicial coherent sheaves on $\Ubf_\bullet$. Thus
\begin{multline}
\Ext^\bullet\big(\Rcal^{\vec{u}_\alpha},\Rcal^{\vec{u}_\beta{}'}\otimes
\Ocal_{\Xscr_k}(-\Dscr_\infty) \big)
 -\Ext^\bullet\big(\imath_\ast(I_\alpha)\otimes\Rcal^{\vec{u}_\alpha},\imath_\ast(I_\beta')\otimes\Rcal^{\vec{u}_\beta{}'}\otimes
 \Ocal_{\Xscr_k}(-\Dscr_\infty) \big)  = \\
 \Ext^\bullet\big(\Rcal^{\vec{u}_\alpha}_{\vert
  \Ubf_\bullet},\Rcal^{\vec{u}_\beta{}'}\otimes
\Ocal_{\Xscr_k}(-\Dscr_\infty)_{\vert \Ubf_\bullet} \big)
 -\Ext^\bullet\big(\imath_\ast(I_\alpha)\otimes\Rcal^{\vec{u}_\alpha}_{\vert
   \Ubf_\bullet},\imath_\ast(I_\beta')\otimes\Rcal^{\vec{u}_\beta{}'}\otimes
 \Ocal_{\Xscr_k}(-\Dscr_\infty)_{\vert \Ubf_\bullet} \big)\ ,
\end{multline}
where for a coherent sheaf $\Gcal$ on $\Xscr_k$ we denote by $\Gcal_{\vert \Ubf_\bullet}$ the corresponding simplicial coherent sheaf on $\Ubf_\bullet$ \cite[Proposition~6.12]{art:olsson2007}.

Recall that $I_\alpha$ and $I_\beta'$ are ideal sheaves of
zero-dimensional subschemes $Z_\alpha$ and $Z_\beta'$ supported at the
$T_t$-fixed points $p_1, \ldots, p_k$ of $X_k$. Hence the restrictions
of $\imath_\ast I_\alpha$ and $\imath_\ast I_\beta'$ on $\Uscr_{j}$
are trivial for $j=k+1,k+2$. For the same reason, the restrictions of
$\imath_\ast I_\alpha$ and $\imath_\ast I_\beta'$ on
$\Uscr_i\times_{\Xscr_k} \Uscr_j$ are also trivial since
$\Uscr_i\times_{\Xscr_k} \Uscr_j\simeq U_i\cap U_j$ for $i,j =1, \ldots, k$. Then for pairwise distinct indices $j_1,\ldots, j_l\in\{1, \ldots, k+2\}$ we get
\begin{equation}
\imath_\ast {I_\alpha}_{\vert \Uscr_{j_1}\times_{\Xscr_k}\Uscr_{j_2}\times_{\Xscr_k}\cdots \times_{\Xscr_k}\Uscr_{j_l}}\simeq {\Ocal_{\Xscr_k}}_{\vert \Uscr_{j_1}\times_{\Xscr_k}\Uscr_{j_2}\times_{\Xscr_k}\cdots \times_{\Xscr_k}\Uscr_{j_l}}\simeq
\imath_\ast {I_\beta'}_{\vert \Uscr_{j_1}\times_{\Xscr_k}\Uscr_{j_2}\times_{\Xscr_k}\cdots \times_{\Xscr_k}\Uscr_{j_l}}
\end{equation}
unless $i=1$ and $j_1=1, \ldots, k$. Then $\big(\imath_\ast
{I_\alpha}_{\vert \Ubf_\bullet}\big)_{\vert \Ubf_n}\simeq
{\Ocal_{\Ubf_\bullet}}_{\vert \Ubf_n}\simeq
\big(\imath_\ast{I_\beta'}_{\vert \Ubf_\bullet}\big)_{\vert
  \Ubf_n}$ for $n\geq
1$. So by using the local-to-global spectral sequence (which
degenerates since $\Ubf$ is a disjoint union of affine spaces) we find
\begin{multline}
\Ext^\bullet\big(\Rcal^{\vec{u}_\alpha},\Rcal^{\vec{u}_\beta{}'}\otimes
\Ocal_{\Xscr_k}(-\Dscr_\infty) \big)
 -\Ext^\bullet\big(\imath_\ast(I_\alpha)\otimes\Rcal^{\vec{u}_\alpha},\imath_\ast(I_\beta')\otimes\Rcal^{\vec{u}_\beta{}'}\otimes
 \Ocal_{\Xscr_k}(-\Dscr_\infty) \big)\\
=\sum_{i=1}^k \ \sum_{j=0}^2 \, (-1)^j\,
\Big(H^0\big(U_i,{\Ocal^j_{\alpha\beta}}_{\vert X_k}
\big)-H^0\big(U_i,{\Ecal_{\alpha\beta}^j}_{\vert X_k} \big)\Big)\ ,
\end{multline}
where
\begin{align}
\Ocal^j_{\alpha\beta}&:=\Ecal
xt^j\big(\Rcal^{\vec{u}_\alpha},\Rcal^{\vec{u}_\beta{}'}\otimes
\Ocal_{\Xscr_k}(-\Dscr_\infty) \big)\ ,\\[4pt]
\Ecal^j_{\alpha\beta}&:=\Ecal
xt^j\big(\imath_\ast(I_\alpha)\otimes\Rcal^{\vec{u}_\alpha},\imath_\ast(I_\beta')\otimes\Rcal^{\vec{u}_\beta{}'}\otimes
\Ocal_{\Xscr_k}(-\Dscr_\infty) \big)\ .
\end{align}
By using the same arguments as in the proof of  \cite[Proposition~5.1]{art:gasparimliu2010}, where $M_{\alpha\beta}(t_1,t_2)$ is computed for framed sheaves on smooth projective toric surfaces, we get
\begin{equation*}
M_{\alpha\beta}(t_1,t_2)=\sum_{i=1}^k\,
\frac{\operatorname{ch}_{T_t}\big({\mathcal R}^{\vec{u}_\beta
{}'}_{\vert U_i}\big)}{\operatorname{ch}_{T_t}\big({\mathcal R}^{\vec
{u}_\alpha}_{\vert U_i}\big)}
\, M_{Y_{\alpha}^{i}, {Y_{\beta}^{i}}{}'}(\chi_1^i,\chi_2^i)\ .
\end{equation*}
The computation of $\operatorname{ch}_{T_t}\big({\mathcal R}^{\vec
{u}}_{\vert U_i}\big)$ for
$i=1, \ldots, k$ and any vector $\vec{u}\in{\mathbb{Z}}^{k-1}$ can
be done by using Lemma~\ref{lem:characterline} and the relation \eqref{eq:tautologicalclasses}.
\endproof

\subsubsection*{Euler class of the Carlsson-Okounkov bundle at a fixed point}

Let us introduce the following notation. Set $a_{\beta\alpha}:=a_\beta'-a_\alpha$, $\vec{v}_{\beta\alpha}:=\vec{v}_\beta{\, '}-\vec{v}_\alpha$ and
\begin{equation}
a_{\beta\alpha}^{(i)}:=a_{\beta\alpha}+(\vec{v}_{\beta\alpha})_{i}\,
\varepsilon_1^{(i)}+(\vec{v}_{\beta\alpha})_{i-1}\, \varepsilon_2^{(i)}
\end{equation}
for $i=1, \ldots, k$ and for $\alpha=1, \ldots, r$, $\beta=1, \ldots, r'$ (we 
set $(\vec{v}_{\beta\alpha})_0=(\vec{v}_{\beta\alpha})_k=0$). Let
$c_{\beta\alpha}$ be the equivalence class of $k\,
(\vec{v}_{\beta\alpha})_{k-1}$ modulo $k$. Set $(C^{-1})^{n 0}=0$ for
$n\in\{1, \ldots, k-1\}$ and $(C^{-1})^{k, c_{\beta\alpha}}=0$. We further define
\begin{align}
m_{Y_\alpha, {Y_\beta}'}(\varepsilon_1, \varepsilon_2,
a_{\beta\alpha}):= &\prod_{s\in Y_\alpha}\,
\big(a_{\beta\alpha}-L_{{Y_\beta}'}(s)\,
\varepsilon_1+(A_{Y_\alpha}(s)+1)\, \varepsilon_2\big) \\ & \times \
\prod_{s' \in
  {Y_\beta}'}\, \big(a_{\beta\alpha}+(L_{Y_\alpha}(s'\,)+1)\,
\varepsilon_1-A_{{Y_\beta}'}(s'\,)\, \varepsilon_2 \big)
\end{align}
for two Young tableaux $Y_\alpha$ and ${Y_\beta}'$, and $\alpha=1,
\ldots, r$, $\beta=1, \ldots, r'$.

The $\tilde{T}$-equivariant Euler class of the vertex contribution $M_{\alpha\beta}(t_1,t_2)$ is
\begin{equation}
\prod_{i=1}^k \, m_{Y_{\alpha}^{i},
  {Y_{\beta}^{i}}'}\big(\varepsilon_1^{(i)},\varepsilon_2^{(i)},
a^{(i)}_{\beta\alpha} \big)\ .
\end{equation}
Now we give the $\tilde{T}$-equivariant Euler class of the edge
contribution. For this, we introduce some more notation. 

If $(\vec{v}_{\beta\alpha})_n-(C^{-1})^{n, c_{\beta\alpha}}> 0$ for $n\in\{1, \ldots, k-1\}$, consider the equation
\small
\begin{multline}\label{eq:cond-index+}
i^2-i\, \Big(\vec{v}_{\beta\alpha}-\sum_{p=1}^{n-1}\,
\big((\vec{v}_{\beta\alpha})_p-(C^{-1})^{p, c_{\beta\alpha}}\big)\,
\vec{e}_p\Big)\cdot C\vec{e}_n +\frac{1}{2}\,
\bigg(\Big(\vec{v}_{\beta\alpha}-\sum_{p=1}^{n-1}\,
\big((\vec{v}_{\beta\alpha})_p-(C^{-1})^{p, c_{\beta\alpha}}\big)\,
\vec{e}_p\Big) \\
\cdot C\Big(\vec{v}_{\beta\alpha}-\sum_{p=1}^{n-1}\,
\big((\vec{v}_{\beta\alpha})_p-(C^{-1})^{p c_{\beta\alpha}}\big)\,
\vec{e}_p\Big) -(C^{-1})^{c_{\beta\alpha}, c_{\beta\alpha}}\bigg)=0\ ,
\end{multline}
\normalsize
and define the set 
\begin{equation}
\splus_n:=\big\{i\in\N\, \big\vert\, i\leq
(\vec{v}_{\beta\alpha})_n-(C^{-1})^{n c_{\beta\alpha}} \mbox{ is a
  solution of Equation \eqref{eq:cond-index+}} \big\}\ .
\end{equation}
Let $\dplus_n=\min(\splus_n)$ if $\splus_n\neq \emptyset$, otherwise $\dplus_n:=(\vec{v}_{\beta\alpha})_n-(C^{-1})^{n, c_{\beta\alpha}}$.

If $(\vec{v}_{\beta\alpha})_n-(C^{-1})^{n, c_{\beta\alpha}}< 0$, consider the equation
\small
\begin{multline}\label{eq:cond-index-}
i^2+i\, \Big(\vec{v}_{\beta\alpha}-\sum_{p=1}^{n-1}\,
\big((\vec{v}_{\beta\alpha})_p-(C^{-1})^{p, c_{\beta\alpha}}\big)\,
\vec{e}_p\Big)\cdot C\vec{e}_n +\frac{1}{2}\,
\bigg(\Big(\vec{v}_{\beta\alpha}-\sum_{p=1}^{n-1}\,
\big((\vec{v}_{\beta\alpha})_p-(C^{-1})^{p, c_{\beta\alpha}}\big)\,
\vec{e}_p\Big) \\
\cdot C\Big(\vec{v}_{\beta\alpha}-\sum_{p=1}^{n-1}\,
\big((\vec{v}_{\beta\alpha})_p-(C^{-1})^{p, c_{\beta\alpha}}\big)\,
\vec{e}_p\Big)-(C^{-1})^{c_{\beta\alpha}, c_{\beta\alpha}}\bigg) =0\ ,
\end{multline}
\normalsize
and define the set 
\begin{equation}
\sminus_n:=\big\{i\in\N\, \big\vert\, i\leq
-(\vec{v}_{\beta\alpha})_n+(C^{-1})^{n c_{\beta\alpha}} \mbox{ is a
  solution of Equation \eqref{eq:cond-index-}} \big\}\ .
\end{equation}
Let $\dminus_n=\min(\sminus_n)$ if $\sminus_n\neq \emptyset$, otherwise $\dminus_n:=-(\vec{v}_{\beta\alpha})_n+(C^{-1})^{n, c_{\beta\alpha}}$. Define $m_{\beta\alpha}$ to be the smallest index $n\in\{1, \ldots, k-1\}$ such that $\splus_n$ or $\sminus_n$ is nonempty, otherwise $m_{\beta\alpha}:=k-1$.

Then the $\tilde{T}$-equivariant Euler class of the edge contribution $L_{\alpha\beta}(t_1, t_2)$ is
\begin{equation}
\prod_{n=1}^{k-1} \,
\ell^{(n)}_{\vec{v}_{\beta\alpha}}\big(\varepsilon_1^{(n)},\varepsilon_2^{(n)},
a_{\beta\alpha} \big)\ ,
\end{equation}
where for fixed $n=1, \ldots, m_{\beta\alpha}$ we set:
\begin{itemize} 
\item[\scriptsize$\blacksquare$] If $(\vec{v}_{\beta\alpha})_n-(C^{-1})^{n, c_{\beta\alpha}}> 0$:
\smallskip
\begin{itemize} 
\item[$\bullet$] For $\delta_{n,c_{\beta\alpha}} -(\vec{v}_{\beta\alpha})_{n+1}+(C^{-1})^{n+1, c_{\beta\alpha}}+2((\vec{v}_{\beta\alpha})_{n}-(C^{-1})^{n, c_{\beta\alpha}}-\dplus_n)\geq 0$:
\begin{multline}
\ell^{(n)}_{\vec{v}_{\beta\alpha}}\big(\varepsilon_1^{(n)},
\varepsilon_2^{(n)}, a_{\beta\alpha} \big)=
\prod_{i=(\vec{v}_{\beta\alpha})_{n}-(C^{-1})^{n,
    c_{\beta\alpha}}-\dplus_n}^{(\vec{v}_{\beta\alpha})_n-(C^{-1})^{n
    , c_{\beta\alpha}}-1}\hspace{4mm} \prod_{j=0}^{2i+\delta_{n,c_{\beta\alpha}} -(\vec{v}_{\beta\alpha})_{n+1}+(C^{-1})^{n+1, c_{\beta\alpha}}}\\ 
\bigg(a_{\beta\alpha}+\Big(i+\left\lfloor\frac{\delta_{n,c_{\beta\alpha}}
    -(\vec{v}_{\beta\alpha})_{n+1}+(C^{-1})^{n+1, c_{\beta\alpha}}}{2}
\right\rfloor\Big)\, \varepsilon_1^{(n)}+j\, \varepsilon_2^{(n)}\bigg)\ .
\end{multline}
\smallskip
\item[$\bullet$]For $2\leq \delta_{n,c_{\beta\alpha}}-(\vec{v}_{\beta\alpha})_{n+1}+(C^{-1})^{n+1, c_{\beta\alpha}}+2((\vec{v}_{\beta\alpha})_n-(C^{-1})^{n, c_{\beta\alpha}})<2 \dplus_n$:
\begin{multline}
\ell^{(n)}_{\vec{v}_{\beta\alpha}}\big(\varepsilon_1^{(n)},
\varepsilon_2^{(n)}, a_{\beta\alpha}\big)=
\prod_{i=(\vec{v}_{\beta\alpha})_n-(C^{-1})^{n, c_{\beta\alpha}}-\dplus_n}^{-\big\lfloor \frac{\delta_{n,c_{\beta\alpha}} -(\vec{v}_{\beta\alpha})_{n+1}+(C^{-1})^{n+1, c_{\beta\alpha}}}{2}\big\rfloor - 1}\hspace{4mm} \prod_{j=1}^{2i- (\delta_{n,c_{\beta\alpha}} -(\vec{v}_{\beta\alpha})_{n+1}+(C^{-1})^{n+1, c_{\beta\alpha}})-1}\\
\shoveright{\bigg(a_{\beta\alpha}+\Big(i-\left\lfloor
    -\frac{\delta_{n,c_{\beta\alpha}}
      -(\vec{v}_{\beta\alpha})_{n+1}+(C^{-1})^{n+1,
        c_{\beta\alpha}}}{2}\right\rfloor\Big)\,
  \varepsilon_1^{(n)}-j\, \varepsilon_2^{(n)}\bigg)^{-1}\times}   \\[4pt]
\shoveleft{ \ \prod_{i=-\big\lfloor \frac{\delta_{n,c_{\beta\alpha}} -(\vec{v}_{\beta\alpha})_{n+1}+(C^{-1})^{n+1, c_{\beta\alpha}}}{2}\big\rfloor}^{2((\vec{v}_{\beta\alpha})_n-(C^{-1})^{n, c_{\beta\alpha}})+\delta_{n,c_{\beta\alpha}} -(\vec{v}_{\beta\alpha})_{n+1}+(C^{-1})^{n+1, c_{\beta\alpha}} -2} \hspace{6mm}\prod_{j=0}^{2i+\delta_{n,c_{\beta\alpha}}-(\vec{v}_{\beta\alpha})_{n+1}+(C^{-1})^{n+1, c_{\beta\alpha}}}} \\[4pt]
  \bigg(a_{\beta\alpha}+\Big(i+\left\lfloor \frac{\delta_{n,c_{\beta\alpha}} -(\vec{v}_{\beta\alpha})_{n+1}+(C^{-1})^{n+1, c_{\beta\alpha}}}{2}\right\rfloor\Big)\, \varepsilon_1^{(n)}-j\, \varepsilon_2^{(n)}\bigg)\ .
\end{multline}
\smallskip
\item[$\bullet$]For $\delta_{n,c_{\beta\alpha}}
  -(\vec{v}_{\beta\alpha})_{n+1}+(C^{-1})^{n+1,
    c_{\beta\alpha}}<2-2((\vec{v}_{\beta\alpha})_n-(C^{-1})^{n , c_{\beta\alpha}})$:
\begin{multline}
\ell^{(n)}_{\vec{v}_{\beta\alpha}}\big(\varepsilon_1^{(n)},
\varepsilon_2^{(n)}, a_{\beta\alpha} \big)=\prod_{i=(\vec{v}_{\beta\alpha})_{n}-(C^{-1})^{n, c_{\beta\alpha}}-\dplus_n}^{(\vec{v}_{\beta\alpha})_n-(C^{-1})^{n, c_{\beta\alpha}}-1}\hspace{0.2cm} \prod_{j=1}^{-2i-\delta_{n,c_{\beta\alpha}} +(\vec{v}_{\beta\alpha})_{n+1}-(C^{-1})^{n+1, c_{\beta\alpha}}-1}\\[4pt]
\bigg(a_{\beta\alpha}+\Big(i-\left\lfloor -
  \frac{\delta_{n,c_{\beta\alpha}}
    -(\vec{v}_{\beta\alpha})_{n+1}+(C^{-1})^{n+1, c_{\beta\alpha}}}{2}
\right\rfloor\Big)\, \varepsilon_1^{(n)}-j\, \varepsilon_2^{(n)}\bigg)^{-1} \ .
\end{multline}
\smallskip
\end{itemize}
\item[\scriptsize$\blacksquare$] If $(\vec{v}_{\beta\alpha})_n-(C^{-1})^{n, c_{\beta\alpha}}= 0$:
\begin{equation}
\qquad\ell^{(n)}_{\vec{v}_{\beta\alpha}}\big(\varepsilon_1^{(n)},
\varepsilon_2^{(n)}, a_{\beta\alpha} \big)=1\ .
\end{equation}
\item[\scriptsize$\blacksquare$] If $(\vec{v}_{\beta\alpha})_n-(C^{-1})^{n, c_{\beta\alpha}}< 0$:
\smallskip
\begin{itemize}
\item[$\bullet$]For $\delta_{n,c_{\beta\alpha}} -(\vec{v}_{\beta\alpha})_{n+1}+(C^{-1})^{n+1, c_{\beta\alpha}}+2(\vec{v}_{\beta\alpha})_n-2(C^{-1})^{n, c_{\beta\alpha}} < 2-2\dminus_n$:
\begin{multline}
\ell^{(n)}_{\vec{v}_{\beta\alpha}}\big(\varepsilon_1^{(n)},
\varepsilon_2^{(n)}, a_{\beta\alpha} \big)=\prod_{i=1-(\vec{v}_{\beta\alpha})_{n}+(C^{-1})^{n, c_{\beta\alpha}}-\dminus_n}^{-(\vec{v}_{\beta\alpha})_{n}+(C^{-1})^{n, c_{\beta\alpha}}}\hspace{6mm} \prod_{j=1}^{2i-(\delta_{n,c_{\beta\alpha}} -(\vec{v}_{\beta\alpha})_{n+1}+(C^{-1})^{n+1, c_{\beta\alpha}})-1}
\\[4pt]
 \bigg(a_{\beta\alpha}-\Big(i+\left\lfloor
   -\frac{\delta_{n,c_{\beta\alpha}}
     -(\vec{v}_{\beta\alpha})_{n+1}+(C^{-1})^{n+1,
       c_{\beta\alpha}}}{2}\right\rfloor\Big)\,
 \varepsilon_1^{(n)}-j\, \varepsilon_2^{(n)}\bigg)\ .
\end{multline}
\item[$\bullet$]For $2-2 \dminus_n\leq \delta_{n,c_{\beta\alpha}} -(\vec{v}_{\beta\alpha})_{n+1}+(C^{-1})^{n+1, c_{\beta\alpha}}+2(\vec{v}_{\beta\alpha})_n-2(C^{-1})^{n, c_{\beta\alpha}}<0$:
\smallskip
\begin{multline}
\ell^{(n)}_{\vec{v}_{\beta\alpha}}\big(\varepsilon_1^{(n)},
\varepsilon_2^{(n)}, a_{\beta\alpha}\big) =
\prod_{i=\big\lfloor \frac{\delta_{n,c_{\beta\alpha}} -(\vec{v}_{\beta\alpha})_{n+1}+(C^{-1})^{n+1, c_{\beta\alpha}}}{2} \big\rfloor+1}^{-(\vec{v}_{\beta\alpha})_{n}+(C^{-1})^{n, c_{\beta\alpha}}}\hspace{6mm} \prod_{j=1}^{2i-(\delta_{n,c_{\beta\alpha}} -(\vec{v}_{\beta\alpha})_{n+1}+(C^{-1})^{n+1, c_{\beta\alpha}})-1 }\\[4pt]
\shoveright{\bigg(a_{\beta\alpha}-\Big(i+\left\lfloor
    -\frac{\delta_{n,c_{\beta\alpha}}
      -(\vec{v}_{\beta\alpha})_{n+1}+(C^{-1})^{n+1,
        c_{\beta\alpha}}}{2}\right\rfloor\Big)\,
  \varepsilon_1^{(n)}-j\, \varepsilon_2^{(n)}\bigg)\times}\\[4pt]
\shoveleft{ \ \prod_{i=1-(\vec{v}_{\beta\alpha})_{n}+(C^{-1})^{n, c_{\beta\alpha}}-\dminus_n}^{\big\lfloor \frac{\delta_{n,c_{\beta\alpha}} -(\vec{v}_{\beta\alpha})_{n+1}+(C^{-1})^{n+1, c_{\beta\alpha}}}{2} \big\rfloor}\hspace{0.2cm} \prod_{j=0}^{-2i+\delta_{n,c_{\beta\alpha}} -(\vec{v}_{\beta\alpha})_{n+1}+(C^{-1})^{n+1, c_{\beta\alpha}}}}\\[4pt]
\bigg(a_{\beta\alpha}+\Big(-i+\left\lfloor
  \frac{\delta_{n,c_{\beta\alpha}}
    -(\vec{v}_{\beta\alpha})_{n+1}+(C^{-1})^{n+1, c_{\beta\alpha}}}{2}
\right\rfloor\Big)\, \varepsilon_1^{(n)}+j\, \varepsilon_2^{(n)}\bigg)^{-1}\ .
\end{multline}
\smallskip
\item[$\bullet$]For $\delta_{n,c_{\beta\alpha}} -(\vec{v}_{\beta\alpha})_{n+1}+(C^{-1})^{n+1, c_{\beta\alpha}}\geq -2(\vec{v}_{\beta\alpha})_n+2(C^{-1})^{n, c_{\beta\alpha}}$:
\begin{multline}
\ell^{(n)}_{\vec{v}_{\beta\alpha}}\big(\varepsilon_1^{(n)},
\varepsilon_2^{(n)}, a_{\beta\alpha} \big) =\prod_{i=1-(\vec{v}_{\beta\alpha})_{n}+(C^{-1})^{n, c_{\beta\alpha}}-\dminus_n}^{-(\vec{v}_{\beta\alpha})_{n}+(C^{-1})^{n, c_{\beta\alpha}}}\hspace{6mm} \prod_{j=0}^{-2i+\delta_{n,c_{\beta\alpha}}-(\vec{v}_{\beta\alpha})_{n+1}+(C^{-1})^{n+1, c_{\beta\alpha}}}\\[4pt]
\bigg(a_{\beta\alpha}+\Big(-i+\left\lfloor
  \frac{\delta_{n,c_{\beta\alpha}}
    -(\vec{v}_{\beta\alpha})_{n+1}+(C^{-1})^{n+1,
      c_{\beta\alpha}}}{2}\right\rfloor\Big)\, \varepsilon_1^{(n)}+j\,
\varepsilon_2^{(n)}\bigg)^{-1}\ .
\end{multline}
\smallskip
\end{itemize}
\end{itemize}
For $n=m_{\beta\alpha}+1, \ldots, k-1$ we set
\begin{equation}
\ell^{(n)}_{\vec{v}_{\beta\alpha}}\big(\varepsilon_1^{(n)},
\varepsilon_2^{(n)}, a_{\beta\alpha} \big)=1\ .
\end{equation}
\begin{remark}
Note that for any fixed $n\in\{1,\dots,k-1\}$,
$d(\beta\alpha)_n^\pm=0$ implies $\ell^{(n)}_{\vec{v}_{\beta\alpha}}
\big(\varepsilon_1^{(n)},\varepsilon_2^{(n)}, a_{\beta\alpha}
\big)=1$.
\end{remark}
Collecting all computations  so far done, the
$\tilde{T}$-equivariant Euler class of the Carlsson-Okounkov bundle
$\Ebf$ at the fixed point ${
  \left([(\Ecal,\phi_\Ecal)]\,,\, [(\Ecal',\phi_{\Ecal'})]\right)}$ is
given by
\begin{empheq}[box=\fbox]{multline}
\eu_{\tilde{T}}\big(\Ebf_{\left([(\Ecal,\phi_\Ecal)]\,,\,
    [(\Ecal',\phi_{\Ecal'})]\right)}\big) =
\prod_{\alpha=1}^r \ \prod_{\beta=1}^{r'} \ \prod_{i=1}^k \, m_{Y_{\alpha}^{i},
  {Y_{\beta}^{i}}'} \big(\varepsilon_1^{(i)},\varepsilon_2^{(i)},
a^{(i)}_{\beta\alpha} \big) \ \prod_{n=1}^{k-1}\,
\ell^{(n)}_{\vec{v}_{\beta\alpha}}
\big(\varepsilon_1^{(n)},\varepsilon_2^{(n)}, a_{\beta\alpha} \big)\ .
\end{empheq}
From a physical viewpoint, this formula expresses the bifundamental weight of an $A_2$ quiver gauge theory, with a $U(r)$ gauge group associated to one node and a $U(r')$ gauge group associated to the other node.

Let $[(\Ecal,\phi_{\Ecal})]$ be a $T$-fixed point of
$\Mcal_{r,\vec{u},\Delta}(\Xscr_k,\Dscr_\infty,\Fcal_\infty^{s,\vec
  w}\, )$ and $(\vec{\boldsymbol{Y}},\vec{\boldsymbol{v}})$ its
corresponding combinatorial data. By Remark
\ref{rem:carlssonokounkovbundle} we can compute the $T$-equivariant
Euler class of the tangent bundle
$T\Mcal_{r,\vec{u},\Delta}(\Xscr_k,\Dscr_\infty,\Fcal_\infty^{s,\vec
  w}\, )$ at $[(\Ecal,\phi_{\Ecal})]$ from this general expression. We get
\begin{equation}
\eu_{T}\big(T_{(\vec{\boldsymbol{Y}},\vec{\boldsymbol{v}})}
\Mcal_{r,\vec{u},\Delta}(\Xscr_k,\Dscr_\infty,\Fcal_\infty^{s,\vec
  w}\, ) \big)=\prod_{\alpha,\beta=1}^r \ \prod_{i=1}^k \,
m_{Y_{\alpha}^{i}, {Y_{\beta}^{i}}}
\big(\varepsilon_1^{(i)},\varepsilon_2^{(i)}, a_{\beta\alpha}^{(i)}
\big) \ \prod_{n=1}^{k-1} \,
\ell^{(n)}_{\vec{v}_{\beta\alpha}}\big(\varepsilon_1^{(n)},\varepsilon_2^{(n)},
a_{\beta\alpha} \big)
\end{equation}
where  $a_{\beta\alpha}=a_\beta-a_\alpha$ and $\vec{v}_{\beta\alpha}=\vec{v}_\beta-\vec{v}_\alpha$.

For $\alpha=1, \ldots, r$ and for a Young tableau $Y_\alpha$ we define
\begin{equation}
m_{Y_\alpha}(\varepsilon_1,\varepsilon_2,a_\alpha):=\prod_{s\in
  Y_\alpha} \big(-L_{Y_\alpha}'(s)\,
\varepsilon_1-A_{Y_\alpha}'(s)\, \varepsilon_2+a_\alpha \big)\ .
\end{equation}
Again, by Remark \ref{rem:carlssonokounkovbundle} we can compute the
$T$-equivariant Euler class of the natural bundle $\Vbf$ at a fixed
point $[(\Ecal,\phi_{\Ecal})]$ from the general character formula. We obtain
\begin{equation}
\eu_{T}\big(\Vbf_{(\vec{\boldsymbol{Y}},\vec{\boldsymbol{v}})}
\big)=\prod_{\alpha=1}^r \ \prod_{i=1}^k \,
m_{Y_{\alpha}^{i}}\big(\varepsilon_1^{(i)},\varepsilon_2^{(i)},
a_\alpha^{(i)} \big) \ \prod_{n=1}^{k-1} \,
\ell^{(n)}_{\vec{v}_{\alpha}}\big(\varepsilon_1^{(n)},\varepsilon_2^{(n)},
a_\alpha \big)\ ,
\end{equation}
where for $i=1, \ldots, k$ we set
\begin{equation}
\vec{a}^{(i)}:=\vec{a}+(\vec{v})_i\,\varepsilon_1^{(i)}
+ (\vec{v})_{i-1}\, \varepsilon_2^{(i)} \ ,
\end{equation}
and $(\vec{v})_i:=\big((\vec{v}_1)_i, \ldots, (\vec{v}_{r})_i \big)$
(we  set $(\vec{v}_\alpha)_0=(\vec{v}_\alpha)_k=0$ for $\alpha=1, \ldots, r$).

\begin{example}
Let $k=2$. Then the Euler class formula becomes
\begin{multline}
\eu_{\tilde{T}}\big(\Ebf_{\left([(\Ecal,\phi_\Ecal)]\,,\,
    [(\Ecal',\phi_{\Ecal'})]\right)}\big) \\ =
\prod_{\alpha=1}^r \ \prod_{\beta=1}^{r'} \, m_{Y_{\alpha}^{1},
  {Y_{\beta}^{1}}'} (2\varepsilon_1,\varepsilon_2-\varepsilon_1,
a_{\beta\alpha}+2v_{\beta\alpha}\,\varepsilon_1 ) \, m_{Y_{\alpha}^{2},
  {Y_{\beta}^{2}}'} (\varepsilon_1-\varepsilon_2,2\varepsilon_2,
a_{\beta\alpha}+2v_{\beta\alpha}\,\varepsilon_2) \\ 
\times \
\ell_{v_{\beta\alpha}}(2\varepsilon_1,\varepsilon_2-\varepsilon_1,
a_{\beta\alpha} \big)\ .
\end{multline}
In this case $c_{\beta\alpha}\in\{0,1\}$ for any $\alpha=1,
\ldots, r$, $\beta=1,
\ldots, r'$, while
$\{v_{\beta\alpha}\}=\frac{\delta_{1,c_{\beta\alpha}}}{2}$ and
$\lfloor v_{\beta\alpha}\rfloor = v_{\beta\alpha}-(C^{-1})^{1,
  c_{\beta\alpha}}$. Since $m=1$, $d(\beta\alpha)_1^+=\lfloor
v_{\beta\alpha}\rfloor$ and $d(\beta\alpha)_1^-=-\lfloor v_{\beta\alpha}\rfloor$,
and we get
\begin{equation*}
\ell_{v_{\beta\alpha}}\big(\varepsilon_1,\varepsilon_2,
a_{\beta\alpha} \big)=\left\{
\begin{array}{cl}
\prod_{i=0}^{\lfloor v_{\beta\alpha}\rfloor-1}\limits\hspace{0.2cm}
\prod_{j=0}^{2i+2\{v_{\beta\alpha}\}}\limits\big(a_{\beta\alpha}+i\,
\varepsilon_1+j\, \varepsilon_2\big) & \mbox{ for } \lfloor v_{\beta\alpha}\rfloor> 0\ ,\\[8pt]
1 & \mbox{ for } \lfloor v_{\beta\alpha}\rfloor= 0\ ,\\[8pt]
\prod_{i=1}^{-\lfloor v_{\beta\alpha}\rfloor}\limits\hspace{0.2cm} \prod_{j=1}^{2i-2\{v_{\beta\alpha}\}-1}\limits
\big(a_{\beta\alpha}+(2\{v_{\beta\alpha}\}-i)\,
\varepsilon_1-j\, \varepsilon_2\big)  & \mbox{ for } \lfloor v_{\beta\alpha}\rfloor< 0 \ .
\end{array}
\right.
\end{equation*}
\normalsize
\end{example}

The explicit formula for $k=3$ is presented in
Appendix~\ref{app:k=3}.

\bigskip\section{BPS correlators and instanton partition functions}\label{sec:instantonpartitionfunctions}

\subsection{$\Ncal=2$ gauge theory}\label{sec:pure}

\subsubsection*{Generating function for correlators of $p$-observables}

For fixed rank $r\geq1$ and
holonomy at infinity $\vec w\in\N^k$, let $\vec{v}\in \frac{1}{k}\, \Z^{k-1}$ be such that $k\, v_{k-1}=
\sum_{i=0}^{k-1}\, i\, w_i \ \bmod{k}$. Define
\begin{multline}\label{eq:deformed-v}
\mathcal{Z}_{\vec{v}}\big(\varepsilon_1, \varepsilon_2, \vec{a}; \qsf,
\vec{\tau}, \vec{t}^{\;(1)}, \ldots, \vec{t}^{\;(k-1)} \big)
:=\sum_{\Delta\in \frac{1}{2r\,k}\, \Z}\, \qsf^{\Delta+\frac{1}{2r}\,
  \vec{v}\cdot C\vec{v}} \\
\times \int_{\Mcal_{r,\vec{u},\Delta}(\Xscr_k,\Dscr_\infty,\Fcal_\infty^{0,\vec
    w}\, )} \,   \exp \bigg(\, \sum_{s=0}^\infty\, \Big(\,
  \sum_{i=1}^{k-1}\, t_s^{(i)} \,
  \big[\ch_{T}(\boldsymbol{\Ecal})/[\Dscr_i]\big]_s+\tau_s \,
  \big[\ch_T(\boldsymbol{\Ecal})/[X_k]\big]_{s-1}\, \Big)\, \bigg)
\end{multline}
where $\boldsymbol{\Ecal}$ is the {universal sheaf}, $\ch_{T}(\boldsymbol{\Ecal})/[\Dscr_i]$ is the \emph{slant product} between $\ch_{T}(\boldsymbol{\Ecal})$ and $[\Dscr_i]$, and the class $\ch_T(\boldsymbol{\Ecal})/[X_k]$ is defined  by localization as
\begin{equation}
\ch_T(\boldsymbol{\Ecal})/[X_k]:=\sum_{i=1}^k \,
\frac{1}{\mathrm{Euler}_{T_t}(T_{p_i} X_k)}\, \imath_{\{p_i\}\times
  \Mcal_{r,\vec{u},\Delta}(\Xscr_k,\Dscr_\infty,\Fcal_\infty^{0,\vec
    w}\, )}^\ast \ch_T(\boldsymbol{\Ecal})\ ;
\end{equation}
here $\imath_{\{p_i\}\times
  \Mcal_{r,\vec{u},\Delta}(\Xscr_k,\Dscr_\infty,\Fcal_\infty^{0,\vec
    w}\, )}$ denotes the inclusion map of $\{p_i\}\times
\Mcal_{r,\vec{u},\Delta}(\Xscr_k,\Dscr_\infty,\Fcal_\infty^{0,\vec
  w}\, )$ in $X_k\times
\Mcal_{r,\vec{u},\Delta}(\Xscr_k,\Dscr_\infty,\Fcal_\infty^{0,\vec
  w}\, )$. The brackets $[-]_s$ indicate the degree $s$ part.
\begin{remark}\label{rem:equivariantcohomology}
Let $\Xscr$ be a topological stack with an action of the Deligne-Mumford torus $T_t$. As explained in  \cite[Section~5]{art:ginotnoohi2012}, there is a well posed notion of $T_t$-equivariant (co)homology theory on $\Xscr$. When $\Xscr$ is a topological space, their definition reduces to Borel's definition of $T_t$-equivariant (co)homology theory on topological spaces. So the slant product is also well-defined for $T_t$-equivariant (co)homology theories on topological stacks.
\end{remark}

\begin{definition}\label{def:deformed-pure}
The \emph{generating function $\Zcal_{X_k}\big(\varepsilon_1, \varepsilon_2, \vec{a}; \qsf,
\vec{\xi}, \vec{\tau}, \vec{t}^{\;(1)}, \ldots, \vec{t}^{\:(k-1)}
\big)$ for correlators of $p$-observables} ($p=0,2$)
of pure $\Ncal=2$ gauge theory on $X_k$ is
\begin{multline}\label{eq:N2deformedpure-1}
\Zcal_{X_k}\big(\varepsilon_1, \varepsilon_2, \vec{a}; \qsf,
\vec{\xi}, \vec{\tau}, \vec{t}^{\;(1)}, \ldots, \vec{t}^{\:(k-1)}
\big) \\
:=\sum_{\stackrel{\scriptstyle \vec{v}\in\frac{1}{k}\, \Z^{k-1}
  }{\scriptstyle k\, v_{k-1}=  \sum_{i=0}^{k-1}\, i\, w_i\bmod{k} }} \,
  \vec{\xi}^{\ \vec{v}}\ \Zcal_{\vec{v}}\big(\varepsilon_1,
  \varepsilon_2, \vec{a}; \qsf, \vec{\tau}, \vec{t}^{\;(1)}, \ldots,
  \vec{t}^{\;(k-1)} \big)\ ,
\end{multline}
where $\vec\xi=(\xi_1,\dots,\xi_{k-1})$ and $\vec\xi\ ^{\vec
  v}:=\prod_{i=1}^{k-1}\, \xi_i^{v_i}$.
\end{definition}

In this subsection we compute explicitly $\Zcal_{X_k}(\varepsilon_1, \varepsilon_2, \vec{a}; \qsf, \vec{\xi},
\vec{\tau}, \vec{t}^{\;(1)}, \ldots, \vec{t}^{\:(k-1)})$. First note that by the localization formula we have
\begin{multline}
\mathcal{Z}_{\vec{v}}\big(\varepsilon_1, \varepsilon_2, \vec{a}; \qsf,
\vec{\tau}, \vec{t}^{\;(1)}, \ldots, \vec{t}^{\;(k-1)} \big) \\
\shoveleft{= \sum_{(\vec{\boldsymbol{Y}},\vec{\boldsymbol{v}})}\,
  \frac{\qsf^{\sum_{\alpha=1}^r\limits\, n_\alpha+\frac{1}{2}\, \sum_{\alpha=1}^r\limits\ \vec{v}_\alpha\cdot C\vec{v}_\alpha}}{\eulerclass}} \\
\times \ \imath_{(\vec{\boldsymbol{Y}},\vec{\boldsymbol{v}})}^\ast\exp
\bigg(\, \sum_{s=0}^\infty\, \Big(\, \sum_{i=1}^{k-1}\, t_s^{(i)} \,
\big[\ch_{T}(\boldsymbol{\Ecal})/[\Dscr_i]\big]_s+\tau_s\,
\big[\ch_T(\boldsymbol{\Ecal})/[X_k]\big]_{s-1}\, \Big)\, \bigg)\ .
\end{multline}

Now we compute $\imath_{(\vec{\boldsymbol{Y}},\vec{\boldsymbol{v}})}^\ast\ch_T(\boldsymbol{\Ecal})/[X_k]$. First note that
\begin{equation}
\imath_{(\vec{\boldsymbol{Y}},\vec{\boldsymbol{v}})}^\ast\ch_T(\boldsymbol{\Ecal})/[X_k]=\sum_{i=1}^k
\, \frac{1}{\varepsilon_1^{(i)}\, \varepsilon_2^{(i)}}\, \imath_{\{p_i\}\times \{(\vec{\boldsymbol{Y}},\vec{\boldsymbol{v}})\}}^\ast \ch_T(\boldsymbol{\Ecal})\ .
\end{equation}
Let us fix an index $i\in \{1, \ldots, k\}$. Let
$\big[(\Ecal,\phi_{\Ecal}) \big]=\big[\bigoplus_{\alpha=1}^r\, (\imath_\ast(I_\alpha)\otimes \Rcal^{C\vec{v}_\alpha},\phi_\alpha)\big]$ be a $T$-fixed point and $(\vec{\boldsymbol{Y}},\vec{\boldsymbol{v}})$ its corresponding combinatorial data. Then
\begin{equation}
\imath_{\{p_i\}\times
  \{(\vec{\boldsymbol{Y}},\vec{\boldsymbol{v}})\}}^\ast
\ch_T(\boldsymbol{\Ecal})=\sum_{\alpha=1}^r\, e_\alpha \, \imath_{p_i}^\ast
\ch_{T_t}\big(\imath_\ast(I_\alpha)\otimes \Rcal^{C\vec{v}_\alpha}
\big)=\sum_{\alpha=1}^r\, e_\alpha\, \ch_{T_t}\big((I_\alpha)_{p_i}
\big) \, \ch_{T_t}\big(\Rcal^{C\vec{v}_\alpha}_{p_i} \big)\ ,
\end{equation}
where $\imath_{p_i}$ denotes the inclusion morphism of the $T_t$-fixed
point $p_i$ into $X_k$.
By \cite[Equation~(4.1)]{book:nakajimayoshioka2004} we have
\begin{equation*}\label{eq:formula4.1}
\operatorname{ch}_{T_t}\big((I_\alpha)_{p_i}
\big)=1-\big(1-\chi_1^i\big)\, \big(1-\chi_2^i\big)\,
\sum_{s\in Y_\alpha^i}\, (\chi_1^i)^{L'(s)}\, (\chi_2^i)^{A'(s)}\ .
\end{equation*}
By Lemma \ref{lem:characterline} and Equation
\eqref{eq:tautologicalclasses}, we have
$\ch_{T_t}(\Rcal^{C\vec{v}_\alpha}_{p_i})=(\chi_1^i)^{(\vec{v}_\alpha)_i}\,
(\chi_2^i)^{(\vec{v}_\alpha)_{i-1}}$. Summing up, we get
\begin{equation}\label{eq:character-surface}
\imath_{(\vec{\boldsymbol{Y}},\vec{\boldsymbol{v}})}^\ast\ch_T(\boldsymbol{\Ecal})/[X_k]=
\sum_{i=1}^k\, \ch_{\vec{Y}^i}\big(\varepsilon_1^{(i)},
\varepsilon_2^{(i)}, \vec{a}^{(i)} \big) \ ,
\end{equation}
where we introduced the notation
\begin{equation}
\operatorname{ch}_{\vec{Y}^i}\big(\varepsilon_1^{(i)}, \varepsilon_2^{(i)},
\vec{a}^{(i)}\big)
:= \sum_{\alpha=1}^r\, \frac{\mathrm{e}^{a^{(i)}_\alpha
}}{\varepsilon_1^{(i)}\, \varepsilon_2^{(i)}}
\, \big(1-\big(1-{\,\mathrm{e}\,}^{\varepsilon_1^{(i)}}\big)\,
\big(1-{\,\mathrm{e}\,}^{\varepsilon_2^{(i)}}\big)\, \sum_{s\in
Y_\alpha^i}\,
{\,\mathrm{e}\,}^{\varepsilon_1^{(i)}\, L'(s)+\varepsilon_2^{(i)}\,
A'(s)}\big)\ .
\end{equation}

\begin{remark}\label{rem:classical}
Note that
\begin{eqnarray*}
\Big[\operatorname{ch}_{\vec{Y}^i}\big(\varepsilon_1^{(i)},
\varepsilon_2^{(i)},
\vec{a}^{(i)}\big)\Big]_s &=&
\frac{1}{\varepsilon_1^{(i)}\, \varepsilon_2^{(i)}\,
(s+2)!}\, \sum_{\alpha=1}^r \, \big(a_\alpha^{(i)}\big)^{s+2} \\
&&{} - \frac{1}{\varepsilon_1^{(i)}\, \varepsilon_2^{(i)}\, (s+2)!}
\,
\sum_{\alpha=1}^r \ \sum_{s\in Y_\alpha^i} \,
\Big(\big(a_\alpha^{(i)} +\varepsilon_1^{(i)}\,
L'(s)+\varepsilon_2^{(i)}\, A'(s) \big)^{s+2} \\
&&{}- \big(a_\alpha^{(i)} +\varepsilon_1^{(i)}\,
(L'(s)+1)+\varepsilon_2^{(i)}\, A'(s)\big)^{s+2} \\
&&{} - \big(a_\alpha^{(i)} +\varepsilon_1^{(i)}\,
L'(s)+\varepsilon_2^{(i)}\, (A'(s)+1)\big)^{s+2} \\
&&{}+ \big(a_\alpha^{(i)} +\varepsilon_1^{(i)}\,
(L'(s)+1)+\varepsilon_2^{(i)}\, (A'(s)+1)\big)^{s+2}\Big) \ .
\end{eqnarray*}
In particular
%
\begin{align}
0=\Big[\operatorname{ch}_{\vec{Y}^i}\big(\varepsilon_1^{(i)},
\varepsilon_2^{(i)},
\vec{a}^{(i)}\big)\Big]_{-1} & =\frac{1}{\varepsilon_1^{(i)}\,
\varepsilon_2^{(i)}}\, \sum_{\alpha=1}^r \, a_\alpha^{(i)} \
,\nonumber\\[4pt]
\label{eq:degreezero}
\Big[\operatorname{ch}_{\vec{Y}^i}\big(\varepsilon_1^{(i)},
\varepsilon_2^{(i)},
\vec{a}^{(i)}\big)\Big]_{0} & =\frac{1}{2\varepsilon_1^{(i)}\,
\varepsilon_2^{(i)}}\, \sum_{\alpha=1}^r \, \big(a_\alpha
^{(i)}\big)^2-\sum_{\alpha=1}^r\, \vert Y_\alpha^i\vert\ .
\end{align}
These formulas will be useful later on.
\end{remark}

Now we compute
$\imath_{(\vec{\boldsymbol{Y}},\vec{\boldsymbol{v}})}^\ast\ch_{T}(\boldsymbol{\Ecal})/[\Dscr_i]$. With
the same conventions as above, we have
\begin{align*}
\operatorname{ch}_T({\mathcal E})=&\sum_{\alpha=1}^r \, e_\alpha\,
\operatorname{ch}_{T_t}\big({\mathcal R}^{C\vec{v}_\alpha}\big)\,
\operatorname{ch}_{T_t}\big(\imath_\ast
I_\alpha\big) \\[4pt] =&
\sum_{\alpha=1}^r \, e_\alpha\, \exp\Big(-\sum\limits
_{j=1}^{k-1} \,
(\vec{v}_\alpha)_j \, [{\mathscr{D}}_j]\Big) \\
&\qquad{}\cdot
\Big(1-\sum_{l=1}^k \, [p_l]\, \big(1-\chi_1^l\big)\,
\big(1-\chi_2^l \big)\ \sum_{s\in Y_\alpha^l}\,
(\chi_1^l)^{L'(s)}\, (\chi_2^l)^{A'(s)}\Big)\ .
\end{align*}
In the following we compute separately the two types of contributions
$\exp\big(-\sum_{j=1}^{k-1}\, (\vec{v}_\alpha)_j \, [\Dscr_j]
\big)/[\Dscr_i]$ and $\exp\big( -\sum_{j=1}^{k-1}\,
(\vec{v}_\alpha)_j\, [\Dscr_j] \big)\, [p_l]/[\Dscr_i]$. For the first quantity we get the expression
\begin{align}
& \exp\Big(-\sum_{j=1}^{k-1}\limits \,
(\vec{v}_\alpha)_j \, [\Dscr_j]\Big) \, \big/\, [\Dscr_i] \ = \ \sum_{m=0}^\infty\, \frac{1}{m!} \, (-1)^m\, \Big(\,
\sum_{j=1}^{k-1}\, (\vec{v}_\alpha)_j \, [\Dscr_j]\, \Big)^m \,
\big/\, [\Dscr_i] \\[4pt] & \qquad\qquad \qquad\qquad\ = \ \sum_{m=0}^\infty \, \frac{1}{m!} \,
(-1)^m\, \int_{X_k}\, \Big(\, \sum_{j=1}^{k-1}\, (\vec{v}_\alpha)_j \,
[D_j]\, \Big)^m \cdot [D_i] \\[4pt]
& \qquad\qquad \qquad\qquad\ = \ \sum_{m=0}^\infty \, \frac{1}{m!} \, (-1)^m \
\sum_{m_1+\cdots+m_{k-1}=m} \, \frac{m!}{m_1!\cdots m_{k-1}!} \,
(\vec{v}_\alpha)_1^{m_1} \cdots (\vec{v}_\alpha)_{k-1}^{m_{k-1}} \\ &
\qquad\qquad \qquad\qquad \qquad\qquad \qquad\qquad \qquad\qquad \times \ \int_{X_k}\, [D_1]^{m_1}\cdots
  [D_i]^{m_i+1}\cdots [D_{k-1}]^{m_{k-1}}\ .
\end{align}
To compute the integral
\begin{equation}\label{eq:integral}
\int_{X_k}\, [D_1]^{m_1}\cdots [D_i]^{m_i+1}\cdots [D_{k-1}]^{m_{k-1}}
\end{equation}
by localization, we need to know the pullback of the class $[D_l]$ at
the fixed points of $X_k$ for each $l=1, \ldots, k-1$. Because of Lemma~\ref{lem:characterline} 
\begin{equation}\label{eq:localizationdivisors}
\imath_{p_i}^\ast[D_l]=\left\{
\begin{array}{l@{\quad}l}
-\varepsilon_1^{(l)} & \mbox{if } i=l \ ,\\ \noalign{\vspace{3pt}}
-\varepsilon_2^{(l+1)} & \mbox{if } i=l+1\ ,\\ \noalign{\vspace{3pt}}
0 & \mbox{otherwise}\ ,
\end{array}
\right.
\end{equation}
the integral \eqref{eq:integral} is nonzero if there exists an index $n\in\{1,\ldots, k-1\}$ such that only the exponent of $[D_n]$ is nonzero or there exists an index $n'\in\{2,\ldots, k-1\}$ such that only the exponents of $[D_{n'-1}]$ and $[D_{n'}]$ are nonzero. Therefore we obtain
\begin{multline*}
\exp\Big(-\sum\limits_{j=1}^{k-1} \,
(\vec{v}_\alpha)_j \, [{\mathscr{D}}_j]\Big) \, \big/\, [{\mathscr
{D}}_i] \\
\shoveleft{=\sum_{m=0}^\infty\, \frac{1}{m!}\, \Big(\,
\sum_{n=1}^{k-1}\, \frac{\left(\vec{v}_\alpha\right)_n^m\,
\big(\varepsilon_1^{(n)}\big)^m\,
\big(-\varepsilon_1^{(n)}\big)^{\delta_{n,i}}}{\varepsilon_1^{(n)}\,
\varepsilon_2^{(n)}}+\sum_{n=1}^{k-1}\,
\frac{\left(\vec{v}_\alpha\right)_n^m\,
\big(\varepsilon_2^{(n+1)}\big)^m\,
\big(-\varepsilon_2^{(n+1)}\big)^{\delta_{n,i}}}{\varepsilon
_1^{(n+1)}\,
\varepsilon_2^{(n+1)}}}\\
\shoveright{+ \sum_{n=2}^{k-1} \ \sum_{l=0}^m \, {m \choose l} \,
\frac{\left(\vec{v}_\alpha\right)_n^l\,
\left(\vec{v}_\alpha\right)_{n-1}^{m-l}\,
\big(\varepsilon_1^{(n)}\big)^l\,
\big(-\varepsilon_1^{(n)}\big)^{\delta_{n,i}}\,
\big(\varepsilon_2^{(n)}\big)^{m-l}\,
\big(-\varepsilon_2^{(n)}\big)^{\delta_{n-1,i}}}{\varepsilon
_1^{(n)}\,
\varepsilon_2^{(n)}}\, \Big)}\\[4pt]
\shoveleft{=\sum_{n=1}^{k-1}\,
\frac{\big(-\varepsilon_1^{(n)}\big)^{\delta_{n,i}}}{\varepsilon
_1^{(n)}\,
\varepsilon_2^{(n)}}\, {\,\mathrm{e}\,}^{\left(\vec{v}_\alpha
\right)_n\,
\varepsilon_1^{(n)}}+\sum_{n=2}^{k}\,
\frac{\big(-\varepsilon_2^{(n)}\big)^{\delta_{n-1,i}}}{\varepsilon
_1^{(n)}\,
\varepsilon_2^{(n)}}\, {\,\mathrm{e}\,}^{\left(\vec{v}_\alpha
\right)_{n-1}\, \varepsilon_2^{(n)}}}\\
+ \sum_{n=2}^{k-1}\,
\frac{\big(-\varepsilon_1^{(n)}\big)^{\delta_{n,i}}\,
\big(-\varepsilon_2^{(n)}\big)^{\delta_{n-1,i}}}{\varepsilon
_1^{(n)}\,
\varepsilon_2^{(n)}}\, {\,\mathrm{e}\,}^{\left(\vec{v}_\alpha
\right)_n\,
\varepsilon_1^{(n)}+\left(\vec{v}_\alpha\right)_{n-1}\,
\varepsilon_2^{(n)}}\ .
\end{multline*}
On the other hand, a similar argument shows that for $l\in\{2, \ldots
, k-2\}$ we get
\begin{multline*}
\exp\Big(-\sum\limits_{j=1}^{k-1} \,
(\vec{v}_\alpha)_j \, [{\mathscr{D}}_j]\Big) \cdot[p_l]\, \big/\,
[{\mathscr{D}}_i]=
\sum_{m=0}^\infty\, \frac{1}{m!}\, (-1)^m\, \Big(\, \sum
_{j=1}^{k-1}\,
(\vec{v}_\alpha)_j \, [{\mathscr{D}}_j]\, \Big)^m \cdot[p_\ell]
\, \big/\,
[{\mathscr{D}}_i] \\[4pt]
=\frac{\big(-\varepsilon_1^{(l)}\big)^{\delta_{l,i}}}{\varepsilon
_1^{(l)}\,
\varepsilon_2^{(l)}}\, {\,\mathrm{e}\,}^{\left(\vec{v}_\alpha
\right)_l\,
\varepsilon_1^{(l)}}+\frac{\big(-\varepsilon_2^{(l)}\big)^{\delta
_{l-1,i}}}{\varepsilon_1^{(l)}\,
\varepsilon_2^{(l)}}\, {\,\mathrm{e}\,}^{\left(\vec{v}_\alpha
\right)_{l-1}\, \varepsilon_2^{(l)}}
+\frac{\big(-\varepsilon_1^{(l)}\big)^{\delta_{l,i}}\,
\big(-\varepsilon_2^{(l)}\big)^{\delta_{l-1,i}}}{\varepsilon
_1^{(l)}\,
\varepsilon_2^{(l)}}\, {\,\mathrm{e}\,}^{\left(\vec{v}_\alpha
\right)_l\,
\varepsilon_1^{(l)}+\left(\vec{v}_\alpha\right)_{l-1}\,
\varepsilon_2^{(l)}}\ ,
\end{multline*}
and
\begin{align*}
\exp\Big(-\sum\limits_{j=1}^{k-1} \,
(\vec{v}_\alpha)_j \, [{\mathscr{D}}_j]\Big) \cdot[p_1]\, \big/\,
[{\mathscr{D}}_i]&=
\frac{\big(-\varepsilon_1^{(1)}\big)^{\delta_{1,i}}}{\varepsilon
_1^{(1)}\,
\varepsilon_2^{(1)}}\, {\,\mathrm{e}\,}^{\left(\vec{v}_\alpha
\right)_1\, \varepsilon_1^{(1)}}\ ,\\[4pt]
\exp\Big(-\sum\limits_{j=1}^{k-1} \,
(\vec{v}_\alpha)_j \, [{\mathscr{D}}_j]\Big) \cdot[p_k]\, \big/\,
[{\mathscr{D}}_i]&=
\frac{\big(-\varepsilon_2^{(k)}\big)^{\delta_{k-1,i}}}{\varepsilon
_1^{(k)}\,
\varepsilon_2^{(k)}}\, {\,\mathrm{e}\,}^{\left(\vec{v}_\alpha
\right)_{k-1}\, \varepsilon_2^{(k)}}\ .
\end{align*}
Summing up, we obtain
%
\begin{eqnarray}
\imath_{(\vec{\boldsymbol{Y}},\vec{\boldsymbol{v}})}^\ast
\operatorname{ch}_{T}(\boldsymbol{{\mathcal E}})/[{\mathscr{D}}_i]
&=& \sum_{l=1}^k\,
\big(-\varepsilon_1^{(l)}\big)^{\delta_{l,i}}\,
\big(-\varepsilon_2^{(l)}\big)^{\delta_{l-1,i}}\,
\operatorname{ch}_{\vec{Y}^l}\big(\varepsilon_1^{(l)},\varepsilon_2^{(l)},
\vec{a}^{(l)} \big) \nonumber\\
&&{} + \sum_{l=2}^{k-1}\,
\Big(\big(-\varepsilon_1^{(l)}\big)^{\delta_{l,i}}\,
\operatorname{ch}_{\vec{Y}^l}\big(\varepsilon_1^{(l)},\varepsilon_2^{(l)},
\vec{a}^{(l)}+\left(\vec{v}\right)_{l-1}\, \varepsilon
_2^{(l)})\nonumber\\
&&{} + \big(-\varepsilon_2^{(l)}\big)^{\delta_{l-1,i}}\,
\operatorname{ch}_{\vec{Y}^l}\big(\varepsilon_1^{(l)},\varepsilon_2^{(l)},
\vec{a}^{(l)}+\left(\vec{v}\right)_{l}\, \varepsilon_1^{(l)}\big
)\Big)
\ .\label{eq:character-divisor}
\end{eqnarray}

By using Equations \eqref{eq:character-surface} and \eqref{eq:character-divisor} we arrive finally at the following result.
\begin{proposition}\label{prop:deformed-pure}
The generating function ${\mathcal Z}_{X_k}(\varepsilon_1, \varepsilon_2,
\vec{a}; \mathsf{q}, \vec{\xi}, \vec{\tau}, \vec{t}^{\;(1)},
\ldots,
\vec{t}^{\;(k-1)})$ assumes the form
\begin{multline}\label{eq:N2deformedpure-2}
{\mathcal Z}_{X_k}\big(\varepsilon_1, \varepsilon_2, \vec{a};
\mathsf{q},
\vec{\xi}, \vec{\tau}, \vec{t}^{\;(1)}, \ldots, \vec{t}^{\;(k-1)}
\big) \\
\shoveleft{=\sum_{\stackrel{\scriptstyle\vec{v}\in\frac{1}{k}\,
{\mathbb{Z}}^{k-1} }{\scriptstyle k\, v_{k-1}= \sum_{i=0}^{k-1}\, i\,
w_i\bmod{k} }}\, \vec{\xi}^{\ \vec{v}} \
\sum_{(\vec{\boldsymbol{Y}},\vec{\boldsymbol{v}})} \,
\begin{aligned}
 \frac{\mathsf{q}^{\sum\limits_{\alpha=1}^r\,
n_\alpha+\frac{1}{2}\, \sum\limits_{\alpha=1}^r\ \vec{v}_\alpha
\cdot C\vec{v}_\alpha}}{\prod\limits_{\alpha,\beta=1}^r\ \prod
\limits_{i=1}^k
\, m_{Y_{\alpha}^{i},
{Y_{\beta}^{i}}}\big(\varepsilon_1^{(i)},\varepsilon_2^{(i)},
a_{\beta\alpha}^{(i)} \big) \ \prod\limits_{n=1}^{k-1} \,
\ell^{(n)}_{\vec{v}_{\beta\alpha}}\big(\varepsilon
_1^{(n)},\varepsilon_2^{(n)},
a_{\beta\alpha} \big)}
\end{aligned}
}\\
\shoveright{ \times \prod_{l=1}^k\, \exp\bigg(\,
\sum_{s=0}^\infty\, \Big(\big(t_s^{(l)}\,
\varepsilon_1^{(l)}+t_s^{(l-1)}\, \varepsilon_2^{(l)}+\tau_s\big)\,
\Big[\operatorname{ch}_{\vec{Y}^l}(\varepsilon_1^{(l)},\varepsilon_2^{(l)},
\vec{a}^{(l)})\Big]_{s-1}\, \Big)\bigg)}\\
\shoveleft{\times\exp\Bigg(\, \sum_{s=0}^\infty\, \bigg(\,
\sum_{l=1}^k\, \Big(\, \sum_{i=1}^{l-2}\,
t_s^{(i)}+\sum_{i=l+1}^{k-1}\, t_s^{(i)}\Big)\, \Big[\operatorname
{ch}_{\vec{Y}^l}\big(\varepsilon_1^{(l)},\varepsilon_2^{(l)}, \vec
{a}^{(l)})\Big]_{s}}\\
+ \sum_{i=1}^{k-1}\, t_s^{(i)} \ \sum_{l=2}^{k-1}\,
\Big(\, \Big[\big(-\varepsilon_1^{(l)}\big)^{\delta_{l,i}}\,
\operatorname{ch}_{\vec{Y}^l}\big(\varepsilon_1^{(l)},\varepsilon_2^{(l)},
\vec{a}^{(l)}+\left(\vec{v}\right)_{l-1}\, \varepsilon_2^{(l)}\big
)\Big]_s \\
+ \Big[\big(-\varepsilon_2^{(l)}\big)^{\delta_{l-1,i}}\,
\operatorname{ch}_{\vec{Y}^l}\big(\varepsilon_1^{(l)},\varepsilon_2^{(l)},
\vec{a}^{(l)}+\left(\vec{v}\right)_{l}\, \varepsilon_1^{(l)}\big)
\Big]_{s}\, \Big)\bigg) \Bigg)
\end{multline}
where we set $t_s^{(0)}=t_s^{(k)}=0$ for any $s$.
\end{proposition}

\begin{example}\label{ex:k2case}
For $k=2$ the generating function ${\mathcal Z}_{X_2}(\varepsilon_1,
\varepsilon_2, \vec{a}; \mathsf{q}, \xi, \vec{\tau}, \vec{t}\ )$ becomes
\begin{multline}
{\mathcal Z}_{X_2}\big(\varepsilon_1, \varepsilon_2, \vec{a};
\mathsf{q}, \xi,
\vec{\tau}, \vec{t} \ \big) \\
\shoveleft{=\sum_{\stackrel{\scriptstyle v\in\frac{1}{2}\, {\mathbb{Z}}
}{\scriptstyle2v= w_1\bmod{2}}}\, \xi^{v} \
\sum_{\vec{\boldsymbol{v}}=(v_1, \ldots, v_r)}\
\frac{\mathsf{q}^{\sum\limits_{\alpha=1}^r\,
v_\alpha^2}}{\prod\limits_{\alpha,\beta=1}^r \,
\ell_{v_{\beta\alpha}}(2\varepsilon_1, \varepsilon_2
-\varepsilon_1, a_{\beta\alpha})}} \\
\quad{}\times{\mathcal Z}_{{\mathbb{C}}^2}\big(2\varepsilon
_1,\varepsilon_2-\varepsilon_1,
\vec{a}+2\varepsilon_1\, \vec{\boldsymbol{v}};
\mathsf{q},\vec{\tau}+2\varepsilon_1 \, \vec{t} \ \big)\,
{\mathcal Z}_{{\mathbb{C}}^2}\big(\varepsilon_1- \varepsilon
_2,2\varepsilon_2,
\vec{a}+2\varepsilon_2\, \vec{\boldsymbol{v}};
\mathsf{q},\vec{\tau}+2\varepsilon_2\, \vec{t} \ \big) ,
\end{multline}
where $\Zcal_{\C^2}(\varepsilon_1,\varepsilon_2, \vec{a}; \qsf,
\vec\tau\,)$ is the \emph{deformed Nekrasov partition function} (the
generating function for $0$-observables) for $U(r)$ gauge theory on $\R^4$ defined in \cite[Section~4.2]{book:nakajimayoshioka2004}:
\begin{equation}
\Zcal_{\C^2}\big(\varepsilon_1,\varepsilon_2, \vec{a}; \qsf, \vec\tau
\, \big):=\sum_{\vec{Y}} \, \frac{\qsf^{\sum_{\alpha=1}^r\limits \, \vert
    Y_\alpha\vert}}{\prod_{\alpha,\beta=1}^r \limits \, m_{Y_\alpha,
    Y_\beta}(\varepsilon_1,\varepsilon_2,a_{\beta\alpha})} \, 
\exp\Big(\, \sum_{s=0}^\infty \, \tau_s\,
\big[\ch_{\vec{Y}}(\varepsilon_1,\varepsilon_2,
\vec{a}\, )\big]_{s-1}\, \Big) \ .
\end{equation}
\end{example}
However, for $k\geq 3$ we see no such nice factorizations of Equation \eqref{eq:N2deformedpure-2}.

\subsubsection*{Instanton partition function}

\begin{definition}
The \emph{instanton partition function} is
\begin{equation}\label{eq:instantonpart-puregaugetheory}
\Zcal^{\mathrm{inst}}_{X_k}\big(\varepsilon_1, \varepsilon_2, \vec{a};
\qsf, \vec{\xi}\ \big):=\Zcal_{X_k}\big(\varepsilon_1, \varepsilon_2,
\vec{a}; \qsf, \vec{\xi},  0, \ldots, 0 \big)\ .
\end{equation}
Here we set $\vec{\tau}=\vec 0$ and $ \vec{t}^{\;(i)}=\vec 0$ for $i=1,\ldots, k-1$.
\end{definition}
By Equation \eqref{eq:N2deformedpure-2} we get
\begin{multline}
\Zcal^{\mathrm{inst}}_{X_k}\big(\varepsilon_1, \varepsilon_2, \vec{a};
\qsf, \vec{\xi} \ \big)\\
= \sum_{\stackrel{\scriptstyle \vec{v}\in\frac{1}{k}\,
    \Z^{k-1}}{\scriptstyle k\, v_{k-1}=  \sum_{i=0}^{k-1}\, i\,
    w_i\bmod{k}} }\, \vec{\xi}^{\ \vec{v}} \
\sum_{(\vec{\boldsymbol{Y}},\vec{\boldsymbol{v}})} \,
\frac{\qsf^{\sum_{\alpha=1}^r\limits\, n_\alpha+\frac{1}{2}\, \sum_{\alpha=1}^r\limits\ \vec{v}_\alpha\cdot C\vec{v}_\alpha}}{\eulerclass}\ .
\end{multline}
Since the \emph{Nekrasov partition function} for pure $U(r)$ gauge
theory on $\R^4$ can be expressed as \cite[Equation~(3.16)]{art:bruzzofucitomoralestanzini2003}
\begin{equation}
\Zcal_{\C^2}^{\mathrm{inst}}(\varepsilon_1, \varepsilon_2, \vec{a};
\qsf):=\sum_{\vec{Y}}\, \frac{\qsf^{\sum_{\alpha=1}^r\limits \, \vert
    Y_\alpha\vert}}{\prod_{\alpha,\beta=1}^r \limits \, m_{Y_\alpha, Y_\beta}(\varepsilon_1,\varepsilon_2,a_{\beta\alpha})}\ ,
\end{equation}
we get
\begin{empheq}[box=\fbox]{multline}
\Zcal^{\mathrm{inst}}_{X_k}\big(\varepsilon_1, \varepsilon_2, \vec{a};
\qsf, \vec{\xi}\ \big)\\
= \sum_{\stackrel{\scriptstyle \vec{v}\in\frac{1}{k}\Z^{k-1}
  }{\scriptstyle k\, v_{k-1}=  \sum_{i=0}^{k-1}\, i\, w_i\bmod{k} }}\,
\vec{\xi}^{\ \vec{v}} \ \sum_{\vec{\boldsymbol{v}}}
\frac{\qsf^{\frac{1}{2}\, \sum_{\alpha=1}^r\limits \,
    \vec{v}_\alpha\cdot C\vec{v}_\alpha}}{\prod_{\alpha,\beta=1}^r
  \limits\ \prod_{n=1}^{k-1}\limits \,
  \ell^{(n)}_{\vec{v}_{\beta\alpha}}\big(\varepsilon_1^{(n)},
  \varepsilon_2^{(n)}, a_{\beta\alpha} \big)} \ \prod_{i=1}^k\,
\Zcal_{\C^2}^{\mathrm{inst}}\big(\varepsilon_1^{(i)},
\varepsilon_2^{(i)}, \vec{a}^{(i)}; \qsf \big)\ .
\end{empheq}

\subsubsection*{Correlators of quadratic $0$-observables}

Let us define
\begin{equation}
\Zcal^{\circ}_{X_k}\big(\varepsilon_1, \varepsilon_2,
\vec{a}; \qsf, \vec{\xi},\tau_1 \big):=\Zcal_{X_k}\big(\varepsilon_1, \varepsilon_2, \vec{a}; \qsf, \vec{\xi}, \,\underline{\vec{\tau}}\,, 0, \ldots, 0\big)\ ,
\end{equation}
where we set $\underline{\vec{\tau}}=(0,-\tau_1,0,\ldots)$ and $
\vec{t}^{\;(i)}=\vec 0$ for $i=1,\ldots, k-1$. From Equations
\eqref{eq:N2deformedpure-2} and \eqref{eq:degreezero} it follows that
\begin{multline}
\Zcal^{\circ}_{X_k}\big(\varepsilon_1, \varepsilon_2, \vec{a};
\qsf,\vec{\xi}, \tau_1 \big) \\
\shoveleft{=\sum_{\stackrel{\scriptstyle \vec{v}\in\frac{1}{k}\, \Z^{k-1}
    }{\scriptstyle k\, v_{k-1}=  \sum_{i=0}^{k-1}\, i\, w_i\bmod{k}
    }} \, \vec{\xi}^{\ \vec{v}} \
    \sum_{(\vec{\boldsymbol{Y}},\vec{\boldsymbol{v}})}\ 
    \frac{\qsf^{\sum_{\alpha=1}^r\limits \, n_\alpha+\frac{1}{2}\,
        \sum_{\alpha=1}^r\limits \, \vec{v}_\alpha\cdot C\vec{v}_\alpha}}{\eulerclass}}\\
\shoveright{\times \ \prod_{l=1}^k\, \exp\Big(-\tau_1\,
  \Big[\ch_{\vec{Y}^l}\big(\varepsilon_1^{(l)},\varepsilon_2^{(l)},
  \vec{a}^{(l)} \big)\Big]_{0}\, \Big)}\\[4pt]
\shoveleft{=\sum_{\stackrel{\scriptstyle \vec{v}\in\frac{1}{k}\, \Z^{k-1}
    }{\scriptstyle k\, v_{k-1}=  \sum_{i=0}^{k-1}\, i\, w_i\bmod{k}
    }} \, \vec{\xi}^{\ \vec{v}} \
    \sum_{(\vec{\boldsymbol{Y}},\vec{\boldsymbol{v}})}\
    \frac{\qsf^{\sum_{\alpha=1}^r\limits \, n_\alpha+\frac{1}{2}\,
        \sum_{\alpha=1}^r\limits \, \vec{v}_\alpha\cdot C\vec{v}_\alpha}}{\eulerclass}}\\
\shoveright{\times \ \prod_{l=1}^k\, \exp\bigg(-\tau_1\, \Big(\,
\frac{1}{2\varepsilon_1^{(l)}\, \varepsilon_2^{(l)}}\,
\sum_{\alpha=1}^r \, \big(a_\alpha^{(l)}\big)^2-\sum_{\alpha=1}^r \,
\vert Y_\alpha^l\vert\, \Big)\bigg)} \\[4pt]
\shoveleft{=\sum_{\stackrel{\scriptstyle \vec{v}\in\frac{1}{k}\, \Z^{k-1}
    }{\scriptstyle k\, v_{k-1}=  \sum_{i=0}^{k-1}\, i\, w_i\bmod{k}
    }} \, \vec{\xi}^{\ \vec{v}} \ \sum_{\vec{\boldsymbol{v}}}\,
  \frac{\qsf^{\frac{1}{2}\, \sum_{\alpha=1}^r\limits\,
      \vec{v}_\alpha\cdot
      C\vec{v}_\alpha}}{\prod_{\alpha,\beta=1}^r\limits \
    \prod_{n=1}^{k-1}\limits \,
    \ell^{(n)}_{\vec{v}_{\beta\alpha}}\big(\varepsilon_1^{(n)},
    \varepsilon_2^{(n)}, a_{\beta\alpha} \big)}} \\
\times \ \prod_{l=1}^k\,
\Zcal_{\C^2}^{\mathrm{cl}}\big(\varepsilon_1^{(l)},
\varepsilon_2^{(l)}, \vec{a}^{(l)}; \tau_1\big)\,
\Zcal_{\C^2}^{\mathrm{inst}} \big(\varepsilon_1^{(l)},
\varepsilon_2^{(l)}, \vec{a}^{(l)}; \qsf_{\mathrm{eff}}\big) \ ,
\end{multline}
where $\qsf_{\mathrm{eff}}:=\qsf \, \e^{\tau_1}$ and we introduced the
\emph{classical partition function} for supersymmetric $U(r)$ gauge theory on
$\mathbb{R}^4$ which is given by
\begin{equation}
\Zcal_{\C^2}^{\mathrm{cl}}(\varepsilon_1, \varepsilon_2, \vec{a};
\tau_1):=\exp\Big(-\frac{\tau_1}{2\varepsilon_1\, \varepsilon_2}\,
\sum_{\alpha=1}^r \, a_\alpha^2\, \Big) \ .
\end{equation}

\begin{proposition}
The partition function $\Zcal^{\circ}_{X_k}(\varepsilon_1,
\varepsilon_2, \vec{a}; \qsf, \tau_1,\vec{\xi}\ )$ factorizes as
\begin{equation}\label{eq:deformedinstanton}
\Zcal^{\circ}_{X_k}\big(\varepsilon_1, \varepsilon_2, \vec{a}; \qsf,
\tau_1,\vec{\xi} \ \big)
=\Zcal_{\C^2}^{\mathrm{cl}}\big(\varepsilon_1, \varepsilon_2, \vec{a};
\tau_1 \big)^{\frac{1}{k}} \ \Zcal^{\circ,\mathrm{inst}}_{X_k}
\big(\varepsilon_1, \varepsilon_2, \vec{a}; \qsf_{\mathrm{eff}},
\vec{\xi}, \tau_1 \big)\ ,
\end{equation}
where
\small
\begin{multline}\label{eq:instantontau}
\Zcal^{\circ,\mathrm{inst}}_{X_k}(\varepsilon_1, \varepsilon_2, \vec{a};\qsf_{\mathrm{eff}}, \vec{\xi}, \tau_1)\\
:= \sum_{\stackrel{\scriptstyle \vec{v}\in\frac{1}{k}\, \Z^{k-1}
    }{\scriptstyle k\, v_{k-1}=  \sum_{i=0}^{k-1}\, i\, w_i\bmod{k}
    }} \, \vec{\xi}^{\ \vec{v}} \ 
\sum_{\vec{\boldsymbol{v}}} \,
\frac{\big(\e^{\tau_1}\big)^{\frac{1}{2}\,
    \sum_{\alpha\neq\beta}\limits \, \vec{v}_\alpha\cdot
    C\vec{v}_\beta}\, \qsf_{\mathrm{eff}}^{\frac{1}{2}\,
    \sum_{\alpha=1}^r\limits \, \vec{v}_\alpha\cdot
    C\vec{v}_\alpha}}{\prod_{\alpha,\beta=1}^r\limits \
  \prod_{n=1}^{k-1} \limits\,
  \ell^{(n)}_{\vec{v}_{\beta\alpha}}\big(\varepsilon_1^{(n)},
  \varepsilon_2^{(n)}, a_{\beta\alpha} \big)} \ \prod_{l=1}^k\,
\Zcal_{\C^2}^{\mathrm{inst}}\big(\varepsilon_1^{(l)},
\varepsilon_2^{(l)}, \vec{a}^{(l)}; \qsf_{\mathrm{eff}} \big)\ .
\end{multline}
\normalsize
\end{proposition}
\proof
From the identities
\begin{equation}
\sum_{l=1}^k\, \frac{1}{\varepsilon_1^{(l)}\,
  \varepsilon_2^{(l)}}=\frac{1}{k\, \varepsilon_1\, \varepsilon_2}
\qquad \mbox{and} \qquad \sum_{l=1}^k\, \frac{(\vec{v}_\alpha)_l\,
  \varepsilon_1^{(l)}+(\vec{v}_\alpha)_{l-1}\,
  \varepsilon_2^{(l)}}{\varepsilon_1^{(l)}\, \varepsilon_2^{(l)}}=0\ ,
\end{equation}
it follows that
\begin{multline}
\sum_{l=1}^k\, \frac{1}{2\varepsilon_1^{(l)}\, \varepsilon_2^{(l)}} \
\sum_{\alpha=1}^r\, \big(a_\alpha^{(l)} \big)^2 \\
=\sum_{\alpha=1}^r\, \frac{a_\alpha^2}{2k\, \varepsilon_1\,
  \varepsilon_2}+\sum_{\alpha,\beta=1}^r \ \sum_{l=1}^k\,
\frac{\big((\vec{v}_\alpha)_l\,
    \varepsilon_1^{(l)}+(\vec{v}_\alpha)_{l-1}\,
    \varepsilon_2^{(l)}\big)\, \big((\vec{v}_\beta)_l\,
    \varepsilon_1^{(l)}+(\vec{v}_\beta)_{l-1}\,
    \varepsilon_2^{(l)}\big)}{2\varepsilon_1^{(l)} \,
    \varepsilon_2^{(l)}} \ .
\end{multline}
By the localization formula (cf.\ Equation \eqref{eq:localizationdivisors}) we get
\begin{multline}
\sum_{l=1}^k \, \frac{\big((\vec{v}_\alpha)_l\,
    \varepsilon_1^{(l)}+(\vec{v}_\alpha)_{l-1}\,
    \varepsilon_2^{(l)}\big)\, \big((\vec{v}_\beta)_l\,
    \varepsilon_1^{(l)}+(\vec{v}_\beta)_{l-1}\,
    \varepsilon_2^{(l)}\big)}{\varepsilon_1^{(l)}\,
    \varepsilon_2^{(l)}} \\
=\Big(\, \sum_{i=1}^{k-1}\, (\vec{v}_\alpha)_i \, [D_i]\, \Big)\cdot
\Big(\, \sum_{i=1}^{k-1}\, (\vec{v}_\beta)_i\, [D_i]\, \Big)=-\vec{v}_\alpha\cdot C\vec{v}_\beta\ .
\end{multline}
Thus
\begin{align}
\prod_{l=1}^k\, \Zcal_{\C^2}^{\mathrm{cl}}\big(\varepsilon_1^{(l)},
\varepsilon_2^{(l)}, \vec{a}^{(l)}; \tau_1
\big) &= \exp\Big(-\sum_{l=1}^k\, \frac{\tau_1}{2\varepsilon_1^{(l)}\,
  \varepsilon_2^{(l)}}\ \sum_{\alpha=1}^r \,
\big(a_\alpha^{(l)} \big)^2\,\Big) \\[4pt]
&= \exp\Big(-\frac{\tau_1}{2k\, \varepsilon_1\, \varepsilon_2} \
\sum_{\alpha=1}^r\, a_\alpha^2\Big)\, \exp\Big(\, \frac{\tau_1}{2}\,
\sum_{\alpha,\beta=1}^r\, \vec{v}_\alpha\cdot C\vec{v}_\beta\, \Big)
\\[4pt] 
&= \Zcal_{\C^2}^{\rm cl}\big(\varepsilon_1, \varepsilon_2, \vec{a};
\tau_1\big)^{\frac{1}{k}}\, \exp\Big(\, \frac{\tau_1}{2}\,
\sum_{\alpha,\beta=1}^r\, \vec{v}_\alpha\cdot C\vec{v}_\beta\, \Big)
\end{align}
and the assertion now follows straightforwardly.
\endproof
\begin{remark}
The partition function $\Zcal^{\circ,
  \mathrm{inst}}_{X_k}(\varepsilon_1, \varepsilon_2,
\vec{a};\qsf_{\mathrm{eff}}, \vec{\xi}, \tau_1)$  for $\tau_1=0$ coincides with
$\Zcal^{\mathrm{inst}}_{X_k}(\varepsilon_1, \varepsilon_2,$ $
\vec{a};\qsf_{\mathrm{eff}}, \vec{\xi}\ )$.
\end{remark}

\subsection{$\Ncal=2^\ast$ gauge theory}

\subsubsection*{Generating function for correlators of $p$-observables}

Let $T_\mu=\C^*$ and
$H^*_{T_\mu}(\mathrm{pt};\mathbb{Q})=\mathbb{Q}[\mu]$. For a
$T$-equivariant locally free sheaf $\boldsymbol G$ of rank $n$ on the moduli
space
$\Mcal_{r,\vec{u},\Delta}(\Xscr_k,\Dscr_\infty,\Fcal_\infty^{0,\vec
  w}\, )$ we define the class
\begin{equation}
{\rm E}_\mu({\boldsymbol G}):=\mu^n+(\crm_1)_T({\boldsymbol G}
)\, \mu^{n-1}+\cdots+(\crm_n)_T({\boldsymbol G} ) \
\in \ H^*_{T\times
  T_\mu}\big(\Mcal_{r,\vec{u},\Delta}(\Xscr_k,\Dscr_\infty,\Fcal_\infty^{0,\vec
  w}\, ) \,; \,\Q \big) \ .
\end{equation}
As previously, let $\vec{v}\in \frac{1}{k}\, \Z^{k-1}$ be such that $k\, v_{k-1}=
\sum_{i=0}^{k-1}\, i\, w_i \ \bmod{k} $ for fixed rank $r$ and
holonomy at infinity $\vec w$. Define
\begin{multline}\label{eq:deformedadjoint-v}
\Zcal^\ast_{\vec{v}}\big(\varepsilon_1, \varepsilon_2, \vec{a}, \mu;
\qsf, \vec{\tau}, \vec{t}^{\;(1)}, \ldots, \vec{t}^{\;(k-1)} \big) \\
\shoveleft{:=\sum_{\Delta\in \frac{1}{2r\, k}\Z} \,
  \qsf^{\Delta+\frac{1}{2r}\, \vec{v}\cdot C\vec{v}} \
  \int_{\Mcal_{r,\vec{u},\Delta}(\Xscr_k,\Dscr_\infty,\Fcal_\infty^{0,\vec
      w}\, )}\,
  {\rm E}_\mu \big(T\Mcal_{r,\vec{u},\Delta}(\Xscr_k,\Dscr_\infty,\Fcal_\infty^{0,\vec
    w}\, )\big)}\\
\times\ \exp \bigg(\, \sum_{s=0}^\infty\, \Big(\, \sum_{i=1}^{k-1}\,
t_s^{(i)} \, \big[\ch_{T}(\boldsymbol{\Ecal})/[\Dscr_i]\big]_s+\tau_s\,
\big[\ch_T(\boldsymbol{\Ecal})/[X_k]\big]_{s-1} \Big) \bigg) \ .
\end{multline}

\begin{definition}\label{def:deformed-adjoint}
The \emph{generating function for correlators of $p$-observables} of $\Ncal=2$ gauge theory on $X_k$ with an adjoint hypermultiplet of mass $\mu$ is
\begin{multline}
\Zcal^{\ast}_{X_k}\big(\varepsilon_1, \varepsilon_2, \vec{a}, \mu; \qsf,
\vec{\xi}, \vec{\tau}, \vec{t}^{\;(1)}, \ldots, \vec{t}^{\:(k-1)}
\big) \\ 
:=\sum_{\stackrel{\scriptstyle \vec{v}\in\frac{1}{k}\, \Z^{k-1}
  }{\scriptstyle k\, v_{k-1}=  \sum_{i=0}^{k-1}\, i\, w_i\bmod{k} }}\,
\vec{\xi}^{\ \vec{v}} \ \Zcal^\ast_{\vec{v}}\big(\varepsilon_1,
\varepsilon_2, \vec{a}, \mu; \qsf, \vec{\tau}, \vec{t}^{\;(1)}, \ldots,
\vec{t}^{\;(k-1)} \big)\ .
\end{multline}
\end{definition}
\begin{proposition}\label{prop:deformed-adjoint}
The generating function $\Zcal^{\ast}_{X_k}(\varepsilon_1,
\varepsilon_2, \vec{a}, \mu; \qsf, \vec{\xi}, \vec{\tau},
\vec{t}^{\;(1)}, \ldots, \vec{t}^{\:(k-1)})$ assumes the form
\begin{multline}\label{eq:N2ast-deformed}
\Zcal^{\ast}_{X_k}\big(\varepsilon_1, \varepsilon_2, \vec{a}, \mu; \qsf,
\vec{\xi}, \vec{\tau}, \vec{t}^{\;(1)}, \ldots, \vec{t}^{\:(k-1)}
\big) \\
\shoveleft{= \sum_{\stackrel{\scriptstyle \vec{v}\in\frac{1}{k}\Z^{k-1}
    }{\scriptstyle k\, v_{k-1}=  \sum_{i=0}^{k-1}\, i\, w_i\bmod{k}}
  }\, \vec{\xi}^{\ \vec{v}} \ 
\sum_{(\vec{\boldsymbol{Y}},\vec{\boldsymbol{v}})} \,
\qsf^{\sum_{\alpha=1}^r\limits\, n_\alpha+\frac{1}{2}\,
  \sum_{\alpha=1}^r\limits \, \vec{v}_\alpha\cdot C\vec{v}_\alpha}}\\
\shoveright{\times\ \frac{\adjointclass}{\eulerclass}}\\
\shoveright{\times\ \prod_{l=1}^k\, \exp\bigg(\sum_{s=0}^\infty\,
  \Big( \big(t_s^{(l)}\, \varepsilon_1^{(l)}+t_s^{(l-1)}\,
  \varepsilon_2^{(l)}+\tau_s\big)\,
  \Big[\ch_{\vec{Y}^l}\big(\varepsilon_1^{(l)},\varepsilon_2^{(l)},
  \vec{a}^{(l)}\big)\Big]_{s-1}\, \Big)\bigg)}\\
\shoveleft{\times\ \exp\Bigg(\, \sum_{s=0}^\infty\, \bigg(\,
  \sum_{l=1}^k\, \Big(\, \sum_{i=1}^{l-2}\,
  t_s^{(i)}+\sum_{i=l+1}^{k-1}\, t_s^{(i)}\Big)\,
  \Big[\ch_{\vec{Y}^l}\,\big(\varepsilon_1^{(l)},\varepsilon_2^{(l)},
  \vec{a}^{(l)}\big)\Big]_{s}}\\
+\, \sum_{i=1}^{k-1}\, t_s^{(i)} \ \sum_{l=2}^{k-1}\,
\Big(\, \Big[\big(-\varepsilon_1^{(l)}\big)^{\delta_{l,i}}\,
\ch_{\vec{Y}^l}\big(\varepsilon_1^{(l)},\varepsilon_2^{(l)},
\vec{a}^{(l)}+\left(\vec{v}_\alpha\right)_{l-1}\, \varepsilon_2^{(l)}
\big)\Big]_s \\
+\, \Big[\big(-\varepsilon_2^{(l)}\big)^{\delta_{l-1,i}}\,
\ch_{\vec{Y}^l}(\varepsilon_1^{(l)},\varepsilon_2^{(l)},
\vec{a}^{(l)}+\left(\vec{v}_\alpha\right)_{\ell}\,
\varepsilon_1^{(l)})\Big]_{s}\, \Big) \bigg)\Bigg) \ .
\end{multline}
\end{proposition}
\proof
By the localization formula we get
\begin{multline}
\Zcal^{\ast}_{X_k}\big(\varepsilon_1, \varepsilon_2, \vec{a}, \mu; \qsf,
\vec{\xi}, \vec{\tau}, \vec{t}^{\;(1)}, \ldots, \vec{t}^{\:(k-1)} \big)\\
\shoveleft{= \sum_{\stackrel{\scriptstyle
      \vec{v}\in\frac{1}{k}\, \Z^{k-1} }{\scriptstyle k\, v_{k-1}=
      \sum_{i=0}^{k-1}\, i\, w_i \bmod{k} }}\, \vec{\xi}^{\ \vec{v}} \
  \sum_{(\vec{\boldsymbol{Y}},\vec{\boldsymbol{v}})} \
  \frac{\qsf^{\sum_{\alpha=1}^r\limits\, n_\alpha+\frac{1}{2}\,
      \sum_{\alpha=1}^r\limits\, \vec{v}_\alpha\cdot
      C\vec{v}_\alpha}}{\eulerclass}}\\
\times \ \imath_{(\vec{\boldsymbol{Y}},\vec{\boldsymbol{v}})}^\ast
{\rm E}_\mu\big(T\Mcal_{r,\vec{u},\Delta}(\Xscr_k,\Dscr_\infty,\Fcal_\infty^{0,\vec
  w}\, )\big) \\
\times\ \imath_{(\vec{\boldsymbol{Y}},\vec{\boldsymbol{v}})}^\ast\exp
\bigg(\, \sum_{s=0}^\infty\, \Big(\, \sum_{i=1}^{k-1}\, t_s^{(i)} \,
\big[\ch_{T}(\boldsymbol{\Ecal})/[\Dscr_i]\big]_s+\tau_s\,
\big[\ch_T(\boldsymbol{\Ecal})/[X_k] \big]_{s-1}\Big)\bigg) \ .
\end{multline}
Note that
\begin{equation}
\imath_{(\vec{\boldsymbol{Y}},\vec{\boldsymbol{v}})}^\ast
{\rm E}_\mu\big(T\Mcal_{r,\vec{u},\Delta}(\Xscr_k,\Dscr_\infty,\Fcal_\infty^{0,\vec
  w}\, )\big)=\sum_{l=0}^d \,\mu^{d-l} \,
(\crm_l)_T\big(T_{(\vec{\boldsymbol{Y}},\vec{\boldsymbol{v}})}\Mcal_{r,\vec{u},\Delta}(\Xscr_k,\Dscr_\infty,\Fcal_\infty^{0,\vec
  w}\, )\big) 
\end{equation}
where $d$ the dimension of the moduli space
$\Mcal_{r,\vec{u},\Delta}(\Xscr_k,\Dscr_\infty,\Fcal_\infty^{0,\vec
  w}\, )$. Since the tangent space at the fixed point
$(\vec{\boldsymbol{Y}},\vec{\boldsymbol{v}})$ as a representation of
the torus $T$ is a direct sum of one-dimensional $T$-modules (see Section \ref{sec:Eulerclasses}), we get
\begin{multline}
\imath_{(\vec{\boldsymbol{Y}},\vec{\boldsymbol{v}})}^\ast
{\rm E}_\mu\big(T\Mcal_{r,\vec{u},\Delta}(\Xscr_k,\Dscr_\infty,\Fcal_\infty^{0,\vec
  w}\, )\big) \\
=\adjointclass\ .
\end{multline}
By using the computations of Section \ref{sec:pure} we then get the assertion.
\endproof

\subsubsection*{Instanton partition function}

\begin{definition}
The \emph{instanton partition function} is
\begin{equation}
\Zcal^{\ast, \mathrm{inst}}_{X_k}\big(\varepsilon_1,
\varepsilon_2, \vec{a}, \mu; \qsf, \vec{\xi}\
\big):=\Zcal^{\ast}_{X_k}\big(\varepsilon_1, \varepsilon_2, \vec{a},
\mu; \qsf, \vec{\xi}, 0, \ldots, 0 \big)\ .
\end{equation}
\end{definition}
By Equation \eqref{eq:N2ast-deformed} and the corresponding expression for the
{Nekrasov partition function} of gauge theory on
$\mathbb{R}^4$ with one adjoint hypermultiplet of mass $\mu$ \cite[Equation~(3.26)]{art:bruzzofucitomoralestanzini2003}
\begin{equation}
\Zcal_{\C^2}^{\ast, \mathrm{inst}}(\varepsilon_1, \varepsilon_2,
\vec{a}, \mu; \qsf):=\sum_{\vec{Y}}\, \qsf^{\sum_{\alpha=1}^r\limits \,
  \vert Y_\alpha\vert} \ \prod_{\alpha,\beta=1}^r\limits\,
\frac{m_{Y_\alpha,
    Y_\beta}(\varepsilon_1,\varepsilon_2,a_{\beta\alpha}+\mu)}{m_{Y_\alpha,
    Y_\beta}(\varepsilon_1,\varepsilon_2,a_{\beta\alpha})} \ ,
\end{equation}
we obtain
\begin{empheq}[box=\fbox]{multline}\label{eq:instantonpart-adjoint}
\Zcal^{\ast, \mathrm{inst}}_{X_k}\big(\varepsilon_1, \varepsilon_2,
\vec{a}, \mu; \qsf, \vec{\xi}\, \big)\\
\shoveleft{=\sum_{\stackrel{\scriptstyle \vec{v}\in\frac{1}{k}\,
      \Z^{k-1} }{\scriptstyle k\, v_{k-1}=  \sum_{i=0}^{k-1}\, i\,
      w_i\bmod{k} }}\, \vec{\xi}^{\ \vec{v}} \
  \sum_{\vec{\boldsymbol{v}}} \, \qsf^{\frac{1}{2}\,
    \sum_{\alpha=1}^r\limits\, \vec{v}_\alpha\cdot C\vec{v}_\alpha}\
  \frac{\prod_{\alpha,\beta=1}^r \limits\ \prod_{n=1}^{k-1}
    \limits\,\ell^{(n)}_{\vec{v}_{\beta\alpha}}\big(\varepsilon_1^{(n)},
    \varepsilon_2^{(n)}, a_{\beta\alpha}+\mu
    \big)}{\prod_{\alpha,\beta=1}^r \limits\
    \prod_{n=1}^{k-1}\limits\,
    \ell^{(n)}_{\vec{v}_{\beta\alpha}}\big(\varepsilon_1^{(n)},
    \varepsilon_2^{(n)}, a_{\beta\alpha} \big)}}\\
\times \ \prod_{i=1}^k\, \Zcal_{\C^2}^{\ast,
  \mathrm{inst}}\big(\varepsilon_1^{(i)}, \varepsilon_2^{(i)},
\vec{a}^{(i)}, \mu; \qsf \big) \ .
\end{empheq}

\subsubsection*{Correlators of quadratic $0$-observables}

As in Section \ref{sec:pure} we define
\begin{equation}
\Zcal^{\ast, \circ}_{X_k}\big(\varepsilon_1, \varepsilon_2, \vec{a},
\mu; \qsf, \vec{\xi}, \tau_1
\big):=\Zcal^{\ast}_{X_k}\big(\varepsilon_1, \varepsilon_2, \vec{a},
\mu; \qsf, \vec{\xi}, \, \underline{\vec{\tau}}\, , 0, \ldots, 0\big)\ .
\end{equation}
From Equation \eqref{eq:N2ast-deformed} it follows that
\begin{multline}
\Zcal^{\ast, \circ}_{X_k}\big(\varepsilon_1, \varepsilon_2, \vec{a},
\mu; \qsf, \vec{\xi}, \tau_1 \big) \\
\shoveleft{=\sum_{\stackrel{\scriptstyle \vec{v}\in\frac{1}{k}\, \Z^{k-1}
    }{\scriptstyle k\, v_{k-1}=  \sum_{i=0}^{k-1}\, i\, w_i\bmod{k}
    }}\, \vec{\xi}^{\ \vec{v}} \
  \sum_{(\vec{\boldsymbol{Y}},\vec{\boldsymbol{v}})} \,
  \qsf^{\sum_{\alpha=1}^r\limits\, n_\alpha+\frac{1}{2}\,
    \sum_{\alpha=1}^r\limits\, \vec{v}_\alpha\cdot C\vec{v}_\alpha} \
  \prod_{l=1}^k\, \exp\Big(-\tau_1\,
  \Big[\ch_{\vec{Y}^l}\big(\varepsilon_1^{(l)},\varepsilon_2^{(l)},
  \vec{a}^{(l)}\big)\Big]_{0}\, \Big)}\\
\times \ \frac{\adjointclass}{\eulerclass}\ .
\end{multline}
As in the case of pure $\Ncal=2$ gauge theory we get
\begin{empheq}[box=\fbox]{equation}
\Zcal^{\ast, \circ}_{X_k}\big(\varepsilon_1, \varepsilon_2, \vec{a},
\mu; \qsf, \vec{\xi}, \tau_1
\big)=\Zcal_{\C^2}^{\mathrm{cl}}\big(\varepsilon_1, \varepsilon_2,
\vec{a};\tau_1 \big)^{\frac{1}{k}} \ \Zcal^{\ast,\circ,
  \mathrm{inst}}_{X_k}\big(\varepsilon_1, \varepsilon_2, \vec{a}, \mu;
\qsf_{\mathrm{eff}}, \vec{\xi}, \tau_1 \big)\ ,
\end{empheq}
where the partition function $\Zcal^{\ast, \circ, \mathrm{inst}}_{X_k}(\varepsilon_1, \varepsilon_2, \vec{a}, \mu; \qsf_{\mathrm{eff}}, \vec{\xi}, \tau_1)$ is defined analogously to \eqref{eq:instantontau}.

\subsection{$\Ncal=4$ gauge theory\label{sec:VW}}

\begin{definition}
The \emph{Vafa-Witten partition function} for $\Ncal=4$ gauge theory on $X_k$ is 
\begin{equation}
\Zcal^{\mathrm{VW}}_{X_k}\big(\qsf, \vec{\xi}\
\big):=\lim_{\mu\to 0}\, \Zcal^{\ast, \mathrm{inst}}_{X_k}
\big(\varepsilon_1, \varepsilon_2, \vec{a}, \mu; \qsf, \vec{\xi}\ \big)
\ .
\end{equation}
\end{definition}
By using our previous computations we get
\begin{align}
\Zcal^{\mathrm{VW}}_{X_k}\big(\qsf, \vec{\xi}\
\big)&=\sum_{\stackrel{\scriptstyle \vec{v}\in\frac{1}{k}\, \Z^{k-1}
  }{ k\, v_{k-1}=\sum_{i=0}^{k-1}\,i\, w_i \bmod{k}}}\, \vec{\xi}^{\
  \vec{v}} \ \sum_{(\vec{\boldsymbol{Y}},\vec{\boldsymbol{v}})} \,
\qsf^{\sum_{\alpha=1}^r\limits \, n_\alpha+\frac{1}{2}\,
  \sum_{\alpha=1}^r\limits \, \vec{v}_\alpha\cdot C\vec{v}_\alpha}\\ \label{eq:limit}
&=\qsf^{\frac{r\, k}{24}} \, \eta(\qsf)^{-r\, k} \,
\sum_{\stackrel{\scriptstyle \vec{v}\in\frac{1}{k}\, \Z^{k-1} }{ k\,
    v_{k-1}= \sum_{i=0}^{k-1}\, i\, w_i\bmod{k}}}\, \vec{\xi}^{\
  \vec{v}} \ \sum_{\vec{\boldsymbol{v}}} \, \qsf^{\frac{1}{2}\,
  \sum_{\alpha=1}^r\limits\, \vec{v}_\alpha\cdot C\vec{v}_\alpha} \ ,
\end{align}
where
\begin{equation}
\eta(\qsf):=\qsf^{\frac{1}{24}}\, \prod_{n=1}^\infty \, \big(1-\qsf^n \big)
\end{equation}
is the \emph{Dedekind $\eta$ function}. 

Let $c\in\{0, 1, \ldots, k-1\}$ be the equivalence class of
$\sum_{i=0}^{k-1}\, i\, w_i$ modulo $k$. Set
$\vec{n}:=\vec{v}-C^{-1}\vec{e}_c$ if $c>0$ and $\vec{n}:=\vec{v}$ otherwise, where $\vec{e}_c$ is the $c$-th coordinate vector of $\Z^{k-1}$. Then $n_l=v_l-(C^{-1})^{l c}$ for any $l=1, \ldots, k-1$. 
Similarly, for any $\alpha=1, \ldots, r$ define
$\vec{n}_\alpha:=\vec{v}_\alpha-C^{-1}\vec{e}_{c_{\alpha}}$ if
$c_\alpha>0$ and $\vec{n}_\alpha:=\vec{v}_\alpha$ otherwise, where
$c_\alpha=i$ if $\sum_{j=0}^{i-1}\, w_j<\alpha\leq\sum_{j=0}^i\, w_j$ for $i=0,1, \ldots, k-1$. Both $\vec{n}$ and $\vec{n}_\alpha$ for $\alpha=1, \ldots, r$ can be regarded as elements in the root lattice $\Qfrak$ of type $A_{k-1}$ (cf.\ Remark \ref{rem:v-y}). Note that $\sum_{\alpha=1}^r\, c_\alpha =c$. Then
Equation \eqref{eq:limit} becomes
\begin{multline}
\Zcal^{\rm VW}_{X_k}\big(\qsf, \vec{\xi}\ \big)\\ =
\qsf^{\frac{r\, k}{24}}\, \eta(\qsf)^{-r\, k} \,
\sum_{\vec{n}\in\Z^{k-1}} \ \sum_{\stackrel{\scriptstyle (\vec{n}_1,
    \ldots, \vec{n}_r) }{\scriptstyle \sum_{\alpha=1}^r\,
    \vec{n}_\alpha=\vec{n}}} \ \prod_{i=1}^{k-1}\,
\xi_i^{\sum_{\alpha=1}^r\limits\,
  ((\vec{n}_\alpha)_i+(C^{-1})^{i,c_\alpha})} \
\qsf^{\frac{1}{2}\, \sum_{\alpha=1}^r\limits\, (\vec{n}_\alpha \cdot
  C\vec{n}_\alpha+2
  (\vec{n}_{\alpha})_{c_\alpha}+\frac{k-c_\alpha}{k}\, c_\alpha)} \ .
\end{multline}

Recall from \cite[Section~14.4]{book:difrancescomathieusenechal1997}
(see also \cite[Appendix~C.2]{art:dijkgraafsulkowski2008} and \cite[Section
5.4]{art:ciraficiszabo2012}) that the
character of an integrable highest weight representation of
$\widehat{\slfrak}(k)$ with highest weight $\widehat{\lambda}$, whose finite part is $\lambda$, at level one is a
combination of string-functions and theta-functions given by
\begin{equation}
\chi^\lambda(\qsf, \zeta):=\eta(\qsf)^{-k+1}\, \sum_{\gamma^\vee \in
  \Qfrak^\vee}\, \qsf^{\frac{1}{2}\,
  \vert\gamma^\vee+\lambda\vert^2}\, \e^{2\pi \ii (\gamma^\vee+\lambda, \zeta)}\ ,
\end{equation}
where $\zeta:=\sum_{i=1}^{k-1}\, z_i\, \gamma_i^\vee$ and
$\gamma_i^\vee$ denotes the $i$-th simple coroot for $i=1, \ldots,
k-1$. Recall that the $d$-th fundamental weight of type $A_{k-1}$ is
\begin{equation}
\omega_d:=\big(\, \underbrace{\mbox{$1-\frac{d}{k}, \ldots,
    1-\frac{d}{k}$}}_{d{\rm -times}}\, ,\mbox{$ -\frac{d}{k}, \ldots,
  -\frac{d}{k}$} \, \big)\in\Z^k
\end{equation}
for $d=1, \ldots, k-1$ and the fundamental weights of type $\hat{A}_{k-1}$ are $\widehat{\omega}_0$ and $\widehat{\omega}_d=\omega_d+\widehat{\omega}_0$ for $d=1, \ldots, k-1$. With $\gamma^\vee:=\sum_{i=1}^{k-1}\,
m_i\, \gamma_i^\vee$ we then get
\begin{align}	
\frac{1}{2}\, \big\vert\gamma^\vee+\omega_d \big\vert^2 &=
\sum_{i=1}^k\, \big(m_i^2-m_i\, m_{i+1} \big) +m_d+\frac{k-d}{2k}\, d\ ,\\[4pt]
\e^{2\pi \ii(\gamma^\vee+\omega_d, \zeta)} &=
 y_1^{2m_1-m_2} \, \cdots\, y_{d-1}^{2 m_{d-1}-m_{d-2}-m_d} \,
y_d^{2m_d-m_{d-1}-m_{d+1}+1} \\ & \qquad \qquad \times \
y_{d+1}^{2m_{d+1}-m_d-m_{d+2}}\, \cdots\, y_{k-1}^{2 m_{k-1}-m_{k-2}}
\ = \ \prod_{i=1}^{k-1}\, \xi_i^{m_i+(C^{-1})^{i,d}}\ ,
\end{align}
where we introduced $y_i:=\e^{2\pi \ii z_i}$ and 
\begin{equation}
\xi_1:=\frac{y_1^2}{y_2}\ ,\quad \xi_2:= \frac{y_2^2}{y_1\, y_3}\
,\quad \ldots\ ,\quad \xi_{k-1}:=\frac{y_{k-1}^2}{y_{k-2}}\ .
\end{equation}
It follows that
\begin{empheq}[box=\fbox]{equation}\label{eq:N4}
\Zcal^{\mathrm{VW}}_{X_k} \big(\qsf, \vec{\xi}\
\big)=\qsf^{\frac{r\, k}{24}}\, \prod_{\alpha=1}^r\,
\frac{\chi^{\widehat{\omega}_{c_\alpha}}(\qsf, \zeta )}{\eta(\qsf)}=
\qsf^{\frac{r\, k}{24}}\, \prod_{j=0}^{k-1}\, \Big(\,
\frac{\chi^{\widehat{\omega}_{j}}(\qsf, \zeta)}{\eta(\qsf)}\, \Big)^{w_j}\ .
\end{empheq}
This expression agrees with the formula obtained in
\cite[Corollary~4.12]{art:fujiiminabe2005} (see also
\cite[Remark~4.2.3]{art:fujiiminabe2005}). For trivial holonomies at
infinity, i.e., when $c_\alpha=0$ for any $\alpha=1, \ldots, r$, which
is equivalent to $w_0=r$ and $w_j=0$ for $j=1, \ldots, k-1$, this
expression agrees with
\cite[Equation~(4.8)]{art:fucitomoralespoghossian2006} and
\cite[Equation~(3.23)]{art:griguoloseminaraszabotanzini2007}. It
displays the Vafa-Witten partition function as a character of the affine
Lie algebra $\widehat{\mathfrak{gl}}(k)_r$, which confirms its
correct modularity properties as implied by S-duality. 

\subsection{Fundamental matter}\label{sec:fundamentalmatter}

\subsubsection*{Generating function for correlators of $p$-observables}

In this section we shall assume that $\Delta\geq 0$. Let $N\leq 2r$ be a positive integer, $T_{N}=(\C^\ast)^N$ and
$H_{T_{N}}^\ast({\rm pt};\mathbb{Q})=\mathbb{Q}[\mu_1, \ldots,$ $
\mu_{N}]$. Let $\Vbf$ be the natural bundle on
$\Mcal_{r,\vec{u},\Delta}(\Xscr_k,\Dscr_\infty,\Fcal_\infty^{0,\vec
  w}\, )$; recall that $\Vbf$ is a $T$-equivariant locally free
sheaf. To define the partition functions and correlators for gauge theories with fundamental matter we consider the class
\begin{equation}
\prod_{s=1}^{N} \, {\rm E}_{\mu_s}(\Vbf) \ \in \ H^\ast_{T\times
  T_N}\big(\Mcal_{r,\vec{u},\Delta}(\Xscr_k,\Dscr_\infty,\Fcal_\infty^{0,\vec
  w}\, ) \,; \,\Q \big)\ .
\end{equation}
It has degree 
\begin{equation}\label{eq:fundamental-condition}
\dim_\C\big(\Mcal_{r,\vec{u},\Delta}(\Xscr_k,\Dscr_\infty,\Fcal_\infty^{0,\vec
  w}\, ) \big)- N\, \rk(\Vbf) \ .
\end{equation}
By using Theorem \ref{thm:moduli} and Proposition
\ref{prop:naturalbundle}, we see that the degree
\eqref{eq:fundamental-condition} is nonnegative if and only if
\begin{equation}\label{eq:fundamental-inequality}
2r\, \Delta - \frac{1}{2}\, \sum_{j=1}^{k-1}\, (C^{-1}\big)^{jj}\,
\vec{w}\cdot\vec{w}(j)-N\, \Big(\Delta +\frac{1}{2r}\, \vec v\cdot
C\vec v - \frac{1}{2} \, \sum_{j=1}^{k-1}\,\big(C^{-1} \big)^{jj}\,
w_j \Big)\geq 0
\end{equation}
for $\vec{v}:=C^{-1}\vec{u}$.
\begin{remark}\label{rem:fundamental-constraint}
We assume   $N\leq 2r$ because this is the constraint on the number
of fundamental matter fields in asymptotically free $\Ncal=2$ gauge
theories on $\R^4$. We impose it in our case because we want to compare
our gauge theories with fundamental matter on $X_k$ to those on
$\R^4$. As an example, let us write   explicitly the inequality
\eqref{eq:fundamental-inequality} for $k=2$, in which case it becomes
\begin{equation}
w_1^2-4v^2\geq 0\ .
\end{equation}
If the fixed holonomy at infinity is trivial, i.e., $\vec{w}=(r, 0)$,
this gives $v=0$; this is the case considered in
\cite{art:bonellimaruyoshitanzini2012}. For arbitrary $k\geq2$, if
$\vec{w}=(r, 0, \ldots, 0)$ and $\vec{v}=\vec 0$, the inequality \eqref{eq:fundamental-inequality} is automatically satisfied when $N\leq 2r$.
\end{remark}
Since $N\leq 2r$, we get
\begin{equation}\label{eq:fundamental-inequality-II}
2r\, \Delta - \frac{1}{2}\, \sum_{j=1}^{k-1}\, \big(C^{-1}\big)^{jj}\,
\vec{w}\cdot\vec{w}(j)-N\,\Big(\Delta +\frac{1}{2r}\, \vec v\cdot
C\vec v - \frac{1}{2}\, \sum_{j=1}^{k-1}\, \big(C^{-1}\big)^{jj}\, w_j
\Big)\geq d_{\vec w}(\vec v\,)
\end{equation}
where we defined
\begin{equation}
d_{\vec{w}}(\vec{v}\,):=r \, \sum_{j=1}^{k-1}\,
\big(C^{-1}\big)^{jj}\, w_j - \frac{1}{2}\, \sum_{j=1}^{k-1}\,
\big(C^{-1} \big)^{jj}\, \vec{w}\cdot\vec{w}(j)-\vec v\cdot C\vec v\ .
\end{equation}
We define the set
\begin{equation}\label{eq:fundamental-set}
\Qfrak_{\vec{w}}:=\Big\{\vec{v}\in \mbox{$\frac{1}{k}$}\, \Z^{k-1}\
\Big\vert\ k\, v_{k-1}= \mbox{$\sum\limits_{i=0}^{k-1}$} \, i\, w_i \ \bmod{k}\;\;
\mbox{and}\;\; d_{\vec{w}}(\vec{v}\, )\geq 0\Big\}\ .
\end{equation}
Then the inequality \eqref{eq:fundamental-inequality} is satisfied for
any choice of discriminant $\Delta$ and for $\vec{v}\in
\Qfrak_{\vec{w}}$. Although the constraint $d_{\vec{w}}(\vec{v}\,)\geq 0$
is stronger than \eqref{eq:fundamental-inequality}, it is natural in light of the discussion in Remark \ref{rem:fundamental-constraint}.

We call the $\Ncal=2$ supersymmetric gauge theory on $X_k$ with $N\leq
2r$ fundamental hypermultiplets \emph{asymptotically free}. If $N=2r$,
we call the associated gauge theory \emph{conformal}; in this case we
consider the subsets 
\begin{equation}
\Qfrak_{\vec{w}}^{\rm conf}:=\Big\{\vec{v}\in \mbox{$\frac{1}{k}$}\, \Z^{k-1}\
\Big\vert\ k\, v_{k-1}= \mbox{$\sum\limits_{i=0}^{k-1}$} \, i\, w_i \ \bmod{k}\;\;
\mbox{and}\;\; d_{\vec{w}}(\vec{v}\, )= 0\Big\}  \ \subset \
\Qfrak_{\vec w}
\end{equation}
for which the degree \eqref{eq:fundamental-condition} is equal to $0$.
In the following we treat only the asymptotically free case. The
conformal case is completely analogous, and simply amounts to
restricting the set $\Qfrak_{\vec w}$ to its subset $\Qfrak_{\vec w}^{\rm conf}$
in all sums below.

For fixed $\vec{v}\in \Qfrak_{\vec{w}}$, we define
\begin{multline}
\Zcal^{N}_{\vec{v}}\big(\varepsilon_1, \varepsilon_2, \vec{a},
\vec{\mu}; \qsf, \vec{\tau}, \vec{t}^{\;(1)}, \ldots, \vec{t}^{\;(k-1)} \big)\\
\shoveleft{:=\sum_{\Delta\in \frac{1}{2rk}\Z}\,
  \qsf^{\Delta+\frac{1}{2r}\, \vec{v}\cdot C\vec{v}} \
  \int_{\Mcal_{r,\vec{u},\Delta}(\Xscr_k,\Dscr_\infty,\Fcal_\infty^{0,\vec
      w}\, )} \ \prod_{s=1}^{N}\, {\rm E}_{\mu_s}(\Vbf)}\\
\times \ \exp \bigg(\, \sum_{s=0}^\infty\, \Big(\, \sum_{i=1}^{k-1}\,
t_s^{(i)} \, \big[\ch_{T}(\boldsymbol{\Ecal})/[\Dscr_i]\big]_s+\tau_s\,
\big[\ch_T(\boldsymbol{\Ecal})/[X_k]\big]_{s-1}\, \Big)\bigg) \ .
\end{multline}
\begin{definition}\label{def:deformed-fund}
The \emph{generating function for correlators of $p$-observables} of $\Ncal=2$ gauge theory on $X_k$ with $N$ fundamental hypermultiplets of masses $\mu_1, \ldots, \mu_{N}$ is
\begin{equation}
\Zcal^{N}_{X_k}\big(\varepsilon_1, \varepsilon_2, \vec{a},\vec{\mu};
\qsf, \vec{\xi}, \vec{\tau}, \vec{t}^{\;(1)}, \ldots,
\vec{t}^{\:(k-1)} \big)
:=\sum_{\vec{v}\in \Qfrak_{\vec{w}}}\, \vec{\xi}^{\ \vec{v}} \
\Zcal^{N}_{\vec{v}} \big(\varepsilon_1, \varepsilon_2, \vec{a},
\vec{\mu} ; \qsf,
\vec{\tau}, \vec{t}^{\;(1)}, \ldots, \vec{t}^{\;(k-1)} \big)\ .
\end{equation}
\end{definition}
\begin{proposition}\label{prop:deformed-fund}
The generating function $\Zcal^{N}_{X_k}(\varepsilon_1, \varepsilon_2, \vec{a},\vec{\mu}; \qsf, \vec{\xi}, \vec{\tau}, \vec{t}^{\;(1)}, \ldots, \vec{t}^{\:(k-1)})$ assumes the form
\begin{multline}\label{eq:N2deformed-fundamental}
\Zcal^{N}_{X_k}\big(\varepsilon_1, \varepsilon_2, \vec{a},\vec{\mu};
\qsf, \vec{\xi}, \vec{\tau}, \vec{t}^{\;(1)}, \ldots,
\vec{t}^{\:(k-1)} \big) \\ 
\shoveleft{=\sum_{\vec{v}\in \Qfrak_{\vec{w}}}\, \vec{\xi}^{\ \vec{v}} \
  \sum_{(\vec{\boldsymbol{Y}},\vec{\boldsymbol{v}})} \,
  \qsf^{\sum_{\alpha=1}^r\limits \, n_\alpha+\frac{1}{2}\, \sum_{\alpha=1}^r\limits\, \vec{v}_\alpha\cdot C\vec{v}_\alpha}}\\
\times \frac{\prod_{s=1}^{N}\limits\,
  \fundclass-s}{\eulerclass}\\
\shoveright{\times\ \prod_{l=1}^k\, \exp\bigg(\sum_{s=0}^\infty\,
  \Big( \big(t_s^{(l)}\, \varepsilon_1^{(l)}+t_s^{(l-1)}\,
  \varepsilon_2^{(l)}+\tau_s\big)\,
  \Big[\ch_{\vec{Y}^l}\big(\varepsilon_1^{(l)},\varepsilon_2^{(l)},
  \vec{a}^{(l)}\big)\Big]_{s-1}\, \Big)\bigg)}\\
\shoveleft{\times\ \exp\Bigg(\, \sum_{s=0}^\infty\, \bigg(\,
  \sum_{l=1}^k\, \Big(\, \sum_{i=1}^{l-2}\,
  t_s^{(i)}+\sum_{i=l+1}^{k-1}\, t_s^{(i)}\Big)\,
  \Big[\ch_{\vec{Y}^l}\,\big(\varepsilon_1^{(l)},\varepsilon_2^{(l)},
  \vec{a}^{(l)}\big)\Big]_{s}}\\
+\, \sum_{i=1}^{k-1}\, t_s^{(i)} \ \sum_{l=2}^{k-1}\,
\Big(\, \Big[\big(-\varepsilon_1^{(l)}\big)^{\delta_{l,i}}\,
\ch_{\vec{Y}^l}\big(\varepsilon_1^{(l)},\varepsilon_2^{(l)},
\vec{a}^{(l)}+\left(\vec{v}_\alpha\right)_{l-1}\, \varepsilon_2^{(l)}
\big)\Big]_s \\
+\, \Big[\big(-\varepsilon_2^{(l)}\big)^{\delta_{l-1,i}}\,
\ch_{\vec{Y}^l}(\varepsilon_1^{(l)},\varepsilon_2^{(l)},
\vec{a}^{(l)}+\left(\vec{v}_\alpha\right)_{\ell}\,
\varepsilon_1^{(l)})\Big]_{s}\, \Big) \bigg)\Bigg) \ .
\end{multline}
\end{proposition}
\proof
The proof is completely analogous to the proof of Proposition
\ref{prop:deformed-adjoint} by noting that
\begin{equation}
\imath_{(\vec{\boldsymbol{Y}},\vec{\boldsymbol{v}})}^\ast
{\rm E}_{\mu_s}(\Vbf)=\sum_{l=0}^{\rk(\Vbf)} \, \mu_s^{\rk(\Vbf)-l} \,
(\crm_l)_T \big(\Vbf_{(\vec{\boldsymbol{Y}},\vec{\boldsymbol{v}})} \big) \ ,
\end{equation}
and since $\Vbf_{(\vec{\boldsymbol{Y}},\vec{\boldsymbol{v}})}$ as a
$T$-module is a direct sum of one-dimensional $T$-modules (see Section
\ref{sec:Eulerclasses}) we get
\begin{equation}
\imath_{(\vec{\boldsymbol{Y}},\vec{\boldsymbol{v}})}^\ast {\rm E}_{\mu_s}(\Vbf)=\fundclass-s\ .
\end{equation}
\endproof

\subsubsection*{Instanton partition function}

\begin{definition}
The \emph{instanton partition function} is
\begin{equation}
\Zcal^{N, \mathrm{inst}}_{X_k}\big(\varepsilon_1, \varepsilon_2,
\vec{a},\vec{\mu}; \qsf, \vec{\xi}\ \big):=\Zcal^{N}_{X_k}
\big(\varepsilon_1, \varepsilon_2, \vec{a},\vec{\mu}; \qsf, \vec{\xi},
0, \ldots, 0 \big)\ .
\end{equation}
\end{definition}
By Equation \eqref{eq:N2deformed-fundamental} and the corresponding
expression for the {Nekrasov partition function} of $U(r)$ gauge theory on $\mathbb{R}^4$ with $N$ fundamental hypermultiplets of masses $\mu_1, \ldots, \mu_{N}$ \cite[Equation~(3.22)]{art:bruzzofucitomoralestanzini2003}
\begin{equation}
\Zcal^{N, \mathrm{inst}}_{\C^2}(\varepsilon_1, \varepsilon_2, \vec{a},
\vec{\mu} ; \qsf):=\sum_{\vec{Y}}\,
\qsf^{\sum_{\alpha=1}^r\limits \, \vert Y_\alpha\vert} \
\frac{\prod_{s=1}^{N}\limits \ \prod_{\alpha=1}^r\limits \, m_{Y_\alpha}(\varepsilon_1,\varepsilon_2,a_\alpha+\mu_s)}{\prod_{\alpha,\beta=1}^r\limits\, m_{Y_\alpha, Y_\beta}(\varepsilon_1,\varepsilon_2,a_{\beta\alpha})}\ ,
\end{equation}
we get
\small
\begin{empheq}[box=\fbox]{multline}
\Zcal^{N, \mathrm{inst}}_{X_k}\big(\varepsilon_1, \varepsilon_2,
\vec{a},\vec{\mu}; \qsf, \vec{\xi} \ \big) \\
=\sum_{\vec{v}\in \Qfrak_{\vec{w}}}\, \vec{\xi}^{\ \vec{v}} \
\sum_{\vec{\boldsymbol{v}}}\, \qsf^{\frac{1}{2}\,
  \sum_{\alpha=1}^r\limits\, \vec{v}_\alpha\cdot C\vec{v}_\alpha} \
\frac{\prod_{s=1}^{N}\limits \ \prod_{\alpha=1}^r\limits \
  \prod_{n=1}^{k-1} \limits
  \ell^{(n)}_{\vec{v}_{\alpha}}\big(\varepsilon_1^{(n)},
  \varepsilon_2^{(n)}, a_\alpha+\mu_s
  \big)}{\prod_{\alpha,\beta=1}^r\limits \ \prod_{n=1}^{k-1}\limits\,
  \ell^{(n)}_{\vec{v}_{\beta\alpha}}\big(\varepsilon_1^{(n)},
  \varepsilon_2^{(n)}, a_{\beta\alpha} \big)} \ 
\prod_{i=1}^k\, \Zcal^{N, \mathrm{inst}}_{\C^2}
\big(\varepsilon_1^{(i)}, \varepsilon_2^{(i)}, \vec{a}^{(i)},\vec{\mu};
\qsf \big)\ .
\end{empheq}
\normalsize

\subsubsection*{Correlators of quadratic $0$-observables}

With notation as before, let us define
\begin{equation}
\Zcal^{N, \circ}_{X_k}\big(\varepsilon_1, \varepsilon_2,
\vec{a},\vec{\mu}; \qsf,\vec{\xi}, \tau_1 \big):=\Zcal^{N}_{X_k}
\big(\varepsilon_1, \varepsilon_2, \vec{a},\vec{\mu}; \qsf, \vec{\xi},
\,\underline{\vec{\tau}} \, , 0, \ldots, 0 \big)\ .
\end{equation}
From Equation \eqref{eq:N2deformed-fundamental} it follows that
\begin{multline}
\Zcal^{N, \circ}_{X_k}\big(\varepsilon_1, \varepsilon_2,
\vec{a},\vec{\mu}; \qsf, \vec{\xi}, \tau_1 \big)=
\sum_{\vec{v}\in \Qfrak_{\vec{w}}}\, \vec{\xi}^{\ \vec{v}} \
\sum_{(\vec{\boldsymbol{Y}},\vec{\boldsymbol{v}})} \,
\qsf^{\sum_{\alpha=1}^r\limits \, n_\alpha+\frac{1}{2}\,
  \sum_{\alpha=1}^r\limits \, \vec{v}_\alpha\cdot C\vec{v}_\alpha}\\
\shoveright{\times \ \frac{\prod_{s=1}^{N}\limits \
    \fundclass-s}{\eulerclass} } \\
\times \ \prod_{l=1}^k\, \exp\Big(-\tau_1\,
  \Big[\ch_{\vec{Y}^l} \big(\varepsilon_1^{(l)},\varepsilon_2^{(l)},
  \vec{a}^{(l)} \big)\Big]_{0}\, \Big)\ .
\end{multline}
As previously we get
\begin{empheq}[box=\fbox]{equation}
\Zcal^{N, \circ}_{X_k}\big(\varepsilon_1, \varepsilon_2,
\vec{a},\vec{\mu}; \qsf, \vec{\xi}, \tau_1 \big)=
\Zcal_{\C^2}^{\mathrm{cl}}\big(\varepsilon_1, \varepsilon_2, \vec{a};
\tau_1 \big)^{\frac{1}{k}} \ \Zcal^{N, \circ, \mathrm{inst}}_{X_k}
\big(\varepsilon_1, \varepsilon_2, \vec{a},\vec{\mu};
\qsf_{\mathrm{eff}}, \vec{\xi}, \tau_1 \big)\ ,
\end{empheq}
where the partition function $\Zcal^{N, \circ,
  \mathrm{inst}}_{X_k}(\varepsilon_1, \varepsilon_2, \vec{a},\vec{\mu};
\qsf_{\mathrm{eff}}, \vec{\xi}, \tau_1)$ is defined analogously to \eqref{eq:instantontau}.

\bigskip\section{Perturbative partition functions}\label{sec:perturbativepartitionfunctions}

\subsection{Perturbative Euler classes}\label{sec:Eulerclasses-pert}

Let $[(\Ecal,\phi_{\Ecal})]$ and $[(\Ecal',\phi_{\Ecal'})]$ be
$T$-fixed points of the respective moduli spaces $\Mcal_{r, \vec{u},
  \Delta}(\Xscr_k, \Dscr_\infty, \Fcal_{\infty}^{s,\vec{w}}\, )$ and
$\Mcal_{r', \vec{u}{}', \Delta'}(\Xscr_k, \Dscr_\infty,
\Fcal_{\infty}^{s,\vec{w}{}'}\, )$. By Equation
\eqref{eq:character-carlssonokounkovbundle} the corresponding character of the
Carlsson-Okounkov bundle is given by
\begin{equation}
\ch_{\tilde{T}} \Ebf_{([(\Ecal,\phi_{\Ecal})]\,,\,
  [(\Ecal',\phi_{\Ecal'})] )}=\sum_{\alpha=1}^r \
\sum_{\beta=1}^{r'}\, e_\beta \, e_\alpha^{\prime\, -1}\,
\big(M_{\alpha\beta}(t_1,t_2)+L_{\alpha\beta}(t_1,t_2) \big)\ ,
\end{equation}
where the vertex contribution
$M_{\alpha\beta}(t_1,t_2)$ is given in Equation
\eqref{eq:vertexcontribution} and the edge contribution $L(t_1,t_2)$
is defined by Equation \eqref{eq:edge} which can be written in the form
\begin{equation} \label{eq:equivEulerchar}
L_{\alpha\beta}(t_1,t_2)= -\chi_{T_t}\big(\bar{X}_k\,,\, {\pi_k}_\ast(\Rcal^{\vec u_\beta{}'-\vec u_\alpha}\otimes \Ocal_{\Xscr_k}(-\Dscr_\infty))\big)\ .
\end{equation} 

\begin{definition}
The \emph{virtual $\tilde{T}$-equivariant Chern character} of the Carlsson-Okounkov bundle at the $\tilde{T}$-fixed point $([(\Ecal,\phi_{\Ecal})], [(\Ecal',\phi_{\Ecal'})])$ is
\begin{equation}
\ch_{\tilde{T}}^{\mathrm{vir}} \Ebf_{([(\Ecal,\phi_{\Ecal})]\,,\,
  [(\Ecal',\phi_{\Ecal'})])}:= \sum_{\alpha=1}^r \
\sum_{\beta=1}^{r'}\, e_\beta \, e_\alpha^{\prime\, -1}\,
\big(M_{\alpha\beta}(t_1,t_2)+L_{\alpha\beta}^\circ(t_1,t_2) \big)\ ,
\end{equation}
where $L_{\alpha\beta}^\circ(t_1,t_2):=-\chi_{T_t}\big(X_k \,,\,
\Rcal^{\vec{u}_\beta{}'-\vec{u}_\alpha}_{\vert X_k} \big)$. 
The \emph{perturbative part} of the $\tilde{T}$-equivariant Chern character of the Carlsson-Okounkov bundle at the $\tilde{T}$-fixed point $([(\Ecal,\phi_{\Ecal})], [(\Ecal',\phi_{\Ecal'})])$ is
\begin{equation}
\ch_{\tilde{T}}^{\mathrm{pert}} \Ebf_{([(\Ecal,\phi_{\Ecal})]\,,\, [(\Ecal',\phi_{\Ecal'})])}:=
\ch_{\tilde{T}}^{\mathrm{vir}} \Ebf_{([(\Ecal,\phi_{\Ecal})]\,,\,
  [(\Ecal',\phi_{\Ecal'})])}-\ch_{\tilde{T}}
\Ebf_{([(\Ecal,\phi_{\Ecal})]\,,\, [(\Ecal',\phi_{\Ecal'})])}\ .
\end{equation}
\end{definition}

By using arguments from the proof of
\cite[Theorem~4.1]{art:brionvergne1997}, one can decompose the
equivariant Euler characteristic in \eqref{eq:equivEulerchar} with
respect to the torus-invariant open subsets of $\bar X_k$ to get
\begin{equation}
L_{\alpha\beta}(t_1,t_2)=-\sum_{\sigma\in\bar{\Sigma}_k(2) } \, \mathsf{S}\Big(\chi_{T_t}^\sigma\big(U_\sigma\,,\, {\pi_k}_\ast(\Rcal^{\vec u_\beta{}'-\vec u_\alpha}\otimes \Ocal_{\Xscr_k}(-\Dscr_\infty))_{\vert U_\sigma}\big)\Big)\ ,
\end{equation}
where $\chi_{T_t}^\sigma$ is the equivariant Euler characteristic over $U_\sigma$ and $\mathsf{S}(f)$ is the \emph{sum} of $f$ for any $f\in \Z[[M]]$ (see \cite[Section~1.3]{art:brionvergne1997} for the definitions of these notions). Since $\pi_k$ is an isomorphism over $X_k$ and similarly
\begin{equation}
L_{\alpha\beta}^\circ (t_1,t_2)=-\sum_{\sigma\in\Sigma_k(2) } \, \mathsf{S}\Big(\chi_{T_t}^\sigma\big(U_\sigma\,,\, {\Rcal^{\vec u_\beta{}'-\vec u_\alpha}}_{\vert U_\sigma}\big)\Big)\ ,
\end{equation}
we get
\begin{multline}\label{eq:perturbative-character}
\ch_{\tilde{T}}^{\mathrm{pert}} \Ebf_{([(\Ecal,\phi_{\Ecal})]\,,\,
  [(\Ecal',\phi_{\Ecal'})])} \\
\shoveleft{=\sum_{\alpha=1}^r \ \sum_{\beta=1}^{r'} \,e_\beta \,
  e_\alpha^{\prime\, -1}\,
  \mathsf{S}\Big(\chi_{T_t}^{\sigma_{\infty,k}}\big(U_{{\infty,k}} \,,\,
  {\pi_k}_\ast(\Rcal^{\vec u_\beta{}'-\vec u_\alpha}\otimes
  \Ocal_{\Xscr_k}(-\Dscr_\infty))_{\vert U_{{\infty,k}}}\big)\Big)} \\
+\, \sum_{\alpha=1}^r \ \sum_{\beta=1}^{r'} \, e_\beta \,e_\alpha^{\prime\, -1}\,
\mathsf{S}\Big(\chi_{T_t}^{\sigma_{\infty,0}}\big(U_{{\infty,0}} \,,\, {\pi_k}_\ast(\Rcal^{\vec u_\beta{}'-\vec u_\alpha}\otimes \Ocal_{\Xscr_k}(-\Dscr_\infty))_{\vert U_{{\infty,0}}}\big)\Big)\ .
\end{multline}

Set $\vec{u}_{\beta\alpha}:=\vec{u}_\beta{}'-\vec{u}_\alpha$ and
define $\vec{v}_{\beta\alpha}:=C^{-1}\vec{u}_{\beta\alpha}$. Let
$c_{\beta\alpha}\in\{0,1,\ldots, k-1\}$ be the equivalence class
modulo $k$ of $k\, (\vec{v}_{\beta\alpha})_{k-1}$. Set
$\vec{z}_{\beta\alpha}:=\vec{u}_{\beta\alpha}-\vec{e}_{c_{\beta\alpha}}$
for $c_{\beta\alpha}>0$ and $\vec{z}_{\beta\alpha}:=\vec{u}_{\beta\alpha}$ otherwise, where $\vec{e}_{c_{\beta\alpha}}$ is the $c_{\beta\alpha}$-th coordinate vector of $\Z^{k-1}$. Define $\vec{s}_{\beta\alpha}:=C^{-1}\vec{z}_{\beta\alpha}\in \Z^{k-1}$. Then 
\begin{equation}
\sum_{i=1}^{k-1}\, (\vec{u}_{\beta\alpha})_i\,
\omega_i=\sum_{j=1}^{k-1}\, (\vec{s}_{\beta\alpha})_j\, \Dscr_j+\omega_{c_{\beta\alpha}}\ ,
\end{equation}
where  we set $\omega_{c_{\beta\alpha}}=0$ if $c_{\beta\alpha}=0$. Therefore
\begin{multline}
\ch_{\tilde{T}}^{\mathrm{pert}} \Ebf_{([(\Ecal,\phi_{\Ecal})]\,,\,
  [(\Ecal',\phi_{\Ecal'})])} \\
\shoveleft{=\sum_{\alpha=1}^r \ \sum_{\beta=1}^{r'}\, e_\beta\,
  e_\alpha^{\prime\, -1}\,
  \mathsf{S}\Big(\chi_{T_t}^{\sigma_{\infty,k}}\big(U_{{\infty,k}}\,,\,
  {\pi_k}_\ast(\Rcal_{c_{\beta\alpha}}\otimes
  \Ocal_{\Xscr_k}(-\Dscr_\infty))_{\vert U_{{\infty,k}}}\big)\Big)} \\
+ \, \sum_{\alpha=1}^r \ \sum_{\beta=1}^{r'}\, e_\beta \,
e_\alpha^{\prime\, -1}\,
\mathsf{S}\Big(\chi_{T_t}^{\sigma_{\infty,0}}\big(U_{{\infty,0}} \,,\, {\pi_k}_\ast(\Rcal_{c_{\beta\alpha}}\otimes \Ocal_{\Xscr_k}(-\Dscr_\infty))_{\vert U_{{\infty,0}}}\big)\Big)\ .
\end{multline}
By using the relations \eqref{eq:tautologicalintegralrelations} and the projection formula for $\pi_k$ we get
\begin{align}
{\pi_k}_\ast\left(\Rcal_{c_{\beta\alpha}}\otimes
  \Ocal_{\Xscr_k}(-\Dscr_\infty)\right)_{\vert U_{{\infty,k}}} & \simeq 
\Ocal_{\bar{X}_k}\big(b_{\beta\alpha}^k \, D_k+ \lfloor
  b_{\beta\alpha}^\infty/k\rfloor\, D_\infty \big)_{\vert
  U_{{\infty,k}}}\ ,\\[4pt] 
{\pi_k}_\ast\left(\Rcal_{c_{\beta\alpha}}\otimes
  \Ocal_{\Xscr_k}(-\Dscr_\infty)\right)_{\vert U_{{\infty,0}}} & \simeq
\Ocal_{\bar{X}_k}\big(b_{\beta\alpha}^0 \, D_0+ \lfloor
b_{\beta\alpha}^\infty/k\rfloor \, D_\infty\big)_{\vert U_{{\infty,0}}}\ ,
\end{align}
where
\begin{align*}
b_{\beta\alpha}^0=&\left\{
\begin{array}{ll}
0 & \mbox{ for } \ c_{\beta\alpha}=0, k-1\ ,\\[5pt]
1 & \mbox{ for } \ 1\leq c_{\beta\alpha}\leq k-2 \ ,
\end{array} \right.
\\[4pt]
b_{\beta\alpha}^k=&\left\{
\begin{array}{ll}
0 & \mbox{ for } \ c_{\beta\alpha}=0\ ,\\[5pt]
c_{\beta\alpha}-1 & \mbox{ for } \ 1\leq c_{\beta\alpha}\leq k-2 \ ,\\[5pt]
1 & \mbox{ for } \ c_{\beta\alpha}=k-1\ ,
\end{array} \right.
\\[4pt]
b_{\beta\alpha}^\infty=&\left\{
\begin{array}{ll}
-1 & \mbox{ for } \ c_{\beta\alpha}=0\ ,\\[5pt]
\tilde{k}\, (c_{\beta\alpha}-2)-1 & \mbox{ for } \ 1\leq c_{\beta\alpha}\leq k-2\ , \\[5pt]
-\tilde{k} -1 & \mbox{ for } \ c_{\beta\alpha}=k-1 \ .
\end{array}
\right.
\end{align*}
By \cite[Section~2.3]{art:brionvergne1997} and the arguments in the proof of \cite[Theorem~4.1]{art:brionvergne1997} we can rewrite the perturbative character \eqref{eq:perturbative-character} as
\begin{multline}
\ch_{\tilde{T}}^{\mathrm{pert}} \Ebf_{([(\Ecal,\phi_{\Ecal})]\,,\,
  [(\Ecal',\phi_{\Ecal'})])} 
=\sum_{\alpha=1}^r \ \sum_{\beta=1}^{r'}\, e_\beta\, e_\alpha^{\prime\, -1}\
\frac{\sum_{m\in M\cap Q(\sigma_{\infty, k}\,;\, b_{\beta\alpha}^k \, D_k+
    \lfloor b_{\beta\alpha}^\infty/k \rfloor \, D_\infty )}\limits\,
  e^m}{\big(1-e^{m^1_{\infty, k}}\big)\, \big(1-e^{m^2_{\infty,
      k}}\big)}\\ + \, \sum_{\alpha=1}^r \ \sum_{\beta=1}^{r'}\,
e_\beta \, e_\alpha^{\prime\, -1} \ \frac{\sum_{m\in M\cap Q(\sigma_{\infty,
      0}\,;\, b_{\beta\alpha}^0 \, D_0+ \lfloor
    b_{\beta\alpha}^\infty/k\rfloor \, D_\infty)}\limits\,
  e^m}{\big(1-e^{m^1_{\infty, 0}} \big)\, \big(1-e^{m^2_{\infty, 0}} \big)}\ ,
\end{multline}
where $m^1_{\infty, k}$ and $m^2_{\infty, k}$ are the \emph{facet normals} to the rays of $\sigma_{\infty, k}$, i.e., the characters orthogonal to the rays of $\sigma_{\infty, k}$, and
\begin{multline}\label{eq:Qk}
Q\big(\sigma_{\infty, k}\,;\,b_{\beta\alpha}^k \, D_k+ \lfloor
b_{\beta\alpha}^\infty/k\rfloor\, D_\infty \big) \\
:=\Big\{x_1 \, m^1_{\infty, k}+x_2 \, m^2_{\infty, k}\ \Big\vert\
x_1,x_2\in \Q \ ,\quad 0\leq
x_1+\mbox{$\frac{b_{\beta\alpha}^k}{\tilde{k}}$} <1 \ , \quad 0\leq
x_2+\mbox{$\frac{\lfloor b_{\beta\alpha}^\infty/k\rfloor}{\tilde{k}}$}
<1\Big\}\ .
\end{multline}
Similarly, for the cone $\sigma_{\infty, 0}$ we denote by $m^1_{\infty, 0}$ and $m^2_{\infty, 0}$ the facet normals to its rays and
\begin{multline}\label{eq:Q0}
Q\big(\sigma_{\infty, 0}\,;\,b_{\beta\alpha}^0 \,D_0+ \left\lfloor
  b_{\beta\alpha}^\infty/k\right\rfloor\, D_\infty \big) \\
:=\Big\{x_1 \, m^1_{\infty, 0}+x_2 \, m^2_{\infty, 0}\ \Big\vert\
x_1,x_2\in \Q \ ,\quad 0\leq
x_1+\mbox{$\frac{b_{\beta\alpha}^0}{\tilde{k}}$} <1 \ , \quad 0\leq
x_2+\mbox{$\frac{\lfloor b_{\beta\alpha}^\infty/k\rfloor}{\tilde{k}}$}
<1\Big\}\ .
\end{multline}
Explicitly, for $\sigma_{\infty, k}$ we have 
\begin{equation*}
m^1_{\infty, k}=\left\{
\begin{array}{ll}
(1-\tilde{k}, - \tilde{k}) &\mbox{for }k\mbox{ even}\ ,\\
(2-k, - k) &\mbox{for }k\mbox{ odd}\ ,
\end{array}\right. \qquad \mbox{and} \qquad 
m^2_{\infty, k}=\left\{
\begin{array}{ll}
(1-2\tilde{k}, -2\tilde{k}) &\mbox{for }k\mbox{ even}\ ,\\
(1-k, -k) &\mbox{for }k\mbox{ odd}\ ,
\end{array}\right.
\end{equation*}
while for $\sigma_{\infty, 0}$ we get
\begin{equation*}
m^1_{\infty, 0}=\left\{
\begin{array}{ll}
(\tilde{k}-1, \tilde{k}) &\mbox{for }k\mbox{ even}\ ,\\
(k-2, k) &\mbox{for }k\mbox{ odd}\ ,
\end{array}\right. \qquad \mbox{and} \qquad 
m^2_{\infty, k}=\left\{
\begin{array}{ll}
(-1,0) &\mbox{for }k\mbox{ even}\ ,\\
(-1,0) &\mbox{for }k\mbox{ odd}\ .
\end{array}\right.
\end{equation*}
To determine the sets \eqref{eq:Qk} and \eqref{eq:Q0} explicitly we need
to consider separately two cases depending on the parity of $k$; they
each have cardinality $\tilde k$.

For even $k$ we get 
\small
\begin{multline}
Q\big(\sigma_{\infty, k}\,;\, b_{\beta\alpha}^k\, D_k+ \left\lfloor
  b_{\beta\alpha}^\infty/k \right\rfloor\, D_\infty \big) \\
\shoveleft{=\bigg\{\Big(i-b_{\beta\alpha}^k-\Big\lfloor\mbox{$\frac{b_{\beta\alpha}^k+\big\lfloor
    \frac{b_{\beta\alpha}^\infty}{k}\big\rfloor}{\tilde{k}} $} \Big\rfloor\,
  (1-2\tilde{k})\,,\, i-b_{\beta\alpha}^k+2\tilde{k}\,
  \Big\lfloor\mbox{$\frac{b_{\beta\alpha}^k+\big\lfloor
    \frac{b_{\beta\alpha}^\infty}{k}\big\rfloor}{\tilde{k}} $} \Big\rfloor\Big)\bigg\}_{i=0,1,\ldots,
    \tilde{k}\, \Big\{\frac{b_{\beta\alpha}^k+\big\lfloor
      \frac{b_{\beta\alpha}^\infty}{k} \big\rfloor}{\tilde{k}}\Big\}}}\\
\bigcup\\
\bigg\{\Big(i-b_{\beta\alpha}^k+\Big(1-\Big\lfloor \mbox{$\frac{b_{\beta\alpha}^k+\big\lfloor
  \frac{b_{\beta\alpha}^\infty}{k}\big\rfloor}{\tilde{k}}$} \Big\rfloor\Big)\,
(1-2\tilde{k})\,,\, i-b_{\beta\alpha}^k-2\tilde{k}\,
\Big(1-\Big\lfloor\mbox{$\frac{b_{\beta\alpha}^k+\big\lfloor
    \frac{b_{\beta\alpha}^\infty}{k}\big\rfloor}{\tilde{k}}$}
\Big\rfloor\Big)\Big)\bigg\}_{i=\tilde{k}\,\Big\{\frac{b_{\beta\alpha}^k+\big\lfloor
    \frac{b_{\beta\alpha}^\infty}{k} \big\rfloor}{\tilde{k}}\Big\}+1, \ldots, \tilde{k}-1}
\end{multline}
\normalsize
and
\small
\begin{equation}
Q\big(\sigma_{\infty, 0}\, ;\, b_{\beta\alpha}^0\, D_0+ \left\lfloor
  b_{\beta\alpha}^\infty/k \right\rfloor\, D_\infty
\big)=\bigg\{\Big(i-b_{\beta\alpha}^0+\Big\lfloor
\mbox{$\frac{b_{\beta\alpha}^0-i}{\tilde{k}}+\frac{1}{\tilde{k}}\,
\big\lfloor
\frac{b_{\beta\alpha}^\infty}{k}\big\rfloor$} \Big\rfloor\,,\,i-b_{\beta\alpha}^0\Big)
\bigg\}_{i=0, 1,\ldots, \tilde{k}-1}\ .
\end{equation}
\normalsize

On the other hand, for odd $k$ we have 
\small
\begin{multline}
Q\big(\sigma_{\infty, k}\, ;\, b_{\beta\alpha}^k\, D_k+ \left\lfloor
 b_{\beta\alpha}^\infty/k \right\rfloor\, D_\infty \big) \\
\shoveleft{=\bigg\{\Big(i-b_{\beta\alpha}^k-\Big\lfloor \mbox{$\frac{2b_{\beta\alpha}^k+\big\lfloor
    \frac{b_{\beta\alpha}^\infty}{k}\big\rfloor}{k} $} \Big\rfloor\, (1-k)
  \,,\, i-b_{\beta\alpha}^k+k\,
  \Big\lfloor\mbox{$ \frac{2b_{\beta\alpha}^k+\big\lfloor
    \frac{b_{\beta\alpha}^\infty}{k}\big\rfloor}{k}$} \Big\rfloor\Big)
  \bigg\}_{i=0,1, \ldots, \Big\lfloor\frac{k}{2}\,
    \Big\{\frac{2b_{\beta\alpha}^k+ \big\lfloor
      \frac{b_{\beta\alpha}^\infty}{k} \big\rfloor}{k}\Big\}\Big\rfloor}}\\
\bigcup\\
 \bigg\{\Big(i-b_{\beta\alpha}^k+\Big(1-\Big\lfloor \mbox{$\frac{2b_{\beta\alpha}^k+\big\lfloor
   \frac{b_{\beta\alpha}^\infty}{k}\big\rfloor}{k}$} \Big\rfloor\Big)\,
 (1-k)\,,\,
 i-b_{\beta\alpha}^k-k\,\Big(1-\Big\lfloor \mbox{$\frac{2b_{\beta\alpha}^k+\big\lfloor
   \frac{b_{\beta\alpha}^\infty}{k}\big\rfloor}{k}$} \Big\rfloor\Big)\Big)\bigg\}_{\stackrel{\scriptstyle
     i=\Big\lfloor\frac{k}{2}\,
   \Big\{\frac{2b_{\beta\alpha}^k+ \big\lfloor
     \frac{b_{\beta\alpha}^\infty}{k} \big\rfloor}{k}\Big\}\Big\rfloor +1,
   \ldots, }{\scriptstyle
   \Big\lfloor\frac{k}{2}+\frac{k}{2}\,
   \Big\{\frac{2b_{\beta\alpha}^k+ \big\lfloor
     \frac{b_{\beta\alpha}^\infty}{k} \big\rfloor}{k}\Big\}\Big\rfloor}
 } \\
 \bigcup\\
 \bigg\{\Big(i-b_{\beta\alpha}^k+\Big(2-\Big\lfloor\mbox{$\frac{2b_{\beta\alpha}^k+\big\lfloor
   \frac{b_{\beta\alpha}^\infty}{k}\big\rfloor}{k}$} \Big\rfloor\Big)\,
 (1-k)\,,\, i-b_{\beta\alpha}^k-k\,
 \Big(2-\Big\lfloor\mbox{$\frac{2b_{\beta\alpha}^k+\big\lfloor
   \frac{b_{\beta\alpha}^\infty}{k}\big\rfloor}{k} $} \Big\rfloor\Big)\Bigg)\bigg\}_{i=\Big\lfloor\frac{k}{2}+\frac{k}{2}\,
   \Big\{\frac{2b_{\beta\alpha}^k+\big\lfloor
     \frac{b_{\beta\alpha}^\infty}{k} \big\rfloor}{k}\Big\}\Big\rfloor+1, \ldots, k-1}
\end{multline}
\normalsize
and
\small
\begin{equation}
Q\big(\sigma_{\infty, 0}\, ;\, b_{\beta\alpha}^0\, D_0+ \left\lfloor
  b_{\beta\alpha}^\infty/k \right\rfloor\, D_\infty
\big)=\bigg\{\Big(i-b_{\beta\alpha}^0+\Big\lfloor \mbox{$
2\frac{b_{\beta\alpha}^0-i}{k}+\frac{1}{k}\, \big\lfloor
\frac{b_{\beta\alpha}^\infty}{k}\big\rfloor $} \Big\rfloor\,,\,i-b_{\beta\alpha}^0\Big)\bigg\}_{i=0,1, \ldots, k-1}\ .
\end{equation}
\normalsize

The corresponding Euler
classes involve infinite products which require regularization. 
 Following \cite[Appendix~A]{art:nekrasovokounkov2006}, we denote by $\Gamma_2(x|-\varepsilon_1,-\varepsilon_2) =
\exp\big(\gamma_{\varepsilon_1,\varepsilon_2}(x)\big)$ the Barnes double gamma-function
which is the double zeta-function
regularization of the infinite product
\begin{equation}
\prod_{i,j=0}^\infty \, \big(x-i\, \varepsilon_1-j\, \varepsilon_2
\big) \ .
\end{equation}

\begin{remark}
The function $\varepsilon_1\, \varepsilon_2\,
\gamma_{\varepsilon_1,\varepsilon_2}(x)$ is analytic near
$\varepsilon_1=\varepsilon_2=0$ with
\begin{equation}
\lim_{\varepsilon_1,\varepsilon_2\to 0}\, \varepsilon_1\, \varepsilon_2\,
\gamma_{\varepsilon_1,\varepsilon_2}(x) = -\frac12\, \log x+\frac34\, x^2 \ .
\end{equation}
This result appears in Section \ref{sec:SWgeometry}.
\end{remark}

\subsection{$\Ncal=2$ gauge theory}\label{sec:N2pure-pert}

Let $\vec{w}:=(w_0, w_1\ldots, w_{k-1})\in \N^k$ and
$r:=\sum_{l=0}^{k-1} \, w_l$. Let us define $\vec{c}:=(c_1, \ldots,
c_r)\in \{0, 1, \ldots, k-1\}^{r}$ such that $c_{\alpha}= i \bmod{k}$
if $\sum_{l=0}^{i-1}\, w_l<\alpha\leq \sum_{l=0}^{i}\, w_l$. For
$\alpha, \beta=1,\ldots, r$, $\alpha\neq \beta$, we define
$c_{\beta\alpha}$ to be the equivalence class modulo $k$ of
$c_\beta-c_\alpha$. Set $a_{\beta\alpha}=a_\beta-a_\alpha$ as before.

For even $k$ define
\small
\begin{multline}
F_{X_k}^{\mathrm{pert}}(\varepsilon_1, \varepsilon_2, \vec{a}\,;
\vec{c}\, )\\
\shoveleft{:=\sum_{ \alpha\neq \beta}\, \bigg(\,
  \sum_{i=0}^{\tilde{k}\,
    \Big\{\frac{b_{\beta\alpha}^k+\big\lfloor
        \frac{b_{\beta\alpha}^\infty}{k}\big\rfloor}{\tilde{k}}
    \Big\}}\gamma_{-\tilde{k}\, \varepsilon_1+\tilde{k}\,
    \varepsilon_2, 2\tilde{k}\,
    \varepsilon_2}\Big(a_{\beta\alpha} +(i-b_{\beta\alpha}^k)\,
  \varepsilon_1+\Big(i-b_{\beta\alpha}^k+2\tilde{k}\,
  \Big\lfloor\mbox{$\frac{b_{\beta\alpha}^k+\big\lfloor
      \frac{b_{\beta\alpha}^\infty}{k}\big\rfloor}{\tilde{k}} $}
  \Big\rfloor\Big)\, \varepsilon_2\Big)}\\[6pt]
+\sum_{i=\tilde{k}\, \Big\{\frac{b_{\beta\alpha}^k+\big\lfloor
    \frac{b_{\beta\alpha}^\infty}{k}\big\rfloor}{\tilde{k}}\Big\}+1}^{\tilde{k}-1}\,
\gamma_{-\tilde{k}\, \varepsilon_1+\tilde{k}\, \varepsilon_2,
  2\tilde{k}\, \varepsilon_2}\Big(a_{\beta\alpha}+
(i-b_{\beta\alpha}^k)\,
\varepsilon_1+\Big(i-b_{\beta\alpha}^k-2\tilde{k}\, \Big(1-\Big\lfloor\mbox{$
  \frac{b_{\beta\alpha}^k+\big\lfloor
    \frac{b_{\beta\alpha}^\infty}{k}\big\rfloor}{\tilde{k}}$}
\Big\rfloor\Big) \Big) \, \varepsilon_2\Big)\\[6pt]
+\, \sum_{i=0}^{\tilde{k}-1}\, \gamma_{\tilde{k}\,
  \varepsilon_1-\tilde{k}\, \varepsilon_2, 2\tilde{k}\,
  \varepsilon_1}\Big(a_{\beta\alpha} +\Big(i-b_{\beta\alpha}^0+k\, \Big\lfloor
\mbox{$ \frac{b_{\beta\alpha}^0-i}{\tilde{k}}+\frac{1}{\tilde{k}}\,
  \big\lfloor \frac{b_{\beta\alpha}^\infty}{k}\big\rfloor $}
\Big\rfloor\Big)\, \varepsilon_1+(i-b_{\beta\alpha}^0)\, \varepsilon_2\Big)\bigg)\ ,
\end{multline}
\normalsize
and for odd $k$
\small
\begin{multline}
F_{X_k}^{\mathrm{pert}}(\varepsilon_1, \varepsilon_2, \vec{a}\, ;
\vec{c}\, )\\
\shoveleft{:=\sum_{\alpha\neq \beta}\, \bigg(\,
  \sum_{i=0}^{\Big\lfloor\frac{k}{2}\,
    \Big\{\frac{2b_{\beta\alpha}^k+\big\lfloor
      \frac{b_{\beta\alpha}^\infty}{k}\big\rfloor}{k}\Big\}\Big\rfloor}\,
  \gamma_{-k\, \varepsilon_1+k\, \varepsilon_2, k\,
    \varepsilon_2}\Big(a_{\beta\alpha} +(i-b_{\beta\alpha}^k)\,
  \varepsilon_1+\Big(i-b_{\beta\alpha}^k+k\,
  \Big\lfloor\mbox{$\frac{2b_{\beta\alpha}^k+\big\lfloor
      \frac{b_{\beta\alpha}^\infty}{k}\big\rfloor}{k} $}
  \Big\rfloor\Big)\, \varepsilon_2\Big)}\\[6pt]
+\sum_{i=\Big\lfloor\frac{k}{2}\,
  \Big\{\frac{2b_{\beta\alpha}^k+\big\lfloor
    \frac{b_{\beta\alpha}^\infty}{k}\big\rfloor}{k}\Big\}\Big\rfloor
  +1}^{\Big\lfloor\frac{k}{2}+\frac{k}{2}\,
  \Big\{\frac{2b_{\beta\alpha}^k+\big\lfloor
    \frac{b_{\beta\alpha}^\infty}{k}\big\rfloor}{k}\Big\}\Big\rfloor}\,
\gamma_{-k\, \varepsilon_1+k\, \varepsilon_2, k\,
  \varepsilon_2}\Big(a_{\beta\alpha} +(i-b_{\beta\alpha}^k)\,
\varepsilon_1+\Big(i-b_{\beta\alpha}^k-k\, \Big(1-\Big\lfloor\mbox{$
  \frac{2b_{\beta\alpha}^k+\big\lfloor
    \frac{b_{\beta\alpha}^\infty}{k}\big\rfloor}{k}$}
\Big\rfloor\Big)\, \varepsilon_2\Big) \\[6pt]
+\sum_{i=\Big\lfloor\frac{k}{2}+\frac{k}{2}\,
  \Big\{\frac{2b_{\beta\alpha}^k+\big\lfloor
    \frac{b_{\beta\alpha}^\infty}{k}\big\rfloor}{k}\Big\}\Big\rfloor+1}^{k-1}\,
\gamma_{-k\, \varepsilon_1+k\, \varepsilon_2, k\,
  \varepsilon_2}\Big(a_{\beta\alpha} +(i-b_{\beta\alpha}^k)\,
\varepsilon_1+\Big(i-b_{\beta\alpha}^k-k\, \Big(2-\Big\lfloor\mbox{$
  \frac{2b_{\beta\alpha}^k+\big\lfloor
    \frac{b_{\beta\alpha}^\infty}{k}\big\rfloor}{k}$}
\Big\rfloor\Big) \Big)\, \varepsilon_2\Big) \\[6pt]
+\, \sum_{i=0}^{k-1}\, \gamma_{-k\, \varepsilon_1+k\, \varepsilon_2,
  k\, \varepsilon_1}\Big(a_{\beta\alpha} +\Big(i-b_{\beta\alpha}^0+k\,
\Big\lfloor \mbox{$ 2\frac{b_{\beta\alpha}^0-i}{k}+\frac{1}{k}\,
  \big\lfloor \frac{b_{\beta\alpha}^\infty}{k}\big\rfloor$}
\Big\rfloor\Big)\, \varepsilon_1+(i-b_{\beta\alpha}^0)\, \varepsilon_2\Big)\bigg)\ .
\end{multline}
\normalsize

\begin{definition}
The \emph{perturbative partition function} for pure $\Ncal=2$ gauge theory on $X_k$ is
\begin{equation}
\Zcal_{X_k}^{\mathrm{pert}}(\varepsilon_1, \varepsilon_2, \vec{a}\,
):=\sum_{\vec{c}}\, \exp\big(-F_{X_k}^{\mathrm{pert}}(\varepsilon_1,
\varepsilon_2, \vec{a}\, ; \vec{c}\, )\big) \ .
\end{equation}
\end{definition}

\begin{example}
For $k=2$ the perturbative partition function becomes
\begin{empheq}[box=\fbox]{multline}
\Zcal_{X_2}^{\mathrm{pert}}(\varepsilon_1, \varepsilon_2, \vec{a}\, ) \\
=\sum_{\vec{c}} \ \prod_{\alpha\neq\beta}\,
\exp\Big(-\gamma_{\varepsilon_2-\varepsilon_1,
  2\varepsilon_2}\big(a_{\beta\alpha}+c_{\beta\alpha}\, (\varepsilon_2-\varepsilon_1)-2\varepsilon_2\big)-
\gamma_{\varepsilon_1-\varepsilon_2,
  2\varepsilon_1}\big(a_{\beta\alpha} -2\varepsilon_1\big) \Big) \ .
\end{empheq}
\end{example}

\subsection{$\Ncal=2^\ast$ gauge theory}

With the same conventions as in Section \ref{sec:N2pure-pert}, define for even $k$
\small
\begin{multline}
F_{X_k}^{\ast, \mathrm{pert}}(\varepsilon_1, \varepsilon_2,
\vec{a}, \mu ; \vec{c}\, ) \\
\shoveleft{:=\sum_{\alpha\neq \beta}\,\Bigg(\, \sum_{i=0}^{\tilde{k}
    \Big\{\frac{b_{\beta\alpha}^k+\big\lfloor
      \frac{b_{\beta\alpha}^\infty}{k}\big\rfloor}{\tilde{k}}\Big\}}\, \bigg(
  \gamma_{-\tilde{k}\, \varepsilon_1+\tilde{k}\, \varepsilon_2,
    2\tilde{k}\,
    \varepsilon_2}\Big(a_{\beta\alpha} +(i-b_{\beta\alpha}^k)\,
  \varepsilon_1+\Big(i-b_{\beta\alpha}^k+2\tilde{k}\,\Big\lfloor\mbox{$\frac{b_{\beta\alpha}^k+\big\lfloor
      \frac{b_{\beta\alpha}^\infty}{k}\big\rfloor}{\tilde{k}}$}\Big\rfloor\Big)\,
  \varepsilon_2\Big)} \\
\shoveright{-\, \gamma_{-\tilde{k}\, \varepsilon_1+\tilde{k}\, \varepsilon_2,
    2\tilde{k}\,
    \varepsilon_2}\Big(\mu+a_{\beta\alpha} +(i-b_{\beta\alpha}^k)\,
  \varepsilon_1+\Big(i-b_{\beta\alpha}^k+2\tilde{k}\,\Big\lfloor\mbox{$\frac{b_{\beta\alpha}^k+\big\lfloor
      \frac{b_{\beta\alpha}^\infty}{k}\big\rfloor}{\tilde{k}}$}\Big\rfloor\Big)\,
  \varepsilon_2\Big) \bigg)}
\\[6pt]
\shoveleft{+\sum_{i=\tilde{k}\,\Big\{\frac{b_{\beta\alpha}^k+\big\lfloor
    \frac{b_{\beta\alpha}^\infty}{k}\big\rfloor}{\tilde{k}}\Big\}+1}^{\tilde{k}-1}\,
\bigg( \gamma_{-\tilde{k}\, \varepsilon_1+\tilde{k}\, \varepsilon_2,
  2\tilde{k}\, \varepsilon_2}\Big(a_{\beta\alpha} +
(i-b_{\beta\alpha}^k)\,
\varepsilon_1+\Big(i-b_{\beta\alpha}^k-2\tilde{k}\,
\Big(1-\Big\lfloor\mbox{$\frac{b_{\beta\alpha}^k+\big\lfloor
    \frac{b_{\beta\alpha}^\infty}{k}\big\rfloor}{\tilde{k}}$}
\Big\rfloor\Big) \Big) \, \varepsilon_2\Big)} \\
\shoveright{-\, \gamma_{-\tilde{k}\, \varepsilon_1+\tilde{k}\, \varepsilon_2,
  2\tilde{k}\, \varepsilon_2}\Big(\mu+ a_{\beta\alpha} +
(i-b_{\beta\alpha}^k)\,
\varepsilon_1+\Big(i-b_{\beta\alpha}^k-2\tilde{k}\,
\Big(1-\Big\lfloor\mbox{$\frac{b_{\beta\alpha}^k+\big\lfloor
    \frac{b_{\beta\alpha}^\infty}{k}\big\rfloor}{\tilde{k}}$}
\Big\rfloor\Big) \Big)\, \varepsilon_2\Big) \bigg)}
\\[6pt]
\shoveleft{+\, \sum_{i=0}^{\tilde{k}-1}\, \bigg(\gamma_{\tilde{k}\,
    \varepsilon_1-\tilde{k}\, \varepsilon_2, 2\tilde{k}\,
    \varepsilon_1}\Big(a_{\beta\alpha}+\Big(i-b_{\beta\alpha}^0+k\,
  \Big\lfloor
  \mbox{$\frac{b_{\beta\alpha}^0-i}{\tilde{k}}+\frac{1}{\tilde{k}}\big\lfloor
    \frac{b_{\beta\alpha}^\infty}{k}\big\rfloor$} \Big\rfloor\Big)\,
  \varepsilon_1+(i-b_{\beta\alpha}^0)\, \varepsilon_2\Big) }\\
-\, \gamma_{\tilde{k}\,
    \varepsilon_1-\tilde{k}\, \varepsilon_2, 2\tilde{k}\,
    \varepsilon_1}\Big(\mu+ a_{\beta\alpha}+\Big(i-b_{\beta\alpha}^0+k\,
  \Big\lfloor
  \mbox{$\frac{b_{\beta\alpha}^0-i}{\tilde{k}}+\frac{1}{\tilde{k}}\big\lfloor
    \frac{b_{\beta\alpha}^\infty}{k}\big\rfloor$} \Big\rfloor\Big)\,
  \varepsilon_1+(i-b_{\beta\alpha}^0)\, \varepsilon_2\Big) \bigg)
  \Bigg) \ ,
\end{multline}
\normalsize
and for odd $k$
\small
\begin{multline}
F_{X_k}^{\ast, \mathrm{pert}}(\varepsilon_1, \varepsilon_2,
\vec{a}, \mu; \vec{c}\, )\\
\shoveleft{:=\sum_{\alpha\neq \beta}\, \Bigg(\,
  \sum_{i=0}^{\Big\lfloor\frac{k}{2}\Big\{\frac{2b_{\beta\alpha}^k+\big\lfloor
      \frac{b_{\beta\alpha}^\infty}{k}\big\rfloor}{k}\Big\}\Big\rfloor}\,
  \bigg( \gamma_{-k\,\varepsilon_1+k\,\varepsilon_2, k\,
    \varepsilon_2}\Big(a_{\beta\alpha} +(i-b_{\beta\alpha}^k)\,
  \varepsilon_1+\Big(i-b_{\beta\alpha}^k+k\,
  \Big\lfloor\mbox{$\frac{2b_{\beta\alpha}^k+\big\lfloor
      \frac{b_{\beta\alpha}^\infty}{k}\big\rfloor}{k}$}\Big\rfloor\Big)\,
  \varepsilon_2\Big)} \\
\shoveright{-\, \gamma_{-k\,\varepsilon_1+k\,\varepsilon_2, k\,
    \varepsilon_2}\Big(\mu+ a_{\beta\alpha} +(i-b_{\beta\alpha}^k)\,
  \varepsilon_1+\Big(i-b_{\beta\alpha}^k+k\,
  \Big\lfloor\mbox{$\frac{2b_{\beta\alpha}^k+\big\lfloor
      \frac{b_{\beta\alpha}^\infty}{k}\big\rfloor}{k}$}\Big\rfloor\Big)\,
  \varepsilon_2\Big) \bigg)} \\[6pt]
\shoveleft{+\sum_{i=\Big\lfloor\frac{k}{2}\,
    \Big\{\frac{2b_{\beta\alpha}^k+\big\lfloor
      \frac{b_{\beta\alpha}^\infty}{k}\big\rfloor}{k}\Big\}\Big\rfloor
    +1}^{\Big\lfloor\frac{k}{2}+\frac{k}{2}\,
    \Big\{\frac{2b_{\beta\alpha}^k+\big\lfloor
      \frac{b_{\beta\alpha}^\infty}{k}\big\rfloor}{k}\Big\}\Big\rfloor}\,
  \bigg( \gamma_{-k\,\varepsilon_1+k\,\varepsilon_2,
    k\,\varepsilon_2}\Big(a_{\beta\alpha} +(i-b_{\beta\alpha}^k)\,
  \varepsilon_1+\Big(i-b_{\beta\alpha}^k-k\,
  \Big(1-\Big\lfloor\mbox{$\frac{2b_{\beta\alpha}^k+\big\lfloor
      \frac{b_{\beta\alpha}^\infty}{k}\big\rfloor}{k}$}
  \Big\rfloor\Big)\Big)\, \varepsilon_2\Big)} \\
\shoveright{-\, \gamma_{-k\,\varepsilon_1+k\,\varepsilon_2,
    k\,\varepsilon_2}\Big(\mu+ a_{\beta\alpha} +(i-b_{\beta\alpha}^k)\,
  \varepsilon_1+\Big(i-b_{\beta\alpha}^k-k\,
  \Big(1-\Big\lfloor\mbox{$\frac{2b_{\beta\alpha}^k+\big\lfloor
      \frac{b_{\beta\alpha}^\infty}{k}\big\rfloor}{k}$}
  \Big\rfloor\Big)\Big)\, \varepsilon_2\Big) \bigg)} \\[6pt]
\shoveleft{+\sum_{i=\Big\lfloor\frac{k}{2}+\frac{k}{2}\Big\{\frac{2b_{\beta\alpha}^k+\big\lfloor
      \frac{b_{\beta\alpha}^\infty}{k}\big\rfloor}{k}\Big\}\Big\rfloor+1}^{k-1}\,
  \bigg(\gamma_{-k\, \varepsilon_1+k\, \varepsilon_2, k\,
    \varepsilon_2}\Big(a_{\beta\alpha} +(i-b_{\beta\alpha}^k)\,
  \varepsilon_1+\Big(i-b_{\beta\alpha}^k-k\,
  \Big(2-\Big\lfloor\mbox{$\frac{2b_{\beta\alpha}^k+\big\lfloor
      \frac{b_{\beta\alpha}^\infty}{k}\big\rfloor}{k}$}
  \Big\rfloor\Big) \Big) \, \varepsilon_2\Big)} \\
\shoveright{-\, \gamma_{-k\, \varepsilon_1+k\, \varepsilon_2, k\,
    \varepsilon_2}\Big(\mu+a_{\beta\alpha} +(i-b_{\beta\alpha}^k)\,
  \varepsilon_1+\Big(i-b_{\beta\alpha}^k-k\,
  \Big(2-\Big\lfloor\mbox{$\frac{2b_{\beta\alpha}^k+\big\lfloor
      \frac{b_{\beta\alpha}^\infty}{k}\big\rfloor}{k}$}
  \Big\rfloor\Big) \Big) \, \varepsilon_2\Big) \bigg)} \\[6pt]
\shoveleft{+\,\sum_{i=0}^{k-1}\, \bigg( \gamma_{-k\, \varepsilon_1+k\,
    \varepsilon_2, k\, \varepsilon_1}\Big(a_{\beta\alpha}
  +\Big(i-b_{\beta\alpha}^0+k\Big\lfloor
  \mbox{$2\frac{b_{\beta\alpha}^0-i}{k}+\frac{1}{k}\big\lfloor
    \frac{b_{\beta\alpha}^\infty}{k}\big\rfloor$} \Big\rfloor\Big)\,
  \varepsilon_1+(i-b_{\beta\alpha}^0)\, \varepsilon_2\Big)} \\
-\, \gamma_{-k\, \varepsilon_1+k\,
    \varepsilon_2, k\, \varepsilon_1}\Big(\mu+a_{\beta\alpha}
  +\Big(i-b_{\beta\alpha}^0+k\Big\lfloor
  \mbox{$2\frac{b_{\beta\alpha}^0-i}{k}+\frac{1}{k}\big\lfloor
    \frac{b_{\beta\alpha}^\infty}{k}\big\rfloor$} \Big\rfloor\Big)\,
  \varepsilon_1+(i-b_{\beta\alpha}^0)\, \varepsilon_2\Big) \bigg)
  \Bigg) \ .
\end{multline}
\normalsize

\begin{definition}
The \emph{perturbative partition function} for $\Ncal=2$ gauge theory on $X_k$ with an adjoint hypermultiplet of mass $\mu$ is
\begin{equation}
\Zcal_{X_k}^{\ast, \mathrm{pert}}(\varepsilon_1, \varepsilon_2,
\vec{a}, \mu):=\sum_{\vec{c}}\, \exp\big(- F_{X_k}^{\ast,
  \mathrm{pert}}\big(\varepsilon_1, \varepsilon_2, \vec{a}, \mu;
\vec{c}\, ) \big)\ .
\end{equation}
\end{definition}

\begin{example}
For $k=2$ the perturbative partition function becomes
\begin{empheq}[box=\fbox]{multline}
\Zcal_{X_2}^{\ast, \mathrm{pert}}(\varepsilon_1,
\varepsilon_2, \vec{a}, \mu) \\
=\sum_{\vec{c}}\, \prod_{\alpha\neq\beta}\,
\frac{\exp\Big(\gamma_{\varepsilon_2-\varepsilon_1,
    2\varepsilon_2}\big(\mu+ a_{\beta\alpha} +c_{\beta\alpha}\,
  (\varepsilon_2-\varepsilon_1)-2\varepsilon_2\big)+
  \gamma_{\varepsilon_1-\varepsilon_2,
    2\varepsilon_1}\big(\mu+a_{\beta\alpha}-2\varepsilon_1\big)\Big)}{\exp\Big(\gamma_{\varepsilon_2-\varepsilon_1,
    2\varepsilon_2}\big(a_{\beta\alpha} +c_{\beta\alpha}\,
  (\varepsilon_2-\varepsilon_1)-2\varepsilon_2\big)+
  \gamma_{\varepsilon_1-\varepsilon_2,
    2\varepsilon_1}\big(a_{\beta\alpha}-2\varepsilon_1\big)\Big)} \ .
\end{empheq}
\end{example}

\subsection{Fundamental matter}\label{sec:perturbative-fundamentalmatter}

With the same notation as before, define for even $k$
\small
\begin{multline}
F_{X_k}^{N, \mathrm{pert}}(\varepsilon_1, \varepsilon_2, \vec{a},
\vec{\mu}\, ; \vec{c}\, )\\
\shoveleft{:=\sum_{ \alpha\neq \beta}\, \bigg(\,
  \sum_{i=0}^{\tilde{k}\,
    \Big\{\frac{b_{\beta\alpha}^k+\big\lfloor
        \frac{b_{\beta\alpha}^\infty}{k}\big\rfloor}{\tilde{k}}
    \Big\}}\gamma_{-\tilde{k}\, \varepsilon_1+\tilde{k}\,
    \varepsilon_2, 2\tilde{k}\,
    \varepsilon_2}\Big(a_{\beta\alpha} +(i-b_{\beta\alpha}^k)\,
  \varepsilon_1+\Big(i-b_{\beta\alpha}^k+2\tilde{k}\,
  \Big\lfloor\mbox{$\frac{b_{\beta\alpha}^k+\big\lfloor
      \frac{b_{\beta\alpha}^\infty}{k}\big\rfloor}{\tilde{k}} $}
  \Big\rfloor\Big)\, \varepsilon_2\Big)}\\[6pt]
+\sum_{i=\tilde{k}\, \Big\{\frac{b_{\beta\alpha}^k+\big\lfloor
    \frac{b_{\beta\alpha}^\infty}{k}\big\rfloor}{\tilde{k}}\Big\}+1}^{\tilde{k}-1}\,
\gamma_{-\tilde{k}\, \varepsilon_1+\tilde{k}\, \varepsilon_2,
  2\tilde{k}\, \varepsilon_2}\Big(a_{\beta\alpha}+
(i-b_{\beta\alpha}^k)\,
\varepsilon_1+\Big(i-b_{\beta\alpha}^k-2\tilde{k}\, \Big(1-\Big\lfloor\mbox{$
  \frac{b_{\beta\alpha}^k+\big\lfloor
    \frac{b_{\beta\alpha}^\infty}{k}\big\rfloor}{\tilde{k}}$}
\Big\rfloor\Big) \Big) \, \varepsilon_2\Big)\\[6pt]
\shoveright{+\, \sum_{i=0}^{\tilde{k}-1}\, \gamma_{\tilde{k}\,
  \varepsilon_1-\tilde{k}\, \varepsilon_2, 2\tilde{k}\,
  \varepsilon_1}\Big(a_{\beta\alpha} +\Big(i-b_{\beta\alpha}^0+k\, \Big\lfloor
\mbox{$ \frac{b_{\beta\alpha}^0-i}{\tilde{k}}+\frac{1}{\tilde{k}}\,
  \big\lfloor \frac{b_{\beta\alpha}^\infty}{k}\big\rfloor $}
\Big\rfloor\Big)\, \varepsilon_1+(i-b_{\beta\alpha}^0)\,
\varepsilon_2\Big) \bigg)} \\[6pt]
\shoveleft{-\sum_{s=1}^N \ \sum_{\alpha=1}^r \, \bigg(\,
  \sum_{i=0}^{\tilde{k}\,
    \Big\{\frac{b_{\alpha}^k+\big\lfloor
        \frac{b_{\alpha}^\infty}{k}\big\rfloor}{\tilde{k}}
    \Big\}}\gamma_{-\tilde{k}\, \varepsilon_1+\tilde{k}\,
    \varepsilon_2, 2\tilde{k}\,
    \varepsilon_2}\Big(\mu_s+ a_{\alpha} +(i-b_{\alpha}^k)\,
  \varepsilon_1+\Big(i-b_{\alpha}^k+2\tilde{k}\,
  \Big\lfloor\mbox{$\frac{b_{\alpha}^k+\big\lfloor
      \frac{b_{\alpha}^\infty}{k}\big\rfloor}{\tilde{k}} $}
  \Big\rfloor\Big)\, \varepsilon_2\Big)}\\[6pt]
+\sum_{i=\tilde{k}\, \Big\{\frac{b_{\alpha}^k+\big\lfloor
    \frac{b_{\alpha}^\infty}{k}\big\rfloor}{\tilde{k}}\Big\}+1}^{\tilde{k}-1}\,
\gamma_{-\tilde{k}\, \varepsilon_1+\tilde{k}\, \varepsilon_2,
  2\tilde{k}\, \varepsilon_2}\Big(\mu_s + a_{\alpha}+
(i-b_{\alpha}^k)\,
\varepsilon_1+\Big(i-b_{\alpha}^k-2\tilde{k}\, \Big(1-\Big\lfloor\mbox{$
  \frac{b_{\alpha}^k+\big\lfloor
    \frac{b_{\alpha}^\infty}{k}\big\rfloor}{\tilde{k}}$}
\Big\rfloor\Big) \Big) \, \varepsilon_2\Big)\\[6pt]
+\, \sum_{i=0}^{\tilde{k}-1}\, \gamma_{\tilde{k}\,
  \varepsilon_1-\tilde{k}\, \varepsilon_2, 2\tilde{k}\,
  \varepsilon_1}\Big(\mu_s+ a_{\alpha} +\Big(i-b_{\alpha}^0+k\, \Big\lfloor
\mbox{$ \frac{b_{\alpha}^0-i}{\tilde{k}}+\frac{1}{\tilde{k}}\,
  \big\lfloor \frac{b_{\alpha}^\infty}{k}\big\rfloor $}
\Big\rfloor\Big)\, \varepsilon_1+(i-b_{\alpha}^0)\,
\varepsilon_2\Big) \bigg) \ ,
\end{multline}
\normalsize
and for odd $k$
\small
\begin{multline}
F_{X_k}^{N, \mathrm{pert}}(\varepsilon_1, \varepsilon_2, \vec{a},
\vec{\mu}\, ; \vec{c}\, )\\
\shoveleft{:=\sum_{\alpha\neq \beta}\, \bigg(\,
  \sum_{i=0}^{\Big\lfloor\frac{k}{2}\,
    \Big\{\frac{2b_{\beta\alpha}^k+\big\lfloor
      \frac{b_{\beta\alpha}^\infty}{k}\big\rfloor}{k}\Big\}\Big\rfloor}\,
  \gamma_{-k\, \varepsilon_1+k\, \varepsilon_2, k\,
    \varepsilon_2}\Big(a_{\beta\alpha} +(i-b_{\beta\alpha}^k)\,
  \varepsilon_1+\Big(i-b_{\beta\alpha}^k+k\,
  \Big\lfloor\mbox{$\frac{2b_{\beta\alpha}^k+\big\lfloor
      \frac{b_{\beta\alpha}^\infty}{k}\big\rfloor}{k} $}
  \Big\rfloor\Big)\, \varepsilon_2\Big)}\\[6pt]
+\sum_{i=\Big\lfloor\frac{k}{2}\,
  \Big\{\frac{2b_{\beta\alpha}^k+\big\lfloor
    \frac{b_{\beta\alpha}^\infty}{k}\big\rfloor}{k}\Big\}\Big\rfloor
  +1}^{\Big\lfloor\frac{k}{2}+\frac{k}{2}\,
  \Big\{\frac{2b_{\beta\alpha}^k+\big\lfloor
    \frac{b_{\beta\alpha}^\infty}{k}\big\rfloor}{k}\Big\}\Big\rfloor}\,
\gamma_{-k\, \varepsilon_1+k\, \varepsilon_2, k\,
  \varepsilon_2}\Big(a_{\beta\alpha} +(i-b_{\beta\alpha}^k)\,
\varepsilon_1+\Big(i-b_{\beta\alpha}^k-k\, \Big(1-\Big\lfloor\mbox{$
  \frac{2b_{\beta\alpha}^k+\big\lfloor
    \frac{b_{\beta\alpha}^\infty}{k}\big\rfloor}{k}$}
\Big\rfloor\Big)\, \varepsilon_2\Big) \\[6pt]
+\sum_{i=\Big\lfloor\frac{k}{2}+\frac{k}{2}\,
  \Big\{\frac{2b_{\beta\alpha}^k+\big\lfloor
    \frac{b_{\beta\alpha}^\infty}{k}\big\rfloor}{k}\Big\}\Big\rfloor+1}^{k-1}\,
\gamma_{-k\, \varepsilon_1+k\, \varepsilon_2, k\,
  \varepsilon_2}\Big(a_{\beta\alpha} +(i-b_{\beta\alpha}^k)\,
\varepsilon_1+\Big(i-b_{\beta\alpha}^k-k\, \Big(2-\Big\lfloor\mbox{$
  \frac{2b_{\beta\alpha}^k+\big\lfloor
    \frac{b_{\beta\alpha}^\infty}{k}\big\rfloor}{k}$}
\Big\rfloor\Big) \Big)\, \varepsilon_2\Big) \\[6pt]
\shoveright{+\, \sum_{i=0}^{k-1}\, \gamma_{-k\, \varepsilon_1+k\, \varepsilon_2,
  k\, \varepsilon_1}\Big(a_{\beta\alpha} +\Big(i-b_{\beta\alpha}^0+k\,
\Big\lfloor \mbox{$ 2\frac{b_{\beta\alpha}^0-i}{k}+\frac{1}{k}\,
  \big\lfloor \frac{b_{\beta\alpha}^\infty}{k}\big\rfloor$}
\Big\rfloor\Big)\, \varepsilon_1+(i-b_{\beta\alpha}^0)\,
\varepsilon_2\Big)\bigg)} \\[6pt]
\shoveleft{-\sum_{s=1}^N \ \sum_{\alpha=1}^r \, \bigg(\,
  \sum_{i=0}^{\Big\lfloor\frac{k}{2}\,
    \Big\{\frac{2b_{\alpha}^k+\big\lfloor
      \frac{b_{\alpha}^\infty}{k}\big\rfloor}{k}\Big\}\Big\rfloor}\,
  \gamma_{-k\, \varepsilon_1+k\, \varepsilon_2, k\,
    \varepsilon_2}\Big(\mu_s+ a_{\alpha} +(i-b_{\alpha}^k)\,
  \varepsilon_1+\Big(i-b_{\alpha}^k+k\,
  \Big\lfloor\mbox{$\frac{2b_{\alpha}^k+\big\lfloor
      \frac{b_{\alpha}^\infty}{k}\big\rfloor}{k} $}
  \Big\rfloor\Big)\, \varepsilon_2\Big)}\\[6pt]
+\sum_{i=\Big\lfloor\frac{k}{2}\,
  \Big\{\frac{2b_{\alpha}^k+\big\lfloor
    \frac{b_{\alpha}^\infty}{k}\big\rfloor}{k}\Big\}\Big\rfloor
  +1}^{\Big\lfloor\frac{k}{2}+\frac{k}{2}\,
  \Big\{\frac{2b_{\alpha}^k+\big\lfloor
    \frac{b_{\alpha}^\infty}{k}\big\rfloor}{k}\Big\}\Big\rfloor}\,
\gamma_{-k\, \varepsilon_1+k\, \varepsilon_2, k\,
  \varepsilon_2}\Big(\mu_s+ a_{\alpha} +(i-b_{\alpha}^k)\,
\varepsilon_1+\Big(i-b_{\alpha}^k-k\, \Big(1-\Big\lfloor\mbox{$
  \frac{2b_{\alpha}^k+\big\lfloor
    \frac{b_{\alpha}^\infty}{k}\big\rfloor}{k}$}
\Big\rfloor\Big)\, \varepsilon_2\Big) \\[6pt]
+\sum_{i=\Big\lfloor\frac{k}{2}+\frac{k}{2}\,
  \Big\{\frac{2b_{\alpha}^k+\big\lfloor
    \frac{b_{\alpha}^\infty}{k}\big\rfloor}{k}\Big\}\Big\rfloor+1}^{k-1}\,
\gamma_{-k\, \varepsilon_1+k\, \varepsilon_2, k\,
  \varepsilon_2}\Big(\mu_s+ a_{\alpha} +(i-b_{\alpha}^k)\,
\varepsilon_1+\Big(i-b_{\alpha}^k-k\, \Big(2-\Big\lfloor\mbox{$
  \frac{2b_{\alpha}^k+\big\lfloor
    \frac{b_{\alpha}^\infty}{k}\big\rfloor}{k}$}
\Big\rfloor\Big) \Big)\, \varepsilon_2\Big) \\[6pt]
+\, \sum_{i=0}^{k-1}\, \gamma_{-k\, \varepsilon_1+k\, \varepsilon_2,
  k\, \varepsilon_1}\Big(\mu_s+ a_{\alpha} +\Big(i-b_{\alpha}^0+k\,
\Big\lfloor \mbox{$ 2\frac{b_{\alpha}^0-i}{k}+\frac{1}{k}\,
  \big\lfloor \frac{b_{\alpha}^\infty}{k}\big\rfloor$}
\Big\rfloor\Big)\, \varepsilon_1+(i-b_{\alpha}^0)\,
\varepsilon_2\Big)\bigg) \ .
\end{multline}
\normalsize

\begin{definition}
The \emph{perturbative partition function} for $\Ncal=2$ gauge theory on $X_k$ with $N$ fundamental hypermultiplets of masses $\mu_1,\dots,\mu_N$ is
\begin{equation}
\Zcal_{X_k}^{N, \mathrm{pert}}(\varepsilon_1, \varepsilon_2, \vec{a},
\vec{\mu}\, ):=\sum_{\vec{c}}\, \exp\big(- F_{X_k}^{N,
  \mathrm{pert}}(\varepsilon_1, \varepsilon_2, \vec{a}, \vec{\mu}\, ;
\vec{c}\, ) \big)\ .
\end{equation}
\end{definition}

\begin{example}
For $k=2$ the perturbative partition function becomes
\begin{empheq}[box=\fbox]{multline}
\Zcal_{X_2}^{N, \mathrm{pert}}(\varepsilon_1, \varepsilon_2, \vec{a},
\vec{\mu}\, ) \\
=\sum_{\vec{c}}\, \frac{\prod\limits_{s=1}^N \ \prod\limits_{\alpha=1}^r \,
  \exp\Big(\gamma_{\varepsilon_2-\varepsilon_1,
    2\varepsilon_2}\big(\mu_s+ a_{\alpha} +c_{\alpha}\,
  (\varepsilon_2-\varepsilon_1)-2\varepsilon_2\big)+
  \gamma_{\varepsilon_1-\varepsilon_2,
    2\varepsilon_1}\big(\mu_s+ a_{\alpha} -2\varepsilon_1\big)\Big)}{\prod\limits_{\alpha\neq\beta}\,
  \exp\Big(\gamma_{\varepsilon_2-\varepsilon_1,
    2\varepsilon_2}\big(a_{\beta\alpha} +c_{\beta\alpha}\,
  (\varepsilon_2-\varepsilon_1)-2\varepsilon_2\big)+
  \gamma_{\varepsilon_1-\varepsilon_2,
    2\varepsilon_1}\big(a_{\beta\alpha} -2\varepsilon_1\big)\Big)} \ .
\end{empheq}
\end{example}

\bigskip\section{Seiberg-Witten geometry\label{sec:SWgeometry}}

\subsection{Instanton prepotentials}\label{sec:instantonprepotential}

Define $F_{\C^2}^{\mathrm{inst}}(\varepsilon_1,\varepsilon_2,\vec{a};
\qsf):=- \varepsilon_1\, \varepsilon_2\, \log \Zcal_{\C^2}^{\mathrm{inst}}(\varepsilon_1, \varepsilon_2, \vec{a}; \qsf)$.
\begin{proposition}[{\cite[Proposition~7.3]{art:nakajimayoshioka2005-I}}]
$F_{\C^2}^{\mathrm{inst}}(\varepsilon_1,\varepsilon_2, \vec{a}; \qsf)$ is analytic in $\varepsilon_1,\varepsilon_2$ near $\varepsilon_1=\varepsilon_2=0$.
\end{proposition}
It is shown by
\cite{art:nakajimayoshioka2005-I,art:nekrasovokounkov2006} that the
instanton part of the \emph{Seiberg-Witten prepotential} of pure
$\Ncal=2$ gauge theory on $\R^4$ is given by
\begin{equation}
\Fcal_{\C^2}^{\mathrm{inst}}(\vec{a}\, ; \qsf):= \lim_{\varepsilon_1,\varepsilon_2\to 0}\, F_{\C^2}^{\mathrm{inst}}(\varepsilon_1,\varepsilon_2,\vec{a}; \qsf)\ .
\end{equation}
This result gives a proof of the \emph{Nekrasov conjecture} for
$\Ncal=2$ gauge theories without matter fields on $\R^4$.

In the case of $\Ncal=2$ gauge theory on the ALE space $X_k$, define
\begin{equation}
F_{X_k}^{\mathrm{inst}}\big(\varepsilon_1,\varepsilon_2,\vec{a}; \qsf,
\vec{\xi} \ \big):=- \tilde{k}\, \varepsilon_1\, \varepsilon_2\, \log
\Zcal_{X_k}^{\mathrm{inst}}\big(\varepsilon_1, \varepsilon_2, \vec{a};
\qsf, \vec{\xi} \ \big)\ .
\end{equation}
The following result is the proof of the analogue of the {Nekrasov
  conjecture} for gauge theory on $X_k$.
\begin{theorem}\label{thm:instantonprepotentialpure}
$F_{X_k}^{\mathrm{inst}}(\varepsilon_1,\varepsilon_2,\vec{a}; \qsf,
\vec{\xi} \ )$ is analytic in $\varepsilon_1,\varepsilon_2$ near $\varepsilon_1=\varepsilon_2=0$ and
\begin{equation}
\Fcal_{X_k}^{\mathrm{inst}}(\vec{a};
\qsf):=\lim_{\varepsilon_1,\varepsilon_2\to 0}\,
F_{X_k}^{\mathrm{inst}}\big(\varepsilon_1,\varepsilon_2,\vec{a};
\qsf, \vec{\xi} \ \big) = \frac{1}{k}\, \Fcal_{\C^2}^{\mathrm{inst}}(\vec{a}; \qsf)\ .
\end{equation}
\end{theorem}
\proof
The proof of the theorem follows the arguments in \cite[Section~5.5]{art:gasparimliu2010}. First note that
\begin{equation}\label{eq:limitblowup}
\prod_{i=1}^k\, \Zcal_{\C^2}^{\mathrm{inst}}\big(\varepsilon_1^{(i)},
\varepsilon_2^{(i)}, \vec{a}^{(i)}; \qsf \big)=\exp\Big(- 
\sum_{i=1}^k\, \frac{F_{\C^2}^{\mathrm{inst}}\big(\varepsilon_1^{(i)},
  \varepsilon_2^{(i)}, \vec{a}^{(i)}; \qsf
  \big)}{\varepsilon_1^{(i)}\, \varepsilon_2^{(i)}}\, \Big) \ .
\end{equation}
Moreover, one finds that
\begin{equation}\label{eq:limitedgefactor}
\lim_{\varepsilon_1,\varepsilon_2\to 0}\,
\frac{1}{\prod_{\alpha,\beta=1}^r\limits \ \prod_{n=1}^{k-1}\limits\,
  \ell^{(n)}_{\vec{v}_{\beta\alpha}}\big(\varepsilon_1^{(n)},
  \varepsilon_2^{(n)}, \vec{a}\, \big)}=\prod_{\alpha\neq \beta}\,
\Big(\, \frac{1}{a_\beta-a_\alpha}\, \Big)^{\frac{1}{2}\,
  (\vec{v}_{\beta\alpha}\cdot
  C\vec{v}_{\beta\alpha}-(C^{-1})^{c_{\beta\alpha} , c_{\beta\alpha}})}\ .
\end{equation}

Given $\vec{\boldsymbol{v}}=(\vec{v}_1,\ldots, \vec{v}_r)$, define
\begin{equation}
F_{\Xscr_k,\vec{\boldsymbol{v}}}^{\mathrm{inst}}(\varepsilon_1,\varepsilon_2,\vec{a}; \qsf):=
\sum_{i=1}^k \,
\frac{F_{\C^2}^{\mathrm{inst}}\big(\varepsilon_1^{(i)},\varepsilon_2^{(i)},
  \vec{a}^{(i)}; \qsf \big)}{\varepsilon_1^{(i)}\,
  \varepsilon_2^{(i)}} + \frac{F_{\C^2}^{\mathrm{inst}}(w,u', \vec{a};
  \qsf)}{w\, u'} + \frac{F_{\C^2}^{\mathrm{inst}}(-w,u'', \vec{a};
  \qsf)}{-w\, u''}\ ,
\end{equation}
where $u',u''$ are the normal weights at the two fixed points on
$\Dscr_\infty$, i.e., the weights of the $T_t$-action on the normal
bundle $\Ncal_{\Dscr_\infty/\Xscr_k}$ at the two points, while $w$ is the
tangent weight at the intersection point with $\Dscr_0$ (thus $-w$ is
the tangent weight at the intersection point with $\Dscr_k$). It
follows that $u'=-\varepsilon_1^{(1)}=-k\, \varepsilon_1$,
$u''=-\varepsilon_2^{(k)}=-k\, \varepsilon_2$, while $w=\tilde{k}\,
\varepsilon_2-\tilde{k}\, \varepsilon_1$.
Arguing along the lines of the proof of
\cite[Lemma~5.16]{art:gasparimliu2010}, we can write
$F_{\Xscr_k,\vec{\boldsymbol{v}}}^{\mathrm{inst}}(\varepsilon_1,\varepsilon_2,\vec{a};
\qsf)$ as an equivariant integral over $\Xscr_k$ of an analytic
function; this is enough to conclude that $F_{\Xscr_k,\vec{\boldsymbol{v}}}^{\mathrm{inst}}(\varepsilon_1,\varepsilon_2,\vec{a}; \qsf)$ is analytic in $\varepsilon_1,\varepsilon_2$ near $\varepsilon_1=\varepsilon_2=0$.
For this argument we use, in addition to the well-posed definition of the equivariant cohomology of a topological stack with an action of a Deligne-Mumford torus noted in Remark \ref{rem:equivariantcohomology}, the functoriality of this construction, in particular with respect to the coarse moduli space morphism \cite[Section 5]{art:ginotnoohi2012}.

By using Equations \eqref{eq:limitblowup} and \eqref{eq:limitedgefactor}, and arguing as in the proof of  \cite[Lemma~5.18]{art:gasparimliu2010}, we find that
\begin{equation}
\log\big(\Zcal_{X_k}^{\mathrm{inst}}(\varepsilon_1,\varepsilon_2,\vec{a};\qsf,\vec{\xi}
\ )\, \Zcal_{\C^2}^{\mathrm{inst}}(w,u',\vec{a};\qsf)\, \Zcal_{\C^2}^{\mathrm{inst}}(-w,u'',\vec{a};\qsf)\big)
\end{equation}
is also analytic in $\varepsilon_1,\varepsilon_2$ near
$\varepsilon_1=\varepsilon_2=0$. Thus the pole of
$\log\big(\Zcal_{X_k}^{\mathrm{inst}}(\varepsilon_1,\varepsilon_2,\vec{a};\qsf,\vec{\xi}
\ )\big)$ at $\varepsilon_1=\varepsilon_2=0$ is the same as the pole of
\begin{equation}
\frac{F_{\C^2}^{\mathrm{inst}}(w,u',\vec{a};\qsf)}{w\, u'} +
\frac{F_{\C^2}^{\mathrm{inst}}(-w,u'',\vec{a};\qsf)}{-w\, u''}\ .
\end{equation}

Now we follow the proof of \cite[Theorem~5.20-(a)]{art:gasparimliu2010}, obtaining that the function
\begin{equation}\label{eq:prepoALE}
-u'\,u''\,
\log\Zcal_{X_k}^{\mathrm{inst}}\big(\varepsilon_1,\varepsilon_2,\vec{a};\qsf,\vec{\xi}
\ ) = \mbox{$\frac{k^2}{\tilde k} $} \, F_{X_k}^{\mathrm{inst}}
\big(\varepsilon_1,\varepsilon_2,\vec{a};\qsf,\vec{\xi} \ \big)
\end{equation}
is analytic in $\varepsilon_1,\varepsilon_2$ near
$\varepsilon_1=\varepsilon_2=0$. By arguing along the lines of the
proof of \cite[Theorem~5.20-(b)]{art:gasparimliu2010} and using the
identity $\tilde{k}\, u' - \tilde{k}\, u''=k\, w$, we find that its limit as $\varepsilon_1,\varepsilon_2\to 0$ is 
\begin{equation}
\mbox{$\frac{k}{\tilde{k}}$} \, \Fcal_{\C^2}^{\mathrm{inst}}(\vec{a};\qsf)\ .
\end{equation}
Multiplying both sides of Equation \eqref{eq:prepoALE} by
$\tilde{k}\,\varepsilon_1\, \varepsilon_2/u'\, u'' = \tilde{k}/k^2$
then implies the statements.
\endproof
\begin{definition}
We call $\Fcal_{X_k}^{\mathrm{inst}}(\vec{a}; \qsf)$ the \emph{instanton prepotential} of pure $\Ncal=2$ gauge theory on $X_k$.
\end{definition}
In addition, one can define
\begin{equation}
F_{X_k}^{\circ}\big(\varepsilon_1,\varepsilon_2, \vec{a}; \qsf,
\tau_1, \vec{\xi} \ \big):=-\tilde{k}\, \varepsilon_1\,
\varepsilon_2\, \log \Zcal_{X_k}^{\circ}\big(\varepsilon_1,
\varepsilon_2, \vec{a}; \qsf, \tau_1, \vec{\xi} \ \big)
\end{equation}
and
\begin{equation}
F_{X_k}^{\circ,\mathrm{inst}}\big(\varepsilon_1,\varepsilon_2, \vec{a};
\qsf, \tau_1, \vec{\xi} \ \big):=- \tilde{k}\, \varepsilon_1\,
\varepsilon_2\, \log \Zcal_{X_k}^{\circ, \mathrm{inst}}\big(\varepsilon_1,
\varepsilon_2, \vec{a}; \qsf, \tau_1, \vec{\xi} \ \big)\ .
\end{equation}
By using Equation \eqref{eq:deformedinstanton} one gets
\begin{equation}
F_{X_k}^{\circ}\big(\varepsilon_1,\varepsilon_2, \vec{a}; \qsf,
\tau_1, \vec{\xi} \ \big)=\frac{\tilde{k}}{2k}\, \tau_1\,
\sum_{\alpha=1}^r\, a_\alpha^2+F_{X_k}^{\circ,
  \mathrm{inst}}\big(\varepsilon_1,\varepsilon_2, \vec{a};
\qsf_{\mathrm{eff}}, \tau_1, \vec{\xi} \ \big)\ .
\end{equation}
\begin{corollary}\label{cor:classinst}
$F_{X_k}^{\circ, \mathrm{inst}}(\varepsilon_1,\varepsilon_2, \vec{a};
\qsf_{\mathrm{eff}}, \tau_1, \vec{\xi} \ )$ and
$F_{X_k}^{\circ}(\varepsilon_1,\varepsilon_2, \vec{a}; \qsf, \tau_1,
\vec{\xi} \ )$ are analytic in $\varepsilon_1,\varepsilon_2$ near
$\varepsilon_1=\varepsilon_2=0$, and
\begin{align}
\lim_{\varepsilon_1,\varepsilon_2\to 0}\,
F_{X_k}^{\circ}\big(\varepsilon_1,\varepsilon_2, \vec{a}; \qsf,
\tau_1, \vec{\xi} \ \big)&=\frac{\tilde{k}}{2k}\, \tau_1
\,\sum_{\alpha=1}^r\, a_\alpha^2+\lim_{\varepsilon_1,\varepsilon_2\to
  0}\, F_{X_k}^{\circ, \mathrm{inst}}\big(\varepsilon_1,\varepsilon_2,
\vec{a}; \qsf_{\mathrm{eff}}, \tau_1, \vec{\xi} \ \big)\\[4pt]
&=\frac1k\, \Big( \, \frac{\tilde{k}\, \tau_1}{2}\, \sum_{\alpha=1}^r\,
a_\alpha^2+ \Fcal_{\C^2}^{\mathrm{inst}}(\vec{a};\qsf_{\mathrm{eff}}) \,
\Big) \ .
\end{align}
\end{corollary}

Now we deal with gauge theories involving matter fields. Define 
\begin{align}
F_{\C^2}^{\ast, \mathrm{inst}}(\varepsilon_1,\varepsilon_2,\vec{a}, \mu;
\qsf) &:=-\varepsilon_1\, \varepsilon_2\, \log \Zcal_{\C^2}^{\ast, \mathrm{inst}}(\varepsilon_1, \varepsilon_2, \vec{a}, \mu; \qsf)\ ,\\[4pt]
F_{\C^2}^{N, \mathrm{inst}}(\varepsilon_1,\varepsilon_2,\vec{a},
\vec{\mu}; \qsf) &:=- \varepsilon_1\, \varepsilon_2\, \log \Zcal_{\C^2}^{N, \mathrm{inst}}(\varepsilon_1, \varepsilon_2, \vec{a}, \vec{\mu}; \qsf)\ .
\end{align}
\begin{proposition}[{\cite{art:nekrasovokounkov2006}}]
$F_{\C^2}^{\ast, \mathrm{inst}}(\varepsilon_1,\varepsilon_2,\vec{a}, \mu; \qsf)$ and $F_{\C^2}^{N, \mathrm{inst}}(\varepsilon_1,\varepsilon_2,\vec{a}, \vec{\mu}; \qsf)$ are analytic in $\varepsilon_1,\varepsilon_2$ near $\varepsilon_1=\varepsilon_2=0$.
\end{proposition}

In \cite{art:nekrasovokounkov2006} it is shown that the {instanton
  part of the Seiberg-Witten prepotential} of $\Ncal=2^\ast$ gauge
theory on $\R^4$ is given by
\begin{equation}
\Fcal_{\C^2}^{\ast, \mathrm{inst}}(\vec{a}, \mu; \qsf):=
\lim_{\varepsilon_1,\varepsilon_2\to 0} \, F_{\C^2}^{\ast, \mathrm{inst}}(\varepsilon_1,\varepsilon_2,\vec{a}, \mu; \qsf)\ ,
\end{equation}
while the {instanton part of the Seiberg-Witten prepotential} of
$\Ncal=2$ gauge theory on $\R^4$ with $N$ fundamental hypermultiplets
is given by
\begin{equation}
\Fcal_{\C^2}^{N, \mathrm{inst}}(\vec{a}, \vec{\mu}; \qsf):=
\lim_{\varepsilon_1,\varepsilon_2\to 0} \, F_{\C^2}^{N, \mathrm{inst}}(\varepsilon_1,\varepsilon_2,\vec{a}, \vec{\mu}; \qsf)\ .
\end{equation}

As shown in \cite[Section~5.5]{art:gasparimliu2010}, the techniques
used for the pure gauge theory work also for gauge theories with matter. Thus we define
\begin{align}
F_{X_k}^{\ast, \mathrm{inst}}\big(\varepsilon_1,\varepsilon_2,\vec{a},
\mu; \qsf, \vec{\xi} \ \big) & :=- \tilde{k}\, \varepsilon_1\,
\varepsilon_2\, \log \Zcal_{X_k}^{\ast,
  \mathrm{inst}}\big(\varepsilon_1, \varepsilon_2, \vec{a}, \mu; \qsf,
\vec{\xi} \ \big)\ ,\\[4pt]
F_{X_k}^{N, \mathrm{inst}}\big(\varepsilon_1,\varepsilon_2,\vec{a},
\vec{\mu}; \qsf, \vec{\xi} \ \big) & :=- \tilde{k}\, \varepsilon_1\,
\varepsilon_2\, \log \Zcal_{X_k}^{N, \mathrm{inst}}\big(\varepsilon_1,
\varepsilon_2, \vec{a}, \vec{\mu}; \qsf, \vec{\xi} \ \big)\ .
\end{align}
\begin{theorem}\label{thm:instantonprepotential-masses}
$ $
\begin{itemize}
\item[(1)] {$\Ncal=2^\ast$ gauge theory:} \ $F_{X_k}^{\ast,
    \mathrm{inst}}(\varepsilon_1,\varepsilon_2,\vec{a}, \mu ; \qsf,
  \vec{\xi} \ )$ is analytic in $\varepsilon_1,\varepsilon_2$ near $\varepsilon_1=\varepsilon_2=0$ and
\begin{equation}
\Fcal_{X_k}^{\ast, \mathrm{inst}}(\vec{a}, \mu;
\qsf):=\lim_{\varepsilon_1,\varepsilon_2\to 0}\, F_{X_k}^{\ast,
  \mathrm{inst}}\big(\varepsilon_1,\varepsilon_2,\vec{a}, \mu; \qsf,
\vec{\xi} \ \big) = \frac{1}{k}\, \Fcal_{\C^2}^{\ast, \mathrm{inst}}(\vec{a}, \mu; \qsf)\ .
\end{equation}
\item[(2)] {Fundamental matter:} \ $F_{X_k}^{N,
    \mathrm{inst}}(\varepsilon_1,\varepsilon_2,\vec{a}, \vec{\mu};
  \qsf, \vec{\xi} \ )$ is analytic in $\varepsilon_1,\varepsilon_2$ near $\varepsilon_1=\varepsilon_2=0$ and
\begin{equation}
\Fcal_{X_k}^{N, \mathrm{inst}}(\vec{a}, \vec{\mu};
\qsf):=\lim_{\varepsilon_1,\varepsilon_2\to 0}\, F_{X_k}^{N,
  \mathrm{inst}}\big(\varepsilon_1,\varepsilon_2,\vec{a}, \vec{\mu};
\qsf, \vec{\xi} \ \big)= \frac{1}{k} \, \Fcal_{\C^2}^{N,
  \mathrm{inst}}(\vec{a}, \vec{\mu}; \qsf) \ .
\end{equation}
\end{itemize}
\end{theorem}
\begin{remark}
Statements analogous to Corollary \ref{cor:classinst} hold as well in
the case of gauge theories with matter fields.
\end{remark}

\subsection{Perturbative prepotentials}\label{sec:perturbativeprepotential}

For any holonomy vector at infinity $\vec c\in\{0,1,\ldots,k-1\}^r$ we have the
  following results.

\begin{theorem} \label{thm:perturbativeprepotential}
$ $
\begin{enumerate}
\item $\Ncal=2$ gauge theory:
\begin{equation}
\lim_{\varepsilon_1,\varepsilon_2\to 0} \, \varepsilon_1 \,
\varepsilon_2 \, F_{X_k}^{\mathrm{pert}}(\varepsilon_1, \varepsilon_2,
\vec{a}; \vec{c}\, )=\frac{1}{k} \,
\Fcal_{\C^2}^{\mathrm{pert}}(\vec{a}\, )\ ,
\end{equation}
where 
\begin{equation}
\Fcal_{\C^2}^{\mathrm{pert}}(\vec{a}):=\sum_{ \alpha\neq \beta}\,
\Big(-\frac{1}{2}\, (a_\alpha-a_\beta)^2\, \log
(a_\alpha-a_\beta)+\frac{3}{4}\, (a_\alpha-a_\beta)^2\Big)
\end{equation}
is the perturbative part of the Seiberg-Witten prepotential of pure $\Ncal=2$ gauge theory on $\R^4$.
\item $\Ncal=2^\ast$ gauge theory:
\begin{equation}
\lim_{\varepsilon_1,\varepsilon_2\to 0} \, \varepsilon_1 \,
\varepsilon_2 \, F_{X_k}^{\ast, \mathrm{pert}}(\varepsilon_1,
\varepsilon_2, \vec{a}, m; \vec{c}\,)=\frac{1}{k} \, \Fcal_{\C^2}^{\ast, \mathrm{pert}}(\vec{a}, \mu)\ ,
\end{equation}
where 
\begin{multline}
\Fcal_{\C^2}^{\ast, \mathrm{pert}}(\vec{a}, \mu) :=  \sum_{ \alpha\neq
  \beta} \, \Big(-\frac{1}{2}\, (a_\alpha-a_\beta)^2\, \log
(a_\alpha-a_\beta)+\frac{3}{4}\, (a_\alpha-a_\beta)^2\\
 +\, \frac{1}{2}\, (a_\alpha-a_\beta+\mu)^2\, \log
 (a_\alpha-a_\beta+\mu)-\frac{3}{4}\, (a_\alpha-a_\beta+\mu)^2\Big)
\end{multline}
is the perturbative part of the Seiberg-Witten prepotential of $\Ncal=2^\ast$ gauge theory on $\R^4$.
\item Fundamental matter: 
\begin{equation}
\lim_{\varepsilon_1,\varepsilon_2\to 0} \, \varepsilon_1 \,
\varepsilon_2 \, F_{X_k}^{N, \mathrm{pert}}(\varepsilon_1,
\varepsilon_2, \vec{a}, \vec{\mu}; \vec{c}\, )= \frac{1}{k} \, \Fcal_{\C^2}^{N, \mathrm{pert}}(\vec{a}, \vec{\mu})\ ,
\end{equation}
where 
\begin{multline}
\Fcal_{\C^2}^{N, \mathrm{pert}}(\vec{a}, \vec{\mu}):=  \sum_{
  \alpha\neq \beta}\, \Big(-\frac{1}{2}\, (a_\alpha-a_\beta)^2\, \log
(a_\alpha-a_\beta)+\frac{3}{4}\, (a_\alpha-a_\beta)^2\Big) \\
 +\, \sum_{s=1}^N \ \sum_{\alpha=1}^r\, \Big(\, \frac{1}{2}\,
 (a_\alpha+\mu_s)^2\, \log (a_\alpha+\mu_s)-\frac{3}{4}\,
 (a_\alpha+\mu_s)^2\, \Big)
\end{multline}
is the perturbative part of the Seiberg-Witten prepotential of $\Ncal=2$ gauge theory on $\R^4$ with $N$ fundamental hypermultiplets.
\end{enumerate}
\end{theorem}
\proof
We prove only (1). The proofs of (2) and (3) are analogous.

Let $k$ be even. Then
\begin{multline}
\lim_{\varepsilon_1,\varepsilon_2\to 0} \, \varepsilon_1\,
\varepsilon_2 \, F_{X_k}^{\mathrm{pert}}(\varepsilon_1, \varepsilon_2,
\vec{a}; \vec{c}\, ) \\
=\lim_{\varepsilon_1,\varepsilon_2\to 0} \ \sum_{\alpha\neq \beta} \,
\Big( \, \tilde{k}\, \varepsilon_1\, \varepsilon_2 \,
\gamma_{-\tilde{k}\, \varepsilon_1+\tilde{k}\, \varepsilon_2,
  2\tilde{k}\, \varepsilon_2}(a_\beta-a_\alpha)+\tilde{k}\,
\varepsilon_1\, \varepsilon_2 \, \gamma_{\tilde{k}\,
  \varepsilon_1-\tilde{k}\, \varepsilon_2, 2\tilde{k}\,
  \varepsilon_1}(a_\beta-a_\alpha)\, \Big) \ .
\end{multline}
By arguing as in the proof of \cite[Theorem 6.8]{art:gasparimliu2010}
one gets the assertion. The case of odd $k$ is similar.
\endproof

\subsection{Blowup equations}\label{sec:blowupequations}

Throughout this subsection we set $k=2$.

\subsubsection*{$\Ncal=2$ gauge theory}

The generating function for correlators of quadratic $2$-observables is defined by
\begin{equation}
\Zcal_{X_2}^\bullet(\varepsilon_1, \varepsilon_2, \vec{a}; \qsf,
{\xi}, t):=\Zcal_{X_2}\big(\varepsilon_1,\varepsilon_2,
\vec{a}; \qsf, \vec{\xi}, \vec{\tau}, \vec{t}\ \big)\ ,
\end{equation}
where $\vec{\tau}=\vec 0$ and $\vec{t}:=(0, -t, 0,
\ldots)$. By Example \ref{ex:k2case}, this partition function becomes
\begin{multline}\label{eq:deformed-t-k2}
\Zcal_{X_2}^\bullet(\varepsilon_1, \varepsilon_2, \vec{a}; \qsf,
\xi, t)= 
\sum_{\stackrel{\scriptstyle v\in\frac{1}{2}\, \Z }{ 2v=
    w_1\bmod{2}}}\, \xi^{v} \ \sum_{\vec{\boldsymbol{v}}=(v_1, \ldots,
  v_r)} \ \frac{\qsf^{\sum_{\alpha=1}^r\limits \,
    v_\alpha^2}}{\prod_{\alpha,\beta=1}^r\limits\,
  \ell_{v_{\beta\alpha}}(2\varepsilon_1, \varepsilon_2-\varepsilon_1,
  a_{\beta\alpha})} \\
\times \ \Zcal_{\C^2}\big(\varepsilon_1^{(1)},\varepsilon_2^{(1)},
\vec{a}^{(1)}; \qsf,\varepsilon_1^{(1)}\, \vec{t} \ \big) \,
\Zcal_{\C^2}\big(\varepsilon_1^{(2)},\varepsilon_2^{(2)},
\vec{a}^{(2)}; \qsf,\varepsilon_2^{(2)}\, \vec{t} \ \big)\ ,
\end{multline}
where $\Zcal_{\C^2}(\varepsilon_1,\varepsilon_2, \vec{a}; \qsf,
\vec\tau)$ is the deformed Nekrasov partition function for $U(r)$ gauge theory on $\R^4$. For $i=1,2$ we have
\begin{multline}
\Zcal_{\C^2}\big(\varepsilon_1^{(i)},\varepsilon_2^{(i)},
\vec{a}^{(i)}; \qsf,\varepsilon_i^{(i)}\, \vec{t} \ \big) \\
=\sum_{\vec{Y}} \, \frac{\qsf^{\sum_{\alpha=1}^r\limits \,\vert
    Y_\alpha\vert}}{\prod_{\alpha,\beta=1}^r\limits\, m_{Y_\alpha,
    Y_\beta}\big(\varepsilon_1^{(i)},\varepsilon_2^{(i)},a_{\beta\alpha}^{(i)}
  \big)} \,
\exp\bigg(-\varepsilon_i^{(i)}\, t\, \Big(\,
\frac{1}{2\varepsilon_1^{(i)}\, \varepsilon_2^{(i)}}\,
\sum_{\alpha=1}^r \, \big(a_\alpha^{(i)} \big)^2-\sum_{\alpha=1}^r \,
\vert Y_\alpha\vert\, \Big)\bigg)\ .
\end{multline}
Since 
\begin{equation}
\frac{1}{2\varepsilon_2^{(1)}}\, \sum_{\alpha=1}^r\,
\big(a_\alpha^{(1)} \big)^2+\frac{1}{2\varepsilon_1^{(2)}}\,
\sum_{\alpha=1}^r\, \big(a_\alpha^{(2)}
\big)^2=\big(\varepsilon_1^{(1)}+\varepsilon_2^{(2)}\big)\,
\sum_{\alpha=1}^r \, v_\alpha^2-2\, \sum_{\alpha=1}^r \, v_\alpha \, a_\alpha\ ,
\end{equation}
Equation \eqref{eq:deformed-t-k2} becomes
\begin{empheq}[box=\fbox]{multline}
\Zcal_{X_2}^\bullet(\varepsilon_1, \varepsilon_2, \vec{a}; \qsf, \xi, t)= 
\sum_{\stackrel{\scriptstyle v\in\frac{1}{2}\, \Z }{\scriptstyle 2v=
    w_1\bmod{2}}} \, \xi^{v} \ \sum_{\vec{\boldsymbol{v}}=(v_1,
  \ldots, v_r)} \, \frac{\big(\qsf
  \,\e^{-2t\,(\varepsilon_1+\varepsilon_2)}
  \big)^{\sum_{\alpha=1}^r\limits\,  v_\alpha^2} \, \e^{2t\,
    \sum_{\alpha=1}^r\limits \, v_\alpha \,
    a_\alpha}}{\prod_{\alpha,\beta=1}^r\limits\,
  \ell_{v_{\beta\alpha}}(2\varepsilon_1, \varepsilon_2-\varepsilon_1,
a_{\beta\alpha})} \\
\times \
\Zcal_{\C^2}^{\mathrm{inst}}\big(2\varepsilon_1,\varepsilon_2-\varepsilon_1,
\vec{a}+2\varepsilon_1\, \vec{\boldsymbol{v}}; \qsf \,
\e^{2\varepsilon_1 \, t}\big)\,
\Zcal_{\C^2}^{\mathrm{inst}}(\varepsilon_1-
\varepsilon_2,2\varepsilon_2 , \vec{a}+2\varepsilon_2\,
\vec{\boldsymbol{v}}; \qsf \, \e^{2\varepsilon_2\, t})\ .
\end{empheq}

We will now compare the generating function $
\Zcal_{X_2}^\bullet(\varepsilon_1, \varepsilon_2, \vec{a}; \qsf, \xi,
t)$ with the original partition function
$\Zcal_{X_2}^{\mathrm{inst}}(\varepsilon_1, \varepsilon_2, \vec{a}; \qsf,
\xi)$ in the {\em low energy limit}, i.e., in the limit where $\varepsilon_1,\varepsilon_2\to 0$, and show that they are related through \emph{theta-functions} on
the genus $r$ Seiberg-Witten curve $\Sigma$ for pure $\Ncal=2$ gauge theory on
$\R^4$. The period matrix $\tau=(\tau_{\alpha\beta})$ of $\Sigma$ is
given by 
\begin{align}
\tau_{\alpha\beta}:= &\frac{1}{2\pi \operatorname{i}}\,
\frac{\partial^2 \Fcal_{\C^2}}{\partial a_\alpha \, \partial
  a_\beta}(\vec{a}; \qsf) \\
=& \frac{1}{2\pi \operatorname{i}} \, \frac{\partial^2
  \Fcal_{\C^2}^{\mathrm{inst}}}{\partial a_\alpha \, \partial
  a_\beta}(\vec{a}; \qsf)-\frac{1}{2\pi \operatorname{i}}\,
\sum_{\alpha'<\beta'}\,
(\delta_{\alpha\alpha'}-\delta_{\alpha\beta'})\,
(\delta_{\beta\alpha'}-\delta_{\beta\beta'})\, \log
(a_{\beta'}-a_{\alpha'})+\frac{\log \qsf}{2\pi \operatorname{i}}\, \delta_{\alpha\beta}\ ,
\end{align}
where $\Fcal_{\C^2}(\vec{a}; \qsf)$ is the total Seiberg-Witten
prepotential defined by Corollary \ref{cor:classinst} and Theorem
\ref{thm:perturbativeprepotential} (here we identify $\qsf:=\e^{2\pi\ii\tau_1}$). We
further define
\begin{equation}
\zeta_\alpha:=-\frac{t}{2\pi
  \operatorname{i}} \, \Big( a_\alpha + \qsf\, 
\frac{\partial^2\Fcal_{\C^2}^{\mathrm{inst}}}{\partial\qsf\, \partial a_\alpha}(\vec{a};
\qsf) \Big) \qquad \mbox{ for } \quad \alpha=1, \ldots, r\ .
\end{equation}

As in {\cite[Appendix~B.1]{book:nakajimayoshioka2004}}, we introduce
basic $Sp(2r,\Z)$ modular forms on the Seiberg-Witten curve $\Sigma$.

\begin{definition}
For $\vec\mu,\vec\nu\in\C^r$, the \emph{theta-function} $\Theta\left[\genfrac{}{}{0pt}{}{\vec{\mu}}{\vec{\nu}}\right](\vec{\zeta}\,\vert\, \tau)$ \emph{with characteristic} $\left[\genfrac{}{}{0pt}{}{\vec{\mu}}{\vec{\nu}}\right]$ is
\begin{multline}
\Theta\left[\genfrac{}{}{0pt}{}{\vec{\mu}}{\vec{\nu}}\right](\vec{\zeta}\,\vert\, \tau)\\
:=\sum_{\vec{\boldsymbol{n}}=(n_1, \ldots, n_r)\in\Z^r} \, \exp
\Big(\pi \operatorname{i}\, \sum_{\alpha,\beta=1}^r \,
(n_\alpha+\mu_\alpha) \, \tau_{\beta\alpha}\,
(n_\beta+\mu_\beta)+2\pi \operatorname{i} \, \sum_{\alpha=1}^r \,
(\zeta_\alpha+\nu_\alpha) \, (n_\alpha+\mu_\alpha) \Big)\ .
\end{multline}
\end{definition}
\begin{theorem}\label{thm:theta-pure}
$\Zcal_{X_2}^\bullet(\varepsilon_1, \varepsilon_2, \vec{a}; \qsf, \xi, t)/\Zcal^{\mathrm{inst}}_{X_2}(\varepsilon_1, \varepsilon_2, \vec{a}; \qsf, \xi)$ is analytic in $\varepsilon_1, \varepsilon_2$ near $\varepsilon_1=\varepsilon_2=0$, and
\begin{multline}
\lim_{\varepsilon_1,\varepsilon_2\to0}\,
\frac{\Zcal_{X_2}^\bullet(\varepsilon_1, \varepsilon_2, \vec{a}; \qsf,
  \xi, t)}{\Zcal^{\mathrm{inst}}_{X_2}(\varepsilon_1, \varepsilon_2,
  \vec{a}; \qsf, \xi)} \\
=\exp\bigg(\Big(\qsf \, \frac{\partial}{\partial \qsf}\Big)^2
\Fcal_{\C^2}^{\mathrm{inst}}(\vec{a}; \qsf)\, t^2+2\pi
\operatorname{i}\, \sum_{\alpha=w_0+1}^r \, \zeta_\alpha \bigg)\,
\frac{\Theta\left[\genfrac{}{}{0pt}{}{0}{C\,\vec{\nu}}\right](C\,\vec{\zeta}+C\,
  \vec\kappa
  \,\vert\,
  C\,\tau)}{\Theta\left[\genfrac{}{}{0pt}{}{0}{C\,\vec{\nu}}\right](C\,
  \vec\kappa
  \,\vert\, C\, \tau)}\ ,
\end{multline}
where $\kappa_\alpha:=\frac1{4\pi\ii}\, \log\xi$ for
$\alpha=1,\dots,r$ and
\begin{equation*}
\nu_\alpha:=\left\{
\begin{array}{ll}
\displaystyle
\sum_{\beta=w_0+1}^r\, \Big(\log (a_{\beta}-a_{\alpha})^2+
\frac{\partial^2 \Fcal_{\C^2}^{\mathrm{inst}}}{\partial a_\alpha
  \, \partial a_\beta}(\vec{a}; \qsf)\Big) -\frac{w_1}{r}\, \log\qsf &
\mbox{ for } \ \alpha=1, \ldots, w_0\ , \\[10pt]
\displaystyle
-2\, \sum_{\beta=1}^{w_0}\, \log
(a_{\beta}-a_{\alpha})+\frac{2w_0}{r}\, \log\qsf & \mbox{ for } \ \alpha=w_0+1, \ldots, r\ .
\end{array}\right.
\end{equation*}
\end{theorem}
\proof
First we show that the value of $\Zcal_{X_2}^\bullet(\varepsilon_1, \varepsilon_2, \vec{a}; \qsf, \xi, t)/\Zcal^{\mathrm{inst}}_{X_2}(\varepsilon_1, \varepsilon_2, \vec{a}; \qsf, \xi)$ at $\varepsilon_1=\varepsilon_2=0$ is
\begin{multline}\label{eq:ratioblowup}
\exp\bigg(\Big(\qsf \, \frac{\partial}{\partial \qsf}\Big)^2
\Fcal_{\C^2}^{\mathrm{inst}}(\vec{a}; \qsf)\, t^2\Big)\,
\sum_{\stackrel{\scriptstyle v\in\frac{1}{2}\, \Z }{ 2v=
    w_1\bmod{2}}} \, \xi^{v} \ \sum_{\vec{\boldsymbol{v}}=(v_1,
  \ldots, v_r)}\, \qsf^{\sum_{\alpha=1}^r\limits \, v_\alpha^2} \
\prod_{\alpha\neq \beta}\, \Big(\,
\frac{1}{a_\beta-a_\alpha}\, \Big)^{v_{\beta\alpha}^2-\frac{1}{4}\, \delta_{1,c_{\beta\alpha}}}\\
\times\ \exp\bigg(\, \sum_{\alpha,\beta=1}^r \, \frac{\partial^2
  \Fcal_{\C^2}^{\mathrm{inst}}}{\partial a_\alpha \, \partial
  a_\beta}(\vec{a}; \qsf) \, v_\alpha\, v_\beta-2\sum_{\alpha=1}^r \,
\Big(\qsf \, \frac{\partial^2 \Fcal_{\C^2}^{\mathrm{inst}}}{\partial
  \qsf \, \partial a_\alpha}(\vec{a}; \qsf)+ a_\alpha\Big)\, v_\alpha
\, t\bigg)\\
\times \ \bigg[\, \sum_{\stackrel{\scriptstyle v\in\frac{1}{2}\, \Z
  }{\scriptstyle 2v=  w_1\bmod{2}}}\, \xi^{v} \
\sum_{\vec{\boldsymbol{v}}=(v_1, \ldots, v_r)}\,
\qsf^{\sum_{\alpha=1}^r\limits\, v_\alpha^2} \ \prod_{\alpha\neq
  \beta} \, \Big(\, \frac{1}{a_\beta-a_\alpha}\,
\Big)^{v_{\beta\alpha}^2-\frac{1}{4}\, \delta_{1,c_{\beta\alpha}}}\\
\times \ \exp\Big(\, \sum_{\alpha,\beta=1}^r \, \frac{\partial^2
  \Fcal_{\C^2}^{\mathrm{inst}}}{\partial a_\alpha\, \partial
  a_\beta}(\vec{a}; \qsf) \, v_\alpha \,v_\beta\,\Big)\, \bigg]^{-1}\ .
\end{multline}
By \cite[Lemma~6.3-(1)]{art:nakajimayoshioka2005-I} we have
\begin{equation}
F_{\C^2}^{\mathrm{inst}}\big(\varepsilon_1^{(1)},\varepsilon_2^{(1)},
\vec{a}; \qsf
\big)=F_{\C^2}^{\mathrm{inst}}\big(\varepsilon_1^{(2)},\varepsilon_2^{(2)},
\vec{a}; \qsf \big)_{\vert \varepsilon_1 \leftrightarrow \varepsilon_2}\ .
\end{equation}
From this symmetry and some algebraic manipulations, one gets
\begin{multline}\label{eq:limit-blowup}
\frac{1}{\varepsilon_2^{(1)}}\, \Big(\,
\frac{F_{\C^2}^{\mathrm{inst}}\big(\varepsilon_1^{(1)},\varepsilon_2^{(1)},
  \vec{a}^{(1)}; \qsf \, \e^{\varepsilon_1^{(1)}\,
    t}\big)-F_{\C^2}^{\mathrm{inst}}\big(\varepsilon_1^{(1)},\varepsilon_2^{(1)},
  \vec{a}; \qsf \big)}{\varepsilon_1^{(1)}} \\
\shoveright{-\,
  \frac{F_{\C^2}^{\mathrm{inst}}\big(\varepsilon_1^{(2)},\varepsilon_2^{(2)},
    \vec{a}^{(2)}; \qsf \, \e^{\varepsilon_2^{(2)}\,
      t}\big)-F_{\C^2}^{\mathrm{inst}}\big(\varepsilon_1^{(2)},\varepsilon_2^{(2)},
    \vec{a}; \qsf \big)}{\varepsilon_2^{(2)}}\Big)_{\vert (\varepsilon_1, \varepsilon_2)=(0,0)}}\\[6pt]
=-\, \Big(\qsf \, \frac{\partial}{\partial \qsf}\Big)^2 \Fcal_{\C^2}^{
  \mathrm{inst}}(\vec{a}; \qsf)\, t^2-\sum_{\alpha,\beta=1}^r \,
\frac{\partial^2 \Fcal_{\C^2}^{\mathrm{inst}}}{\partial a_\alpha
  \, \partial a_\beta}(\vec{a}; \qsf) \, v_\alpha \,
v_\beta+2\sum_{\alpha=1}^r \qsf \, \frac{\partial^2
  \Fcal_{\C^2}^{\mathrm{inst}}}{\partial \qsf \, \partial
  a_\alpha}(\vec{a}; \qsf) \, v_\alpha \, t \ .
\end{multline}
One computes the numerator of Equation \eqref{eq:ratioblowup} by
dividing $\Zcal_{X_2}^\bullet(\varepsilon_1, \varepsilon_2, \vec{a};
\qsf, \xi, t)$ by the product of partition functions
$\Zcal_{\C^2}^{\mathrm{inst}}\big(\varepsilon_1^{(1)},\varepsilon_2^{(1)},
\vec{a}; \qsf \big)\,
\Zcal_{\C^2}^{\mathrm{inst}}\big(\varepsilon_1^{(2)},\varepsilon_2^{(2)},
\vec{a}; \qsf \big)$ and by using Equation \eqref{eq:limit-blowup}. Similarly, one gets the denominator of Equation \eqref{eq:ratioblowup}.

The rewriting of Equation \eqref{eq:ratioblowup} in terms of
theta-functions follows from straightforward but tedious algebraic manipulations, after setting
\begin{equation*}
v_\alpha=\left\{
\begin{array}{ll}
n_\alpha & \mbox{for } \ \alpha=1, \ldots, w_0\ ,\\
n_\alpha+\frac{1}{2} & \mbox{for } \ \alpha=w_0+1,\ldots, r\ ,
\end{array}\right.
\end{equation*}
where $n_\alpha\in\Z$, and using the identity 
\begin{equation}
 \sum_{ \alpha\neq \beta} \, (v_\alpha-v_\beta)^2=r\,
 \sum_{\alpha=1}^r \, v_\alpha^2-v^2 \ .
\end{equation}
\endproof
\begin{remark}
Our formula for the generating function $\Zcal_{X_2}^\bullet(\varepsilon_1, \varepsilon_2, \vec{a};
\qsf, \xi, t)$ resembles that of
\cite[Equation~(6.10)]{art:nakajimayoshioka2005-I}. The expression in
Theorem \ref{thm:theta-pure} resembles the formula given in
\cite[Theorem~8.1]{art:nakajimayoshioka2005-I}. In particular, note
that when $w_0=r$ and $w_1=0$, i.e., when the holonomy at infinity is trivial, we get
$\vec{\nu}=\vec 0$. 

Recall that by Theorem
\ref{thm:instantonprepotentialpure}, the value of
$F^{\mathrm{inst}}_{X_2}(\varepsilon_1, \varepsilon_2, \vec{a}; \qsf,
\xi)$ at $(\varepsilon_1,\varepsilon_2)=(0,0)$ is 
\begin{equation}
\Fcal^{\mathrm{inst}}_{X_2}(\vec{a}; \qsf)=\mbox{$\frac{1}{2}$} \, \Fcal_{\C^2}^{\mathrm{inst}}(\vec{a}; \qsf)\ .
\end{equation}
In this sense the blowup equation of Theorem \ref{thm:theta-pure}
relates the generating function $\Zcal_{X_2}^\bullet(\varepsilon_1,
\varepsilon_2, \vec{a}; \qsf, \xi, t)$ to the partition function
$\Zcal^{\mathrm{inst}}_{\C^2}(\varepsilon_1, \varepsilon_2, \vec{a};
\qsf)^{\frac{1}{2}}$ which in the low energy limit is identified as
the partition function of $\Ncal=2$ gauge theory on the quotient
$\C^2/\Z_2$.
\end{remark}

\subsubsection*{$\Ncal=2^\ast$ gauge theory}

The generating function $\Zcal_{X_2}^{\ast,\bullet}(\varepsilon_1, \varepsilon_2, \vec{a}, \mu;
\qsf, {\xi}, t)$ for correlators of quadratic $2$-observables
is the generating function
$\Zcal_{X_2}^{\ast}(\varepsilon_1,$ $\varepsilon_2, \vec{a}, \mu; \qsf,
{\xi}, \vec{\tau}, \vec{t} \ )$ introduced in Definition
\ref{def:deformed-adjoint} specialized at $\vec{\tau}=\vec 0$ and
$\vec{t}:=(0, -t, 0, \ldots)$. By using similar arguments as above we get
\begin{empheq}[box=\fbox]{multline}
\Zcal_{X_2}^{\ast,\bullet}(\varepsilon_1, \varepsilon_2, \vec{a}, \mu; \qsf, \xi, t) \\
\shoveleft{=\sum_{\stackrel{\scriptstyle v\in\frac{1}{2}\Z
    }{\scriptstyle 2v=  w_1\bmod{2}}} \, \xi^{v} \
  \sum_{\vec{\boldsymbol{v}}=(v_1, \ldots, v_r)} \, \big(\qsf \,
  \e^{-2t\, (\varepsilon_1+\varepsilon_2)})^{\sum_{\alpha=1}^r\limits
    \, v_\alpha^2} \, \e^{2t\, \sum_{\alpha=1}^r\limits\, v_\alpha\,
    a_\alpha} \ 
\frac{\prod\limits_{\alpha,\beta=1}^r\,
  \ell_{v_{\beta\alpha}}(2\varepsilon_1, \varepsilon_2-\varepsilon_1,
  a_{\beta\alpha}+\mu)}{\prod\limits_{\alpha,\beta=1}^r \,
  \ell_{v_{\beta\alpha}}(2\varepsilon_1, \varepsilon_2-\varepsilon_1 , a_{\beta\alpha})}}\\
\times\Zcal_{\C^2}^{\ast,
  \mathrm{inst}}\big(2\varepsilon_1,\varepsilon_2-\varepsilon_1,
\vec{a}+2\varepsilon_1\, \vec{\boldsymbol{v}}, \mu ; \qsf \,
\e^{2\varepsilon_1\, t}\big)\, \Zcal_{\C^2}^{\ast,
  \mathrm{inst}}\big(\varepsilon_1-\varepsilon_2 ,2\varepsilon_2 ,
\vec{a}+2\varepsilon_2\, \vec{\boldsymbol{v}}, \mu; \qsf \,
\e^{2\varepsilon_2\, t}\big) \ .
\end{empheq}
The period matrix $\tau=(\tau_{\alpha\beta})$ of the genus $r$
Seiberg-Witten curve for $\Ncal=2^\ast$ gauge theory on $\R^4$ is
given by
\begin{multline}
\tau_{\alpha\beta}:= \frac{1}{2\pi \operatorname{i}}\,
 \frac{\partial^2 \Fcal_{\C^2}^{\ast,
    \mathrm{inst}}}{\partial a_\alpha \, \partial a_\beta}(\vec{a}, \mu; \qsf)\\
+\, \frac{1}{2\pi \operatorname{i}}\,
\sum_{\alpha'<\beta'}(\delta_{\alpha\alpha'}-\delta_{\alpha\beta'})\,
(\delta_{\beta\alpha'}-\delta_{\beta\beta'})\, \log \Big(\,
\frac{a_{\beta'}-a_{\alpha'}+\mu}{a_{\beta'}-a_{\alpha'}}\, \Big)
+\frac{\log\qsf}{2\pi \operatorname{i}}\, \delta_{\alpha\beta}\ .
\end{multline}
We define
\begin{equation}
\zeta_\alpha:=-\frac{t}{2\pi
  \operatorname{i}} \, \Big( a_\alpha+\frac{\partial^2\Fcal_{\C^2}^{\ast,
    \mathrm{inst}}}{\partial\qsf\, \partial a_\alpha}(\vec{a}, \mu;
\qsf) \Big) \ .
\end{equation}
By using the same arguments as in the proof of Theorem
\ref{thm:theta-pure} we arrive at the following result.
\begin{theorem}\label{thm:theta-adjoint}
$\Zcal_{X_2}^{\ast,\bullet}(\varepsilon_1, \varepsilon_2, \vec{a},\mu ;
\qsf, \xi, t)/\Zcal^{\ast, \mathrm{inst}}_{X_2}(\varepsilon_1,
\varepsilon_2, \vec{a},\mu; \qsf, \xi)$ is analytic in $\varepsilon_1,
\varepsilon_2$ near $\varepsilon_1=\varepsilon_2=0$, and
\begin{multline}
\lim_{\varepsilon_1,\varepsilon_2\to0}\,
\frac{\Zcal_{X_2}^{\ast,\bullet}(\varepsilon_1, \varepsilon_2,
  \vec{a},\mu ;
\qsf, \xi, t)}{\Zcal^{\ast, \mathrm{inst}}_{X_2}(\varepsilon_1,
\varepsilon_2, \vec{a},\mu ; \qsf, \xi)} \\
= \exp\bigg(\Big(\qsf\,
\frac{\partial}{\partial \qsf}\Big)^2 \Fcal_{\C^2}^{\ast,
  \mathrm{inst}}(\vec{a}, \mu; \qsf)\, t^2+2\pi \operatorname{i}\,
\sum_{\alpha=w_0+1}^r \, \zeta_\alpha \Big)\,
\frac{\Theta\left[\genfrac{}{}{0pt}{}{0}{C\,
      \vec{\nu}}\right](C\,\vec{\zeta}+C\, \vec\kappa
  \,\vert\,
  C\,\tau)}{\Theta\left[\genfrac{}{}{0pt}{}{0}{C\,
      \vec{\nu}}\right](C\, \vec\kappa\,\vert\,
  C\, \tau)}\ ,
\end{multline}
where $\kappa_\alpha:=\frac1{4\pi\ii}\, \log\xi$ for
$\alpha=1,\dots,r$ and
\begin{equation*}
\nu_\alpha:=\left\{
\begin{array}{ll}
\displaystyle
\sum_{\beta=w_0+1}^r\, \Big( \log
\big((a_{\beta}-a_{\alpha})^2-\mu^2\big) + \frac{\partial^2
  \Fcal_{\C^2}^{\ast, \mathrm{inst}}}{\partial a_\alpha \, \partial
  a_\beta}(\vec{a}, \mu; \qsf) \Big) &
\\ \displaystyle  \hfill\vadjust{}-\frac{w_1}{r}\, \log\qsf & \mbox{ for } \ \alpha=1, \ldots, w_0\ , \\[10pt]
\displaystyle
-\sum_{\beta=1}^{w_0}\, \log
\big((a_{\beta}-a_{\alpha})^2-\mu^2\big)+2\frac{w_0}{r}\, \log\qsf &
\mbox{ for } \ \alpha=w_0+1, \ldots, r\ .
\end{array}\right.
\end{equation*}
\end{theorem}

\subsubsection*{Fundamental matter}

We consider only the asymptotically free case, i.e., $N\leq 2r$. The conformal case, i.e., $N=2r$, can be treated in a similar way.

The generating function for quadratic $2$-observables
$\Zcal_{X_2}^{N,\bullet}(\varepsilon_1, \varepsilon_2, \vec{a},
\vec{\mu}; \qsf, {\xi}, t)$ is the generating function
$\Zcal_{X_2}^{N}(\varepsilon_1,$ $\varepsilon_2, \vec{a}, \vec{\mu};
\qsf, \vec{\xi}, \vec{\tau}, \vec{t} \ )$ introduced in Definition
\ref{def:deformed-fund} specialized at $\vec{\tau}=\vec 0$ and
$\vec{t}:=(0, -t, 0, \ldots)$. By using similar arguments as above we get
\begin{empheq}[box=\fbox]{multline}
\Zcal_{X_2}^{N,\bullet}(\varepsilon_1, \varepsilon_2, \vec{a},\vec \mu; \qsf, \xi, t) \\
\shoveleft{=\sum_{\vec{v}\in \Qfrak_{\vec{w}}}\, \xi^{v} \
  \sum_{\vec{\boldsymbol{v}}=(v_1, \ldots, v_r)} \,
  \big(\qsf \,
  \e^{-2t\, (\varepsilon_1+\varepsilon_2)} \big)^{\sum_{\alpha=1}^r\limits\, v_\alpha^2} \, \e^{2t\,
    \sum_{\alpha=1}^r\limits\, v_\alpha \, a_\alpha} \ 
\frac{\prod\limits_{s=1}^N \ \prod\limits_{\alpha=1}^r\,
  \ell_{v_{\alpha}}(2\varepsilon_1, \varepsilon_2-\varepsilon_1 ,
  a_{\alpha}+\mu_s)}{\prod\limits_{\alpha,\beta=1}^r \, \ell_{v_{\beta\alpha}}(2\varepsilon_1, \varepsilon_2-\varepsilon_1, a_{\beta\alpha})}}\\
\times \ \Zcal_{\C^2}^{N,
  \mathrm{inst}}\big(2\varepsilon_1,\varepsilon_2-\varepsilon_1,
\vec{a}+2\varepsilon_1\, \vec{\boldsymbol{v}}, \vec\mu ; \qsf \,
\e^{2\varepsilon_1\, t} \big)\, \Zcal_{\C^2}^{N,
  \mathrm{inst}}\big(\varepsilon_1-\varepsilon_2,2\varepsilon_2,
\vec{a}+2\varepsilon_2\, \vec{\boldsymbol{v}}, \vec\mu ; \qsf \,
\e^{2\varepsilon_2\,t} \big)\ .
\end{empheq}
The period matrix $\tau=(\tau_{\alpha\beta})$ of the genus $r$ Seiberg-Witten
curve for $\Ncal=2$ gauge theory on $\R^4$ with $N$ fundamental
hypermultiplets is given by
\begin{multline}
\tau_{\alpha\beta}:= \frac{1}{2\pi \operatorname{i}}\,
 \frac{\partial^2 \Fcal_{\C^2}^{N,
    \mathrm{inst}}}{\partial a_\alpha \, \partial a_\beta}(\vec{a},
\vec{\mu}; \qsf)-\frac{1}{2\pi \operatorname{i}}\,
\sum_{\alpha'<\beta'}\,
(\delta_{\alpha\alpha'}-\delta_{\alpha\beta'})\,
(\delta_{\beta\alpha'}-\delta_{\beta\beta'})\, \log(a_{\beta'}-a_{\alpha'})\\
+\, \frac{1}{2\pi \operatorname{i}}\, \sum_{s=1}^N\,
\delta_{\alpha\beta}\, \log(a_{\beta}+\mu_s)+\frac{\log\qsf}{2\pi \operatorname{i}}\ .
\end{multline}
We define
\begin{equation}
\zeta_\alpha:=-\frac{t}{2\pi
  \operatorname{i}} \, \Big( a_\alpha+\frac{\partial^2\Fcal_{\C^2}^{N,
    \mathrm{inst}}}{\partial\qsf\, \partial a_\alpha}(\vec{a}, \vec\mu;
\qsf) \Big) \ .
\end{equation}
Let us introduce the set
\begin{equation}
\tilde\Qfrak_{\vec{w}}:=\Big\{(n_1, \ldots,
n_r)\in\Z^r\,\Big\vert\,\Big(\, \mbox{$\sum\limits_{\alpha=1}^r\,
  n_\alpha$}\, \Big)^2+w_1\, \mbox{$\sum\limits_{\alpha=1}^r$}\,
n_\alpha-\frac{3}{4} \, w_1^2\leq 0 \Big\}\ .
\end{equation}
\begin{definition}
For $\vec\mu,\vec\nu\in\C^r$ the \emph{modified} theta function $\tilde \Theta\left[\genfrac{}{}{0pt}{}{\vec{\mu}}{\vec{\nu}}\right](\vec{\zeta}\,\vert\, \tau)$ \emph{with characteristic} $\left[\genfrac{}{}{0pt}{}{\vec{\mu}}{\vec{\nu}}\right]$ is
\begin{multline}
\tilde
\Theta\left[\genfrac{}{}{0pt}{}{\vec{\mu}}{\vec{\nu}}\right](\vec{\zeta}\,\vert\,
\tau) \\
:=\sum_{\vec{\boldsymbol{n}}=(n_1, \ldots, n_r)\in
  \tilde\Qfrak_{\vec{w}}} \, \exp \Big(\pi \operatorname{i}\,
\sum_{\alpha,\beta=1}^r \, (n_\alpha+\mu_\alpha) \,
\tau_{\beta\alpha}\, (n_\beta+\mu_\beta)+2\pi \operatorname{i} \,
\sum_{\alpha=1}^r \, (\zeta_\alpha+\nu_\alpha) \, (n_\alpha+\mu_\alpha) \Big)\ .
\end{multline}
\end{definition}
By using the same arguments as in the proof of Theorem
\ref{thm:theta-pure} we then arrive at the following result.
\begin{theorem}\label{thm:theta-fund}
$\Zcal_{X_2}^{N,\bullet}(\varepsilon_1, \varepsilon_2, \vec{a},\vec\mu ;
\qsf, \xi, t)/\Zcal^{N, \mathrm{inst}}_{X_2}(\varepsilon_1,
\varepsilon_2, \vec{a},\vec\mu; \qsf, \xi)$ is analytic in $\varepsilon_1,
\varepsilon_2$ near $\varepsilon_1=\varepsilon_2=0$, and
\begin{multline}
\lim_{\varepsilon_1,\varepsilon_2\to0}\,
\frac{\Zcal_{X_2}^{N,\bullet}(\varepsilon_1, \varepsilon_2,
  \vec{a},\vec\mu ;
\qsf, \xi, t)}{\Zcal^{N, \mathrm{inst}}_{X_2}(\varepsilon_1,
\varepsilon_2, \vec{a},\vec\mu ; \qsf, \xi)} \\
= \exp\bigg(\Big(\qsf\,
\frac{\partial}{\partial \qsf}\Big)^2 \Fcal_{\C^2}^{N,
  \mathrm{inst}}(\vec{a}, \vec\mu; \qsf)\, t^2+2\pi \operatorname{i}\,
\sum_{\alpha=w_0+1}^r \, \zeta_\alpha \Big)\,
\frac{\tilde\Theta\left[\genfrac{}{}{0pt}{}{0}{C\,
      \vec{\nu}}\right](C\,\vec{\zeta}+C\, \vec\kappa
  \,\vert\,
  C\,\tau)}{\tilde\Theta\left[\genfrac{}{}{0pt}{}{0}{C\,
      \vec{\nu}}\right](C\, \vec\kappa\,\vert\,
  C\, \tau)}\ ,
\end{multline}
where $\kappa_\alpha:=\frac1{4\pi\ii}\, \log\xi$ for
$\alpha=1,\dots,r$ and
\begin{equation*}
\nu_\alpha:=\left\{
\begin{array}{ll}
\displaystyle
\sum_{\beta=w_0+1}^r\,\Big( \log(a_{\beta}-a_{\alpha})^2 
+ \frac{\partial^2 \Fcal_{\C^2}^{N,
    \mathrm{inst}}}{\partial a_\alpha \, \partial a_\beta}(\vec{a},
\vec{\mu}; \qsf) \Big) -\frac{w_1}{r}\log\qsf & \mbox{ for } \ \alpha=1, \ldots, w_0\ , \\[10pt]
\displaystyle
-2\, \sum_{\beta=1}^{w_0}\, \log(a_{\beta}-a_{\alpha})+\frac{2w_0}{r}\log\qsf
+ \, \sum_{s=1}^N\, \log(a_{\alpha}+\mu_s) & \mbox{ for } \ \alpha=w_0+1, \ldots, r\ .
\end{array}\right.
\end{equation*}
\end{theorem}

\begin{remark}
When the holonomy at infinity is trivial, i.e., $w_1=0$, our formula for the generating function
$\Zcal_{X_2}^{N,\bullet}(\varepsilon_1, \varepsilon_2, \vec{a}, \vec\mu;
\qsf, \xi, t)$ is analogous to 
\cite[Equation~(5.3)]{art:gottschenakajimayoshioka2011}, while the expression in
Theorem \ref{thm:theta-fund} is analogous to
\cite[Equation (5.5)]{art:gottschenakajimayoshioka2011}.

Theorems \ref{thm:theta-pure}, \ref{thm:theta-adjoint} and
\ref{thm:theta-fund} are the analogues for $\Ncal=2$ gauge theories of the representation of the
Vafa-Witten partition function in terms of modular forms. For the
$\Ncal=4$ gauge theory considered in Section \ref{sec:VW}, whose prepotential is $\frac{1}{4\pi\ii}\,
\log\qsf\, \sum_{\alpha=1}^r a_\alpha^2$, the
$SL(2,\Z)$ monodromies act on the coupling $\frac{1}{2\pi\ii}\,
\log\qsf$, while for the $\Ncal=2$ gauge theories considered in this
section the $Sp(2r,\Z)$
monodromies act on the periods of the Seiberg-Witten curve, which
determine the low energy
effective gauge couplings, twisted by the intersection form (Cartan
matrix) $C$.
\end{remark}

\bigskip
\appendix

\section{Rank and dimension formulas}\label{sec:rankdimension}

In this appendix  we compute the rank of the natural bundle $\Vbf$ and
the dimension of the moduli space $\Mcal_{r,\vec{u},\Delta}(\Xscr_k,
\Dscr_\infty, \Fcal^{s,\vec{w}}_\infty\, )$  by using the
To\"en-Riemann-Roch theorem \cite{art:toen1999,
  phd:toen1999}. Additional details can be found in \cite[Appendix A]{art:bryancadmanyoung2012}.

\begin{theorem}\label{thm:rank}
The rank of the natural bundle $\Vbf$  of the moduli space
$\Mcal_{r,\vec{u},\Delta}(\Xscr_k, \Dscr_\infty,
\Fcal^{0,\vec{w}}_\infty\, )$ is given by
\begin{equation}\label{eq:rankformula}
\rk(\Vbf) = \Delta +\frac{1}{2r}\, \vec v\cdot C\vec v - \frac{1}{2}
\, \sum_{j=1}^{k-1}\, (C^{-1})^{jj}\, w_j\ .
\end{equation}
\end{theorem}

\begin{corollary}\label{cor:dimension} 
The dimension of the moduli space
$\Mcal_{r,\vec{u},\Delta}(\Xscr_k, \Dscr_\infty,
\Fcal^{s,\vec{w}}_\infty\, )$ 
is
\begin{equation}\label{eq:dimensionformula}
\dim_{\C} \big(\Mcal_{r,\vec{u},\Delta}(\Xscr_k, \Dscr_\infty,
\Fcal^{s,\vec{w}}_\infty\, ) \big) = 2r\, \Delta - \frac{1}{2}\,
\sum_{j=1}^{k-1}\, (C^{-1})^{jj}\, \vec{w}\cdot\vec{w}(j)\ ,
\end{equation}
where for $j=1, \ldots, k-1$ the vector $\vec{w}(j)$ is $(w_j, \ldots, w_{k-1}, w_0, w_1, \ldots, w_{j-1})$.
\end{corollary}
Before attacking the proofs of these results, that we shall give in Section \ref{sec:eulerchar}, in Section \ref{sec:inertiastack} we study the inertia stack of $\Xscr_k$, and in Section \ref{sec:topologicalinv}  we compute some topological invariants of $\Xscr_k$ and $\Dscr_\infty$.

\subsection{Inertia stack}\label{sec:inertiastack}

In this subsection we will compute the inertia stack
$\mathcal{I}(\Xscr_k)$ of $\Xscr_k$. This is a fundamental ingredient
in the application of the To\"en-Riemann-Roch theorem.

\subsubsection{Characterization of the stacky points $p_0$ and $p_\infty$}

We give a characterization of the stacky points $p_0$ and $p_\infty$
of $\Dscr_\infty\subset\Xscr_k$ as trivial gerbes over a point. We
also characterize their Picard groups and the restrictions to them of the generators of the Picard group of $\Dscr_\infty$.

\begin{lemma}
Both stacks $p_0$ and $p_\infty$ are isomorphic to $\Bscr\mu_{k\,
  \tilde{k}}=[{\rm pt}/\mu_{k\,\tilde{k}}]$. At the gerbe structure level, the maps between the banding groups $\mu_k$ of $\Dscr_\infty$ and $\mu_{k\,\tilde{k}}$ of $p_0$, $p_\infty$ are given by
\begin{equation*}
\omega\in\mu_k \ \longmapsto \ \omega^{\pm\, \tilde{k}}\in\mu_{k\, \tilde{k}}\ ,
\end{equation*}
where we take the minus sign for $p_0$ and the plus sign for $p_\infty$.
\end{lemma}
\proof
Consider the cone $\sigma_{\infty,k+2}$. We can compute  the
quotient stacky fan $\stackyfank/\sigma_{\infty,k+2}$. First note that
$N(\sigma_{\infty,k+2})\simeq\Z^2/(\Z v_0\oplus k\, \Z
v_\infty)\simeq\Z_{k\, \tilde{k}}$, and the quotient map $N\to
N(\sigma_{\infty,k+2})$ sends $a\, e_1+b\, e_2$ to $a\bmod{k} \,
\tilde{k}$. The quotient fan
$\bar{\Sigma}_k/\sigma_{\infty,k+2}\subset
N(\sigma_{\infty,k+2})_\Q=0$ is just $\{0\}$. Thus
$\stackyfank/\sigma_{\infty,k+2}=(\Z_{k\, \tilde{k}}, 0, 0)$ and
therefore $p_0$ is the trivial $\mu_{k\, \tilde{k}}$- gerbe
$\Bscr\mu_{k\, \tilde{k}}:=[{\rm pt}/\mu_{k\tilde{k}}]$ over a point
$\rm pt$. 

The quotient map $N\to N(\sigma_{\infty,k+2})$ factorizes through the
quotient map $N(\rho_\infty)\to N(\sigma_{\infty,k+2})$ which is given by
\begin{equation*}
(c,d) \ \longmapsto \ \left\{ 
\begin{array}{ll}
c-d\, \tilde{k}\bmod{k} \, \tilde{k}\; & \mbox{ for even $k$}\ , \\[8pt]
c \, \frac{k-1}{2} - d \, \tilde{k} \bmod{k} \, \tilde{k}\; & \mbox{ for odd $k$}\ .
\end{array}\right.
\end{equation*}
The induced map between the torsion subgroups $\Z_k\to\Z_{k\, \tilde{k}}$ is multiplication by $-\tilde{k}$, and the map between the banding groups of $\Dscr_\infty$ and $p_0$ is given by
\begin{equation*}
\mu_k\simeq\mathrm{Hom}\big(N(\rho_\infty)_{\mathrm{tor}},\C^\ast\big)
\ \xrightarrow{(-)^{-\tilde{k}}} \
\mathrm{Hom}\big(N(\sigma_{\infty,k+2})_{\mathrm{tor}},\C^\ast
\big)\simeq\mu_{k\, \tilde{k}}\ .
\end{equation*}
For $p_\infty$ one can argue similarly. 
\endproof
Now we give a characterization of line bundles over $p_0$ and
$p_\infty$, regarded as trivial $\mu_{k\, \tilde{k}}$-gerbes over a point. 
\begin{lemma}\label{lem:pic}
The Picard group $\Pic(p_0)$ (resp.\ $\Pic(p_\infty)$) of $p_0$
(resp.\ $p_\infty$) is generated by the line bundle $\Lcal_{p_0}$
(resp.\ $\Lcal_{p_\infty}$) corresponding to the character
$\chi\colon\omega\in\mu_{k\, \tilde{k}}\to\omega\in\C^\ast$. In
particular, $\Pic(p_0)\simeq \Pic(p_\infty)\simeq \Z_{k\,
  \tilde{k}}$. The restrictions of the generators $\Lcal_1,\Lcal_2$ of
$\Pic(\Dscr_\infty)$ to $p_0,p_\infty$ are given by
\begin{equation*}
\begin{array}{cclccl}
{\Lcal_1}_{\vert p_0} & \simeq & \Lcal_{p_0} \ ; & {\Lcal_1}_{\vert p_\infty} & \simeq & \Lcal_{p_\infty}\ ,\\[8pt]
{\Lcal_2}_{\vert p_0} & \simeq  & \left\{
\begin{array}{ll}
\Lcal_{p_0}^{\otimes \tilde{k}} & \mbox{ for even $k$}\ , \\
& \\
\Lcal_{p_0}^{\otimes 2k} & \mbox{ for odd $k$}\ ;
\end{array}
\right. & {\Lcal_2}_{\vert p_\infty} & \simeq & \left\{
\begin{array}{ll}
\Lcal_{p_\infty}^{\otimes -\tilde{k}} & \mbox{ for even $k$}\ , \\
& \\
\Lcal_{p_\infty}^{\otimes -2k} & \mbox{ for odd $k$}\ .
\end{array}\right.
\end{array}
\end{equation*}
\end{lemma}
\proof
Consider the cone $\sigma_{\infty,0}$. By arguing as in the
proof of Lemma \ref{lem:restrictionDinfty}, one finds that the
restrictions of the line bundles on $\Xscr_k$ are given by
\begin{equation*}
\Ocal_{\Xscr_k}(\Dscr_0)_{\vert p_0}\simeq\left\{
\begin{array}{ll}
\Lcal_{p_0}^{\otimes k}\; & \mbox{ for even $k$}\ , \\
& \\
\Lcal_{p_0}^{\otimes 2k}\; & \mbox{ for odd $k$}\ .
\end{array}\right.
\qquad \mbox{and} \qquad
\Ocal_{\Xscr_k}(\Dscr_\infty)_{\vert p_0}\simeq\Lcal_{p_0}\ .
\end{equation*}
Hence we obtain
\begin{equation*}
{\Lcal_1}_{\vert p_0}\simeq \Lcal_{p_0} \qquad \mbox{and} \qquad {\Lcal_2}_{\vert p_0}\simeq\left\{
\begin{array}{ll}
\Lcal_{p_0}^{\otimes \tilde{k}}\; & \mbox{ for even $k$}\ , \\
& \\
\Lcal_{p_0}^{\otimes 2k}\; & \mbox{ for odd $k$}\ .
\end{array}\right.
\end{equation*}
For $p_\infty$ one proceeds similarly. 
\endproof

\subsubsection{Characterization of the inertia stack $\mathcal{I}(\Xscr_k)$}

By Theorem \ref{thm:inertiatoric}, we have
\begin{equation*}
\mathcal{I}(\Xscr_k)=\bigsqcup_{v\in\boxx(\stackyfank)}\,
\Xscr\big(\stackyfank/\sigma(\bar{v}) \big)\ .
\end{equation*}
One can show that the cardinality of $\boxx(\stackyfank)$ is $k\,
(2\tilde{k}-1)$, and its elements are classified as follows. The
element $0$ belongs to $\boxx(\sigma)$ for every two-dimensional cone
$\sigma\in\{\sigma_1, \ldots,
\sigma_k,\sigma_{\infty,0},\sigma_{\infty,k}\}$. Its corresponding
minimal cone is $\{0\}\in\Sigma(0)$. Thus
$\Xscr(\stackyfan/\{0\})\simeq \Xscr_k$. Moreover,
$\boxx(\stackyfank)$ contains $k-1$ elements of the form $v_\infty,
2v_\infty, \ldots, (k-1)\, v_\infty$ which belong to
$\rho_\infty\setminus 0$, so that their corresponding minimal cone is
$\rho_\infty$. Thus for $g_i\in G_{\stackyfank}$ corresponding to $i\,
v_\infty$ for $i=1, \ldots, k-1$, we have an isomorphism $\kappa_i\colon [Z^{g_i}_{\bar{\Sigma}_k}/G_{\stackyfank}]\xrightarrow{\sim} \Xscr(\stackyfank/\rho_\infty)=\Dscr_\infty$. 

Let $i=1, \ldots, k-1$. In the following we denote by
$\Dscr_\infty^i$ the substack
$[Z^{g_i}_{\bar{\Sigma}_k}/G_{\stackyfank}]\subset
\mathcal{I}(\Xscr_k)$. After fixing a primitive $k$-th root of unity $\omega$, it is easy to see that the element $g_i$ is $(1,\ldots, 1, \omega^i)\in G_{\stackyfank}$. Then for a scheme $S$, the objects of $\Dscr_\infty^i(S)$ are pairs of the form $(x, g_i)$, where $x$ is an object of $\Dscr_\infty(S)$. The case $i=0$ is excluded because the pairs $(x, 1)$ with $x\in \Dscr_\infty$ are in $\Xscr_k\subset \IXkscr$. The group of automorphisms of $(x,g_i)$ is $\mu_k$ and the inclusion of $\mu_k$ into $G_{\stackyfank}$ is given by the map
\begin{equation*}
\gamma_k^i\, \colon \, \omega\in\mu_k \ \longmapsto \ g_i\in
G_{\stackyfank}= (\mathbb{C}^\ast)^k\ .
\end{equation*}
The isomorphism $\kappa_i$ implies $\imath\circ\varphi_k^i = \gamma_k^i$,
where the maps $\varphi_k^i$ and $\imath$ are given by
\begin{itemize}
\item for even $k$:
\begin{align}
\varphi_k^i\, \colon\, \omega\in\mu_k &\longmapsto
\big(\omega^i,\omega^{i\, \tilde{k}}\big)\in\mathbb{C}^\ast\times \mu_k\ ,\\[4pt]
\imath \, \colon\, (t, \omega)\in\mathbb{C}^\ast\times \mu_k
&\longmapsto \big(1, \ldots, 1, t^{\tilde{k}}\, \omega^{-1},t \big)\in(\mathbb{C}^\ast)^k\ ;
\end{align}
\item for odd $k$:
\begin{align}
\varphi_k^i\, \colon\, \omega\in\mu_k & \longmapsto  \big(\omega^i,1
\big)\in\mathbb{C}^\ast\times \mu_k\ ,\\[4pt] 
\imath \, \colon\, (t, \omega)\in\mathbb{C}^\ast\times \mu_k &
\longmapsto \big(1, \ldots, 1, t^{k}\, \omega^{\frac{k-1}{2}},t \big)\in(\mathbb{C}^\ast)^k\ .
\end{align}
\end{itemize}

The set $\boxx(\stackyfan)$ also contains $k\, \tilde{k}$ elements
belonging to $\sigma_{\infty,k}$. Among them, there are exactly $k$
elements which belong to $\rho_\infty$; their minimal cone is
$\rho_\infty$. The minimal cone of the remaining $k\, \tilde{k}-k$
elements is $\sigma_{\infty,k}$. The corresponding group elements are
$h_j=(1, \ldots, 1, \eta^{2j\, \tilde{k}},\eta^j)\in G_{\stackyfank}$
for $j=0, \ldots, k\, \tilde{k}-1$, where $\eta$ is a primitive
$k\,\tilde{k}$-th root of unity. For $\tilde{k}\mid j$ we have
$h_j=g_{j/\tilde{k}}$ and therefore
$\big[Z^{h_j}_{\bar{\Sigma}_k}/G_{\stackyfank} \big]\simeq
\Dscr_\infty^{j/\tilde{k}}$. So from now on we consider only elements
$h_j$ with $j=1, \ldots, \tilde{k}k-1$, $\tilde{k}\nmid j$. Then for
any $h_j$ we have an isomorphism $\kappa^\infty_j\colon
\big[Z^{h_j}_{\bar{\Sigma}_k}/G_{\stackyfank} \big]\xrightarrow{\sim}\Xscr(\stackyfank/\sigma_{\infty,k})=p_\infty$. Let
$j=1, \ldots, k\, \tilde{k}-1$, $\tilde{k}\nmid j$. Denote by
$p_\infty^j$ the substack
$\big[Z^{h_j}_{\bar{\Sigma}_k}/G_{\stackyfank} \big]\subset
\mathcal{I}(\Xscr_k)$. Then for a scheme $S$, the objects of
$p_\infty^j(S)$ are pairs of the form $(y, h_j)$, where $y\in
p_\infty(S)$. The group of automorphisms of $(y,h_j)$ is $\mu_{k\,
  \tilde{k}}$ and the inclusion of $\mu_{k\, \tilde{k}}$ into $G_{\stackyfank}$ is given by the map
\begin{equation*}
\gamma^{k,\infty}_j\, \colon\, \eta\in \mu_{k\, \tilde{k}} \ \longmapsto
\ h_j\in G_{\stackyfank}=(\mathbb{C}^\ast)^k\ .
\end{equation*}
The isomorphism $\kappa^\infty_j$ implies that $\imath\circ\jmath^\infty\circ \varphi^{k,\infty}_j = \gamma^{k,\infty}_j$,
where the maps $\varphi^{k,\infty}_j$ and $\jmath^\infty$ are given by
\begin{itemize}
\item for even $k$:
\begin{align}
\varphi^{k,\infty}_j\, \colon\, \eta\in\mu_{k\, \tilde{k}} &
\longmapsto  \eta^j\in\mu_{k\, \tilde{k}}\ ,\\[4pt]
\jmath^\infty \, \colon\, \eta\in \mu_{k\, \tilde{k}} & \longmapsto 
\big(\eta, \eta^{-\tilde{k}} \big)\in\mathbb{C}^\ast\times \mu_k\ ;
\end{align}
\item for odd $k$:
\begin{align}
\varphi^{k,\infty}_j\, \colon\, \eta\in\mu_{k\, \tilde{k}} &
\longmapsto  \eta^j\in\mu_{k\, \tilde{k}}\ ,\\[4pt]
\jmath^\infty \, \colon\, \eta\in \mu_{k\, \tilde{k}} & \longmapsto 
\big(\eta, \eta^{-2k} \big)\in \mathbb{C}^\ast\times \mu_k\ .
\end{align}
\end{itemize}

In a similar way, we obtain substacks $p_0^j\subset
\mathcal{I}(\Xscr_k)$ associated to $f_j=(1,\ldots, 1, \eta^j)\in
G_{\stackyfank}$, which are isomorphic to $p_0$, where $\eta$ is a
primitive $k\, \tilde{k}$-th root of unity and $j=1,\ldots, k\,
\tilde{k}-1$, $\tilde{k}\nmid j$. Therefore we obtain an isomorphism $\imath\circ j^0\circ \varphi^{k,0}_j = \gamma^{k,0}_j$,
where the maps $\varphi^{k,0}_j$ and $\jmath^0$ are given by
\begin{itemize}
\item for even $k$:
\begin{align}
\varphi^{k,0}_j\, \colon\, \eta\in\mu_{k\, \tilde{k}} & \longmapsto 
\eta^j\in\mu_{k\, \tilde{k}}\ ,\\[4pt]
\jmath^0 \,\colon\, \eta\in \mu_{k\, \tilde{k}} & \longmapsto 
\big(\eta, \eta^{\tilde{k}} \big)\in\mathbb{C}^\ast\times \mu_k\ ;
\end{align}
\item for odd $k$:
\begin{align}
\varphi^{k,0}_j\, \colon \, \eta\in\mu_{k\, \tilde{k}} & \longmapsto 
\eta^j\in\mu_{k\, \tilde{k}}\ ,\\[4pt]
\jmath^0 \,\colon\, \eta\in \mu_{k\, \tilde{k}} & \longmapsto 
\big(\eta, \eta^{2k} \big)\in \mathbb{C}^\ast\times \mu_k\ .
\end{align}
\end{itemize}

Thus we can write the inertia stack as
\begin{equation}\label{eq:inertiastack}
\IXkscr=\Xscr_k \sqcup \bigg( \, \bigsqcup_{i=1}^{k-1}\,
\Dscr_\infty^i \, \bigg) \sqcup \bigg(\,
\bigsqcup_{\stackrel{\scriptstyle j=1}{\scriptstyle \tilde{k}\nmid
    j}}^{k\, \tilde{k}-1} \, p_0^j \, \bigg) \sqcup \bigg(\,
\bigsqcup_{\stackrel{\scriptstyle j=1}{\scriptstyle \tilde{k}\nmid
    j}}^{k\, \tilde{k}-1} \, p_\infty^j \, \bigg)\ .
\end{equation}

\subsection{Topological invariants}\label{sec:topologicalinv}

We compute the integrals of Chern classes of the tangent bundles to $\Xscr_k$ and $\Dscr_\infty$.
The  canonical bundle of $\Dscr_\infty$ is
$\Kcal_{\Dscr_\infty}\simeq\Ocal_{\Dscr_\infty} (-p_0-p_\infty)$;
this can be regarded as a generalization of the analogous result for
toric varieties \cite[Theorem~8.2.3]{book:coxlittleschenck2011} (see also \cite{art:kawamata2006}). By Corollary \ref{cor:decomposition} we obtain $\Kcal_{\Dscr_\infty}\simeq \Lcal_1^{\otimes -2\tilde{k}}$. Let $\Tcal_{\Dscr_\infty}$ denote the tangent sheaf to $\Dscr_\infty$.
By Lemma \ref{lem:degree} one has 
\begin{equation}\label{eq:intc1dinfty}
\int_{\Dscr_\infty} \, \crm_1 (\Tcal_{\Dscr_\infty}) =
\int_{\Dscr_\infty} \, \crm_1 \big(\Ocal_{\Dscr_\infty} (p_0+p_\infty)
\big) = \frac{2}{k\, \tilde{k}}\ .
\end{equation}

By applying \cite[Theorem~3.4]{art:tseng2011}, which enables us to compute the integral of total Chern class $\crm(\Tcal_{\IXkscr})$, we have
\begin{equation*}
\int_{\IXkscr} \, \crm(\Tcal_{\IXkscr}) = \int_{\bar{X}_k} \, \crm^{\rm
  SM}(\bar{X}_k) = e(\bar{X}_k) = \big\vert\bar\Sigma_k(2) \big\vert = k+2\ ,
\end{equation*}
where $\crm^{\rm SM}(\bar{X}_k)$ denotes the Chern-Schwartz-MacPherson
class. The second equality comes from
\cite{art:parusinskipragacz1995}, while the third equality comes from \cite[Theorem~12.3.9]{book:coxlittleschenck2011}. On the other hand, by the decomposition \eqref{eq:inertiastack} of the inertia stack $\IXkscr$ we have
\begin{equation*}
\int_{\IXkscr} \, \crm(\Tcal_{\IXkscr}) = \int_{\Xscr_k} \,
\crm_2(\Tcal_{\Xscr_k}) + (k-1)\, \int_{\Dscr_\infty} \,
\crm_1(\Tcal_{\Dscr_\infty}) + k\, (\tilde{k}-1)\, \int_{p_0} \, 1 +
k\, (\tilde{k}-1)\, \int_{p_\infty} \, 1\ .
\end{equation*}
Recall that the order of the stabilizers of $p_0$ and $p_\infty$ is
$k\, \tilde{k}$, so that $\int_{p_0}\, 1 = \frac{1}{k\, \tilde{k}}\,
\int_{\rm pt} \, 1 = \frac{1}{k\, \tilde{k}}$ where $\rm pt$ is
understood to be the one-point scheme which we may regard as the coarse moduli space of $p_0$. For $p_\infty$ one obtains the same result, so that
\begin{equation}\label{eq:intc2xk}
\int_{\Xscr_k}\, \crm_2(\Tcal_{\Xscr_k}) = k + \frac{2}{k\, \tilde{k}}\ .
\end{equation}

\subsection{Euler characteristic}\label{sec:eulerchar}

In this subsection we collect the results described so far and compute
all the ingredients needed to prove Theorem \ref{thm:rank}. As
explained in the proof of Proposition \ref{prop:naturalbundle}, the
rank of the natural bundle $\Vbf$ is given by the Euler characteristic
as $-\chi(\Xscr_k, \Ecal\otimes \Ocal_{\Xscr_k}(-\Dscr_\infty))$,
where $\Ecal$ is the underlying sheaf of a point
$[(\Ecal,\phi_{\Ecal})]\in\Mcal_{r,\vec{u},\Delta}(\Xscr_k,\Dscr_\infty,\Fcal_\infty^{0,\vec
  w}\, )$.

By using the To\"en-Riemann-Roch theorem we have 
\begin{equation*}
\chi\big(\Xscr_k,  \Ecal\otimes \Ocal_{\Xscr_k}(-\Dscr_\infty)
\big)=\int_{\IXkscr} \, \frac{ \ch \big(\rho(\varpi^\ast(\Ecal\otimes
  \Ocal_{\Xscr_k}(-\Dscr_\infty))) \big)}{ \ch
  \big(\rho(\lambda_{-1}(\Ncal^\vee)) \big) } \cdot \mathrm{Td}(\Tcal_{\IXkscr})\ ,
\end{equation*}
where $\Ncal$ is the normal bundle to the local immersion $\varpi\colon \IXkscr\to\Xscr_k$.

Using the decomposition \eqref{eq:inertiastack}, the integral over the
inertia stack becomes a sum of four terms given by
\begin{align}
A  := & \int_{\Xscr_k} \, \ch \big(
\Ecal\otimes\Ocal_{\Xscr_k}(-\Dscr_\infty) \big) \cdot \mathrm{Td}(\Xscr_k)\ , \\[4pt]
B  := & \sum_{i=1}^{k-1} \ \int_{\Dscr_\infty^i} \, \frac{\ch \Big(
  \rho \big((\Ecal\otimes\Ocal_{\Xscr_k}(-\Dscr_\infty) \big)_{\vert
  \Dscr_\infty^i} \big) \Big)}{\ch \Big( \rho \big( \lambda_{-1}
(\Ncal_{\Dscr_\infty^i / \Xscr_k}^\vee) \big) \Big)} \cdot \mathrm{Td}(\Dscr_\infty)\ , \\[4pt]
C  := & \sum_{\stackrel{\scriptstyle i=1}{\scriptstyle \tilde{k}\nmid
    i}}^{k\,\tilde k-1} \ \int_{p_0^i} \, \frac{\ch \Big( \rho \big((
\Ecal\otimes\Ocal_{\Xscr_k}(-\Dscr_\infty) \big)_{\vert p_0^i} \big)
\Big)}{\ch \Big( \rho \big( \lambda_{-1} (\Ncal_{p_0^i /
  \Xscr_k}^\vee) \big) \Big)} \cdot \mathrm{Td}(p_0)\ , \\[4pt]
D  := & \sum_{\stackrel{\scriptstyle i=1}{\scriptstyle \tilde{k}\nmid
    i}}^{k\,\tilde k-1} \ \int_{p_\infty^i} \, \frac{\ch \Big( \rho \big((
\Ecal\otimes\Ocal_{\Xscr_k}(-\Dscr_\infty) \big)_{\vert p_\infty^i} \big)
\Big)}{\ch \Big( \rho \big( \lambda_{-1} (\Ncal_{p_\infty^i /
  \Xscr_k}^\vee) \big) \Big)} \cdot \mathrm{Td}(p_\infty) \ .
\end{align}
We compute each part separately. In the ensuing calculations we make
use of a number of identities involving sums of complex roots of
unity, which are proven in Appendix \ref{ap:rootsunity} below.

\subsubsection*{Computation of $A$}
\begin{align}
A =& \int_{\Xscr_k} \, \big(\ch_2(\Ecal) - \crm_1(\Ecal)\cdot
[\Dscr_\infty] +\crm_1(\Ecal)\cdot \mathrm{Td}_1 (\Xscr_k) ) \\ & +\,
r\, \int_{\Xscr_k} \, \Big(\mathrm{Td}_2 (\Xscr_k) + \frac{1}{2} \,
[\Dscr_\infty]^2 - [\Dscr_\infty] \cdot \mathrm{Td}_1 (\Xscr_k)
\Big) = \int_{\Xscr_k} \, \ch_2(\Ecal)  +r\, \frac{k^2\, \tilde{k}^2
  + 4\tilde{k}^2 -6\tilde{k} +1}{12\,k\, \tilde{k}^2} 
\end{align}
since $\crm_1(\Ecal)\cdot [\Dscr_\infty] = \crm_1(\Ecal)\cdot
\mathrm{Td}_1 (\Xscr_k) =0$; the evaluation of the second integral
follows from Equation \eqref{eq:intc2xk}, Proposition \ref{prop:intersections} and the adjunction formula  \cite[Theorem 3.8]{art:nironi2008-II}.
 
\subsubsection*{Computation of $B$}
We first need to compute $\rho \big( ( \Ecal\otimes \Ocal_{\Xscr_k}
(-\Dscr_\infty) )_{\vert \Dscr_\infty^i} \big)$. For a fixed value of
$i$, the homomorphism $\rho$ sends the K-theory class $[\Gcal]$ of a
vector bundle $\Gcal$ on $\Dscr_\infty^i$ to $\sum_m\, \omega^{i \, m}\, [\Gcal_m]$ where $\Gcal_m$ is the $m$-th summand in the decomposition of $\Gcal$ with respect to the action of $\omega^i$. Note that $(\Ecal\otimes \Ocal_{\Xscr_k} (-\Dscr_\infty))_{\vert \Dscr_\infty} \simeq \Fcal_\infty^{0,\vec{w}} \otimes \Lcal_1^{\otimes -1}\simeq \bigoplus_{j=0}^{k-1}\, \Ocal_{\Dscr_\infty}(-1,j)^{\oplus w_j}$. Then
\begin{equation*}
\rho \big( \left( \Ecal\otimes \Ocal_{\Xscr_k} (-\Dscr_\infty)
\right)_{\vert \Dscr_\infty^i} \big) =\sum_{j=0}^{k-1} \, w_j \,
\rho\big(\Ocal_{\Dscr_\infty}(-1,j) \big)\ .
\end{equation*}

\begin{lemma}\label{lem:restrictionpoint} 
For fixed $i=1,\dots,k-1$ one has
\begin{equation*} 
\rho\big(\Ocal_{\Dscr_\infty}(-1,j) \big) = \left\{
\begin{array}{ll}
\omega^{i\, (\tilde{k}\, j-1)}\big[ \Ocal_{\Dscr_\infty}(-1,j) \big] \qquad & \mbox{for even $k$}\ ,\\[8pt]
\omega^{-i}\, \big[ \Ocal_{\Dscr_\infty}(-1,j) \big] \qquad & \mbox{for odd $k$}\ .
\end{array} \right.
\end{equation*}
\end{lemma}
\proof 
Suppose $k$ is even. Recall that $\Ocal_{\Dscr_\infty}(-1,j) \simeq
\Lcal_1^{\otimes -1} \otimes \Lcal_2^{\otimes j}$ corresponds to the
character $\chi^{(-1,j)} \colon$ $ (t,\omega) \in \C^\ast \times \mu_k
\mapsto t^{-1}\, \omega^j \in \C^\ast$. The element
$\rho(\Ocal_{\Dscr_\infty}(-1,j))$ is computed with respect to the map
$\varphi^i_k\colon\omega \in \mu_k \mapsto (\omega^i, \omega^{i\,
  \tilde{k}}) \in \C^\ast \times \mu_k$, where $\omega$ is a primitive
$k$-th root of unity. So the composition of the latter map with
$\chi^{(-1,j)}$ gives $\omega \in \mu_k \mapsto \omega^{i\,
  (\tilde{k}\, j-1)} \in \C^\ast$.

For odd $k$ one has a similar result. In that case the map $\varphi_k^i$ is given by $\omega \in \mu_k \mapsto (\omega^i, 1) \in \C^\ast \times \mu_k$, which by composition with $\chi^{(-1,j)}$ yields $\omega \in \mu_k \mapsto \omega^{-i} \in \C^\ast$.
\endproof

By Lemma \ref{lem:restrictionpoint}, we obtain
\begin{equation*}
\ch \Big( \rho \big( \left( \Ecal\otimes \Ocal_{\Xscr_k} (-\Dscr_\infty) \right)_{\vert \Dscr_\infty^i} \big) \Big) = \left\{
\begin{array}{ll}
\sum_{j=0}^{k-1}\limits\, w_j \, \omega^{i\, (\tilde{k}\, j-1)} \, \ch\big( \Ocal_{\Dscr_\infty}(-1,j) \big) \; & \mbox{for even $k$}\ ,\\[8pt]
\sum_{j=0}^{k-1}\limits\, w_j\,  \omega^{-i} \, \ch\big( \Ocal_{\Dscr_\infty}(-1,j) \big) \; & \mbox{for odd $k$}\ .
\end{array}
\right.
\end{equation*}
The normal bundle $\Ncal_{\Dscr_\infty/\Xscr_k}$ is isomorphic to $\Ocal_{\Xscr_k}(\Dscr_\infty)_{\vert \Dscr_\infty}\simeq \Lcal_1$ (cf.\ Lemma \ref{lem:restrictionDinfty}). Thus again by Lemma \ref{lem:restrictionpoint} we get
\begin{equation*}
\ch \Big( \rho \big( \lambda_{-1} \Ncal_{\Dscr_\infty^i/\Xscr_k}^\vee
\big) \Big) = \ch \big( 1 - \rho( \Lcal_1^{\otimes -1} ) \big) = 1 -
\omega^{-i} \, \ch ( \Lcal_1^{\otimes -1} ) = 1 - \omega^{-i} \, \big(
1-\crm_1(\Lcal_1) \big) \ .
\end{equation*}
Set $s_j=\tilde{k}\, j-1$ if $k$ is even and $s_j=-1$ if $k$ is odd. We obtain
\begin{align}
B  =  & \sum_{i=1}^{k-1} \ \int_{\Dscr_\infty} \, \Big(\,
\sum_{j=0}^{k-1} w_j \, \omega^{i\, s_j} \, \big( 1 + \crm_1(
\Ocal_{\Dscr_\infty}(-1,j) ) \big) \,\Big) \\
& \qquad \qquad \qquad \cdot \ \Big(\, \frac{1}{1 - \omega^{-i}} -
\frac{\omega^{-i}}{(1 - \omega^{-i})^2} \, \crm_1(\Lcal_1) \, \Big) \cdot \big( 1 + \mathrm{Td}_1 (\Dscr_\infty) \big) \\[4pt]
 = & \sum_{j=0}^{k-1} \, w_j \ \sum_{i=1}^{k-1} \, \frac{\omega^{i\,
     s_j}}{1 - \omega^{-i}} \, \Big(\, \frac{1}{k\, \tilde{k}} -
 \frac{1}{k\, \tilde{k}^2} - \frac{1}{k\, \tilde{k}^2} \,
 \frac{\omega^{-i}}{1 - \omega^{-i}} \, \Big) \ ,
\end{align}
where the last equality follows from Equation \eqref{eq:intc1dinfty}
and Lemma \ref{lem:degree}. Now we use the identity  \cite{art:bryan2012}
\begin{equation}\label{eq:bryanidentity}
\frac{1}{k} \, \sum_{i=1}^{k-1} \, \frac{\omega^{i\, s}}{1 - \omega^{-i}} = \left\lfloor \frac{s}{k} \right\rfloor - \frac{s}{k} + \frac{k-1}{2k}\ ,
\end{equation}
together with the fact that
\begin{equation*}
\left\lfloor \frac{s_j}{k} \right\rfloor - \frac{s_j}{k} = \left\{
\begin{array}{ll}
\frac{1}{k} - 1 \qquad & \mbox{for odd $k$, or even $k$ and even $j$}
\ , \\[3pt]
\frac{1}{k} - \frac{1}{2} \qquad & \mbox{for even $k$ and odd $j$} \ ,
\end{array}
\right.
\end{equation*}
and $\sum_{j=0}^{k-1}\, w_j = r$. Thus for odd $k$ we easily get  
\begin{equation*}
\sum_{j=0}^{k-1} \, w_j \ \sum_{i=1}^{k-1} \, \frac{\omega^{i\,
    s_j}}{1 - \omega^{-i}} \, \Big(\, \frac{1}{k\, \tilde{k}} -
\frac{1}{k\, \tilde{k}^2} \Big) =\frac{\tilde{k} - 1}{\tilde{k}^2} \,
\sum_{j=0}^{k-1}\, w_j \,\Big(\, \left\lfloor \frac{s_j}{k}
\right\rfloor - \frac{s_j}{k} + \frac{k-1}{2k} \, \Big) = r\,
-\frac{(k-1)^2}{2k^3}\ .
\end{equation*}
In the case of even $k$, define the natural numbers $r_e = \sum_{j \;
  {\rm even}}\, w_j$ and $r_o = \sum_{j \; {\rm odd}}\, w_j$; then $r
= r_e + r_o$ and we have
\begin{equation*}
 \frac{\tilde{k} - 1}{\tilde{k}^2} \, \sum_{j=0}^{k-1} \, w_j \,
 \Big(\, \left\lfloor \frac{s_j}{k} \right\rfloor - \frac{s_j}{k} +
 \frac{k-1}{2k} \, \Big) =  r_e \, \frac{(\tilde{k}-1)\, (1-k)}{2k\,
   \tilde{k}^2}+  r_o \, \frac{\tilde{k} - 1}{2k\, \tilde{k}^2}\ .
\end{equation*}
Thus 
\begin{equation*}
B = \left\{
\begin{array}{ll} \displaystyle
-r\, \frac{(k-1)^2}{2k^3} -\frac{r}{k^3}\, \sum_{i=1}^{k-1} \, \frac{\omega^{-2i}}{(1 - \omega^{-i})^2} \; & \mbox{ for odd $k$}\ , \\[5pt]
\displaystyle
r_e\, \frac{(\tilde{k}-1)\, (1-k)}{2k\, \tilde{k}^2}+r_o\,
\frac{\tilde{k} - 1}{2k\, \tilde{k}^2} -\frac{1}{k\, \tilde{k}^2} \,
\sum_{j=0}^{k-1} \, w_j \ \sum_{i=1}^{k-1}\, \frac{\omega^{i\,
    (\tilde{k}\, j-2)}}{(1 - \omega^{-i})^2} \; & \mbox{ for even $k$}\ .
\end{array}
\right.
\end{equation*}
Using Lemma \ref{lem:simpleidentity}, with $s=k-2$ for odd $k$, and
for even $k$ and even $j$, and with $s=\tilde{k}-2$ for even $k$ and
odd $j$, we get
\begin{equation*}
B= \left\{
\begin{array}{ll} \displaystyle
-\frac{5k^2 - 6k + 1}{12k^3}\, r \qquad & \mbox{for odd $k$}\ , \\[10pt] \displaystyle
-\frac{2k^2 - 3k + 1}{3k^3}\, r + \frac{r_o}{2k} \qquad & \mbox{for even $k$} \ .
\end{array}
\right.
\end{equation*}

\subsubsection*{Computation of $C$ and $D$} 

We shall need the following characterization.
\begin{lemma}\label{lem:conormal0}
The conormal bundles to $p_0$  and  $p_\infty$ in $\Xscr_k$ are given by
\begin{equation*}
\Ncal_{p_0/\Xscr_k}^\vee \simeq \Lcal_{p_0}^{\otimes - 2 \tilde{k}}
\oplus \Lcal_{p_0}^{\otimes -1}\qquad \mbox{and} \qquad
\Ncal_{p_\infty/\Xscr_k}^\vee \simeq \Lcal_{p_\infty}^{\otimes - 2 \tilde{k}} \oplus \Lcal_{p_\infty}^{\otimes -1}\ .
\end{equation*}
\end{lemma}
\proof 
Since $p_0$ is a zero-dimensional substack in $\Xscr_k$, its tangent
bundle is trivial and so $\Ncal_{p_0/\Xscr_k}^\vee \simeq
{\Tcal_{\Xscr_k}}_{\vert p_0}^\vee$. The divisors $\Dscr_0$ and
$\Dscr_\infty$, which intersect in $p_0$, are normal crossing and so the tangent bundle splits as 
\begin{equation*}
{\Tcal_{\Xscr_k}}_{\vert p_0} \simeq {\Tcal_{\Dscr_0}}_{\vert p_0} \oplus {\Tcal_{\Dscr_\infty}}_{\vert p_0}\ .
\end{equation*}
By the adjunction formula \cite[Theorem~3.8]{art:nironi2008-II}, we obtain
\begin{align}
{\Tcal_{\Dscr_0}}_{\vert p_0} \simeq & \big(\Kcal_{\Dscr_0}^\vee
\big)_{\vert p_0} \\[4pt] \simeq & \big( (\Kcal_{\Xscr_k} \otimes
\Ocal_{\Xscr_k}(\Dscr_0))_{\vert \Dscr_\infty}^\vee \big)_{\vert p_0} \\[4pt]
\simeq & \bigg( {\Ocal_{\Xscr_k}\Big(-\mbox{$\sum\limits_{i=1,
      \ldots, k, \infty} $}\, \Dscr_i\Big)}_{\vert \Dscr_\infty}^\vee \bigg)_{\vert p_0} \simeq {\Lcal_1}_{\vert p_0} \simeq \Lcal_{p_0}\ .
\end{align}
Since the canonical line bundle $\Kcal_{\Dscr_\infty}$ is isomorphic
to $\Ocal_{\Dscr_\infty}(-p_0-p_\infty)$, we get
${\Tcal_{\Dscr_\infty}}_{\vert p_0} \simeq (\Ocal_{\Dscr_\infty}(p_0 +
p_\infty))_{\vert p_0} \simeq {\Lcal_1^{\otimes 2 \tilde{k}}}_{\vert
  p_0} \simeq \Lcal_{p_0}^{\otimes 2 \tilde{k}}$ and the first
isomorphism is proved. The second isomorphism is proven in the same way.
\endproof
By Lemma \ref{lem:conormal0} we obtain
\begin{equation*}
\ch_0 \Big( \rho \big( \lambda_{-1} (\Ncal_{p_0^i/\Xscr_k}^\vee) \big)
\Big)  = \big(1 - \eta^{-i}\big)\, \big(1 - \eta^{-2i\, \tilde{k}}
\big) = \ch_0 \Big( \rho \big( \lambda_{-1}
(\Ncal_{p_\infty^i/\Xscr_k}^\vee) \big) \Big) \ .
\end{equation*}
By Lemma \ref{lem:pic} we have
\begin{multline}
\big( \Ecal\otimes \Ocal_{\Xscr_k}(-\Dscr_\infty) \big)_{\vert p_0}
\simeq   \big( \Fcal_\infty^{0, \vec{w}}\otimes \Lcal_1^{\otimes -1}
\big)_{\vert p_0}  \simeq \Big( \, \bigoplus_{j=0}^{k-1} \,
\Ocal_{\Dscr_\infty}(-1,j)^{\oplus w_j} \, \Big)_{\vert p_0} \simeq 
 \bigoplus_{j=0}^{k-1}\, \big( \Lcal_{p_0}^{\otimes \tilde{k}\, j - 1}
 \big)^{\oplus w_j} \ ,
\end{multline}
and similarly
\begin{equation*}
\big(\Ecal\otimes \Ocal_{\Xscr_k}(-\Dscr_\infty) \big)_{\vert
  p_\infty} \simeq \bigoplus_{j=0}^{k-1} \, \big(
\Lcal_{p_\infty}^{\otimes -\tilde{k}\, j - 1} \big)^{\oplus w_j}\ .
\end{equation*}
Thus we have
\begin{equation}\label{eq:C+D}
C+D=\frac{1}{k\, \tilde{k}} \, \sum_{j=0}^{k-1} \, w_j \
\sum_{\stackrel{\scriptstyle i=1}{\scriptstyle \tilde{k} \nmid
    i}}^{k\,\tilde k-1}\, \frac{\eta^{i\, (\tilde{k}\,
    j-1)}+\eta^{-i\, (\tilde{k}\, j+1)}}{\big(1 - \eta^{-i}\big)\,
  \big(1 - \eta^{-2i\, \tilde{k}}\big)}\ .
\end{equation}

We now have to distinguish two cases.

\subsubsection*{Odd $k$}\label{ssec:CDkodd}

By Lemma \ref{prelimLemma3} we have
\begin{align}
C+D=&
\frac1{k^2}\, \sum_{j=0}^{k-1}\, w_j\ \sum_{i=1}^{k-1}\, \frac{\omega^{i\,j}+\omega^{-i\,j}}{1-\omega^{-2i}} \ \sum_{l=0}^{k-1}\, \frac1{\eta^i\, \omega^{l}-1} \\[4pt] 
=& \frac1{4k}\, \sum_{j=0}^{k-1}\, w_j \ \sum_{i=1}^{k-1}\, (\omega^{i\,j}+\omega^{-i\,j})\, \Big(\, \frac{3-\omega^i}{(1-\omega^{-i})^2}+ \frac{\omega^{2i}}{1+\omega^i}\, \Big) \ .
\end{align}
It is convenient to separate the contributions from $j=0$ and $j\geq1$ in the above sum; we denote them by $(C+D)_0$ and $(C+D)_>$, respectively. By Lemmas \ref{prelimLemma1} and  \ref{lem:simpleidentity} with $s=k$ and $s=1$ we find
\begin{equation*}
(C+D)_0=-\frac{k^2-1}{12k}\, w_0\ ,
\end{equation*}
while by Lemmas \ref{prelimLemma2} and \ref{lem:simpleidentity} with $s=j, s=k-j$ and $s=j+1, s=-j+1$ we get
\begin{equation}
(C+D)_> = \sum_{j=1}^{k-1}\, w_j \,\Big(\, \frac{j\, (k-j)}{2k} -
\frac{k^2-1}{12k} \, \Big) \ .
\end{equation}
Thus
\begin{equation}\label{eq:CDdimensionodd}
C + D = - \frac{k^2-1}{12k}\, r+ \sum_{j=1}^{k-1}\, \frac{j\, (k-j)}{2k}\, w_j\ .
\end{equation}

\subsubsection*{Even $k$}\label{ssec:CDkeven}

By Lemma \ref{prelimLemma4} we have
\begin{align}
C+D =& \frac{2}{k^2}\, \sum_{j=0}^{k-1}\, w_j \ \sum_{i=1}^{\tilde k-1}\, \frac{\omega^{i\,j}+\omega^{-i\,j}}{1-\omega^{-2i}} \ \sum_{l=0}^{k-1}\, \frac{(-1)^{l\, j}}{\eta^i\, \omega^{l}-1} \\[4pt]
=& \frac{2}{k}\, \sum_{j\; {\rm even}}\, w_j \
\sum_{i=1}^{\tilde{k}-1} \, \frac{\omega^{i\, (j-2)} + \omega^{i\, (-j-2)}}{(1-\omega^{-2i})^2}+
\frac{2}{k}\, \sum_{j\; {\rm odd}}\, w_j \ \sum_{i=1}^{\tilde{k}-1} \,
\frac{\omega^{i\, (j-1)} + \omega^{i\, (-j-1)}}{(1-\omega^{-2i})^2}\ .
\end{align}
Here we use Lemma \ref{lem:simpleidentity}, with $\tilde{k}$ instead
of $k$ and $\omega^2$ (a $\tilde{k}$-th root of unity) instead of
$\omega$. For even $j$, we use $s=j/2-1$ and $s=\tilde{k}-j/2-1$,
while for odd $j$ we use $s=(j-1)/2$ and $s= \tilde{k}-(j+1)/2$. This gives
\begin{equation}\label{eq:CDdimensioneven}
C+D = - \frac{k^2-1}{12k}\, r + \sum_{j=1}^{k-1}\, \frac{j\, (k-j)}{2k}\, w_j + \frac{r_o - r_e}{4k} \ .
\end{equation}

Thus Equation \eqref{eq:C+D} becomes for odd $k$
\begin{equation*}
C+D = -\frac{k^2-1}{12k}\, r + \frac{1}{2}\, \sum_{j=1}^{k-1}\,
(C^{-1})^{jj}\, w_j\ ,
\end{equation*}
while for even $k$
\begin{equation*}
C+D = -\frac{k^2-4}{12k}\, r - \frac{r_o}{2k} + \frac{1}{2} \,
\sum_{j=1}^{k-1}\, (C^{-1})^{jj}\, w_j\ .
\end{equation*}
 
\proof[Proof of Theorem \ref{thm:rank}.] 
One just needs to sum the contributions $A$, $B$, $C$ and $D$ and change sign.
\endproof

\subsection{Dimension formula}

The proof of  Corollary \ref{cor:dimension} requires a Riemann-Roch formula for the  Euler characteristic  
\begin{equation*} 
\chi(\Ecal,\Gcal) = \sum_{i} \, (-1)^i\, \dim_\C
\Ext^i(\Ecal,\Gcal) \ ,
\end{equation*}
where $\Ecal$ and $\Gcal$ are coherent sheaves on $\Xscr_k$. One can prove (as in \cite[Lemma~6.1.3]{book:huybrechtslehn2010}) that
\begin{equation*} 
\chi(\Ecal,\Gcal) =  \int_{\IXkscr} \,
\frac{\ch^\vee\big(\rho(\varpi^\ast\Ecal)\big)\cdot
  \ch\big(\rho(\varpi^\ast\Gcal) \big)}
{\ch\big(\rho(\lambda_{-1} \mathcal N^\vee) \big)} \cdot
 \operatorname{Td}(\Xscr_k) \ ,
\end{equation*}
where for an element $x\otimes\omega\in K(\IXkscr)\otimes\mu_{\infty}$ one defines
\begin{equation*} 
\ch^\vee(x\otimes\omega) := \sum_i \, (-1)^i \, \ch_i(x)\otimes\omega^{-1}\ .
\end{equation*}
Here $\mu_{\infty}\subset \C$ is the group of all roots of unity.  If
$\Ecal$ is locally free then $\ch^\vee(\rho(\varpi^\ast\Ecal))= \ch(\rho(\varpi^\ast\Ecal^\vee)).$

By Proposition \ref{prop:vanishing} and the proof of Theorem \ref{thm:moduli} we have
\begin{equation*}
\dim_{\mathbb{C}} \Mcal_{r,\vec{u},\Delta}\big(\Xscr_k, \Dscr_\infty,
\Fcal^{s,\vec{w}}_\infty\, \big) = -\chi\big(\Ecal,\Ecal\otimes
\Ocal_{\Xscr_k}(-\Dscr_\infty) \big) \ ,
\end{equation*}
where $[(\Ecal,\phi_{\Ecal})]$ is a point in the moduli space. By the
same argument as in Section \ref{sec:eulerchar} we obtain the required formula.

\begin{example}
As we saw in Section \ref{sec:rankonecase},
$\mathscr{M}_{1,\vec{u},n}(\Xscr_k, \Dscr_\infty,
\Ocal_{\Dscr_\infty}(s,i))$ is isomorphic to the Hilbert scheme ${\rm
  Hilb}^n(X_k)$ of $\Delta=n$ points on $X_k$. In the rank one
case one has $\vec w\cdot \vec w(j)=0$ for all $j\geq1$ and the
dimension formula
\eqref{eq:dimensionformula} agrees with the dimension of ${\rm
  Hilb}^n(X_k)$. Likewise, when $w_i=r$ for some
$i\in\{0,1,\dots,k-1\}$ and $w_j=0$ for all $j\neq i$, the dimension
of the moduli space is $2\, r\,\Delta$.
 
For rank $r\geq2$ we can give another formulation of the dimension
formula \eqref{eq:dimensionformula}. Let us consider a locally free
sheaf on $\Xscr_k$ of the form
\begin{equation}\label{eq:example}
\Ecal:=\bigoplus_{\alpha=1}^r\, \Mcal_\alpha=\bigoplus_{\alpha=1}^r \
\bigotimes_{i=1}^{k-1}\, \mathcal{R}_i^{\otimes
  u_{\alpha,i}}\otimes\Ocal_{\Xscr_k}(s \, \Dscr_\infty)
\end{equation}
for some choice of integers $u_{\alpha,j}$. The discriminant of $\Ecal$ is
\begin{equation*}
\Delta(\Ecal) = \frac{r-1}{2r}\, \sum_{\alpha=1}^{r} \ \sum_{i,j=1}^{k-1}
\, u_{\alpha,i}\, (C^{-1})^{ij}\,u_{\alpha,j}-
\frac1{2r} \, \sum_{ \alpha\neq \beta} \ \sum_{i,j=1}^{k-1}\, u_{\alpha,i}\,
(C^{-1})^{ij}\,u_{\beta,j} \ ,
\end{equation*}
and by Lemma \ref{lem:restrictionDinfty} and Corollary
\ref{cor:restrictiontautological} we have
\begin{equation*}
\Mcal_{\alpha\vert\Dscr_\infty} \simeq \Ocal_{\Dscr_\infty}\Big(s \,,\,
\mbox{$\sum\limits_{i=1}^{k-1}$}\, i \,u_{\alpha,i} \Big) \ .
\end{equation*}
Hence $\Ecal$ has an induced framing $\phi_{\Ecal}$ to
$\Fcal_\infty^{s,\vec w}$, where $\vec w = (w_0,w_1, \dots,w_{k-1})$ with
\begin{equation*}
w_l :=\#\Big\{\alpha\,\Big\vert\, \mbox{$\sum\limits_{i=1}^{k-1}$}\, i\,u_{\alpha,i} =
l \ \bmod{k}\Big\}\ .
\end{equation*}
Then $[(\Ecal,\phi_{\Ecal})]$ is a point of the moduli space
$\Mcal_{r,\vec{u},\Delta}(\Xscr_k, \Dscr_\infty,
\Fcal^{s,\vec{w}}_\infty\, )$, where 
\begin{equation*}
\vec{u}:=\Big(\, \mbox{$\sum\limits_{\alpha=1}^r\, u_{\alpha, 1}\,,\,
  \ldots\,,\, \sum\limits_{\alpha=1}^r u_{ \alpha , k-1}$}\, \Big) \ .
\end{equation*}
By the arguments of Section \ref{sec:eulerchar}, we find that the
dimension of the moduli space (cf. Equation
\eqref{eq:dimensionformula}) is
\begin{equation*}
(r-1) \, \sum_{\alpha=1}^{r} \ \sum_{i,j=1}^{k-1}
\, u_{\alpha,i}\, (C^{-1})^{ij}\,u_{\alpha,j}-
\sum_{ \alpha\neq \beta} \ \sum_{i,j=1}^{k-1}\, u_{\alpha,i}\,
(C^{-1})^{ij}\,u_{\beta,j}-\frac{1}{2}\, \sum_{j=1}^{k-1}\,
(C^{-1})^{jj}\, \vec{w}\cdot\vec{w}(j)\ .
\end{equation*}
For $k=2$ this formula reduces to 
\begin{equation*}
\frac{r}{2}\, \sum_{\alpha=1}^r u_\alpha^2-\frac{1}{2}\,
\big(u^2+w_0\, w_1 \big)\ ,
\end{equation*}
where $u=\sum_{\alpha=1}^r \, u_\alpha$. By considering explicitly all possible parities of $w_0$ and $w_1$ one can check that this number is an integer.
\end{example}

\bigskip

\section{Summation identities for complex roots of unity}\label{ap:rootsunity}
In this appendix we collect a number of identities that were used in
Appendix \ref{sec:rankdimension}. In the following $\omega$ will
denote a complex $k$-th root of unity, and $\eta$ a complex $k\,
\tilde{k}$-th root of unity.

\begin{lemma}\label{lem:simpleidentity} For $s=1,\ldots,k$, we have
\begin{equation*}
\sum_{i=1}^{k-1}\, \frac{\omega^{s\, i}}{(1 -
  \omega^{-i})^2}=-\frac{(k-5)\, (k-1)}{12}+\frac{s\, (k-2-s)}{2}\ .
\end{equation*}
\end{lemma}
\proof
By the same arguments as in the proof of Equation \eqref{eq:bryanidentity} in \cite{art:bryan2012}, one can prove the identity
\begin{equation}\label{eq:identitybryan}
\sum_{i=1}^{k-1}\, \frac{\omega^{s\, i}}{(1-\omega^{-i})^2} = \sum_{ l
  =0}^{k-1}\, (s- l )\, \Big( \, - \left\{\frac{s- l }{k}\right\} +
\frac{k-1}{2k} \, \Big)\ .
\end{equation}
From this identity we get
\begin{align}
\sum_{i=1}^{k-1}\, \frac{\omega^{s\, i}}{(1 - \omega^{-i})^2}=& -
\sum_{ l =0}^{s-1}\, (s- l )\, \frac{s- l }{k} - \sum_{ l
  =s+1}^{k-1}\, (s- l )\, \frac{k-s- l }{k} +\frac{k-1}{2k}\, \sum_{ l
  =0}^{k-1}\, (s- l )\\[4pt] 
=&-\sum_{m=1}^{k-1}\, \frac{m^2}{k} - \sum_{m=s+1}^{k-1}\, m +
\frac{k-1}{2}\, \Big(\, s-\frac{k-1}{2} \, \Big)\ .
\end{align}
By standard algebraic manipulations one then gets the assertion.
\endproof

We now prove various identities which were used for the calculation of
the $C$ and $D$ contributions in Section \ref{sec:eulerchar}.
We divide these identities according to the parity of $k$. We start with odd $k$.
\begin{lemma}
For any fixed $1\leq i\leq k-1$ and $x\in \C\setminus\mu_k$, we have
\begin{equation}
\prod_{\stackrel{\scriptstyle j=1}{\scriptstyle j\neq i}}^{k-1}\,
(x-\omega^j)=-\sum_{n=0}^{k-2}\, x^n\ \sum_{ l =1}^{n+1}\,
\omega^{- l \,i} \qquad \mbox{and} \qquad
\sum_{i=1}^{k-1}\, \frac1{x-\omega^i}= \frac{\sum\limits_{n=0}^{k-2}\,
  (n+1)\, x^n}{\sum\limits_{n=0}^{k-1}\, x^n} \ .
\end{equation}
\label{prelimLemma-main}\end{lemma}
\proof
We start by noting that, as pointed out in \cite{art:bryan2012}, one has  $\sum_{i=0}^{k-1}\, \omega^{s\,i}=k$ if $s= \, 0
\bmod{k}$ and $\sum_{i=0}^{k-1}\, \omega^{s\,i}=0$ otherwise; moreover,
\begin{equation}
\prod_{\stackrel{\scriptstyle j=1}{\scriptstyle j\neq i}}^{k-1}\,
(x-\omega^j)=\frac1{x-\omega^i}\, \sum_{n=0}^{k-1}\, x^n =:
\sum_{n=0}^{k-2}\, c_n\, x^n \ .
\end{equation}
The polynomial coefficients $n!\,c_n$ can be obtained by differentiating
the second expression $n$ times with respect to $x$ at $x=0$, and it is
straightforward to prove by induction that
\begin{equation}
c_n=-\sum_{ l =1}^{n+1}\, \omega^{- l \, i} \ .
\end{equation}
Note in particular that $c_{k-2}=-\sum_{ l =1}^{k-1}\, \omega^{- l \,
  i}=1$ as expected.

For the second identity, we write
\begin{equation}
\sum_{i=1}^{k-1}\, \frac1{x-\omega^i}=\frac{\sum\limits_{i=1}^{k-1}\
  \prod\limits_{\stackrel{\scriptstyle j=1}{\scriptstyle j\neq
      i}}^{k-1}\, \, (x-\omega^j)}{\prod\limits_{i=1}^{k-1}\,
  (x-\omega^i)} \ .
\end{equation}
From above we have
\begin{equation}
\prod_{i=1}^{k-1}\, (x-\omega^i)=\frac{x^k-1}{x-1}= \sum_{n=0}^{k-1}\,
x^n
\end{equation}
and
\begin{equation}
\sum\limits_{i=1}^{k-1}\
  \prod\limits_{\stackrel{\scriptstyle j=1}{\scriptstyle j\neq
      i}}^{k-1}\, \, (x-\omega^j) = -\sum_{n=0}^{k-2}\, x^n \
  \sum_{ l =1}^{n+1}\ \sum_{i=1}^{k-1}\, \omega^{- l \,j}=
  \sum_{n=0}^{k-2}\, (n+1)\, x^n \ ,
\end{equation}
and the result follows.
\endproof

\begin{lemma}
$\qquad \displaystyle{\sum_{i=1}^{k-1}\,
  \frac{\omega^{2i}}{1+\omega^i} =-\frac{k+1}2 \ .}$
\label{prelimLemma1}\end{lemma}
\proof
Since $k$ is odd, setting $x=-1$ in Lemma~\ref{prelimLemma-main} gives
\begin{equation}
\prod_{i=1}^{k-1}\, (1+\omega^i) = \sum_{n=0}^{k-1}\, (-1)^n=1
\end{equation}
and
\begin{align}
\sum_{i=1}^{k-1}\, \omega^{2i} \ \prod\limits_{\stackrel{\scriptstyle j=1}{\scriptstyle j\neq i}}^{k-1}\, \, (1+\omega^j) &= \sum_{n=0}^{k-2}\, (-1)^n\ \sum_{ l =1}^{n+1}\ \sum_{i=1}^{k-1}\, \omega^{-( l -2)\, i} \\[4pt]
&=-1-\sum_{n=1}^{\frac{k-1}2}\, 2n+\sum_{n=1}^{\frac{k-1}2}\, (2n-1) \
= \ -\frac{k+1}2 \ .
\end{align}
\endproof

\begin{lemma}
For any $1\leq j\leq k-1$, one has
\begin{equation}
\sum_{i=1}^{k-1}\, \frac{\omega^{i\, (j+2)}+\omega^{-i\,
    (j-2) }}{1+\omega^i} = -1 \ .
\end{equation}
\label{prelimLemma2}\end{lemma}
\proof
Putting $x=-1$ in Lemma~\ref{prelimLemma-main} again we find
\begin{align}
\sum_{i=1}^{k-1}\, \frac{\omega^{i\, (j+2)}}{1+\omega^i} &=
\sum_{n=1}^{k-2}\, (-1)^n\ \sum_{ l =1}^{n+1} \ \sum_{i=1}^{k-1}\,
\omega^{-i\, ( l -j-2)} \\[4pt] &=\sum_{n=0}^j\, (-1)^{n+1}\,
(n+1)+\sum_{n=j+1}^{k-1}\, (-1)^n\, (k-1)+\sum_{n=j+1}^{k-1}\,
(-1)^{n+1}\, n \ .
\end{align}
For odd $j$ this gives
\begin{equation}
\frac{j+1}2+(k-1)-\frac{k-1}2-\frac{j+1}2=\frac{k-1}2
\end{equation}
while for even $j$ we get
\begin{equation}
-\frac j2-1-\frac{k-1}2+\frac j2=-\frac{k+1}2 \ .
\end{equation}
Now the sum
\begin{equation}
\sum_{i=1}^{k-1}\, \frac{\omega^{-i\, (j-2)}}{1+\omega^i} =
\sum_{i=1}^{k-1}\, \frac{\omega^{i\, (k-j+2)}}{1+\omega^i}
\end{equation}
is computed in an identical way by just replacing $j$ with
$k-j$. Since $j$ and $k-j$ have opposite parity, we get the result.
\endproof

\begin{lemma}
Let $\eta$ be a $k$-th root of $\omega$ and $1\leq
i\leq k-1$. Then
\begin{equation}
\frac1k\, \sum_{j=0}^{k-1}\, \frac1{\eta^i\,
  \omega^j-1}=\frac1{\omega^i-1} \ .
\end{equation}
\label{prelimLemma3}\end{lemma}
\proof
Using Lemma~\ref{prelimLemma-main} with $x=\eta^{-i}$ we compute
\begin{align}
\frac{1}{k}\, \sum_{j=0}^{k-1}\, \frac{1}{\eta^i\, \omega^j-1} &= -\frac{\eta^{-i}}{k}\, \Bigg(\, \frac{1}{\eta^{-i}-1}+\frac{\sum\limits_{n=0}^{k-2}\, (n+1)\, \eta^{-i\,n}}{\sum\limits_{n=0}^{k-1}\, \eta^{-i\,n}}\, \Bigg) \\[4pt]
&= -\frac{\eta^{-i}}k\, \frac{2\eta^{-i\,(k-1)}+ (k-2)\,
  \eta^{-i\,(k-1)}}{\eta^{-i\,k}-1} \ = \ \frac{1}{\omega^i-1} \ .
\end{align}
\endproof

Now we consider the case of even $k$.

\begin{lemma} Let $\eta$ be a $\tilde k$-th root of $\omega$ and $1\leq i\leq \tilde k-1$. Then
\begin{equation} 
\frac1k\, \sum_{j=0}^{k-1}\, \frac1{\eta^i\, \omega^j-1} =
\frac1{\omega^{2i}-1} \qquad \mbox{and} \qquad
\frac1k\, \sum_{j=0}^{k-1}\, \frac{(-1)^j}{\eta^i\, \omega^j-1} =
\frac{\omega^i}{\omega^{2i}-1} \ .
\end{equation}
\label{prelimLemma4}\end{lemma}
\proof
The first identity follows exactly as in the proof of
Lemma~\ref{prelimLemma3}, except that now $\eta^k=\omega^2$. For the
second identity, we proceed   as in the proof of
Lemma~\ref{prelimLemma3}. Using $\omega^{\tilde k}=-1$ we first compute
\begin{align}
\sum_{j=1}^{k-1}\, (-1)^j \ \sum_{n=0}^{k-2}\, \eta^{-i\,n}\ \sum_{ l =1}^{n+1}\, \omega^{- l \,j} &= -\sum_{n=0}^{\tilde k-2}\, \eta^{-i\,n}\, (n+1)+\sum_{n=\tilde k-1}^{k-2}\, \eta^{-i\,n}\, (k-1)-\sum_{n=\tilde k-1}^{k-2}\, \eta^{-i\,n}\, n \\[4pt] 
&= -\eta^i\, \sum_{n=0}^{\tilde k-1}\, \big(n\, (1+\omega^{-i})-\tilde k\, \omega^{-i}\big)\, \eta^{-i\, n} \ .
\end{align}
Using also
\begin{equation}
\sum_{n=0}^{k-1}\, \eta^{-i\,n}=\Big(\, \sum_{n=0}^{\tilde
k-1}+\sum_{n=\tilde k}^{2\tilde k-1}\, \Big)\, \eta^{-i\,n} =
(1+\omega^{-i})\, \sum_{n=0}^{\tilde k-1}\, \eta^{-i\,n} \ ,
\end{equation}
we arrive at
\begin{align}
\frac1k\, \sum_{j=0}^{k-1}\, \frac{(-1)^j}{\eta^i\, \omega^j-1} &= -\frac{\eta^{-i}}{k}\, \Bigg(\, \frac1{\eta^{-i}-1}+\frac{\eta^i}{1+\omega^{-i}}\, \frac{\sum\limits_{n=0}^{\tilde k-1}\, \big(n\, (1+\omega^{-i})-\tilde k\, \omega^{-i}\big)\, \eta^{-i\,   n}}{\sum\limits_{n=0}^{\tilde k-1}\, \eta^{-i\,n}}\, \Bigg) \\[4pt]
&= -\frac{\eta^{-i}}{k}\, \frac{-\tilde k\, \omega^{-i}\,
  \omega^{-i}-1)\, \eta^i+(1+\omega^{-i})\, \tilde k\, \eta^i\,
  \omega^{-i}}{(\omega^{-i}-1)\, (1+\omega^{-i})} \ = \ \frac{\omega^i}{\omega^{2i}-1} \ .
\end{align}
\endproof

\bigskip\section{Edge contribution}\label{sec:edgecontribution}

In this appendix we compute the equivariant Euler characteristic
\begin{equation*}
\chi_{T_t}\big(\Rcal^{\vec u}\otimes \Ocal_{\Xscr_k}(-\Dscr_\infty)
\big):= \chi_{T_t}\big(\Xscr_k\,,\,\Rcal^{\vec u}\otimes
\Ocal_{\Xscr_k}(-\Dscr_\infty)\big)
\end{equation*} 
for a vector $\vec{u}\in \Z^{k-1}$. 

Set $\vec{v}:=C^{-1}\vec{u}$. Then $\vec{v}\in\frac{1}{k}\, \Z^{k-1}$. 
If $c$ is the equivalence class of $\sum_{j=1}^{k-1}\, j\, u_j$ modulo
$k$, then for $ l=1, \ldots, k-1$ we have
\begin{equation}\label{eq:congruences}
k\, v_{ l} = - l \, c \bmod{k}\ .
\end{equation}
As a consequence, a component $v_ l$ is integral if and only if every
component is (cf. Remark \ref{rem:v-y}). Define
$\vec{z}:=\vec{u}-\vec{e}_{c}$ if $c>0$ and $\vec{z}:=\vec{u}$ otherwise, where $\vec{e}_{c}$ is the $c$-th coordinate vector of $\mathbb{Z}^{k-1}$. Set also $\vec{s}:=C^{-1}\vec{z}$. Then by construction $\vec{s}\in\Z^{k-1}$.

\subsection{Preliminaries}

\begin{lemma}\label{lem:exactseqtautological}
Given a vector $\vec{u}\in\Z^{k-1}$, for every $ l=1, \ldots, k-1$ there is an exact sequence
\begin{equation}\label{eq:exactsequencei}
0 \ \longrightarrow \ \Rcal^{\vec{u}+C\vec{e}_ l} \ \longrightarrow \
\Rcal^{\vec{u}} \ \longrightarrow \ \Rcal^{\vec{u}}_{\vert \Dscr_ l} \
\longrightarrow \ 0
\end{equation}
where $C$ is the Cartan matrix of type $A_{k-1}$.
\end{lemma}
\proof 
Fix $ l=1, \ldots, k-1$ and consider the short exact sequence
\begin{equation}
0 \ \longrightarrow \ \Ocal_{\Xscr_k}(-\Dscr_ l) \ \longrightarrow \
\Ocal_{\Xscr_k} \ \longrightarrow \ \Ocal_{\Dscr_ l} \ \longrightarrow
\ 0\ .
\end{equation}
We obtain the assertion simply by  tensoring this sequence with
$\Rcal^{\vec{u}}$; we need only prove that $\Rcal^{\vec{u}}\otimes
\Ocal_{\Xscr_k}(-\Dscr_ l) = \Rcal^{\vec{u}+C\vec{e}_ l}$. By
definition $\Rcal^{\vec{u}}=\Ocal_{\Xscr_k}\big(\sum_{i=1}^{k-1}\, u_i
\, \omega_i \big)$, and we have
\begin{equation}
\sum_{i=1}^{k-1} u_i \, \omega_i-\Dscr_ l=
\sum_{i=1}^{k-1} \, u_i \, \omega_i+\sum_{i=1}^{k-1}\, C_{ l i}\,
\omega_i=\sum_{i=1}^{k-1}\, \big(u_i+(C\vec{e}_ l)_i \big)\, \omega_i \ .
\end{equation}
\endproof
By Equation \eqref{eq:tautologicalintegralrelations} we get the following result.
\begin{lemma}\label{lem:restrictions} 
Let $\vec{u}\in\Z^{k-1}$ and $ l=1,\ldots, k-1$. Then $\Rcal^{\vec{u}}_{\vert \Dscr_ l} \simeq \Ocal_{\Dscr_ l}(u_ l)$, 
where $\Ocal_{\Dscr_ l}(1):={\pi_k}_{\vert \Dscr_ l}^\ast \Ocal_{D_ l}(1)$.
\end{lemma}

\subsection{Iterative procedure}

Let us recall that by the arguments in Appendix
\ref{sec:eulerchar} the dimension of $H^1(\Rcal^{\vec{u}}\otimes
\Ocal_{\Xscr_k}(-\Dscr_\infty))$ is $\frac{1}{2}\, (\vec{v}\cdot C\vec{v}-(C^{-1})^{c c})$, where we set  $(C^{-1})^{c c}=0$ if $c=0$. 

First, let us assume that $s_ l \geq 0$ for every $ l=1, \ldots, k-1$. Consider the equation
\begin{equation}\label{eq:condition-index+}
\frac{C_{ l  l}}{2}\, i^2- i\, \Big(\vec{v}-\sum_{p=1}^{ l-1}\, s_p\,
\vec{e}_p\Big)\cdot C\vec{e}_ l+\frac{1}{2}\,
\bigg(\Big(\vec{v}-\sum_{p=1}^{ l-1}\, s_p\, \vec{e}_p\Big)\cdot
C\Big(\vec{v}-\sum_{p=1}^{ l-1}\, s_p\, \vec{e}_p\Big)-\big(C^{-1} \big)^{c c}\bigg)=0\ ,
\end{equation}
and define the set 
\begin{equation}
S_ l^+:=\{i\in\N\, \vert\, i\leq s_ l \ \mbox{ is a solution of Equation \eqref{eq:condition-index+}}\}\ .
\end{equation}
Let $m$ be the smallest integer $ l\in\{1, \ldots, k-1\}$ such that $S_ l^+$ is nonempty; if all sets are empty, let $m:=k-1$.

Now we can compute $\chi_{T_t}(\Rcal^{\vec{u}}\otimes \Ocal_{\Xscr_k}(-\Dscr_\infty))$ by using Lemma \ref{lem:exactseqtautological}. Recall that by convention $\Rcal_0:=\Ocal_{\Xscr_k}$. 

Let $d_1^+:=\min(S_1^+)$ if $S_1^+\neq \emptyset$, otherwise
$d_1^+:=s_1$. By using the exact sequence \eqref{eq:exactsequencei}
$d_1^+$ times for $ l=1$, we obtain 
\begin{multline}
\chi_{T_t}\big(\Rcal^{\vec{u}}\otimes \Ocal_{\Xscr_k}(-\Dscr_\infty)
\big)= \chi_{T_t}\big(\Rcal^{\vec{u}-C\vec{e}_1} \otimes
\Ocal_{\Xscr_k}(-\Dscr_\infty) \big) -
\chi_{T_t} \big(\Rcal^{\vec{u}-C\vec{e}_1} \otimes
\Ocal_{\Xscr_k}(-\Dscr_\infty)_{\vert \Dscr_1} \big) \\[4pt]
\vdots \\[4pt]
= \chi_{T_t}\big(\Rcal^{\vec{u}-d_1^+\, C\vec{e}_1} \otimes
\Ocal_{\Xscr_k}(-\Dscr_\infty) \big) - \sum_{i=1}^{d_1^+}\,
\chi_{T_t}\big(\Rcal^{\vec{u}-i\,C \vec{e}_1} \otimes
\Ocal_{\Xscr_k}(-\Dscr_\infty)_{\vert \Dscr_1} \big)\ .
\end{multline}
If $m=1$, we conclude the inductive procedure. If not, set $d_2^+:=
\min(S_2^+)$ if $S_2^+\neq \emptyset$, otherwise $d_2^+:=s_2$. Now we do   $d_2^+$ steps more with the sequence \eqref{eq:exactsequencei} for $ l=2$ and obtain
 
\begin{multline}
\chi_{T_t}\big(\Rcal^{\vec{u}} \otimes
\Ocal_{\Xscr_k}(-\Dscr_\infty) \big) \\
\shoveleft{= \chi_{T_t}\big(\Rcal^{\vec{u}-d_1^+\, C\vec{e}_1 -
    C\vec{e}_2} \otimes \Ocal_{\Xscr_k}(-\Dscr_\infty) \big) -
  \chi_{T_t}\big(\Rcal^{\vec{u}-d_1^+\, C\vec{e}_1 - C\vec{e}_2} \otimes \Ocal_{\Xscr_k}(-\Dscr_\infty)_{\vert \Dscr_2}) }\\
\shoveright{- \, \sum_{i=1}^{d_1^+}\,\chi_{T_t}\big(\Rcal^{\vec{u}-i\,
    C\vec{e}_1} \otimes \Ocal_{\Xscr_k}(-\Dscr_\infty)_{\vert \Dscr_1}
  \big)}\\[4pt]
 \vdots \\[4pt]
\shoveleft{= \chi_{T_t}\big(\Rcal^{\vec{u}-C(d_1^+\, \vec{e}_1+d_2^+\,
    \vec{e}_2)} \otimes \Ocal_{\Xscr_k}(-\Dscr_\infty) \big)
 - \sum_{i=1}^{d_1^+}\, \chi_{T_t}\big(\Rcal^{\vec{u}-i\, C\vec{e}_1}
 \otimes \Ocal_{\Xscr_k}(-\Dscr_\infty)_{\vert \Dscr_1} \big)}\\
- \sum_{i=1}^{d_2^+}\, \chi_{T_t}\big(\Rcal^{\vec{u}-C(d_1^+\,
  \vec{e}_1+i\, \vec{e}_2)} \otimes
\Ocal_{\Xscr_k}(-\Dscr_\infty)_{\vert \Dscr_2} \big)\ .
\end{multline}

Iterating this procedure $d_ l^+$ times using the sequence \eqref{eq:exactsequencei}
for $ l=3, \ldots, m$, where $d_ l^+:= \min(S_ l^+)$ if
$S_ l^+\neq \emptyset$ and $d_ l^+:=s_ l$ otherwise, we get
\begin{multline}
\chi_{T_t}\big(\Rcal^{\vec{u}}\otimes \Ocal_{\Xscr_k}(-\Dscr_\infty) \big) \\
\shoveleft{ =\chi_{T_t}\big(\Rcal^{\vec{u}-\sum_{p=1}^m\limits\, d_p^+
    \, C \vec{e}_p} \otimes \Ocal_{\Xscr_k}(-\Dscr_\infty) \big) -
  \sum_{ l=1}^{m}  \ \sum_{i=1}^{d_ l^+} \, \chi_{T_t}\big(\Rcal^{
    \vec{u}-C(\, \sum_{p=1}^{ l-1}\limits \, s_p \, \vec{e}_p + i\, \vec{e}_ l ) } \otimes \Ocal_{\Xscr_k}(-\Dscr_\infty)_{\vert \Dscr_ l} \big)}\ .
\end{multline}
Note that if $d_m^+\in S_m^+$, then
$\chi_{T_t}\big(\Rcal^{\vec{u}-\sum_{p=1}^m \, d_p^+ \, C \vec{e}_p}
\otimes \Ocal_{\Xscr_k}(-\Dscr_\infty) \big)=0$ by Equation
\eqref{eq:condition-index+} (as then the left-hand side of the equation
is exactly the dimension of $H^1\big(\Xscr_k,
\Rcal^{\vec{u}-\sum_{p=1}^m\, d_p^+ \, C  \vec{e}_p} \otimes \Ocal_{\Xscr_k}(-\Dscr_\infty)\big)$). Otherwise
\begin{multline}
\chi_{T_t}\big(\Rcal^{\vec{u}-\sum_{p=1}^m\limits\, d_p^+ \, C
  \vec{e}_p} \otimes \Ocal_{\Xscr_k}(-\Dscr_\infty)\big) \\
=
\chi_{T_t}\big(\Rcal_c\otimes\Rcal^{\vec{z}-\sum_{p=1}^{k-1}\limits\,
  s_p\, C \vec{e}_p} \otimes \Ocal_{\Xscr_k}(-\Dscr_\infty)
\big)=\chi_{T_t}\big(\Rcal_{c} \otimes \Ocal_{\Xscr_k}(-\Dscr_\infty) \big)\ ,
\end{multline}
and the last characteristic is equal to zero thanks to the computations in Section \ref{sec:eulerchar}.

When $s_ l<0$ for $ l=1,\ldots, k-1$ consider the equation
\begin{equation}\label{eq:condition-index-}
\frac{C_{ l  l}}{2}\, i^2+i\, \Big(\vec{v}-\sum_{p=1}^{ l-1}\, s_p\,
\vec{e}_p\Big)\cdot C\vec{e}_ l +\frac{1}{2}\,
\bigg(\Big(\vec{v}-\sum_{p=1}^{ l-1}\, s_p\, \vec{e}_p\Big)\cdot
C\Big(\vec{v}-\sum_{p=1}^{ l-1}\, s_p\, \vec{e}_p\Big)-\big(C^{-1} \big)^{c c}\bigg)=0\ ,
\end{equation}
and define the set 
\begin{equation}
S_ l^-:=\{i\in\N\, \vert\, i\leq -s_ l \ \mbox{ is a solution of Equation \eqref{eq:condition-index-}}\}\ .
\end{equation}
One can follow the procedure just described by using the short exact
sequence \eqref{eq:exactsequencei} $d_ l^-$ times, where $d_ l^-:=
\min(S_ l^-)$ if $S_ l^-\neq \emptyset$ and $d_ l^-:=-s_ l$ otherwise. In this case, one exchanges the roles played by the left and middle terms of the sequence. 

In general, let $m$ be the smallest  integer $ l\in\{1, \ldots, k-1\}$ such that $S_ l^+$ or $S_ l^-$ is nonempty; if all these sets are empty,  let $m:=k-1$. Define for $ l=1, \ldots, m$ the \emph{$L$-factors} as  
\begin{equation*}
L^{( l)}\big(\chi_1^{( l)},\chi_2^{( l)} \big):= \left\{
\begin{array}{ll}
- \sum_{i=1}^{d_ l^+}\limits\, \chi_{T_t}\big(\Rcal^{ \vec{u}-C(\,
  \sum_{p=1}^{ l-1}\limits\, s_p \, \vec{e}_p + i \, \vec{e}_ l) } \otimes \Ocal_{\Xscr_k}(-\Dscr_\infty)_{\vert \Dscr_ l} \big) & \mbox{ for }s_ l\geq 0 \ ,\\[8pt]
\sum_{i=0}^{d_ l^- -1}\limits\, \chi_{T_t}\big(\Rcal^{ \vec{u}-C(\,
  \sum_{p=1}^{ l-1}\limits s_p \, \vec{e}_p - i \, \vec{e}_ l ) } \otimes \Ocal_{\Xscr_k}(-\Dscr_\infty)_{\vert \Dscr_ l} \big) & \mbox{ for }s_ l<0\ .
\end{array}
\right.
\end{equation*}
Then we obtain
\begin{equation}
\chi_{T_t}\big(\Rcal^{\vec{u}}\otimes \Ocal_{\Xscr_k}(-\Dscr_\infty)
\big)  =\sum_{ l=1}^{m}\, L^{( l)}\big(\chi_1^{( l)},\chi_2^{( l)} \big)\ .
\end{equation}
It remains just to compute the $L$-factors. By Lemma \ref{lem:restrictions}, we have
\begin{equation}
\Rcal^{ \vec{u}-C(\, \sum_{p=1}^{ l-1}\limits\, s_p\, \vec{e}_p \pm
  i\, \vec{e}_ l) } \otimes \Ocal_{\Xscr_k}(-\Dscr_\infty)_{\vert
  \Dscr_ l} \simeq  \Ocal_{\Dscr_ l}(\delta_{ l,c} + z_ l + s_{
  l-1} \mp 2i)\ .
\end{equation}
Recalling that $z_l = (C\vec{s}\, )_l = -s_{l-1} +2s_l -s_{l+1}$ for $ l=1, \ldots, k-1$, we can rewrite the $L$-factors as
\begin{equation}\label{eq:gammafactorscrestricted}
L^{( l)}\big(\chi_1^{( l)},\chi_2^{( l)} \big)=\left\{
\begin{array}{ll}
- \sum_{i=s_ l-d_ l^+}^{s_ l-1}\limits\, \chi_{T_t}\big( \Ocal_{\Dscr_ l}(\delta_{ l,c} -s_{ l+1} +2i) \big)\quad & \mbox{for }s_ l> 0\ ,\\[8pt]
\qquad \qquad \qquad \qquad 0 & \mbox{for } s_ l=0\ ,\\[8pt]
\sum_{i=1-s_ l-d_ l^-}^{-s_ l}\limits\, \chi_{T_t}\big( \Ocal_{\Dscr_ l}(\delta_{ l,c} -s_{ l+1} -2i) \big)\quad & \mbox{for }s_ l<0\ .
\end{array}
\right.
\end{equation}

\subsection{Characters of restrictions and $L$-factors}\label{sec:edgecontribution-result}

Here we choose the $T_t$-equivariant structure on the line bundles
$\Ocal_{\Dscr_ l}(a)$ given by the isomorphism
\begin{equation}\label{eq:secondiso}
\Ocal_{\Dscr_ l}(a) \simeq  \Ocal_{\Xscr_k}\big( -\left\lfloor
  \mbox{$\frac{a}{2}$} \right\rfloor \, \Dscr_ l +
2\left\{\mbox{$\frac{a}{2}$} \right\} \, \Dscr_{ l+1} \big)_{\vert \Dscr_ l}\ .
\end{equation}
\begin{theorem}\label{thm:restrictions}
Fix $ l\in\{1, \ldots, k-1\}$. For $a\geq 0$ we have
\begin{equation}
\chi_{T_t}\big(\Ocal_{\Dscr_ l}(a) \big) = \big( \chi_1^ l
\big)^{\lfloor \frac{a}{2} \rfloor}\, \sum_{j=0}^{a}\, \big(
\chi_2^ l \big)^j \qquad \mbox{and} \qquad
\chi_{T_t}\big(\Ocal_{\Dscr_ l}(-a) \big) = -\big( \chi_1^ l
\big)^{-\lfloor \frac{a}{2} \rfloor} \, \sum_{j=1}^{a-1}\, \big( \chi_2^ l \big)^{-j}\ .
\end{equation}
\end{theorem}
\proof
Let $a\geq0$ and consider the short exact sequence
\begin{multline}
0 \ \longrightarrow \ \Ocal_{\Xscr_k}\left(\left( -\left\lfloor
      \mbox{$\frac{a}{2}$} \right\rfloor - 1 \right)\, \Dscr_ l +
  2\left\{\mbox{$\frac{a}{2}$} \right\} \, \Dscr_{ l+1} \right) \\
\longrightarrow \ \Ocal_{\Xscr_k}\left( -\left\lfloor
    \mbox{$\frac{a}{2}$} \right\rfloor\, \Dscr_ l +
  2\left\{\mbox{$\frac{a}{2}$} \right\} \, \Dscr_{ l+1} \right) \
\longrightarrow \ \Ocal_{\Xscr_k}\left( -\left\lfloor
    \mbox{$\frac{a}{2}$} \right\rfloor \, \Dscr_ l + 2\left\{
    \mbox{$\frac{a}{2}$} \right\} \, \Dscr_{ l+1} \right)_{\vert
  \Dscr_ l} \ \longrightarrow \ 0\ .
\end{multline}
Then for the Euler characteristic we have
\begin{multline}
\chi_{T_t}\left( \Ocal_{\Xscr_k}\left( -\left\lfloor
      \mbox{$\frac{a}{2}$} \right\rfloor \, \Dscr_ l +
    2\left\{\mbox{$\frac{a}{2}$}\right\} \, \Dscr_{ l+1} \right)_{\vert \Dscr_ l} \right) \\
= \chi_{T_t}\Big( \Ocal_{\Xscr_k}\left( -\left\lfloor
    \mbox{$\frac{a}{2}$} \right\rfloor \, \Dscr_ l +
  2\left\{\mbox{$\frac{a}{2}$}\right\} \, \Dscr_{ l+1} \right) \Big) -
\chi_{T_t}\Big( \Ocal_{\Xscr_k}\left(\left( -\left\lfloor
      \mbox{$\frac{a}{2}$} \right\rfloor - 1 \right)\, \Dscr_ l +
  2\left\{\mbox{$\frac{a}{2}$}\right\}\, \Dscr_{ l+1} \right) \Big) \\
=\chi_{T_t}\Big( \Ocal_{\bar{X}_k}\left( -\left\lfloor
    \mbox{$\frac{a}{2}$} \right\rfloor\, D_ l +
  2\left\{\mbox{$\frac{a}{2}$}\right\}\, D_{ l+1} \right) \Big) -
\chi_{T_t}\Big( \Ocal_{\bar{X}_k}\left(\left( -\left\lfloor
      \mbox{$\frac{a}{2}$} \right\rfloor - 1 \right) \, D_ l + 2\left\{\mbox{$\frac{a}{2}$}\right\}\, D_{ l+1} \right) \Big)\ ,
\end{multline}
where the last equality follows from the fact that the pushforward ${\pi_k}_\ast$ preserves the equivariant decomposition of the cohomology groups. To complete the proof it is sufficient to compute for $m\geq0$
\begin{equation}\label{eq:case-1}
\chi_{T_t}\big(\Ocal_{\bar{X}_k}(-m\, D_ l) \big) - \chi_{T_t}\big(
\Ocal_{\bar{X}_k}(-(m+1)\, D_ l) \big)
\end{equation}
which corresponds to the case $\{\frac{a}{2}\}=0$, and
\begin{equation}\label{eq:case-2}
\chi_{T_t}\big( \Ocal_{\bar{X}_k}(D_{ l+1}-m\, D_ l) \big) -
\chi_{T_t}\big( \Ocal_{\bar{X}_k}(D_{ l+1}-(m+1)\, D_ l) \big)
\end{equation}
which corresponds to the case $\{\frac{a}{2}\}=\frac{1}{2}$.

By \cite[Proposition 9.1.6]{book:coxlittleschenck2011}, it is easy to
verify that the zeroth and second cohomology groups that appear in
Equation \eqref{eq:case-1} vanish. To compute the first cohomology
groups in \eqref{eq:case-1} it is enough to count the integer points on
the line of direction $( l-1, l)$ between the points $(-( l-2)\, m,-(
l-1)\, m)$ and $( l \, m,( l+1)\, m)$. We easily get 
\begin{multline}
\chi_{T_t}\big(\Ocal_{{X}_k}(-m\, D_ l) \big) - \chi_{T_t}\big(
\Ocal_{{X}_k}(-(m+1)\, D_ l) \big) \\
= \sum_{j=0}^{2m} \, T_1^{-(
  l-2)\, m+j\, ( l-1)} \, T_2^{-( l-1)\, m+j\, l}
= \big(\chi_1^ l \big)^m \, \sum_{j=0}^{2m} \, \big(\chi_2^ l \big)^j\ ,
\end{multline}
where the last equality follows from the expression
\eqref{eq:variablerelations} for the variables $T_1$ and $T_2$ introduced in Section \ref{sec:minimalresolution}. 

In the same way, Equation \eqref{eq:case-2} becomes
\begin{multline}
\chi_{T_t}\big( \Ocal_{{X}_k}(D_{ l+1}-m\, D_ l) \big) -
\chi_{T_t}\big( \Ocal_{{X}_k}(D_{ l+1}-(m+1)\, D_ l) \big) \\ =
\sum_{j=0}^{2m+1}\, T_1^{-( l-2)\, m+j\, ( l-1)} \, T_2^{-( l-1)\,
  m+j\, l} 
 = \big(\chi_1^ l\big)^m \, \sum_{j=0}^{2m+1} \, \big(\chi_2^ l \big)^j\ .
\end{multline}
For $a<0$ one argues in a similar way.
\endproof

Now we use Theorem \ref{thm:restrictions} to compute the explicit
expressions for the $L$-factors $L^{( l)}\big(\chi_1^{( l)},\chi_2^{(
  l)} \big)$ in \eqref{eq:gammafactorscrestricted}. Fix $ l\in\{1,
\ldots, k-1\}$ and  set $(C^{-1})^{ l, 0}=0$. We also 
set $(C^{-1})^{k, c}=0$. Then we get:
\begin{itemize} \setlength{\itemsep}{4mm}
\item[\scriptsize$\blacksquare$] For $v_ l-(C^{-1})^{ l c}> 0$:
\smallskip
\begin{itemize} \setlength{\itemsep}{0.8cm}
\item[$\bullet$]For $\delta_{ l,c} -v_{ l+1}+(C^{-1})^{ l+1, c}+2(v_{ l}-(C^{-1})^{ l c}-d_ l^+)\geq 0$:
\begin{equation}
L^{( l)}\big(\chi_1^{( l)},\chi_2^{( l)} \big) =- \sum_{i=v_{
    l}-(C^{-1})^{ l c}-d_ l^+}^{v_ l-(C^{-1})^{ l c}-1}\hspace{0.2cm}
\sum_{j=0}^{2i+\delta_{ l,c} -v_{ l+1}+(C^{-1})^{ l+1, c}}\,
\big(\chi_1^ l \big)^{i+\big\lfloor\frac{\delta_{ l,c} -v_{
      l+1}+(C^{-1})^{l+1, c}}{2} \big\rfloor}\, \big(\chi_2^ l \big)^j\ .
\end{equation}
\item[$\bullet$]For $2\leq \delta_{ l,c}-v_{ l+1}+(C^{-1})^{ l+1, c}+2(v_ l-(C^{-1})^{ l c})<2 d_ l^+$:
\begin{multline}
L^{( l)}\big(\chi_1^{( l)},\chi_2^{( l)} \big) \\
\shoveleft{=\sum_{i=v_ l-(C^{-1})^{ l c}-d_ l^+}^{-\big\lfloor
    \frac{\delta_{ l,c} -v_{ l+1}+(C^{-1})^{ l+1, c}}{2}\big\rfloor -
    1}\hspace{6mm} \sum_{j=1}^{2i- (\delta_{ l,c} -v_{
      l+1}+(C^{-1})^{ l+1, c})-1} \, \big(\chi_1^ l
  \big)^{i-\big\lfloor -\frac{\delta_{ l,c} -v_{ l+1}+(C^{-1})^{ l+1,
        c}}{2}\big\rfloor}\, \big(\chi_2^ l \big)^{-j}} \\[4pt]
 - \, \sum_{i=-\big\lfloor \frac{\delta_{ l,c} -v_{ l+1}+(C^{-1})^{
       l+1, c}}{2}\big\rfloor}^{2(v_ l-(C^{-1})^{ l c})+\delta_{ l,c}
   -v_{ l+1}+(C^{-1})^{ l+1, c} -2}
 \hspace{6mm}\sum_{j=0}^{2i+\delta_{ l,c}-v_{ l+1}+(C^{-1})^{ l+1,
     c}}\, \big(\chi_1^ l \big)^{i+\big\lfloor \frac{\delta_{ l,c}
     -v_{ l+1}+(C^{-1})^{ l+1, c}}{2}\big\rfloor} \, \big(\chi_2^ l \big)^j\ .
\end{multline}
\item[$\bullet$]For $\delta_{ l,c} -v_{ l+1}+(C^{-1})^{l+1, c}<2-2(v_ l-(C^{-1})^{ l c})$:
\begin{multline}
L^{( l)}\big(\chi_1^{( l)},\chi_2^{( l)} \big) \\
=\sum_{i=v_{
    l}-(C^{-1})^{ l c}-d_ l^+}^{v_ l-(C^{-1})^{ l c}-1}\hspace{0.2cm}
\sum_{j=1}^{-2i-\delta_{ l,c} +v_{ l+1}-(C^{-1})^{l+1,c}-1} \,
\big(\chi_1^ l \big)^{i-\big\lfloor - \frac{\delta_{ l,c} -v_{
      l+1}+(C^{-1})^{l+1,c}}{2} \big\rfloor} \, \big(\chi_2^ l \big)^{-j} \ .
\end{multline}
\end{itemize}
\item[\scriptsize$\blacksquare$] For $v_ l-(C^{-1})^{ l c}= 0$:
  $\qquad L^{( l)} \big(\chi_1^{( l)},\chi_2^{( l)} \big)=0$.
\item[\scriptsize$\blacksquare$] For $v_ l-(C^{-1})^{ l c}< 0$:
\smallskip
\begin{itemize}\setlength{\itemsep}{0.5cm}
\item[$\bullet$]For $\delta_{ l,c} -v_{ l+1}+(C^{-1})^{l+1,c} +2v_{ l}-2(C^{-1})^{ l c} < 2-2 d_ l^{-}$:
\begin{multline}
L^{( l)}\big(\chi_1^{( l)},\chi_2^{( l)} \big) \\
=-\sum_{i=1-v_{
    l}+(C^{-1})^{ l c}-d_ l^-}^{-v_{ l}+(C^{-1})^{ l c}}\hspace{0.2cm}
\sum_{j=1}^{2i-(\delta_{ l,c} -v_{ l+1}+(C^{-1})^{l+1,c})-1} \,
\big(\chi_1^ l \big)^{-i-\big\lfloor -\frac{\delta_{ l,c} -v_{
      l+1}+(C^{-1})^{l+1,c}}{2}\big\rfloor} \, \big(\chi_2^ l \big)^{-j}\ .
\end{multline}
\item[$\bullet$]For $2-2 d_ l^{-}\leq \delta_{ l,c} -v_{ l+1}+(C^{-1})^{l+1,c}+2v_{ l}-2(C^{-1})^{ l c}<0$:
\begin{multline}
L^{( l)}\big(\chi_1^{( l)},\chi_2^{( l)} \big) \\
\shoveleft{=\sum_{i=1-v_{ l}+(C^{-1})^{ l c}-d_ l^-}^{\big\lfloor
    \frac{\delta_{ l,c} -v_{ l+1}+(C^{-1})^{l+1,c}}{2}
    \big\rfloor}\hspace{0.2cm} \sum_{j=0}^{-2i+\delta_{ l,c} -v_{
      l+1}+(C^{-1})^{l+1,c}} \, \big(\chi_1^ l \big)^{-i+\big\lfloor
    \frac{\delta_{ l,c} -v_{ l+1}+(C^{-1})^{l+1,c}}{2} \big\rfloor} \,
  \big(\chi_2^ l \big)^j}\\
-\, \sum_{i=\big\lfloor \frac{\delta_{ l,c} -v_{
      l+1}+(C^{-1})^{l+1,c}}{2} \big\rfloor+1}^{-v_{ l}+(C^{-1})^{ l
    c}}\hspace{0.2cm} \sum_{j=1}^{2i-(\delta_{ l,c} -v_{
    l+1}+(C^{-1})^{l+1,c})-1 } \, \big(\chi_1^ l \big)^{-i-
  \big\lfloor -\frac{\delta_{ l,c} -v_{
      l+1}+(C^{-1})^{l+1,c}}{2}\big\rfloor} \, \big(\chi_2^ l \big)^{-j}\ .
\end{multline}
\item[$\bullet$]For $\delta_{ l,c} -v_{ l+1}+(C^{-1})^{l+1,c}\geq -2v_{ l}+2(C^{-1})^{ l c}$:
\begin{equation}
L^{( l)}\big(\chi_1^{( l)},\chi_2^{( l)} \big) =\sum_{i=1-v_{
    l}+(C^{-1})^{ l c}-d_ l^-}^{-v_{ l}+(C^{-1})^{ l c}}\hspace{0.2cm}
\sum_{j=0}^{-2i+\delta_{ l,c}-v_{ l+1}+(C^{-1})^{l+1,c}} \,
\big(\chi_1^ l \big)^{-i+\big\lfloor \frac{\delta_{ l,c} -v_{
      l+1}+(C^{-1})^{l+1,c}}{2}\big\rfloor}\, \big(\chi_2^ l \big)^{j}\ .
\end{equation}
\end{itemize}
\end{itemize}

\begin{remark}\label{rem:edge-choice}
The procedure we described in this appendix involves some
choices. In particular, one can consider other realizations of
$\Ocal_{\Dscr_l}(a)$ as a restriction of a line bundle on $\Xscr_k$
for $l=1, \ldots, k-1$. We chose the realization in Equation
\eqref{eq:secondiso} because it resembles similar choices made in
\cite[Section~4.3]{art:bruzzopoghossiantanzini2011} and in the proof
of \cite[Theorem~3.4]{art:nakajimayoshioka2005-I}. Different
realizations of the line bundles $\Ocal_{\Dscr_l}(a)$ yield equivalent results in equivariant K-theory.
\end{remark}

\section{Gauge theory on $X_3$\label{app:k=3}}

In this appendix we present some explicit calculations for $k=3$.

\subsection{Euler classes}

For $k=3$, the $\tilde T$-equivariant Euler class of the
Carlsson-Okounkov bundle from Section \ref{sec:Eulerclasses} becomes
\begin{equation}
\eu_{\tilde{T}}\big(\Ebf_{\left([(\Ecal,\phi_\Ecal)] \,,\, [(\Ecal',\phi_{\Ecal'})]\right)}\big) =
\prod_{\alpha=1}^r \ \prod_{\beta=1}^{r'} \ \prod_{i=1}^3 \,
m_{Y_{\alpha}^{i}, {Y_{\beta}^{i}}'}
\big(\varepsilon_1^{(i)},\varepsilon_2^{(i)}, a^{(i)}_{\beta\alpha}
\big) \ \prod_{n=1}^{2}\,
\ell^{(n)}_{\vec{v}_{\beta\alpha}}\big(\varepsilon_1^{(n)},\varepsilon_2^{(n)},
a_{\beta\alpha} \big)\ .
\end{equation}
In this case $c_{\beta\alpha}\in\{0,1,2\}$ for any $\alpha=1, \ldots,
r$, $\beta=1,\ldots,r'$. By Equation \eqref{eq:congruences} we have
$(C^{-1})^{1, c_{\beta\alpha}}=\{(\vec{v}_{\beta\alpha})_1\}$ and
$(C^{-1})^{2,c_{\beta\alpha}}=\{(\vec{v}_{\beta\alpha})_2\}$. Hence
the edge factor
$\ell^{(1)}_{\vec{v}_{\beta\alpha}}\big(\varepsilon_1^{(1)},\varepsilon_2^{(1)},
a_{\beta\alpha} \big)$ assumes the following form:
\begin{itemize}
\item[\scriptsize$\blacksquare$] For $\lfloor (\vec{v}_{\beta\alpha})_1\rfloor> 0$:
\smallskip
\begin{itemize} 
\item[$\bullet$]For $\delta_{1,c_{\beta\alpha}} -\lfloor (\vec{v}_{\beta\alpha})_2\rfloor+2(\lfloor (\vec{v}_{\beta\alpha})_1\rfloor-\dplus_1)\geq 0$:
\small
\begin{multline}
\ell^{(1)}_{\vec{v}_{\beta\alpha}}\big(\varepsilon_1^{(1)},
\varepsilon_2^{(1)}, a_{\beta\alpha} \big) \\
=\prod_{i=\lfloor (\vec{v}_{\beta\alpha})_1\rfloor-\dplus_1}^{\lfloor (\vec{v}_{\beta\alpha})_1\rfloor-1}\hspace{0.2cm} \prod_{j=0}^{2i+\delta_{1,c_{\beta\alpha}} -\lfloor (\vec{v}_{\beta\alpha})_2\rfloor}
\,
\bigg(a_{\beta\alpha}+\Big(i+\left\lfloor\frac{\delta_{1,c_{\beta\alpha}}
    -\lfloor (\vec{v}_{\beta\alpha})_2\rfloor}{2} \right\rfloor\Big)\,
\varepsilon_1^{(1)}+j\, \varepsilon_2^{(1)}\bigg)\ .
\end{multline}
\normalsize
\item[$\bullet$]For $2\leq \delta_{1,c_{\beta\alpha}}-\lfloor (\vec{v}_{\beta\alpha})_2\rfloor+2\lfloor (\vec{v}_{\beta\alpha})_1\rfloor<2 \dplus_1$:
\small
\begin{multline}
\ell^{(1)}_{\vec{v}_{\beta\alpha}}\big(\varepsilon_1^{(1)},
\varepsilon_2^{(1)}, a_{\beta\alpha} \big) \\
\shoveleft{=\prod_{i=\lfloor
    (\vec{v}_{\beta\alpha})_1\rfloor-\dplus_1}^{-\big\lfloor
    \frac{\delta_{1,c_{\beta\alpha}} -\lfloor
      (\vec{v}_{\beta\alpha})_2\rfloor}{2}\big\rfloor -
    1}\hspace{0.2cm} \prod_{j=1}^{2i- (\delta_{1,c_{\beta\alpha}}
    -\lfloor (\vec{v}_{\beta\alpha})_2\rfloor)-1} \,
\bigg(a_{\beta\alpha}+\Big(i-\left\lfloor
  -\frac{\delta_{1,c_{\beta\alpha}} -\lfloor
    (\vec{v}_{\beta\alpha})_2\rfloor}{2}\right\rfloor\Big)\,
\varepsilon_1^{(1)}-j\, \varepsilon_2^{(1)}\bigg)^{-1}} \\[4pt]
\times \ \prod_{i=-\big\lfloor \frac{\delta_{1,c_{\beta\alpha}} -\lfloor (\vec{v}_{\beta\alpha})_2\rfloor}{2}\big\rfloor}^{2(\lfloor \vec{v}_{\beta\alpha})_1\rfloor+\delta_{1,c_{\beta\alpha}} -\lfloor (\vec{v}_{\beta\alpha})_2\rfloor -2} \hspace{0.2cm}\prod_{j=0}^{2i+\delta_{1,c_{\beta\alpha}}-\lfloor (\vec{v}_{\beta\alpha})_2\rfloor}
  \, \bigg(a_{\beta\alpha}+\Big(i+\left\lfloor \frac{\delta_{1,c_{\beta\alpha}} -\lfloor (\vec{v}_{\beta\alpha})_2\rfloor}{2}\right\rfloor\Big)\, \varepsilon_1^{(1)}-j\, \varepsilon_2^{(1)}\bigg)\ .
\end{multline}
\normalsize
\item[$\bullet$]For $\delta_{1,c_{\beta\alpha}} -\lfloor (\vec{v}_{\beta\alpha})_2\rfloor<2-2\lfloor (\vec{v}_{\beta\alpha})_1\rfloor$:
\small
\begin{multline}
\ell^{(1)}_{\vec{v}_{\beta\alpha}}\big(\varepsilon_1^{(1)},
\varepsilon_2^{(1)}, a_{\beta\alpha} \big) \\
=\prod_{i=\lfloor (\vec{v}_{\beta\alpha})_1\rfloor-\dplus_1}^{\lfloor (\vec{v}_{\beta\alpha})_1\rfloor-1}\hspace{0.2cm} \prod_{j=1}^{-2i-\delta_{1,c_{\beta\alpha}} +\lfloor (\vec{v}_{\beta\alpha})_2\rfloor-1}
\, \bigg(a_{\beta\alpha}+\Big(i-\left\lfloor -
  \frac{\delta_{1,c_{\beta\alpha}} -\lfloor
    (\vec{v}_{\beta\alpha})_2\rfloor}{2} \right\rfloor\Big)\,
\varepsilon_1^{(1)}-j\, \varepsilon_2^{(1)}\bigg)^{-1} \ .
\end{multline}
\normalsize
\end{itemize}
\item[\scriptsize$\blacksquare$] For $\lfloor (\vec{v}_{\beta\alpha})_1\rfloor= 0$:
\begin{equation}
\ell^{(1)}_{\vec{v}_{\beta\alpha}}\big(\varepsilon_1^{(1)},
\varepsilon_2^{(1)}, a_{\beta\alpha} \big)=1\ .
\end{equation}
\item[\scriptsize$\blacksquare$] For $\lfloor (\vec{v}_{\beta\alpha})_1\rfloor< 0$:
\smallskip
\begin{itemize}
\item[$\bullet$]For $\delta_{1,c_{\beta\alpha}} -\lfloor (\vec{v}_{\beta\alpha})_2\rfloor+2\lfloor (\vec{v}_{\beta\alpha})_1\rfloor < 2-2\dminus_1$:
\small
\begin{multline}
\ell^{(1)}_{\vec{v}_{\beta\alpha}}\big(\varepsilon_1^{(1)},
\varepsilon_2^{(1)}, a_{\beta\alpha} \big) \\
=\prod_{i=1-\lfloor (\vec{v}_{\beta\alpha})_1\rfloor-\dminus_1}^{-\lfloor (\vec{v}_{\beta\alpha})_1\rfloor}\hspace{0.2cm} \prod_{j=1}^{2i-(\delta_{1,c_{\beta\alpha}} -\lfloor (\vec{v}_{\beta\alpha})_2\rfloor)-1}
 \bigg(a_{\beta\alpha}-\Big(i+\left\lfloor
   -\frac{\delta_{1,c_{\beta\alpha}} -\lfloor
     (\vec{v}_{\beta\alpha})_2\rfloor}{2}\right\rfloor\Big)\,
 \varepsilon_1^{(1)}-j\, \varepsilon_2^{(1)}\bigg)\ .
\end{multline}
\normalsize
\item[$\bullet$]For $2-2 \dminus_1\leq \delta_{1,c_{\beta\alpha}} -\lfloor (\vec{v}_{\beta\alpha})_2\rfloor+2\lfloor (\vec{v}_{\beta\alpha})_1\rfloor<0$:
\small
\begin{multline}
\ell^{(1)}_{\vec{v}_{\beta\alpha}}\big(\varepsilon_1^{(1)},
\varepsilon_2^{(1)}, a_{\beta\alpha} \big) \\
\shoveleft{=\prod_{i=\big\lfloor \frac{\delta_{1,c_{\beta\alpha}} -\lfloor (\vec{v}_{\beta\alpha})_2\rfloor}{2} \big\rfloor+1}^{-\lfloor (\vec{v}_{\beta\alpha})_1\rfloor}\hspace{0.2cm} \prod_{j=1}^{2i-(\delta_{1,c_{\beta\alpha}} -\lfloor (\vec{v}_{\beta\alpha})_2\rfloor)-1 }
\bigg(a_{\beta\alpha}-\Big(i+\left\lfloor
  -\frac{\delta_{1,c_{\beta\alpha}} -\lfloor
    (\vec{v}_{\beta\alpha})_2\rfloor}{2}\right\rfloor\Big)\,
\varepsilon_1^{(1)}-j\, \varepsilon_2^{(1)}\bigg)}\\
\times \ \prod_{i=1-\lfloor (\vec{v}_{\beta\alpha})_1\rfloor-\dminus_1}^{\big\lfloor \frac{\delta_{1,c_{\beta\alpha}} -\lfloor (\vec{v}_{\beta\alpha})_2\rfloor}{2} \big\rfloor}\hspace{0.2cm} \prod_{j=0}^{-2i+\delta_{1,c_{\beta\alpha}} -\lfloor (\vec{v}_{\beta\alpha})_2\rfloor}
\, \bigg(a_{\beta\alpha}+\Big(-i+\left\lfloor
  \frac{\delta_{1,c_{\beta\alpha}} -\lfloor
    (\vec{v}_{\beta\alpha})_2\rfloor}{2} \right\rfloor\Big)\,
\varepsilon_1^{(1)}+j\, \varepsilon_2^{(1)}\bigg)^{-1}\ .
\end{multline}
\normalsize
\item[$\bullet$]For $\delta_{1,c_{\beta\alpha}} -\lfloor (\vec{v}_{\beta\alpha})_2\rfloor\geq -2\lfloor (\vec{v}_{\beta\alpha})_1\rfloor$:
\small
\begin{multline}
\ell^{(1)}_{\vec{v}_{\beta\alpha}}\big(\varepsilon_1^{(1)},
\varepsilon_2^{(1)}, a_{\beta\alpha} \big) \\
=\prod_{i=1-\lfloor (\vec{v}_{\beta\alpha})_1\rfloor-\dminus_1}^{-\lfloor (\vec{v}_{\beta\alpha})_1\rfloor}\hspace{0.2cm} \prod_{j=0}^{-2i+\delta_{1,c_{\beta\alpha}}-\lfloor (\vec{v}_{\beta\alpha})_2\rfloor}
\, \bigg(a_{\beta\alpha}+\Big(-i+\left\lfloor
  \frac{\delta_{1,c_{\beta\alpha}} -\lfloor
    (\vec{v}_{\beta\alpha})_2\rfloor}{2}\right\rfloor\Big)\,
\varepsilon_1^{(1)}+j\, \varepsilon_2^{(1)}\bigg)^{-1}\ .
\end{multline}
\normalsize
\end{itemize}
\end{itemize}
If $m_{\beta\alpha}=1$, then
\begin{equation}
\ell^{(2)}_{\vec{v}_{\beta\alpha}}\big(\varepsilon_1^{(2)},
\varepsilon_2^{(2)}, a_{\beta\alpha} \big)=1\ ,
\end{equation}
otherwise
\small
\begin{multline*}
\ell^{(2)}_{\vec{v}_{\beta\alpha}}\big(\varepsilon_1^{(2)},\varepsilon_2^{(2)},a_{\beta\alpha}
\big) \\
=\left\{
\begin{array}{cl}
\prod_{i=\lfloor (\vec{v}_{\beta\alpha})_2\rfloor-\dplus_2}^{\lfloor
  (\vec{v}_{\beta\alpha})_2\rfloor-1}\limits\hspace{0.2cm}
\prod_{j=0}^{2i+\delta_{2,c_{\beta\alpha}}}\limits\,
\big(a_{\beta\alpha}+i\, \varepsilon_1^{(2)}+j\, \varepsilon_2^{(2)}\big) & \mbox{for } \lfloor (\vec{v}_{\beta\alpha})_2\rfloor> 0\ ,\\[8pt]
1 & \mbox{for } \lfloor (\vec{v}_{\beta\alpha})_2\rfloor= 0\ ,\\[8pt]
\prod_{i=1-\lfloor
  (\vec{v}_{\beta\alpha})_2\rfloor-\dminus_2}^{-\lfloor
  (\vec{v}_{\beta\alpha})_{2}\rfloor}\limits\hspace{0.2cm}
\prod_{j=1}^{2i-\delta_{2,c_{\beta\alpha}}-1}\limits\,
\big(a_{\beta\alpha}-\big\lfloor
-\frac{\delta_{2,c_{\beta\alpha}}}{2}\big\rfloor\,
\varepsilon_1^{(2)}-(i\, \varepsilon_1^{(2)}+j\, \varepsilon_2^{(2)})\big) & \mbox{for } \lfloor (\vec{v}_{\beta\alpha})_2\rfloor<0\ .
\end{array}
\right.
\end{multline*}
\normalsize

\subsection{$U(2)$ gauge theory}

In this subsection we provide some explicit computations for the instanton
partition function \eqref{eq:instantonpart-puregaugetheory} of pure
$\Ncal=2$ gauge theory on $X_3$ for rank $r=2$. For $k=3$ the $U(2)$ partition function \eqref{eq:instantonpart-puregaugetheory} becomes
\begin{equation}
\Zcal^{\mathrm{inst}}_{X_3}(\varepsilon_1, \varepsilon_2, a_1,a_2 ;
\qsf, \xi_1,\xi_2 )= \sum_{\stackrel{\scriptstyle
    \vec{v}=(\upsilon_1,\upsilon_2)\in\frac{1}{3}\, \Z^2
  }{\scriptstyle 3\upsilon_2=  w_1+2w_2 \bmod{3}}}\,
  \xi_1^{\upsilon_1}\, \xi_2^{\upsilon_2} \
  \sum_{\vec{v}_1+\vec{v}_2=\vec{v}}\,
  \Zcal_{\vec{v}_1,\vec{v}_2}(\varepsilon_1, \varepsilon_2, a_1,a_2 ; \qsf) \ ,
\end{equation}
where
\begin{equation}
\Zcal_{\vec{v}_1,\vec{v}_2}(\varepsilon_1, \varepsilon_2, a_1,a_2;
\qsf):=\frac{\qsf^{\frac{1}{2}\, \sum_{\alpha=1}^2\limits\,
    \vec{v}_\alpha\cdot
    C\vec{v}_\alpha}}{\prod_{\alpha,\beta=1}^2\limits \
  \prod_{l=1}^{2}\limits\,
  \ell^{(l)}_{\vec{v}_{\beta\alpha}}\big(\varepsilon_1^{(l)},
  \varepsilon_2^{(l)}, a_{\beta\alpha} \big)} \ \prod_{i=1}^3\,
\Zcal_{\C^2}^{\mathrm{inst}} \big(\varepsilon_1^{(i)},
\varepsilon_2^{(i)}, a_1^{(i)},a_2^{(i)} ; \qsf \big) \ .
\end{equation}

We compute the lowest order of the expansion of this partition
function in $\qsf$ for $-1\leq u_1, u_2\leq 1$, where $(u_1,u_2):=C\vec{v}$. Let $I_{\vec{v}}^{(w_0,w_1,w_2)}$ be the set consisting of pairs $(\vec{v}_1, \vec{v}_2)$ such that $ \vec{v}_1+\vec{v}_2=\vec{v}$ and
\begin{equation}
3(\vec{v}_1)_i=-i\, l \bmod{3} \qquad\mbox{for } \ i=1,2 
\end{equation}
if $\sum_{j=0}^{l-1}\, w_j<1\leq \sum_{j=0}^{l}\, w_j$ for $l=0, 1,
2$. We denote by $\min\big(I_{\vec{v}}^{(w_0,w_1,w_2)} \big)$ the set
consisting of the vectors which provide the smallest power of $\qsf$. Set $a=\frac{a_1-a_2}{2}$.

\subsubsection*{{\scriptsize$\blacksquare$} $(w_0, w_1, w_2)=(2, 0, 0)$.} We need to consider the
cases $\vec{v}=(0,0)$, $(1,1)$ and $(-1,-1)$.
\begin{itemize}[leftmargin=0.2cm]
\item For $\vec{v}=(0,0)$ the set $\min\big(I_{(0,0)}^{(2,0,0)}\big)$
  is $\big\{\big((0,0)\,,\, (0,0)\big) \big\}$ and
\begin{equation}
\sum_{(\vec{v}_1, \vec{v}_2)\in I_{(0,0)}^{(2,0,0)}}\, \Zcal_{\vec{v}_1,\vec{v}_2}(\varepsilon_1, \varepsilon_2, {a}; \qsf)=1+\cdots\ .
\end{equation}
\item For $\vec{v}=(1,1)$ the set $\min\big(I_{(1,1)}^{(2,0,0)} \big)$
  is $\big\{\big((0,0)\,,\, (1,1) \big), \big((1,1)\,,\,(0,0) \big)\big\}$ and
\begin{equation}
\sum_{(\vec{v}_1, \vec{v}_2)\in I_{(1,1)}^{(2,0,0)}}\, \Zcal_{\vec{v}_1,\vec{v}_2}(\varepsilon_1, \varepsilon_2, a; \qsf)=
\qsf\, \Big(\, -\frac{1}{2a\, (2a-\varepsilon_1-\varepsilon_2)}-\frac{1}{2a\, (2a+\varepsilon_1+\varepsilon_2)}\, \Big)+\cdots\ .
\end{equation}
\item For $\vec{v}=(-1,-1)$ the set $\min\big(I_{(-1,-1)}^{(2,0,0)}
  \big)$ is $\big\{\big((0,0)\,,\, (-1,-1) \big),
  \big((-1,-1)\,,\,(0,0) \big) \big\}$ and
\begin{equation}
\sum_{(\vec{v}_1, \vec{v}_2)\in I_{(-1,-1)}^{(2,0,0)}}\, \Zcal_{\vec{v}_1,\vec{v}_2}(\varepsilon_1, \varepsilon_2, a; \qsf)=
\qsf\, \Big(\, -\frac{1}{2a\, (2a-\varepsilon_1-\varepsilon_2)}-\frac{1}{2a\, (2a+\varepsilon_1+\varepsilon_2)}\, \Big)+\cdots\ .
\end{equation}
\end{itemize}

\subsubsection*{{\scriptsize$\blacksquare$}  $(w_0, w_1, w_2)=(0, 2, 0)$.} We need to consider the cases $\vec{v}=(\frac{1}{3},\frac{2}{3})$, $(-\frac{2}{3},-\frac{1}{3})$ and $(\frac{1}{3},-\frac{1}{3})$.
\begin{itemize}[leftmargin=0.2cm]
\item For $\vec{v}=(\frac{1}{3},\frac{2}{3})$ the set
  $\min\big(I_{(\frac{1}{3},\frac{2}{3})}^{(0,2, 0)} \big)$ is
  $\big\{\big((\frac{2}{3},\frac{1}{3})\,,\,
  (-\frac{1}{3},\frac{1}{3}) \big),
  \big((-\frac{1}{3},\frac{1}{3})\,,\, (\frac{2}{3},\frac{1}{3})\big) \big\}$ and 
\begin{equation}
\sum_{(\vec{v}_1, \vec{v}_2)\in
  I_{(\frac{1}{3},\frac{2}{3})}^{(0,2,0)}}\, \Zcal_{\vec{v}_1,\vec{v}_2}(\varepsilon_1, \varepsilon_2, a; \qsf)=
\qsf^{\frac{2}{3}}\, \Big(\, -\frac{1}{2a\, (2a-\varepsilon_1-\varepsilon_2)}-\frac{1}{2a\, (2a+\varepsilon_1+\varepsilon_2)}\, \Big)+\cdots\ .
\end{equation}
\item For $\vec{v}=(-\frac{2}{3},-\frac{1}{3})$ the set
  $\min\big(I_{(-\frac{2}{3},-\frac{1}{3})}^{(0,2, 0)} \big)$ is
  $\big\{\big((-\frac{1}{3},\frac{1}{3})\,,\,
  (-\frac{1}{3},-\frac{2}{3})\big),
  \big((-\frac{1}{3},-\frac{2}{3})\,,\,
  (-\frac{1}{3},\frac{1}{3})\big) \big\}$ and 
\begin{equation}
\sum_{(\vec{v}_1, \vec{v}_2)\in
  I_{(-\frac{2}{3},-\frac{1}{3})}^{(0,2,0)}}\, \Zcal_{\vec{v}_1,\vec{v}_2}(\varepsilon_1, \varepsilon_2, a; \qsf)=
\qsf^{\frac{2}{3}}\, \Big(\, -\frac{1}{2a\, (2a-\varepsilon_1-\varepsilon_2)}-\frac{1}{2a\, (2a+\varepsilon_1+\varepsilon_2)}\, \Big)+\cdots\ .
\end{equation}
\item For $\vec{v}=(\frac{1}{3},-\frac{1}{3})$ the set
  $\min\big(I_{(\frac{1}{3},-\frac{1}{3})}^{(0,2, 0)} \big)$ is
  $\big\{\big((\frac{2}{3},\frac{1}{3})\,,\,
  (-\frac{1}{3},-\frac{2}{3})
  \big),\big((-\frac{1}{3},-\frac{2}{3})\,,\,
  (\frac{2}{3},\frac{1}{3}) \big) \big\}$ and 
\begin{equation}
\sum_{(\vec{v}_1, \vec{v}_2)\in
  I_{(\frac{1}{3},-\frac{1}{3})}^{(0,2,0)}}\, \Zcal_{\vec{v}_1,\vec{v}_2}(\varepsilon_1, \varepsilon_2, a; \qsf)=
\qsf^{\frac{2}{3}}\, \Big(\, -\frac{1}{2a\, (2a-\varepsilon_1-\varepsilon_2)}-\frac{1}{2a\, (2a+\varepsilon_1+\varepsilon_2)}\, \Big)+\cdots\ .
\end{equation}
\end{itemize}

\subsubsection*{{\scriptsize$\blacksquare$} $(w_0, w_1, w_2)=(0, 0, 2)$.} We need to consider separately the cases $\vec{v}=(\frac{2}{3},\frac{1}{3})$, $(-\frac{1}{3},-\frac{2}{3})$ and $(-\frac{1}{3},\frac{1}{3})$.
\begin{itemize}[leftmargin=0.2cm]
\item For $\vec{v}=(\frac{2}{3},\frac{1}{3})$ the set
  $\min\big(I_{(\frac{2}{3},\frac{1}{3})}^{(0,0,2)} \big)$ is
  $\big\{\big((\frac{1}{3},\frac{2}{3})\,,\,
  (\frac{1}{3},-\frac{1}{3})
  \big),\big((\frac{1}{3},-\frac{1}{3})\,,\,
  (\frac{1}{3},\frac{2}{3})\big) \big\}$ and 
\begin{equation}
\sum_{(\vec{v}_1, \vec{v}_2)\in
  I_{(\frac{2}{3},\frac{1}{3})}^{(0,0,2)}}\, \Zcal_{\vec{v}_1,\vec{v}_2}(\varepsilon_1, \varepsilon_2, a; \qsf)=
\qsf^{\frac{2}{3}}\, \Big(\, -\frac{1}{2a\, (2a-\varepsilon_1-\varepsilon_2)}-\frac{1}{2a\, (2a+\varepsilon_1+\varepsilon_2)}\, \Big)+\cdots\ .
\end{equation}
\item For $\vec{v}=(-\frac{1}{3},-\frac{2}{3})$ the set
  $\min\big(I_{(-\frac{1}{3},-\frac{2}{3})}^{(0,0,2)} \big)$ is
  $\big\{\big((-\frac{2}{3},-\frac{1}{3})\,,\,
  (\frac{1}{3},-\frac{1}{3})
  \big),\big((\frac{1}{3},-\frac{1}{3})\,,\,
  (-\frac{2}{3},-\frac{1}{3})\big) \big\}$ and 
\begin{equation}
\sum_{(\vec{v}_1, \vec{v}_2)\in
  I_{(-\frac{1}{3},-\frac{2}{3})}^{(0,0,2)}}\, \Zcal_{\vec{v}_1,\vec{v}_2}(\varepsilon_1, \varepsilon_2, a; \qsf)=
\qsf^{\frac{2}{3}}\, \Big(\, -\frac{1}{2a\, (2a-\varepsilon_1-\varepsilon_2)}-\frac{1}{2a\, (2a+\varepsilon_1+\varepsilon_2)}\, \Big)+\cdots\ .
\end{equation}
\item For $\vec{v}=(-\frac{1}{3},\frac{1}{3})$ the set
  $\min\big(I_{(-\frac{1}{3},\frac{1}{3})}^{(0,0,2)} \big)$ is
  $\big\{\big((\frac{1}{3},\frac{2}{3})\,,\,
  (-\frac{2}{3},-\frac{1}{3})
  \big),\big((-\frac{2}{3},-\frac{1}{3})\,,\,
  (\frac{1}{3},\frac{2}{3})\big) \big\}$ and 
\begin{equation}
\sum_{(\vec{v}_1, \vec{v}_2)\in
  I_{(-\frac{1}{3},\frac{1}{3})}^{(0,0,2)}}\, \Zcal_{\vec{v}_1,\vec{v}_2}(\varepsilon_1, \varepsilon_2, a; \qsf)=
\qsf^{\frac{2}{3}}\, \Big(\, -\frac{1}{2a\, (2a-\varepsilon_1-\varepsilon_2)}-\frac{1}{2a\, (2a+\varepsilon_1+\varepsilon_2)}\, \Big)+\cdots\ .
\end{equation}
\end{itemize}

\subsubsection*{{\scriptsize$\blacksquare$} $(w_0, w_1, w_2)=(1, 1, 0)$.} We need to consider separately the cases $\vec{v}=(\frac{2}{3}, \frac{1}{3})$, $(-\frac{1}{3}, \frac{1}{3})$ and $(-\frac{1}{3}, -\frac{2}{3})$.
\begin{itemize}[leftmargin=0.2cm]
\item For $\vec{v}=(\frac{2}{3},\frac{1}{3})$ the set
  $\min\big(I_{(\frac{2}{3},\frac{1}{3})}^{(1,1,0)} \big)$ is
  $\big\{\big((0,0)\,,\, (\frac{2}{3},\frac{1}{3})\big) \big\}$ and
\begin{equation}
\sum_{(\vec{v}_1, \vec{v}_2)\in
  I_{(\frac{2}{3},\frac{1}{3})}^{(1,1,0)}}\, \Zcal_{\vec{v}_1,\vec{v}_2}(\varepsilon_1, \varepsilon_2, a; \qsf)=\qsf^{\frac{1}{3}}+\cdots\ .
\end{equation}
\item For $\vec{v}=(-\frac{1}{3},\frac{1}{3})$ the set
  $\min\big(I_{(-\frac{1}{3},\frac{1}{3})}^{(1,1,0)} \big)$ is
  $\big\{\big((0,0)\,,\, (-\frac{1}{3},\frac{1}{3})\big) \big\}$ and
\begin{equation}
\sum_{(\vec{v}_1, \vec{v}_2)\in
  I_{(-\frac{1}{3},\frac{1}{3})}^{(1,1,0)}}\, \Zcal_{\vec{v}_1,\vec{v}_2}(\varepsilon_1, \varepsilon_2, a; \qsf)=\qsf^{\frac{1}{3}}+\cdots\ .
\end{equation}
\item For $\vec{v}=(-\frac{1}{3},-\frac{2}{3})$ the set
  $\min\big(I_{(-\frac{1}{3},-\frac{2}{3})}^{(1,1,0)} \big)$ is
  $\big\{\big((0,0)\,,\, (-\frac{1}{3},-\frac{2}{3})\big) \big\}$ and
\begin{equation}
\sum_{(\vec{v}_1, \vec{v}_2)\in
  I_{(-\frac{1}{3},-\frac{2}{3})}^{(1,1,0)}}\, \Zcal_{\vec{v}_1,\vec{v}_2}(\varepsilon_1, \varepsilon_2, a; \qsf)=
\qsf^{\frac{1}{3}}+\cdots\ .
\end{equation}
\end{itemize}

\subsubsection*{{\scriptsize$\blacksquare$}  $(w_0, w_1, w_2)=(1, 0, 1)$.} We need to consider separately the cases $\vec{v}=(\frac{1}{3}, \frac{2}{3})$, $(\frac{1}{3}, -\frac{1}{3})$ and $(-\frac{2}{3}, -\frac{1}{3})$.
\begin{itemize}[leftmargin=0.2cm]
\item For $\vec{v}=(\frac{1}{3}, \frac{2}{3})$ the set
  $\min\big(I_{(\frac{1}{3}, \frac{2}{3})}^{(1,0,1)} \big)$ is
  $\big\{\big((0,0)\,,\, (\frac{1}{3},\frac{2}{3})\big) \big\}$ and
\begin{equation}
\sum_{(\vec{v}_1, \vec{v}_2)\in I_{(\frac{1}{3},
    \frac{2}{3})}^{(1,0,1)}}\, \Zcal_{\vec{v}_1,\vec{v}_2}(\varepsilon_1, \varepsilon_2, a; \qsf)=
\qsf^{\frac{1}{3}}+\cdots\ .
\end{equation}
\item For $\vec{v}=(\frac{1}{3}, -\frac{1}{3})$ the set
  $\min\big(I_{(\frac{1}{3}, -\frac{1}{3})}^{(1,0,1)} \big)$ is
  $\big\{\big((0,0)\,,\, (\frac{1}{3},-\frac{1}{3})\big) \big\}$ and
\begin{equation}
\sum_{(\vec{v}_1, \vec{v}_2)\in I_{(\frac{1}{3},
    -\frac{1}{3})}^{(1,0,1)}}\, \Zcal_{\vec{v}_1,\vec{v}_2}(\varepsilon_1, \varepsilon_2, a; \qsf)=\qsf^{\frac{1}{3}}+\cdots\ .
\end{equation}
\item For $\vec{v}=(-\frac{2}{3}, -\frac{1}{3})$ the set
  $\min\big(I_{(-\frac{2}{3}, -\frac{1}{3})}^{(1,0,1)} \big)$ is
  $\big\{\big((0,0)\,,\, (-\frac{2}{3},-\frac{1}{3})\big) \big\}$ and
\begin{equation}
\sum_{(\vec{v}_1, \vec{v}_2)\in I_{(-\frac{2}{3},
    -\frac{1}{3})}^{(1,0,1)}}\, \Zcal_{\vec{v}_1,\vec{v}_2}(\varepsilon_1, \varepsilon_2, a; \qsf)=\qsf^{\frac{1}{3}}+\cdots\ .
\end{equation}
\end{itemize}

\subsubsection*{{\scriptsize$\blacksquare$}  $(w_0, w_1, w_2)=(0, 1, 1)$.} We need to consider   the cases $\vec{v}=(0,0)$, $(1,1)$ and $(-1,-1)$.
\begin{itemize}[leftmargin=0.2cm]
\item For $\vec{v}=(0,0)$ the set $\min\big(I_{(0,0)}^{(0,1,1)} \big)$
  is $\big\{ \big((\frac{2}{3},\frac{1}{3})\,,\,
  (-\frac{2}{3},-\frac{1}{3})\big), \big((-\frac{1}{3},
  -\frac{2}{3})\,,\,(\frac{1}{3},\frac{2}{3})
  \big),\big((-\frac{1}{3},\frac{1}{3})\,,\,
  (\frac{1}{3},-\frac{1}{3})\big) \big\}$ and 
\begin{multline}
\sum_{(\vec{v}_1, \vec{v}_2)\in I_{(0,0)}^{(0,1,1)}}\, \Zcal_{\vec{v}_1,\vec{v}_2}(\varepsilon_1, \varepsilon_2, a; \qsf)=
\qsf^{\frac{2}{3}}\, \Big(\, -\frac{1}{2a\,
  (2a+\varepsilon_1+\varepsilon_2)}-\frac{1}{(2a-3\varepsilon_2)\,
  (2a-\varepsilon_1-\varepsilon_2)} \\
-\, \frac{1}{(2a-\varepsilon_1-\varepsilon_2)\, (2a+\varepsilon_1+\varepsilon_2)}\, \Big)+\cdots\ .
\end{multline}
\item For $\vec{v}=(1,1)$ the set $\min\big(I_{(1,1)}^{(0,1,1)} \big)$
  is $\big\{\big((\frac{2}{3},\frac{1}{3})\,,\,
  (\frac{1}{3},\frac{2}{3})\big) \big\}$ and 
\begin{equation}
\sum_{(\vec{v}_1, \vec{v}_2)\in I_{(1,1)}^{(0,1,1)}}\, \Zcal_{\vec{v}_1,\vec{v}_2}(\varepsilon_1, \varepsilon_2, a; \qsf)=\qsf^{\frac{2}{3}}+\cdots \ .
\end{equation}
\item For $\vec{v}=(-1,-1)$  the set $\min\big(I_{(-1,-1)}^{(0,1,1)}
  \big)$ is $\big\{\big((-\frac{1}{3},-\frac{2}{3})\,,\,
  (-\frac{2}{3},-\frac{1}{3})\big) \big\}$ and 
\begin{equation}
\sum_{(\vec{v}_1, \vec{v}_2)\in I_{(-1,-1)}^{(0,1,1)}}\, \Zcal_{\vec{v}_1,\vec{v}_2}(\varepsilon_1, \varepsilon_2, a; \qsf)=\qsf^{\frac{2}{3}}+\cdots \ .
\end{equation}
\end{itemize}

These results coincide in all but one case with the computations in
\cite[Appendix~C]{art:bonellimaruyoshitanziniyagi2012}, which are
claimed to coincide with those of
\cite{art:fucitomoralespoghossian2004}. In particular, the coefficient
of $\qsf^{\frac{2}{3}}$ in the expansion of the function $\sum_{(\vec{v}_1, \vec{v}_2)\in
  I_{(0,0)}^{(0,1,1)}}\, \Zcal_{\vec{v}_1,\vec{v}_2}(\varepsilon_1,
\varepsilon_2, {a}; \qsf)$ does not coincide with the
corresponding coefficient of $Z^{(0,0)}(1,2)$ (to compare these two sets
of results, one needs to set $C\vec{v}=-\big(c_1^{(1)}, c_1^{(2)}
\big)$). As pointed out in Remark \ref{rem:edge-choice}, our edge
contributions depend on some choices; since different choices give
equivalent expressions, we believe that one can find the correct
choices in the procedure described in Appendix
\ref{sec:edgecontribution}  which yields agreement with the results of \cite {art:bonellimaruyoshitanziniyagi2012}. 

Since our partition functions are defined by means of integrals over
moduli spaces of framed sheaves on $\Xscr_k$, it seems more natural
to focus on a comparison between our moduli spaces and the
moduli spaces of $\Z_k$-equivariant framed sheaves on the projective
plane $\PP^2$; these moduli spaces may be regarded as a geometric foundation for the 
computations of the partition functions  done in
\cite{art:fucitomoralespoghossian2004}. We will address this problem in
a future work.

\end{document}